\documentclass[10pt,a4paper]{article}
\usepackage[utf8]{inputenc}
\usepackage[T1]{fontenc}
\usepackage[english]{babel}
\usepackage{amsmath,amssymb,mathtools,array,booktabs,longtable,bm}
\usepackage[a4paper,top=18mm,bottom=20mm,left=23mm,right=23mm]{geometry}
\usepackage{fancyhdr,graphicx,xcolor,enumitem,docmute,etoolbox,hyperref,tikz}
\usetikzlibrary{calc}
\hypersetup{colorlinks=false,pdfborder={0 0 0},hypertexnames=false,
pdftitle={New Symbolic Procedures in the Study of Dynamical Systems: English Translation of the 1995 Doctoral Thesis},
pdfauthor={Daniel Condurache},
pdfsubject={1995 doctoral thesis -- conforming English translation and audited LaTeX reconstruction, 2026},
pdfkeywords={symbolic representations, dynamical systems, algebras, mechanics, non-inertial frames}}
\newcommand{\eqprefix}{I.1}
\newcommand{\eqn}[2]{%
  \ifstrequal{#1}{*}{\begin{equation*}#2\end{equation*}}%
  {\begin{equation}\tag{\eqprefix.#1}#2\end{equation}}}
\newcommand{\alg}[1]{\bm{#1}}
\newcommand{\vct}[1]{\boldsymbol{#1}}
\newcommand{\onehat}{\widehat{1}}\newcommand{\ihat}{\widehat{\imath}}
\newcommand{\ehat}{\widehat{\varepsilon}}\newcommand{\epshat}{\widehat{\varepsilon}}
\newcommand{\what}{\widehat w}\newcommand{\vhat}{\widehat v}
\newcommand{\xor}{\mathbin{\oplus}}
\newcommand{\wh}[1]{\widehat{#1}}\newcommand{\W}{\mathcal W}
\newcommand{\R}{\mathbb R}\newcommand{\C}{\mathbb C}\newcommand{\N}{\mathbb N}
\newcommand{\z}{\widehat 0}\newcommand{\one}{\widehat 1}\newcommand{\zero}{\widehat 0}
\newcommand{\A}{\bigotimes C^{\,n+1}}\newcommand{\FT}{\widehat{\mathcal F}}
\newcommand{\eps}{\widehat\varepsilon}\newcommand{\jhat}{\widehat j}
\newcommand{\Lh}{\widehat{\mathcal L}}
\newcommand{\ep}{\widehat e_p}\newcommand{\mat}[1]{\widetilde{#1}}
\newcommand{\tens}[1]{\widetilde{#1}}
\newcommand{\figeight}{%
\begin{tikzpicture}[scale=.72,>=stealth]
  \coordinate (O) at (0,0); \coordinate (A) at (1.05,3.2); \coordinate (B) at (4.2,.75);
  \coordinate (C) at ($(A)!.55!(B)$);
  \draw[->,thick] (O)--(A) node[midway,left] {$\vct r_1$};
  \draw[->,thick] (O)--(B) node[midway,below] {$\vct r_2$};
  \draw[->,thick] (O)--(C) node[midway,right] {$\vct r_c$};
  \draw[->,thick] (A)--(C) node[midway,above] {$\vct F_{12}$};
  \draw[->,thick] (B)--(C) node[midway,below] {$\vct F_{21}$};
  \fill (A) circle (1.3pt) (B) circle (1.3pt) (C) circle (1.3pt);
  \node[above] at (A) {$P_1(m_1)$}; \node[right] at (B) {$P_2(m_2)$};
  \node[above right] at (C) {$C$}; \node[below left] at (O) {$O$};
\end{tikzpicture}}

\newcommand{\fignine}{%
\begin{tikzpicture}[scale=.72,>=stealth]
  \coordinate (O) at (0,0); \coordinate (A) at (.9,3.3); \coordinate (B) at (4.15,.65);
  \coordinate (C) at ($(A)!.55!(B)$);
  \draw[->,thick] (O)--(A) node[midway,left] {$\vct r_1^{0}$};
  \draw[->,thick] (O)--(B) node[midway,below] {$\vct r_2^{0}$};
  \draw[->,thick] (O)--(C) node[midway,right] {$\vct r_c^{0}$};
  \draw (A)--(B);
  \draw[->] (A)--++(-.65,.75) node[left] {$\vct v_1^{0}$};
  \draw[->] (B)--++(1,.45) node[right] {$\vct v_2^{0}$};
  \draw[->] (C)--++(.15,1.0) node[right] {$\vct v_c^{0}$};
  \fill (A) circle (1.3pt) (B) circle (1.3pt) (C) circle (1.3pt);
  \node[above] at (A) {$P_1^{0}$}; \node[below] at (B) {$P_2^{0}$};
  \node[right] at (C) {$C^0$}; \node[below left] at (O) {$O$};
\end{tikzpicture}}

\newcommand{\figten}{%
\begin{tikzpicture}[scale=.72,>=stealth]
  \coordinate (O) at (0,0); \coordinate (C) at (0,3.3);
  \draw[dashed] (O)--(-1.65,3.15) arc[start angle=180,end angle=360,x radius=1.65,y radius=.48]--cycle;
  \draw[->,thick] (O)--(1.15,3.48) node[right] {$\vct H(t)$};
  \draw[->] (O)--(0,4.15) node[above] {$\vct\omega$};
  \draw[->] (1.5,3.25) arc[start angle=10,end angle=105,x radius=.55,y radius=.22];
  \node[left] at (-1.8,3.25) {$\pi$}; \node[below] at (O) {$O$};
\end{tikzpicture}}

\newcommand{\figeleven}{%
\begin{tikzpicture}[scale=.72,>=stealth]
  \coordinate (O) at (0,0); \coordinate (A) at (.3,3.1); \coordinate (B) at (3.6,.55);
  \draw[->,thick] (O)--(A) node[midway,left] {$\vct r_1^{0}$};
  \draw[->,thick] (O)--(B) node[midway,below] {$\vct r_2^{0}$};
  \draw[->,thick] (A)--(B) node[midway,right] {$\vct r^{0}$};
  \draw[->] (A)--++(-.8,.7) node[left] {$\vct v_1^{0}$};
  \draw[->] (B)--++(1,.55) node[right] {$\vct v_2^{0}$};
  \node[below left] at (O) {$O$}; \node[above] at (A) {$P_1^{0}$}; \node[below] at (B) {$P_2^{0}$};
\end{tikzpicture}}

\newcommand{\figtwelve}{%
\begin{tikzpicture}[scale=.72,>=stealth]
  \coordinate (O) at (0,0); \coordinate (A) at (.9,3.1); \coordinate (B) at (3.5,.7);
  \draw[->,thick] (O)--(A) node[midway,left] {$\vct r_1^{0}$};
  \draw[->,thick] (O)--(B) node[midway,below] {$\vct r_2^{0}$};
  \draw[->,thick] (A)--(B) node[midway,right] {$\vct r^{0}$};
  \draw[->] (A)--++(-1,.65) node[left] {$\vct v_1^{0}+\vct\omega_0\!\times\!\vct r_1^{0}$};
  \draw[->] (B)--++(1,.65) node[right] {$\vct v_2^{0}+\vct\omega_0\!\times\!\vct r_2^{0}$};
  \node[below left] at (O) {$O$}; \node[above] at (A) {$P_1^{0}$}; \node[below] at (B) {$P_2^{0}$};
\end{tikzpicture}}

\newcommand{\figthirteen}{\figten}
\newcommand{\figfourteen}{%
\begin{tikzpicture}[scale=.72,>=stealth]
  \coordinate (O) at (0,0);
  \draw[dashed] (O)--(-1.55,3.05) arc[start angle=180,end angle=360,x radius=1.55,y radius=.46]--cycle;
  \draw[->,thick] (O)--(.95,3.35) node[right] {$\vct L_0(t)$};
  \draw[->] (O)--(0,4.05) node[above] {$\vct\omega$};
  \draw[->] (1.45,3.12) arc[start angle=5,end angle=100,x radius=.52,y radius=.20];
  \node[left] at (-1.7,3.15) {$\pi$}; \node[below] at (O) {$O$};
\end{tikzpicture}}

\newcommand{\figfifteen}{%
\begin{tikzpicture}[scale=.78,>=stealth]
  \coordinate (O) at (0,0); \coordinate (A) at (.7,3.4); \coordinate (B) at (4.2,.55);
  \coordinate (C) at ($(A)!.58!(B)$);
  \draw[->,thick] (O)--(A) node[midway,left] {$\vct r_1$};
  \draw[->,thick] (O)--(B) node[midway,below] {$\vct r_2$};
  \draw[->,thick] (O)--(C) node[midway,right] {$\vct r_c$};
  \draw[->] (C)--(A) node[midway,left] {$\vct r'_1$};
  \draw[->] (C)--(B) node[midway,right] {$\vct r'_2$};
  \fill (A) circle (1.3pt) (B) circle (1.3pt) (C) circle (1.3pt);
  \node[above] at (A) {$P_1(m_1)$}; \node[right] at (B) {$P_2(m_2)$};
  \node[above right] at (C) {$C$}; \node[below left] at (O) {$O$};
\end{tikzpicture}}

\newcommand{\setchapterheader}[1]{\fancyhf{}\fancyhead[L]{\thepage}\fancyhead[R]{#1}}

\makeatletter
\let\oldsection\section
\let\oldsubsection\subsection
\let\oldsubsubsection\subsubsection
\newcommand{\sectionstar}[1]{\oldsection*{#1}\addcontentsline{toc}{section}{#1}}
\newcommand{\subsectionstar}[1]{\oldsubsection*{#1}\addcontentsline{toc}{subsection}{#1}}
\newcommand{\subsubsectionstar}[1]{\oldsubsubsection*{#1}\addcontentsline{toc}{subsubsection}{#1}}
\newcommand{\fourthsectionstar}[1]{%
  \par\medskip\noindent\textbf{#1}\par\smallskip
  \phantomsection\addcontentsline{toc}{paragraph}{#1}}
\renewcommand{\section}{\@ifstar{\sectionstar}{\oldsection}}
\renewcommand{\subsection}{\@ifstar{\subsectionstar}{\oldsubsection}}
\renewcommand{\subsubsection}{\@ifstar{\subsubsectionstar}{\oldsubsubsection}}
\makeatother

\newcommand{\chaptertoc}[2]{%
  \clearpage\phantomsection\addcontentsline{toc}{part}{#1. #2}%
  \setchapterheader{#1. #2}}

\begin{document}
\begin{titlepage}
\centering
\vspace*{18mm}
{\Large ``GHEORGHE ASACHI'' TECHNICAL UNIVERSITY OF IAȘI\par}
\vspace{7mm}
{\large Faculty of Machine Manufacturing\par}
\vfill
{\Large\bfseries DANIEL CONDURACHE\par}
\vspace{12mm}
{\Huge\bfseries New Symbolic Procedures\par}
\vspace{3mm}
{\Huge\bfseries in the Study of Dynamical Systems\par}
\vspace{14mm}
{\Large Doctoral Thesis, 1995\par}
\vspace{4mm}
{\large Conforming English Translation and Audited \LaTeX{} Reconstruction, 2026\par}
\vfill
{\large Iași\par}
\end{titlepage}

\pagenumbering{roman}
\section*{Editorial Note to the 2026 English Edition}

\subsection*{Status and provenance}

This document is a conforming English translation and an audited \LaTeX{}
reconstruction of the doctoral thesis \emph{Noi procedee simbolice în studiul
sistemelor dinamice}, defended by Daniel Condurache at the ``Gheorghe Asachi''
Technical University of Iași in 1995. It is not a revised 2026 research
monograph. The mathematical results, hypotheses, equation numbering,
bibliography, applications, and historical statements reproduce the 1995
thesis. The intervention made in 2026 is limited to translation, typographical
reconstruction, recovery of the figures, and structural verification. No
post-1995 mathematical result has been inserted into the thesis text.

The original work was produced in a period in which international scientific
communication in Romania depended largely on printed journals, exchanged
offprints, and locally available books. The present edition makes the full
text searchable and internationally accessible while preserving its historical
state. Consequently, statements concerning novelty or the literature are to
be read as statements made in 1995, not as claims based on a 2026 literature
survey.

The reconstruction was audited module by module against the conforming
Romanian source. The order and structure of the displayed mathematics,
equation identifiers, internal numerical citations, and all 164 bibliographic
records were preserved. Bibliographic titles remain in their original
publication languages.

\subsection*{Contemporary abstract}

The thesis develops symbolic procedures for the analysis of dynamical systems
by representing spaces of scalar- and vector-valued signals in finite-
dimensional real algebras. Its algebraic core includes explicit bases,
orthogonal idempotents, zero divisors, and exponential forms for direct-product
complex algebras, a generalized Sobrero algebra, Walsh algebras, and polynomial
algebras. Binary, Gray-code, and Walsh--Hadamard indexing are used to organize
the idempotent channels and the associated spectral components.

On this structure the thesis constructs symbolic representations of continuous
and discrete signals, a procedure for nonlinear-system identification, and
algebra-valued Fourier- and Laplace-type integral transforms together with
their inversion and energy relations. The explicit idempotent decomposition
provides the computational channel structure needed to define, evaluate, and
invert these transforms. The framework is then extended to vector-valued
functions through matrix and tensor symbolic representations. Applications
include motion in non-inertial frames, charged-particle motion in electromagnetic
fields, the Foucault pendulum, central positional-force fields, and the two-body
problem in arbitrarily rotating frames. The document is published in this
historically conforming form because the conjunction of explicit algebraic
decomposition, signal interpretation, integral transforms, and dynamical
applications remains relevant to current work on multicomplex and related
finite-dimensional commutative algebras.

\noindent\textbf{Keywords.}
Symbolic representation; dynamical systems; finite-dimensional commutative
algebras; orthogonal idempotents; zero divisors; Gray code; Walsh--Hadamard
matrices; Fourier transform; Laplace transform; non-inertial mechanics.

\subsection*{Authorship and rights statement}

Daniel Condurache is the author of the 1995 doctoral thesis and of this
conforming English edition. The reconstructed figures derive from the
author's thesis copy. No third-party edition or later publisher-held text has
been incorporated. Distribution of this version is authorized by the author
under the license selected for its repository deposit.

\clearpage
{\small\let\section\oldsection\tableofcontents}
\clearpage
\pagenumbering{arabic}

\chaptertoc{CHAPTER I}{Introduction}
\renewcommand{\eqprefix}{I.1.1}
\thispagestyle{empty}
\begin{center}{\Large\bfseries CHAPTER I}\par\vspace{4mm}
{\large\bfseries Introduction}\end{center}
\section*{I.1. Algebraic and Geometric Procedures in the Study of Dynamical Systems}
\subsection*{I.1.1. Algebraic and Operational Procedures}

The theory of dynamical systems has its origins in Count Louis Lagrange's
(1736--1813) studies of the motion of mechanical systems in the neighborhood
of equilibrium configurations.  Within this framework, the concepts of natural
frequency, damping factor, free and forced motion, resonance, and modulation
began to take shape.

Heaviside operational calculus is a symbolic procedure whose elegance has
fascinated generations of researchers.  Introduced heuristically in 1887 for
electrical circuits, it denoted $D=d/dt$ and treated $(d^n/dt^n)f(t)$ as
$D^nf$, thereby transforming the equation
\eqn{1}{\sum_{k=0}^{n}a_k\frac{d^kx}{dt^k}=f(t)}
into $P(D)x=f(t)$, where $P(D)=\sum_{k=0}^{n}a_kD^k$, and formally yielded
$x=f(t)/P(D)$.  Its meaning was recovered by decomposing $[P(D)]^{-1}$ into
partial fractions and identifying the resulting terms with elementary
differential equations.  Oliver Heaviside (1850--1925) did not insist on a
formal justification, but the procedure anticipated the productive idea of
algebraizing a differential equation and determining its solution symbolically.

A symbolic representation due to Fresnel [18] associates with the harmonic
signal $s(t)=A\sin(\omega t+\varphi)$, $A>0$, $\varphi\in[0,2\pi)$, first a
phasor and then a complex number:
\eqn{2}{s=A\sin(\omega t+\varphi)\xmapsto{\mathcal F}
\widehat s=Ae^{j\varphi}.}

The representation is $\mathbb R$-linear and has the property
\eqn{3}{\dot s=\omega A\sin\left(\omega t+\varphi+\frac\pi2\right)
\xmapsto{\mathcal F}j\omega Ae^{j\varphi}.}
If $f\in S_\omega$ in (1), there exists a solution $x\in S_\omega$ that can be
obtained algebraically.  Setting $\widehat f=\mathcal F(f)$ and
$\widehat x=\mathcal F(x)$ gives
\eqn{4}{\sum_{k=0}^{n}a_k(j\omega)^k\widehat x=\widehat f,}
with the algebraic solution
\eqn{5}{\widehat x=\frac{\widehat f}{\displaystyle\sum_{k=0}^{n}a_k(j\omega)^k}.}
The solution is a harmonic oscillation of the same angular frequency, with
amplitude and phase given by the modulus and argument of the complex number
(5).  A vanishing denominator corresponds to resonance with an unbounded
response.  This procedure marks the beginning of frequency-domain methods in
the analysis of dynamical systems.

Frequency-domain methods also apply to periodic disturbances by expanding them
in Fourier series [23] and treating each harmonic separately.  Linearity makes
it possible to recover the response series, although this series may converge
slowly or even diverge [23], [105].

The heuristic interpretation of pulses as periodic functions with period
$T\to\infty$ led to the Fourier transform.  With a signal $f\in L_1$ satisfying
$\int_{-\infty}^{\infty}|f(t)|dt<\infty$, one associates
\eqn{6}{\mathcal Ff(j\omega)=\int_{-\infty}^{\infty}f(t)e^{-j\omega t}\,dt,
\qquad\omega\in\mathbb R,\quad j^2=-1.}

The mapping (6) is $\mathbb R$-linear and preserves the property
\eqn{7}{\mathcal F\dot f(j\omega)=j\omega\,\mathcal Ff(j\omega).}
Under the conditions specified in [151], [23], the signal is recovered by
\eqn{8}{\frac{f(t+0)+f(t-0)}2
=\frac1{2\pi}\int_{-\infty}^{\infty}\mathcal Ff(j\omega)e^{j\omega t}\,d\omega.}
This formula expresses the synthesis of the signal in terms of harmonic signals
[28].  Property (7) algebraizes differential equations and leads to the transfer
function of a linear time-invariant system, namely the ratio between the
transforms of the output and input quantities:
\eqn{9}{H(j\omega)=\frac{\mathcal Fe(j\omega)}{\mathcal Fi(j\omega)}.}
The transfer function depends only on the characteristics of the system and can
be determined experimentally, thereby opening the way to dynamical-system
identification [61], [143], [139].

For signals that do not possess a Fourier transform in the classical sense,
the Laplace transform is introduced.  It associates a complex-valued function
with an original signal $f\in\Omega$ that is locally integrable, supported on
$[0,+\infty)$, and piecewise continuous:
\eqn{10}{\mathcal Lf(s)=\int_0^{\infty}f(t)e^{-st}\,dt,\qquad s\in\mathbb C.}
The mapping is $\mathbb R$-linear and satisfies
\eqn{11}{\mathcal L\dot f(s)=s\mathcal Lf(s)-f(0).}

Under the conditions specified in [23], [151], the original signal is recovered
from
\eqn{12}{\frac{f(t+0)+f(t-0)}2
=\frac1{2\pi j}\int_{\delta-j\infty}^{\delta+j\infty}
\mathcal Lf(s)e^{st}\,ds,\qquad\delta>a_f,}
where $a_f$ is the abscissa of convergence.  The formula expresses synthesis in
terms of harmonic signals whose amplitudes are modulated by an exponential
signal [23].

For a linear time-invariant dynamical system, the Laplace transfer function is
\eqn{13}{H(s)=\frac{\mathcal Le(s)}{\mathcal Li(s)}.}
It depends only on the system characteristics, and its expression is obtained
through elementary algebraic calculations.  The Laplace transform brings the
theory of complex functions into dynamical-system theory; the holomorphy and
regularity properties of the transfer function acquire particular significance
[151], with practical applications [152], [158].

Constant, step, or polynomial signals do not always admit classical transforms,
which led to generalized functions, or distributions [151], [156].  Functional
analysis makes it possible to differentiate discontinuous functions and to
differentiate convergent series of distributions indefinitely.  The notion of
a signal is thus extended to tempered distributions $\mathcal S'$, and the
Fourier transform becomes an isomorphism of these spaces.

For linear time-invariant dynamical systems, the input--output relation has the
convolution form
\eqn{14}{x*h=f.}
Here $f$ is the input signal, $x$ is the output signal, and $h$ is the system's
response distribution, also called its weighting function, obtained for the
Dirac input $\delta$.  The Laplace transfer function is the transform of the
weighting function.

Particularly noteworthy is the essentially algebraic construction of the space
of distributions due to the Polish school led by J. Mikusiński [57];
distributions are the elements of the field of fractions of the

J. Mikusiński's construction [57] regards distributions as elements of the field
of fractions of the integral domain of functions, with convolution as the
product.  The absence of zero divisors follows from Titchmarsh's theorem: the
convolution of two functions vanishes if and only if at least one of them
vanishes.

This approach justified Heaviside-type symbolic procedures and made it possible
to represent dynamical systems by operators between signal spaces [156], [116].
Their geometric treatment as Hilbert spaces with a resolution of the identity,
together with spectral theory, led to optimal Karhunen--Loève-type transforms
[118].

The sampling theorem [9], [151] permits the recovery of a band-limited continuous
signal from its uniformly spaced samples.  Discrete signals became objects in
their own right, while digital filters are modeled by finite-difference equations
[73], [77].  These equations are algebraized by the $Z$-transform [151], [26].
The discrete Fourier transform and FFT algorithms enable real-time analysis,
with applications to system identification and diagnosis [61], [43].

For nonlinear systems [22], [63]--[68], consider the equation
\eqn{15}{\sum_{k=0}^{n}a_k\frac{d^kx}{dt^k}
+\sum_{p=2}^{m}b_px^p=f(t),\qquad a_k,b_p\in\mathbb R.}
Michael Fliess's method uses formal languages and automata.  For the alphabet
$X=\{x_0\}$, let $X^*$ be the free monoid, with $1$ as the empty word.  The
algebra of formal series $\mathbb R\langle\!\langle X\rangle\!\rangle$ consists
of the elements
\eqn{16}{s=\sum_{w\in X^*}(s,w)w,\qquad(s,w)\in\mathbb R.}
The operations are
\eqn{17}{s_1+s_2=\sum_{w\in X^*}[(s_1,w)+(s_2,w)]w.}

\eqn{18}{s_1s_2=\sum_{w\in X^*}
\left\{\sum_{u\cdot v=w}(s_1,u)(s_2,v)\right\}w,}
\eqn{19}{\lambda s=\sum_{w\in X^*}[\lambda(s,w)]w.}
The Hurwitz product is introduced on this algebra by
\eqn{20}{s_1\otimes s_2=\sum_{\substack{w_1\in X^*\\w_2\in X^*}}
(s_1,w_1)(s_2,w_2)(w_1\otimes w_2),}
where
\eqn{21}{x_0^k\otimes x_0^{n-k}=\binom nkx_0^n,
\qquad k\leq n,\quad k,n\in\mathbb N.}
The series $s$ is invertible if there exists an $s^{-1}$ such that
$ss^{-1}=s^{-1}s=1$.  If $(s,1)\neq0$, write
$s=(s,1)(1-s')$, $(s',1)=0$, and
\[
s^{-1}=\frac1{(s,1)}(1-s')^{-1}=\frac1{(s,1)}\sum_{k\geq0}(s')^k.
\]

Let $\mathcal A$ be the set of real-valued functions analytic at the origin,
\eqn{22}{f(t)=\sum_{k=0}^{\infty}f_k\frac{t^k}{k!},\qquad f_k\in\mathbb R.}
Define the transform
\eqn{23}{\mathcal F:\mathcal A\to\mathbb R\langle\!\langle X\rangle\!\rangle,
\qquad\mathcal Ff=\sum_{k=0}^{\infty}f_kx_0^k.}
The transform is $\mathbb R$-linear and has the properties

\eqn{24}{\mathcal F\!\left(\int_0^t f(\tau)\,d\tau\right)=x_0\mathcal Ff,}
\eqn{25}{\mathcal F(f_1f_2)=\mathcal Ff_1\otimes\mathcal Ff_2.}
The transform can be derived from the Laplace transform [22]:
\eqn{26}{\mathcal Ff=\left.s\mathcal Lf(s)\right|_{s=x_0^{-1}}.}
For example,
$\mathcal F\sigma(t)=1$, $\mathcal F(t^n/n!)=x_0^n$, and
$\mathcal F(\cos\omega t)=(1-\omega^2x_0^2)^{-1}$.

For zero initial conditions, equation (15) is equivalent to the integral equation
\begin{equation}\tag{\eqprefix.27}\begin{aligned}
a_nx+a_{n-1}\int_0^t x(\tau_1)d\tau_1+\cdots
&+a_0\int_0^t d\tau_n\int_0^{\tau_n}\!\cdots
\int_0^{\tau_2}x(\tau_1)d\tau_1\\
&+\sum_{p=2}^{m}b_p\int_0^t d\tau_n\cdots
\int_0^{\tau_2}x^p(\tau_1)d\tau_1\\
&=\int_0^t d\tau_n\cdots\int_0^{\tau_2}f(\tau_1)d\tau_1.
\end{aligned}\end{equation}
Setting $g=\mathcal Fx$, equations (24)--(25) yield
\eqn{28}{\left(\sum_{k=0}^{n}a_kx_0^{n-k}\right)g
+x_0^n\sum_{p=2}^{m}b_pg^{\otimes p}=x_0^{n-1}\mathcal Ff.}
The transform of the response therefore satisfies an algebraic equation of
degree $m$ in the algebra of formal series.  The completeness of this algebra
allows a fixed-point theorem to be applied and the solution to be constructed
iteratively as
\eqn{29}{g=g_1+g_2+\cdots+g_n+\cdots.}

The terms of the iteration are
\begin{equation}\tag{\eqprefix.30}\begin{aligned}
g_1&=\left(\sum_{k=0}^{n}a_kx_0^{n-k}\right)^{-1}x_0^{n-1}\mathcal Ff,\\
g_n&=\left(\sum_{k=0}^{n}a_kx_0^{n-k}\right)^{-1}x_0^n
\sum_{p=2}^{m}b_p\!\sum_{n_1+\cdots+n_p=n}
g_{n_1}\otimes\cdots\otimes g_{n_p}.
\end{aligned}\end{equation}
The algorithm can be implemented in symbolic languages such as PL/I, REDUCE,
or MACSYMA.  Can and Ünal [22] define the transfer function of system (15) as
the transform of its response to the unit step.  Since $\mathcal F\sigma=1$,
it satisfies
\eqn{31}{\left(\sum_{k=0}^{n}a_kx_0^{n-k}\right)H
+x_0^n\sum_{p=2}^{m}b_pH^{\otimes p}=x_0^{n-1}.}
If $H$ is the transfer function, the transform of the response is
\eqn{32}{\mathcal Fx=H\cdot\mathcal Ff.}
Fliess-series calculus generalizes Heaviside symbolic calculus to nonlinear
systems.

\renewcommand{\eqprefix}{I.1.2}
\subsection*{I.1.2. Symbolic Representations}
Let $V$ and $A$ be two real vector spaces, each endowed with a linear unary
operation:
\eqn{1}{\begin{gathered}
``\dot{\ }'':V\to V,\\
\dot{(\alpha a+\beta b)}=\alpha\dot a+\beta\dot b,
\qquad a,b\in V,\quad\alpha,\beta\in\mathbb R,
\end{gathered}}
\eqn{2}{\begin{gathered}
``*'':A\to A,\\
(\mu x+\lambda y)^*=\mu x^*+\lambda y^*,
\qquad x,y\in A,\quad\mu,\lambda\in\mathbb R.
\end{gathered}}

\renewcommand{\eqprefix}{I.1.2}

The operations introduced above are linear endomorphisms of the spaces $V$ and
$A$, respectively.

\noindent\textbf{Definition I.1.2.} The vector space $V$ is represented
symbolically, with respect to the operator ``$\dot{\ }$'', on the space $A$
equipped with the operator ``$*$'', if there exists a homomorphism
$\Psi:V\to A$ such that
\eqn{3}{\Psi(\dot f)=\Psi(f)^*,\qquad\forall f\in V.}
If $\Psi$ is an isomorphism, the representation is exact.  The space $A$ is
the representation space, and ``$*$'' is the operator of the symbolic
representation.

For $p\in\mathbb R[X]$, $p(X)=X^n+a_1X^{n-1}+\cdots+a_n$, define
\eqn{4}{p(\dot{\ })f=f^{(n)}+a_1f^{(n-1)}+\cdots+a_{n-1}\dot f+a_nf,}
\eqn{5}{p(*)w=w^{[n]}+a_1w^{[n-1]}+\cdots+a_{n-1}w^*+a_nw.}
By induction, equation (3) gives
\eqn{6}{\Psi[p(\dot{\ })f]=p(*)\Psi(f),\qquad f\in V, p\in\mathbb R[X].}
This transfers linear operator equations from $V$ to $A$, a procedure called
the \emph{algebraization of the operator equation}.

\noindent\textbf{Examples.}

$1^\circ$ For the space of infinitely differentiable real-valued functions,
$\dot f=df/dt$, and an associative, commutative, unital real algebra in which
$w^*=d\,w$ for a fixed $d$, one obtains Alfred Braier's representation [18],
[20], used in vibration mechanics.

$2^\circ$ For the space $V_3$ of free vectors, let
\[
\dot{\vct f}=\vct\omega\times\vct f,
\qquad [\vct f]^*=\widetilde\omega[\vct f],
\]
where $\widetilde\omega$ is the skew-symmetric matrix associated with
$\vct\omega$.  The mapping
\[
\Psi:V_3\to\mathcal V_3,\qquad
\Psi(\vct f)=[f_1,f_2,f_3]^T
\]
is an isomorphism and satisfies
\[
\Psi(\vct\omega\times\vct f)=\widetilde\omega[\vct f],
\]
and therefore defines the exact symbolic representation proposed by
N. Irimiciuc [88], [89].

$3^\circ$ Let $L_1$ be the space of Lebesgue-integrable functions on
$\mathbb R$ that are differentiable and whose derivative also belongs to
$L_1$, and let
\[
\mathcal F=\{F:\mathbb R\to\mathbb C\mid F=F(j\omega),\ j^2=-1\},
\qquad F^*=j\omega F.
\]
The Fourier transform
\[
\Psi:L_1\to\mathcal F,\qquad
\Psi(f)=\int_{-\infty}^{\infty}f(t)e^{-j\omega t}\,dt
\]
is a homomorphism and satisfies
\[
\Psi(\dot f)=j\omega\Psi(f),
\]
and is therefore a symbolic representation.

$4^\circ$ Let $\Omega$ be the space of original functions that are locally
integrable, supported on $[0,+\infty)$, and piecewise continuous.  For
\[
\mathcal L=\{F:\mathbb C\to\mathbb C\mid F=F(s),\ s=\delta+j\omega\},
\qquad F^*=sF-f(0),
\]
the Laplace transform
\[
\Psi:\Omega\to\mathcal L,\qquad
\Psi(f)=\int_0^{\infty}f(t)e^{-st}\,dt
\]
satisfies $\Psi(\dot f)=s\Psi(f)-f(0)$ and is a symbolic representation.

$5^\circ$ Consider the space of discrete signals
\[
S_d=\{f:\mathbb Z\to\mathbb R\mid\operatorname{supp}f\subseteq[0,+\infty),
\ f=f(n)\},
\]
with the right-shift operator
\[
\dot f(n)=f(n+1).
\]
Let $\mathcal L_Z$ be the space of Laurent series convergent in a circular
annulus:
\[
\mathcal L_Z=\left\{L(z)=\sum_{n=-\infty}^{\infty}a_nz^{-n},\quad
a_n\in\mathbb R,\quad r<|z|<R\right\},
\]
with the operator $L^*(z)=zL(z)$.  The $Z$-transform [151]
\[
\Psi:S_d\to\mathcal L_Z,\qquad
\Psi(f)=\sum_{n=-\infty}^{\infty}f(n)z^{-n}
\]
has the representation property associated with signal shifting.

The $Z$-transform satisfies
\[
\Psi[f(n+1)]=z\,\Psi[f(n)],
\]
and is therefore a symbolic representation in the sense of Definition I.1.2.
The preceding examples illustrate the consistency of the definition.  This
framework encompasses almost all classical procedures used in the analysis of
linear dynamical systems and allows new symbolic representations to be
introduced.

\section*{I.2. Scope and Organization of the Thesis}
The thesis begins with the symbolic representation of a vector space endowed
with an $\mathbb R$-linear unary operation on another vector space equipped
with an endomorphism.  The representation is a homomorphism that intertwines
the two unary operations.  The transfer of operator equations to a more
convenient space is called the \emph{algebraization of the operator equation}.

Chapter II solves the problem of symbolic representation on finite-dimensional,
associative, commutative, unital real algebras.  It determines the spaces that
can be represented for a fixed algebra and operator, as well as the algebras
and operators capable of representing the solutions of a linear operator
equation.  Representations of continuous [20] and discrete [38], [39] signals
are included.

Chapter III studies the structural properties of direct products of complex
algebras $\bigotimes C^{n+1}$, generalized Sobrero algebras, and Walsh
algebras: zero divisors, orthogonal idempotents, nilpotent elements, and
exponential forms.  Generating polynomials and cyclic bases are determined;
some of these results were published in [43], [47].

Chapter IV is devoted to applications.  It characterizes the signals
representable on $\bigotimes C^{n+1}$ and on the generalized Sobrero algebra,
as well as the distinguished elements of these algebras.

A representation of continuous, periodic, integrable signals on the Walsh
algebra is introduced, leading to a parametric-identification procedure for
nonlinear systems with analytic nonlinearities.  An algebraic procedure is
also proposed for determining the response of systems to representable
excitations, using transmittance in the sense of symbolic representation.

Chapter V studies integral transforms on finite-dimensional real algebras,
introduced by the author between 1982 and 1987 [40]--[42], [45].  A
Fourier-type transform on $\bigotimes C^{n+1}$ and a Laplace-type transform on
the generalized Sobrero algebra are analyzed.  Their monogenic character is
proved, and they are extended to the class of distributions.

Chapter VI addresses the symbolic representation of real vector-valued
functions.  The matrix representation introduced in [46] algebraizes Cauchy
problems for motion in uniformly rotating non-inertial frames and in gyroscopic
fields.  Exact vector solutions are obtained for motion in a uniform
gravitational field, for electric charges in crossed electric and magnetic
fields, and for the Foucault pendulum; the results generalize Larmor's theorem
[49], [50].

The tensor representation extends the method to nonuniform rotations about a
fixed direction, providing exact solutions and first integrals for uniform and
central positional fields, as well as for the motion of charges in a slowly
varying electromagnetic field.

The final extension allows arbitrary rotations of the non-inertial frame.  New
first integrals in central fields, generalizations of the Foucault pendulum and
Larmor's theorem, and a complete vector treatment of the two-body problem in
an arbitrarily rotating non-inertial frame are obtained.

\chaptertoc{CHAPTER II}{Symbolic Representations on \texorpdfstring{$\mathbb R$}{R}-Algebras}
\renewcommand{\eqprefix}{II.1}

\thispagestyle{empty}\vspace*{28mm}
\begin{center}\Huge\bfseries CHAPTER II\par\vspace{16mm}
\LARGE Symbolic Representations on $\mathbb R$-Algebras
\end{center}
\vspace{25mm}\section*{II.1. General Results}

Let $A_n$ be a finite-dimensional, commutative, associative, unital
$\mathbb R$-algebra of order $n\in\mathbb N$ [2], [54].  If
$B=\{\alg\varepsilon_j\}_{j=1}^n$ is a basis, every $\alg w\in A_n$ has a
unique expansion of the form
\eqn{1}{\alg w=\sum_{i=1}^{n}w_i\alg\varepsilon_i,\qquad w_i\in\mathbb R.}
Multiplication is determined by the multiplication table of the basis,
\eqn{2}{\alg\varepsilon_i\alg\varepsilon_j=
\sum_{k=1}^{n}\gamma_{ijk}\alg\varepsilon_k,\qquad
\gamma_{ijk}\in\mathbb R.}
The algebra $A_n$ is canonically isomorphic to a commutative algebra of square
matrices.  The element $\alg w$ corresponds to
\eqn{3}{[\alg w]=(a_{jk})_{j,k=1}^{n},}
where
\eqn{4}{a_{jk}=\sum_{i=1}^{n}w_i\gamma_{ijk}.}
For the fixed element
\eqn{5}{\alg d=\sum_{i=1}^{n}d_i\alg\varepsilon_i}
define the linear operator
\eqn{6}{* :A_n\to A_n,\qquad \alg w^*=\alg d\alg w.}
Let $V$ be a real vector space and let $\dot{\ }:V\to V$ be an endomorphism.
We study exact symbolic representations, that is, linear bijections
\eqn{7}{\Psi:V\to A_n}
that satisfy
\eqn{8}{\Psi(\dot f)=\alg d\Psi(f),\qquad f\in V.}

\noindent\textbf{Remarks.}

$1^\circ$ Since $\Psi$ is an isomorphism, $V$ is finite-dimensional and
$\dim_{\mathbb R}V=n$.

$2^\circ$ The operator $\alg d$, if it exists, is unique.  Indeed, if
$\alg d'$ also satisfied
\eqn{9}{\Psi(\dot f)=\alg d'\Psi(f),\qquad f\in V,}
then $(\alg d-\alg d')\Psi(f)=\alg0$.  Choosing $f$ such that
$\Psi(f)=\alg e$, the identity of the algebra, gives $\alg d=\alg d'$.

Two fundamental problems arise: determining the spaces $V$ that can be
represented on a given algebra with a fixed operator, and determining an
algebra and operator capable of representing a given vector space.

\subsection*{II.1.1. Determination of the Vector Spaces Symbolically Representable on a Given $\mathbb R$-Algebra}
\noindent\textbf{Theorem II.1.1.} On the algebra $A_n$, with the fixed
operator $\alg d\in A_n$, one can symbolically represent the solution space of
\eqn{10}{p(\dot{\ })f=0,}
where $p\in\mathbb R[X]$ is the characteristic polynomial of $\alg d$, provided
that the solution space has dimension $n$.

\noindent\textit{Proof.} For the representation (7), the fundamental relation becomes
\eqn{11}{\Psi[p(\dot{\ })f]=p(\alg d)\Psi(f).}
By the Cayley--Hamilton theorem, the matrix $[\alg d]$ has the polynomial
$p(\lambda)=\det([\alg d]-\lambda I_n)$, and the canonical isomorphism gives
\eqn{12}{p(\alg d)=\alg0.}
Equation (11) then yields
\eqn{13}{\Psi[p(\dot{\ })f]=\alg0.}
The bijectivity of $\Psi$ gives
\eqn{14}{p(\dot{\ })f=0.}

To describe $\Psi$, set $f_j=\Psi^{-1}(\alg\varepsilon_j)$.  Equations (2),
(5), and (8) give
\[
\dot f_j=\Psi^{-1}(\alg d\alg\varepsilon_j)
=\sum_{k=1}^{n}\alpha_{jk}f_k,
\]
where
\eqn{15}{\alpha_{jk}=\sum_{i=1}^{n}d_i\gamma_{ijk}.}
The elements represented by the basis $B$ therefore satisfy the operator system
\eqn{16}{\dot f_j=\sum_{k=1}^{n}\alpha_{jk}f_k,\qquad j=1,\ldots,n.}
The exact representation exists if and only if system (16) has $n$ linearly
independent solutions, called a fundamental system.  Since
$[\alg d]=(\alpha_{jk})$, the system takes the matrix form
\eqn{17}{[\dot f]=[\alg d]\,[f].}
If $\{f_j\}_{j=1}^{n}$ is a fundamental system, the general solution of (10)
and its symbol are
\[
f=\sum_{j=1}^{n}c_jf_j,\qquad
\alg f=\sum_{j=1}^{n}c_j\alg\varepsilon_j.
\]
Conversely, the element
\eqn{18}{\alg w=\sum_{j=1}^{n}w_j\alg\varepsilon_j}
represents the solution
\eqn{19}{w=\sum_{j=1}^{n}w_jf_j.}

\noindent\textbf{Remark.} The characteristic polynomial of $\alg d$ has degree $n$:
\eqn{20}{p(\lambda)=(-1)^n\lambda^n+a_1\lambda^{n-1}+\cdots+a_n,
\qquad p(\alg d)=\alg0.}
It is proved in [54], [140] that
\eqn{21}{a_n=\det[\alg d].}
If $\det[\alg d]=0$, then $\alg d$ is a zero divisor, because
\eqn{22}{\alg d\bigl[(-1)^n\alg d^{\,n-1}+a_1\alg d^{\,n-2}
+\cdots+a_{n-1}\alg e\bigr]=\alg0.}
If there exists a polynomial $p_1\in\mathbb R[X]$, $\deg p_1=n_1<n$, such that
$p_1(\alg d)=\alg0$, then the characteristic polynomial is not minimal and the
representable elements satisfy
\eqn{23}{p_1(\dot{\ })f=0.}

\noindent\textbf{Corollary II.1.1.} If the characteristic polynomial of
$\alg d$ is not minimal, the elements representable on $A_n$ with operator
$\alg d$ satisfy
\eqn{24}{p_1(\dot{\ })f=0,\qquad \deg p_1=n_1<n.}
The minimal polynomial is one of the divisors of the characteristic polynomial.

\subsection*{II.1.2. Determination of an Algebra and Operator Capable of Symbolically Representing a Given Vector Space}

\renewcommand{\eqprefix}{II.2}
\phantomsection\addcontentsline{toc}{section}{II.2. Symbolic Representation of Continuous Signals on $\mathbb R$-Algebras}

\noindent\textbf{Completion of the proof of Theorem II.1.2.}
For $p(X)=X^n+a_1X^{n-1}+\cdots+a_n$, the polynomial algebra $P_n$ has the
cyclic basis $B=\{\alg v^0,\alg v,\ldots,\alg v^{n-1}\}$ and the relation
$p(\alg v)=\alg0$.  The operator $\alg d=\alg v$ has the Frobenius matrix
\[
[\alg d]=\begin{bmatrix}
0&1&0&\cdots&0\\0&0&1&\cdots&0\\
\vdots&\vdots&\vdots&\ddots&\vdots\\0&0&0&\cdots&1\\
-a_n&-a_{n-1}&-a_{n-2}&\cdots&-a_1
\end{bmatrix}.
\]
By Theorem II.1.1, every solution of $p(\dot{\ })f=0$ is represented on $P_n$
with the operator $\alg v$.

\noindent\textbf{II.2. Symbolic Representation of Continuous Signals.}
For $V=C^\infty(I)$ and $\dot f=df/dt$, the solutions of
\[
f^{(n)}+a_1f^{(n-1)}+\cdots+a_{n-1}\dot f+a_nf=0
\]
are represented on an algebra whose operator $\alg d$ has the same
characteristic polynomial.  The fundamental system satisfies
$[\dot f]=[\alg d][f]$, with solution
\[
[f]=e^{[\alg d]t}[c],\qquad
e^{[\alg d]t}=\sum_{k\ge0}\frac{[\alg d]^kt^k}{k!}.
\]
In the algebra $A_n$, define
$e^{\alg d t}=\sum_{k\ge0}\alg d^{\,k}t^k/k!$; then
$\frac d{dt}e^{\alg d t}=\alg d e^{\alg d t}$.  The following formulas hold:
\[
e^{\alg d t}=\alg u\cos t+\alg d\sin t\quad(\alg d^2=-\alg u),\qquad
e^{\alg d t}=\alg u\cosh t+\alg d\sinh t\quad(\alg d^2=\alg u),
\]
and, if $\alg d^{p+1}=\alg0$,
$e^{\alg d t}=\sum_{k=0}^{p}\alg d^{\,k}t^k/k!$.

\noindent\textbf{Theorem II.2.3.}
If the minimal polynomial of $\alg d$ has degree $n$ and
\[
e^{\alg d t}=\sum_{j=1}^nf_j(t)\alg\varepsilon_j,
\]
then the functions $f_j$ form a fundamental system.  With
$\alpha_{jk}=\sum_i d_i\gamma_{ijk}$, differentiation of the expansion gives
\[
\dot f_k=\sum_{j=1}^n\alpha_{jk}f_j.
\]

The differential system obtained above can be written in matrix form as
\eqn{18}{[\dot f]=[\alg d]^T[f].}
The functions $f_k$, $k=1,\ldots,n$, also satisfy differential equation
II.2.(3), because the characteristic polynomials of $[\alg d]$ and
$[\alg d]^T$ coincide [54].

To prove the linear independence of the functions $f_j$, consider the matrix
\eqn{19}{[W]=
\begin{bmatrix}
f_1&f_2&\cdots&f_n\\
\dot f_1&\dot f_2&\cdots&\dot f_n\\
\vdots&\vdots&\ddots&\vdots\\
f_1^{(n-1)}&f_2^{(n-1)}&\cdots&f_n^{(n-1)}
\end{bmatrix}.}
Its determinant is the Wronskian of the system [8].  It is sufficient to show
that $\det[W]_{t=0}\ne0$.  Differentiating relation II.2.(14) $k$ times gives
\eqn{20}{\alg d^{\,k}e^{\alg d t}=
f_1^{(k)}\alg\varepsilon_1+\cdots+
f_n^{(k)}\alg\varepsilon_n,\quad k=0,\ldots,n-1.}
At $t=0$, this yields
\eqn{21}{\alg d^{\,k}=
f_1^{(k)}(0)\alg\varepsilon_1+\cdots+
f_n^{(k)}(0)\alg\varepsilon_n.}
With the notation
\eqn{22}{[\alg d]_{\mathrm c}=
[\alg d^0,\alg d^1,\ldots,\alg d^{\,n-1}]^T,}
\eqn{23}{[\alg\varepsilon]=
[\alg\varepsilon_1,\alg\varepsilon_2,\ldots,\alg\varepsilon_n]^T,}
relations (21) become
\eqn{24}{[\alg d]_{\mathrm c}=[W]_{t=0}[\alg\varepsilon].}
If $\det[W]_{t=0}=0$, the powers in (22) would be linearly dependent, contrary
to the fact that the minimal polynomial of $\alg d$ has degree $n$.  Hence the
functions $f_j$ form a fundamental system.

\noindent\textbf{Remarks.}

$1^\circ$ If the basis $B=\{\alg\varepsilon_k\}_{k=1}^n$ satisfies Euler-type
relations such as II.2.(10)--(13), then $e^{\alg d t}$ and the functions $f_k$
can be computed efficiently.

$2^\circ$ Once a fundamental system $[f]=[f_1,\ldots,f_n]^T$ is known, the
representation mapping
\[
\Psi:V\longrightarrow A_n
\]
is uniquely specified by the values $\Psi(f_k)$.  For a fixed basis
$B=\{\alg\varepsilon_k\}_{k=1}^n$, there are $n!$ representations that are
notationally distinct but equivalent, corresponding to the bijections between
$\{f_1,\ldots,f_n\}$ and $\{\alg\varepsilon_1,\ldots,\alg\varepsilon_n\}$.

\renewcommand{\eqprefix}{II.3}
\section*{II.3. Symbolic Representation of Discrete Signals on $\mathbb R$-Algebras}
Let $S_d$ be the space of real discrete signals:
\eqn{1}{S_d=\{f:\mathbb Z\to\mathbb R\mid f=f(k)\}.}
The right-shift operator is
\eqn{2}{\dot f\stackrel{\mathrm{def}}{=}f(k+1).}

\noindent\textbf{Theorem II.3.1.} On the algebra $A_n$, with operator
$\alg d\in A_n$, one can symbolically represent the solution space of the
difference equation
\eqn{3}{f(k+n)+a_1f(k+n-1)+\cdots+a_{n-1}f(k+1)+a_nf(k)=0,}
whose characteristic polynomial coincides with that of $\alg d$.

\noindent\textbf{Theorem II.3.2.} The solution space of the linear difference equation
\eqn{4}{f(k+n)+a_1f(k+n-1)+\cdots+a_{n-1}f(k+1)+a_nf(k)=0}
is symbolically represented on the polynomial algebra generated by the
characteristic polynomial of the equation, with the element generating the
cyclic basis as representation operator.

To obtain $\Psi$, a fundamental system must be determined for
\eqn{5}{[f(k+1)]=[\alg d]\,[f(k)].}
The general solution is
\eqn{6}{[f(k)]=[\alg d]^k[c],}
where $[c]$ is a column of $n$ real constants.

\noindent\textbf{Theorem II.3.3.} If $\alg d\in A_n$ has a minimal polynomial
of degree $n$ and
\eqn{7}{\alg d^{\,k}=f_1(k)\alg\varepsilon_1+\cdots+
f_n(k)\alg\varepsilon_n,\qquad k\in\mathbb Z,}
then the functions $f_i$, $i=1,\ldots,n$, form a fundamental system for
difference equation (3).

\noindent\textit{Proof.} Equation (7) gives
\eqn{8}{\alg d^{\,k+1}=\sum_{p=1}^nf_p(k+1)\alg\varepsilon_p.}
On the other hand,
\[
\alg d^{\,k+1}=\alg d\sum_{j=1}^nf_j(k)\alg\varepsilon_j
=\sum_{j=1}^n f_j(k)\sum_{i=1}^nd_i
\sum_{p=1}^n\gamma_{ijp}\alg\varepsilon_p.
\]
With $\alpha_{jp}=\sum_{i=1}^nd_i\gamma_{ijp}$, one obtains
\eqn{9}{\alg d^{\,k+1}=
\sum_{p=1}^n\left(\sum_{j=1}^nf_j(k)\alpha_{jp}\right)
\alg\varepsilon_p.}

Comparing (8) and (9), the functions $f_p$ satisfy the difference system
\eqn{10}{f_p(k+1)=\sum_{j=1}^n\alpha_{jp}f_j(k),
\qquad p=1,\ldots,n.}
Since $[\alg d]=(\alpha_{ij})$, the system can be written in matrix form as
\eqn{11}{[f(k+1)]=[\alg d]^T[f(k)].}
The functions $f_p$ satisfy equation (3), because $[\alg d]$ and
$[\alg d]^T$ have the same characteristic polynomial, which is identical to
that of $\alg d$.  Linear independence is proved as in Theorem II.2.3, using
the discrete counterpart of the Wronskian [73] and the fact that the minimal
polynomial has degree $n$.

\noindent\textbf{Remarks.}

$1^\circ$ If the operator admits the exponential representation
\eqn{12}{\alg d=e^{\alg\omega},\qquad \alg\omega\in A_n,}
its powers and the fundamental functions can be determined from the
Euler-type relations of the algebra.

\renewcommand{\eqprefix}{II.3}

If the operator $\alg d$ can be written as $\alg d=e^{\alg\omega}$, then,
under the stated conditions, the discrete signals in II.3.(7) are determined
using the Euler-type relations II.2.(10)--(13).  Such a representation is
always possible if $\alg d$ is not a zero divisor in $A_n$ [140].

\noindent $2^\circ$ After a fundamental system
$[f]=[f_1,\ldots,f_n]^T$ has been determined, the representation
\eqn{13}{\Psi:S_d\longrightarrow A_n}
is specified by the values $\Psi(f_p)$, that is, by a bijection
\eqn{14}{\Psi:\{f_1,\ldots,f_n\}\longrightarrow
\{\alg\varepsilon_1,\ldots,\alg\varepsilon_n\}.}
For fixed bases, there are $n!$ such representations, notationally distinct
but equivalent.

Discrete signals can be obtained by uniform sampling [23], [151] of continuous
signals.  If
\eqn{15}{f:\mathbb R\to\mathbb R,\qquad f=f(t),}
then the sampled signal is the restriction
\eqn{16}{f_d:\mathbb Z\to\mathbb R,\qquad f_d(k)=f(k).}

\noindent\textbf{Theorem II.3.4.} If the continuous signal $f$ is represented
on $A_n$, with operator $\alg d$, by $\alg w\in A_n$, then the discrete signal
$f_d$ obtained by uniform sampling is represented on the same algebra by the
same element $\alg w$, with $e^{\alg d}$ as operator.

\noindent\textit{Proof.} Let $B=\{\alg\varepsilon_k\}_{k=1}^n$ be a basis and
let $V$ be the solution space of
\eqn{17}{f^{(n)}+a_1f^{(n-1)}+\cdots+a_{n-1}\dot f+a_nf=0,}
whose characteristic polynomial is that of $\alg d$.

The representation
\eqn{18}{\Psi:V\longrightarrow A_n}
is specified by $\Psi(f_k)=\alg\varepsilon_k$, where the functions $f_k$ are
the coefficients of the expansion in Theorem II.2.3:
\eqn{19}{e^{\alg d t}=\sum_{i=1}^nf_i(t)\alg\varepsilon_i.}
If the signal $f$ is represented by
\eqn{20}{\alg w=\sum_{i=1}^nw_i\alg\varepsilon_i,}
then
\eqn{21}{f(t)=\sum_{i=1}^nw_if_i(t),\qquad t\in\mathbb R.}
Sampling gives
\eqn{22}{f_d(k)=\sum_{i=1}^nw_if_i(k),\qquad k\in\mathbb Z.}
From (19), for $t=k\in\mathbb Z$,
\eqn{23}{(e^{\alg d})^k=\sum_{i=1}^nf_i(k)\alg\varepsilon_i.}
Theorem II.3.3 allows the representation
\eqn{24}{\Psi_d:V_d\longrightarrow A_n,}
to be defined with operator $e^{\alg d}$ and
$\Psi_d[f_i(k)]=\alg\varepsilon_i$.

By linearity and equation (22),
\eqn{25}{\Psi_d[f_d(k)]
=\sum_{i=1}^nw_i\Psi_d[f_i(k)]
=\sum_{i=1}^nw_i\alg\varepsilon_i=\alg w.}
Therefore, the discrete signal is represented by the same element of the
algebra, with $e^{\alg d}$ as shift operator.

\noindent\textbf{Remark.} In (24), $V_d$ is the solution space of
\eqn{26}{f(k+n)+b_1f(k+n-1)+\cdots+b_{n-1}f(k+1)+b_nf(k)=0,}
where
\[
q(\lambda)=\lambda^n+b_1\lambda^{n-1}+\cdots+b_n
=\det\!\left(\lambda I_n-e^{[\alg d]}\right).
\]

\noindent\textbf{Corollary II.3.5.} The class of discrete signals obtained by
uniformly sampling the solutions of
\eqn{27}{f^{(n)}+a_1f^{(n-1)}+\cdots+a_{n-1}\dot f+a_nf=0}
is the solution space of the difference equation
\eqn{28}{f(k+n)+b_1f(k+n-1)+\cdots+b_{n-1}f(k+1)+b_nf(k)=0,}
whose characteristic polynomial is
$q(\lambda)=\det(\lambda I_n-e^D)$, where
\eqn{29}{D=\begin{bmatrix}
0&1&0&\cdots&0\\
0&0&1&\cdots&0\\
\vdots&\vdots&\vdots&\ddots&\vdots\\
0&0&0&\cdots&1\\
-a_n&-a_{n-1}&-a_{n-2}&\cdots&-a_1
\end{bmatrix}.}

\noindent\textit{Proof.} By Theorem II.2.2, the solutions of equation (27) are
represented on the polynomial algebra generated by
\eqn{30}{p(\lambda)=\lambda^n+a_1\lambda^{n-1}
+\cdots+a_{n-1}\lambda+a_n.}
The algebra has the cyclic basis
$B=\{\alg v^0,\alg v,\ldots,\alg v^{\,n-1}\}$, with
$p(\alg v)=\alg0$, and its representation operator is $\alg v$.  Its canonical
matrix is
\eqn{31}{[\alg v]=\begin{bmatrix}
0&1&0&\cdots&0\\
0&0&1&\cdots&0\\
\vdots&\vdots&\vdots&\ddots&\vdots\\
0&0&0&\cdots&1\\
-a_n&-a_{n-1}&-a_{n-2}&\cdots&-a_1
\end{bmatrix}.}
By Theorem II.3.4, the discrete signals obtained by sampling are represented on
the same algebra, with operator $e^{\alg v}$.  By Theorem II.3.1, they satisfy
a linear homogeneous difference equation with constant coefficients, whose
characteristic polynomial is
\eqn{32}{q(\lambda)=
\det\!\left(\lambda I_n-[e^{\alg v}]\right)
=\det\!\left(\lambda I_n-e^{[\alg v]}\right).}

\chaptertoc{CHAPTER III}{Particular \texorpdfstring{$\mathbb R$}{R}-Algebras: Structural Properties}
\phantomsection\addcontentsline{toc}{section}{III.1. Algebra $\bigotimes C^{n+1}$}
\renewcommand{\eqprefix}{III.1}

\thispagestyle{empty}
\begin{center}
{\Large\bfseries CHAPTER III}\par
\vspace{4mm}
{\large\bfseries Particular $R$-Algebras: Structural Properties}\par
\vspace{5mm}
{\large\bfseries III.1. The Algebra $\bigotimes C^{n+1}$}\par
\vspace{4mm}
{\large\bfseries III.1.1. Preliminaries: Bases in $\bigotimes C^{n+1}$}\addcontentsline{toc}{subsection}{III.1.1. Preliminaries: Bases}
\end{center}

Let $C_k$, $k=\overline{0,n}$, be $n+1$ real algebras of order $2$,
isomorphic to the complex field $C$. If $\{\onehat_k,\ihat_k\}$ is a basis
of $C_k$, its multiplication table is
\eqn{1}{
\begin{aligned}
\onehat_k\ihat_k&=\ihat_k\onehat_k=\ihat_k,\\
\onehat_k^{2}&=-\ihat_k^{2}=\onehat_k.
\end{aligned}}

The algebras $C_k$, $k\in\{0,1,\ldots,n\}$, are associative, commutative,
and unital. Consider the direct-product algebra of these $(n+1)$ algebras,
denoted by
\eqn{2}{
\bigotimes C^{n+1}=C_0\otimes C_1\otimes\cdots\otimes C_n .}
Algebra (2) has order $2^{n+1}$ and is commutative, associative, and unital,
with identity element
\[
\onehat\overset{\mathrm{not}}{=}
\onehat_0\onehat_1\cdots\onehat_n.
\]

The natural basis elements of $\bigotimes C^{n+1}$ are obtained as formal
products of basis elements of the algebras $C_k$, $k=\overline{0,n}$
[2, 54]. The multiplication table of the direct-product algebra follows from
the tables (1) of its component algebras and from associativity and
commutativity.

Using the notation
\eqn{3}{
\ihat_k^{p_k}=\begin{cases}
\onehat_k, & p_k=0,\\
\ihat_k, & p_k=1,
\end{cases}}
the basis elements of algebra (2) can be written as
\eqn{4}{
\ehat_p=\ihat_0^{p_0}\ihat_1^{p_1}\cdots\ihat_n^{p_n},
\qquad p=\overline{0,2^{n+1}-1}.}

Here $p=(p_0p_1\ldots p_n)_2$ is the natural binary representation of $p$
with $(n+1)$ bits.

In (4), the identity of the algebra is
$\ehat_0=\onehat_0\onehat_1\cdots\onehat_n\overset{\mathrm{not}}{=}\onehat$.
In what follows, the identities of the algebras $C_k$ will be omitted from
(4). Thus, we write $\ehat_{2^n}=\ihat_1$ instead of
$\ehat_{2^n}=\onehat_0\ihat_1\onehat_2\cdots\onehat_n$, and
$\ehat_1=\ihat_n$ instead of
$\ehat_1=\onehat_0\onehat_1\cdots\onehat_{n-1}\ihat_n$.
This convention simplifies the notation.

The multiplication table of algebra (2) follows from (4). If we set
\eqn{5}{
\ehat_q=\ihat_0^{q_0}\ihat_1^{q_1}\cdots\ihat_n^{q_n},
\qquad q\in\{0,1,\ldots,2^{n+1}-1\},}
where $q=(q_0q_1\ldots q_n)_2$, then (4) and (5) give
\eqn{6}{
\ehat_p\ehat_q=
\ihat_0^{p_0+q_0}\ihat_1^{p_1+q_1}\cdots
\ihat_n^{p_n+q_n},\quad 
p,q\in\{0,1,\ldots,2^{n+1}-1\}.}
Since $p_k,q_k\in\{0,1\}$, $k=\overline{0,n}$, equation (1) yields
\eqn{7}{
\ehat_p\ehat_q=(-1)^{\langle p,q\rangle}\ehat_{p\xor q},
\quad p,q\in\{0,1,\ldots,2^{n+1}-1\},}
where $\langle p,q\rangle=p_0q_0+p_1q_1+\cdots+p_nq_n$, and $p\xor q$
denotes the bitwise modulo-$2$ sum of $p$ and $q$. Equation (7) shows that
algebra (2) is monomial [54] and that
\eqn{8}{
\ehat_p^{2}=(-1)^{\langle p,p\rangle}\ehat_0,
\quad p\in\{0,1,\ldots,2^{n+1}-1\}.}
Since $\ehat_0=\onehat$, the basis elements fall into two classes:
\eqn{9}{
\ehat_p^{2}\in\{\onehat,-\onehat\},
\quad p\in\{0,1,\ldots,2^{n+1}-1\}.}
The elements whose square is the identity, respectively minus the identity,
have indices containing an even, respectively odd, number of bits equal to
one. Each class contains $2^n$ elements.

For reindexing, we use Gray coding with bit reversal [1]:
\eqn{10}{
g:\{0,1,\ldots,2^{n+1}-1\}\longrightarrow
\{0,1,\ldots,2^{n+1}-1\},}
\eqn{11}{g(k)=(\bar k_0\bar k_1\ldots\bar k_n)_2,}
where
\eqn{12}{
\bar k_p=\begin{cases}
k_n, & p=0,\\
k_{n-p}\xor k_{n-p+1}, & p\in\{1,2,\ldots,n\}.
\end{cases}}

Here $k=(k_0k_1\ldots k_n)_2$ is the natural binary representation of $k$
with $n+1$ bits. The mapping (10) is bijective and has the property that the
first $2^n$ numbers, $\{0,1,\ldots,2^n-1\}$, contain an even number of bits
equal to one, whereas the numbers $\{2^n,\ldots,2^{n+1}-1\}$ contain an odd
number of bits equal to one.

Now set
\eqn{13}{
\widehat v_k=\ehat_{g(k)}=
\ihat_0^{\bar k_0}\ihat_1^{\bar k_1}\cdots\ihat_n^{\bar k_n},
\qquad k\in\{0,1,\ldots,2^{n+1}-1\}.}
The set $\{\widehat v_k\}_{k=0}^{2^{n+1}-1}$ is a basis of
$\bigotimes C^{n+1}$ with the property
\eqn{14}{
\widehat v_k^{2}=\begin{cases}
\onehat, & k\in\{0,1,\ldots,2^n-1\},\\
-\onehat, & k\in\{2^n,2^n+1,\ldots,2^{n+1}-1\}.
\end{cases}}
In what follows, write
\[
\widehat u_k\overset{\mathrm{not}}{=}\widehat v_k,
\qquad
\widehat j_k\overset{\mathrm{not}}{=}\widehat v_{2^n+k},
\qquad k\in\{0,1,\ldots,2^n-1\}.
\]
In accordance with (14),
\[
\widehat u_k^{2}=-\widehat j_k^{2}=\onehat,
\qquad k\in\{0,1,\ldots,2^n-1\}.
\]

\begin{center}
{\large\bfseries III.1.2. Structural Properties}\addcontentsline{toc}{subsection}{III.1.2. Structural Properties}
\end{center}

We now present several properties concerning the structure of the algebra
$\bigotimes C^{n+1}$.

\medskip
\noindent\textbf{Theorem III.1.2.1.} The elements
\eqn{15}{
\widehat e_p=\frac{1}{2^n}
\prod_{k=1}^{n}\left[\onehat-(-1)^{p_k}\ihat_0\ihat_k\right],
\qquad p=(p_1p_2\ldots p_n)_2\in\{0,1,\ldots,2^n-1\},}
constitute $2^n$ orthogonal idempotents in $\bigotimes C^{n+1}$; that is,
they satisfy
\eqn{16}{
\widehat e_p\widehat e_q=\begin{cases}
\widehat e_p, & p=q,\\
\widehat 0, & p\ne q,
\end{cases}
\qquad p,q\in\{0,1,\ldots,2^n-1\}.}

\noindent\textit{Proof.}

Let $p=(p_1,p_2,\ldots,p_n)_2$ and $q=(q_1,q_2,\ldots,q_n)_2$ be two
numbers in $\{0,1,\ldots,2^n-1\}$. Using (15), compute the products
$\widehat e_p\widehat e_q$:
\eqn{17}{
\begin{aligned}
\widehat e_p\widehat e_q
&=\frac{1}{2^{n+1}}\prod_{k=1}^{n}
\bigl[\onehat-(-1)^{p_k}\ihat_0\ihat_k\bigr]
\bigl[\onehat-(-1)^{q_k}\ihat_0\ihat_k\bigr]\\
&=\frac{1}{2^{n+1}}\prod_{k=1}^{n}
\bigl[\onehat-(-1)^{p_k}\ihat_0\ihat_k
-(-1)^{q_k}\ihat_0\ihat_k+(-1)^{p_k+q_k}\onehat\bigr].
\end{aligned}}
If $p=q$, then $p_k=q_k$ for every $k\in\{1,2,\ldots,n\}$, and (17) gives
\[
\widehat e_p\widehat e_q
=\frac{1}{2^{n+1}}\prod_{k=1}^{n}
\bigl[2\onehat-2(-1)^{p_k}\ihat_0\ihat_k\bigr]
=\widehat e_p.
\]
If $p\ne q$, there exists a $k\in\{1,2,\ldots,n\}$ such that
$p_k\ne q_k$. It follows from (17) that
\[
\widehat e_p\widehat e_q=\widehat 0,
\qquad p,q\in\{0,1,\ldots,2^n-1\},\quad p\ne q.
\]
Relations (16) also show that the elements $\widehat e_p$,
$p=\overline{0,2^n-1}$, are zero divisors in $\bigotimes C^{n+1}$.

\medskip
\noindent\textbf{Remark.}
The real algebra $\bigotimes C^{n+1}$ contains a subalgebra isomorphic to
the complex field $C$, namely
\[
C=\{\lambda\onehat+\mu\ihat_0\mid \lambda,\mu\in\mathbb R\}.
\]
Accordingly, $\bigotimes C^{n+1}$ may be regarded as a $C$-algebra, denoted
by $\bigotimes C^n$. A basis of $\bigotimes C^n$ is obtained by retaining
from (4) the $2^n$ elements corresponding to $p_0=0$. In this $C$-algebra,
the elements (15) are $2^n$ orthogonal idempotents. A standard property [54]
gives
\eqn{18}{
\widehat e_0+\widehat e_1+\cdots+\widehat e_{2^n-1}=\onehat.}
This raises the question whether the set
$B=\{\widehat e_p,\ihat_0\widehat e_p\}_{p=0}^{2^n-1}$ is a basis of
$\bigotimes C^{n+1}$. Equation (18) gives
\eqn{19}{
\ihat_k=\sum_{p=0}^{2^n-1}\ihat_k\widehat e_p,
\qquad k=\overline{0,n}.}

Computing $\ihat_k\widehat e_p$ from (15), one obtains
\eqn{20}{\ihat_k\widehat e_p=(-1)^{p_k}\ihat_0\widehat e_p,
\quad k\in\{0,1,\ldots,n\},\ p\in\{0,1,\ldots,2^n-1\}.}
Equations (19) and (20) yield
\eqn{21}{\ihat_k=\sum_{p=0}^{2^n-1}(-1)^{p_k}\ihat_0\widehat e_p,
\qquad k\in\{0,1,\ldots,n\}.}
Using (13) and (21), it follows that
\eqn{22}{\widehat v_k=
\sum_{p=0}^{2^n-1}(-1)^{p_0k_0+p_1k_1+\cdots+p_nk_n}
\ihat_0^{k_0+k_1+\cdots+k_n}\widehat e_p,
\quad k\in\{0,1,\ldots,2^{n+1}-1\}.}
Since $\ihat_0^2=-\onehat$, we have
\eqn{23}{\ihat_0^{k_0+\cdots+k_n}=\begin{cases}
\sigma(k)\onehat,&k\in\{0,1,\ldots,2^n-1\},\\
\sigma(k)\ihat_0,&k\in\{2^n,\ldots,2^{n+1}-1\},
\end{cases}}
where
\eqn{24}{\sigma(k)=(-1)^{\left\lfloor (k_0+k_1+\cdots+k_n)/2\right\rfloor},
\qquad k=(k_0k_1\ldots k_n)_{\mathrm{Gray}}.}
Therefore,
\eqn{25}{\widehat u_k=\widehat v_k=
\sum_{p=0}^{2^n-1}(-1)^{p_0k_0+\cdots+p_nk_n}\sigma(k)\widehat e_p,
\quad k\in\{0,1,\ldots,2^n-1\},}
while for the second half of the basis,
\eqn{26}{\widehat j_{k-2^n}=\widehat v_k=
\sum_{p=0}^{2^n-1}(-1)^{p_0k_0+\cdots+p_nk_n}
\sigma(k)\ihat_0\widehat e_p,
\quad k\in\{2^n,\ldots,2^{n+1}-1\}.}
Relations (25) and (26) admit a convenient matrix form.

Introduce the column matrices
\eqn{27}{\boldsymbol v=
[\widehat v_0,\widehat v_1,\ldots,\widehat v_{2^{n+1}-1}]^T,}
\eqn{28}{\boldsymbol e=
[\widehat e_0,\ldots,\widehat e_{2^n-1},
\ihat_0\widehat e_0,\ldots,\ihat_0\widehat e_{2^n-1}]^T.}
Consider the square matrix of order $2^n$
\eqn{29}{H_w=(h_{p,k})_{p,k=0}^{2^n-1},\qquad
h_{p,k}=(-1)^{p_0\bar k_0+p_1\bar k_1+\cdots+p_n\bar k_n},}
which is the Walsh matrix [1] of order $2^n$ and satisfies
\eqn{30}{H_w^{-1}=\frac1{2^n}H_w^T=\frac1{2^n}H_w.}
Introduce the matrices
\eqn{31}{\sigma=\operatorname{diag}\bigl(\sigma(0),\sigma(1),\ldots,
\sigma(2^n-1)\bigr),}
and the antidiagonal matrix
\eqn{32}{\sigma'=
\begin{bmatrix}
0&\cdots&0&\sigma(2^n)\\
0&\cdots&\sigma(2^n+1)&0\\
\vdots&\ddots&\vdots&\vdots\\
\sigma(2^{n+1}-1)&0&\cdots&0
\end{bmatrix}.}
Relations (25) and (26) can be combined by means of
\eqn{33}{S=\begin{bmatrix}\sigma H_w&0\\0&\sigma'H_w\end{bmatrix},}
where $0$ is the zero matrix of order $2^n$, in the form
\eqn{34}{\boldsymbol v=S\boldsymbol e.}

Multiplication of the Walsh matrix by $\sigma$ and $\sigma'$ permutes its
rows and may change their signs, since $\sigma(k)\in\{-1,1\}$. Equations
(30)--(33) show that $S$ is invertible and
\eqn{35}{S^{-1}=\frac1{2^n}S^T.}
From (34),
\eqn{36}{\boldsymbol e=\frac1{2^n}S^T\boldsymbol v,}
which proves that $B=\{\widehat e_p,\ihat_0\widehat e_p\}_{p=0}^{2^n-1}$
is a basis of $\bigotimes C^{n+1}$. An arbitrary element has the expansion
\eqn{37}{\widehat w=\sum_{k=0}^{2^{n+1}-1}w_k\widehat v_k,
\qquad w_k\in\mathbb R,}
and the unique expansion in the basis $B$,
\eqn{38}{\widehat w=\sum_{k=0}^{2^n-1}
(a_k\widehat e_k+b_k\ihat_0\widehat e_k),}
or
\eqn{39}{\widehat w=\sum_{k=0}^{2^n-1}
(a_k\onehat+\ihat_0b_k)\widehat e_k.}
The coefficients are determined by
\eqn{40}{[a_0,a_1,\ldots,a_{2^n-1}]^T
=(\sigma H_w)^T[w_0,w_1,\ldots,w_{2^n-1}]^T,}
\eqn{41}{[b_0,b_1,\ldots,b_{2^n-1}]^T
=(\sigma'H_w)^T[w_{2^n},w_{2^n+1},\ldots,w_{2^{n+1}-1}]^T.}
Conversely,
\eqn{42}{[w_0,\ldots,w_{2^n-1}]^T
=\frac1{2^n}(\sigma H_w)^T[a_0,\ldots,a_{2^n-1}]^T,}
\eqn{43}{[w_{2^n},\ldots,w_{2^{n+1}-1}]^T
=\frac1{2^n}(\sigma'H_w)^T[b_0,\ldots,b_{2^n-1}]^T.}

\noindent\textbf{Remark.} The passage from the ``direct-product'' basis (13)
to the orthogonal-idempotent basis (36) is equivalent to applying the discrete
Walsh transform twice [1, 76], in accordance with (40), (41). The procedure
is directly implementable on a computer by means of fast Walsh-transform
algorithms [1, 76, 151].

Equation (39), together with the orthogonality of
$\{\widehat e_p\}_{p=0}^{2^n-1}$, gives the following theorem.

\medskip\noindent\textbf{Theorem III.1.2.2.}
The algebra $\bigotimes C^{n+1}$ is the direct sum of $2^n$ algebras
isomorphic to the complex field; it is polynomial and semisimple.

\smallskip\noindent\textit{Proof.}
The direct-sum property follows from (39) and (16). Since the algebra is the
direct sum of $2^n$ semisimple polynomial algebras with generating polynomial
$P(X)=X^2+1$, the assertion follows [54].

\begin{center}
{\large\bfseries III.1.3. Exponential Form of the Elements of
$\bigotimes C^{n+1}$: Zero Divisors}\addcontentsline{toc}{subsection}{III.1.3. Exponential Form}
\end{center}
We say that $\widehat w$ admits an exponential form if there exists
$\widehat u\in\bigotimes C^{n+1}$ such that
\eqn{45}{\widehat w=e^{\widehat u},}
where
\eqn{46}{e^{\widehat u}=\sum_{k=0}^{\infty}\frac{\widehat u^{k}}{k!}.}
In the basis $B$, write
\eqn{47}{\widehat w=\sum_{k=0}^{2^n-1}
(a_k\onehat+\ihat_0b_k)\widehat e_k,
\qquad a_k,b_k\in\mathbb R,}
and set
\eqn{48}{a_k\onehat+\ihat_0b_k
=r_k\onehat\,e^{\ihat_0\varphi_k},\qquad
r_k\ge0,\quad \varphi_k\in[0,2\pi).}

The parameters in (48) satisfy
\eqn{49}{r_k=\sqrt{a_k^2+b_k^2},}
\eqn{50}{r_k\cos\varphi_k=a_k,\qquad r_k\sin\varphi_k=b_k,
\qquad k=\overline{0,2^n-1}.}
If $r_k>0$ for every $k$, then
\eqn{51}{\widehat w=\sum_{k=0}^{2^n-1}
r_ke^{\ihat_0\varphi_k}\widehat e_k
=\sum_{k=0}^{2^n-1}e^{\ln r_k\onehat+\ihat_0\varphi_k}\widehat e_k.}
By orthogonality of the idempotents, (51) becomes
\eqn{52}{\widehat w=\exp\!\left[
\sum_{k=0}^{2^n-1}(\ln r_k\onehat+\ihat_0\varphi_k)\widehat e_k
\right],\qquad r_k>0.}
This is the exponential form (45), with
\eqn{53}{\widehat u=\sum_{k=0}^{2^n-1}\left[
\ln\sqrt{a_k^2+b_k^2}\,\widehat e_k
+\operatorname{atan2}(b_k,a_k)\ihat_0\widehat e_k\right],
\qquad a_k^2+b_k^2>0.}
The angle in (53) is the unique angle in $[0,2\pi)$ satisfying (50).
The zero divisors in $\bigotimes C^{n+1}$ are characterized by
\eqn{54}{\widehat w\text{ is a zero divisor}
\iff\prod_{k=0}^{2^n-1}r_k
=\prod_{k=0}^{2^n-1}\sqrt{a_k^2+b_k^2}=0.}

\fancyhead[R]{Chapter III -- Particular $R$-Algebras: Structural Properties}

Therefore, element (47) is a zero divisor in $\bigotimes C^{n+1}$ if there
exists a $k\in\{0,1,\ldots,2^n-1\}$ such that
\eqn{55}{a_k=b_k=0.}
An element $\what$ written as in (47) and satisfying (55) will be called a
\emph{zero divisor of rank $k$}. Equations (47) and (16) give
\eqn{56}{\what\,\epshat_k=\widehat 0.}
Elements that are not zero divisors are invertible. For (47), the inverse is
\eqn{57}{\what^{-1}=\sum_{k=0}^{2^n-1}
\frac{a_k\onehat-\ihat_0b_k}{a_k^2+b_k^2}\,\epshat_k.}
The elements $\xi_k=a_k\onehat+\ihat_0b_k$ satisfy the equations
\eqn{58}{p_k(\xi_k)=\widehat0,\qquad
p_k(X)=X^2-2a_kX+a_k^2+b_k^2.}
By the properties (16) of the orthogonal idempotents, element (47) satisfies
\eqn{59}{p(\what)=\widehat0,\qquad
p(X)=\prod_{k=0}^{2^n-1}(X^2-2a_kX+a_k^2+b_k^2).}
If $\what$ is not a zero divisor, polynomial (59) is its minimal polynomial
in $\bigotimes C^{n+1}$. If $\what$ is a zero divisor of rank $k$, the
minimal polynomial has degree $2^{n+1}-2$ and is
\eqn{60}{q(\what)=\widehat0,\qquad
q(X)=\prod_{\substack{p=0\\p\ne k}}^{2^n-1}
(X^2-2a_pX+a_p^2+b_p^2).}
For ranks $k$ and $p$, the degree decreases to $2^{n+1}-4$, and so forth.
There are no nilpotent elements in $\bigotimes C^{n+1}$.

If an arbitrary element is written in basis (13), then
\eqn{61}{\what=\sum_{p=0}^{2^{n+1}-1}w_p\vhat_p,
\qquad w_p\in\mathbb R.}

The components of element (61) in the basis
$\{\epshat_p,\ihat_0\epshat_p\}_{p=0}^{2^n-1}$ are computed from (40),
(41). Componentwise,
\eqn{62}{a_k=\sum_{p=0}^{2^n-1}(-1)^{\langle p,k\rangle}
\sigma(p)w_p,\qquad k\in\{0,1,\ldots,2^n-1\},}
\eqn{63}{b_k=\sum_{p=0}^{2^n-1}(-1)^{\langle p,k\rangle}
\sigma(2^n+p)w_{2^n+p}.}
We have set
\eqn{64}{\langle p,k\rangle=p_1k_1+\cdots+p_nk_n,}
where $p=(p_0p_1\ldots p_n)_2$ is written in natural binary code, whereas
$k=(k_0k_1\ldots k_n)_{\mathrm{Gray}}$ is written in Gray code with bit
reversal. The term $p_0k_0$ does not occur because $p_0=0$ for
$p\in\{0,\ldots,2^n-1\}$.

The zero divisors of rank $k$ among the elements (61) are characterized by
\eqn{65}{\sum_{p=0}^{2^n-1}(-1)^{\langle p,k\rangle}\sigma(p)w_p
=\sum_{p=0}^{2^n-1}(-1)^{\langle p,k\rangle}
\sigma(2^n+p)w_{2^n+p}=0.}
If $\what$ is not a zero divisor, it admits the exponential form
\eqn{66}{\what=e^{\widehat u},}
where
\eqn{67}{\widehat u=\sum_{k=0}^{2^{n+1}-1}u_k\vhat_k,
\qquad u_k\in\mathbb R.}
The coefficients are computed from
\eqn{68}{[u_0,\ldots,u_{2^n-1}]^T
=\frac1{2^n}\sigma H_w[\ln r_0,\ldots,\ln r_{2^n-1}]^T,}
\eqn{69}{[u_{2^n},\ldots,u_{2^{n+1}-1}]^T
=\frac1{2^n}\sigma'H_w[\varphi_0,\ldots,\varphi_{2^n-1}]^T,}
with
\eqn{70}{r_k=\Biggl(\left[
\sum_{p=0}^{2^n-1}(-1)^{\langle p,k\rangle}\sigma(p)w_p
\right]^2+\left[
\sum_{p=0}^{2^n-1}(-1)^{\langle p,k\rangle}
\sigma(2^n+p)w_{2^n+p}\right]^2\Biggr)^{1/2}.}
The corresponding phase is
\eqn{71}{\varphi_k=\operatorname{atan2}\!\left(
\sum_{p=0}^{2^n-1}(-1)^{\langle p,k\rangle}\sigma(2^n+p)w_{2^n+p},
\sum_{p=0}^{2^n-1}(-1)^{\langle p,k\rangle}\sigma(p)w_p
\right).}

Relation (66) may also be written as
\eqn{72}{\what=\widehat\rho\,e^{\widehat\theta},}
with
\eqn{73}{\widehat\rho=\sum_{k=0}^{2^n-1}\rho_k\widehat u_k,
\qquad \rho_k\in\mathbb R,}
\eqn{74}{\widehat\theta=\sum_{k=0}^{2^n-1}\theta_k\widehat j_k,
\qquad \theta_k\in[0,2\pi),}
where
\eqn{75}{[\rho_0,\rho_1,\ldots,\rho_{2^n-1}]^T
=\frac1{2^n}\sigma H_w[r_0,r_1,\ldots,r_{2^n-1}]^T,}
\eqn{76}{[\theta_0,\theta_1,\ldots,\theta_{2^n-1}]^T
=\frac1{2^n}\sigma'H_w[\varphi_0,\varphi_1,\ldots,
\varphi_{2^n-1}]^T.}
The numbers $r_k$ and $\varphi_k$ are given by (70), (71). Representation
(72) will be called the \emph{semiexponential form} of $\what$. Zero divisors,
which do not admit the exponential form (66), can be expressed in
semiexponential form.

If, in (72),
\eqn{77}{\widehat\rho=\rho\onehat,\qquad\rho\in\mathbb R,}
then the element
\eqn{78}{\what=\rho e^{\widehat\theta}}
is called an \emph{element of special type}.

\medskip\noindent\textbf{Theorem III.1.3.1.}
The product of two elements of special type is an element of special type.
Every element of $\bigotimes C^{n+1}$ is the sum of at most $2^n$ elements
of special type.

\smallskip\noindent\textit{Proof.}
Let two elements of special type be
\eqn{79}{\widehat w_1=\rho_1\exp\!\left(
\sum_{k=0}^{2^n-1}\alpha_k\widehat j_k\right),}

\eqn{80}{\widehat w_2=\rho_2\exp\!\left(
\sum_{k=0}^{2^n-1}\beta_k\widehat j_k\right),
\quad \rho_1,\rho_2>0,\quad \alpha_k,\beta_k\in[0,2\pi).}
Their product is
\eqn{81}{\widehat w_1\widehat w_2=\rho\exp\!\left(
\sum_{k=0}^{2^n-1}\gamma_k\widehat j_k\right),}
where
\eqn{82}{\rho=\rho_1\rho_2>0,
\qquad\gamma_k=\alpha_k+\beta_k\pmod{2\pi}.}
By definition, (81) is an element of special type.

Now let an arbitrary element be written in semiexponential form as
\eqn{83}{\what=\widehat\rho\,e^{\widehat\theta},}
where
\eqn{84}{\widehat\rho=\sum_{k=0}^{2^n-1}\rho_k\widehat u_k,
\qquad\rho_k\in\mathbb R,}
\eqn{85}{\widehat\theta=\sum_{k=0}^{2^n-1}\theta_k\widehat j_k,
\qquad\theta_k\in[0,2\pi).}
We express the basis elements $\widehat u_k$ in exponential form. For
$k=(k_0k_1\ldots k_n)_{\mathrm{Gray}}$,
\eqn{86}{\widehat u_k=\ihat_0^{k_0}\ihat_1^{k_1}
\cdots\ihat_n^{k_n},}
\eqn{87}{k=(k_0k_1\ldots k_n)_{\mathrm{Gray}}.}
Since
\eqn{88}{\ihat_p=e^{\ihat_p\pi/2},
\qquad \ihat_p^2=-\onehat,\quad p=\overline{0,n},}
it follows that
\eqn{89}{\widehat u_k=\exp\!\left(
\sum_{p=0}^{n}\ihat_p k_p\frac{\pi}{2}\right).}
Substitution of (89) into (83) gives the required decomposition as a sum of
elements of special type.

More precisely,
\eqn{90}{\what=\sum_{k=0}^{2^n-1}\rho_k
\exp\!\left(\sum_{p=0}^{n}\ihat_p k_p\frac{\pi}{2}\right)
\exp\widehat\theta.}
Since $\ihat_p\in\{\widehat j_0,\widehat j_1,\ldots,
\widehat j_{2^n-1}\}$, the properties of the exponential give
\eqn{91}{\what=\sum_{k=0}^{2^n-1}\rho_k e^{\widehat\alpha_k},}
where
\eqn{92}{\widehat\alpha_k=\widehat\theta+
\sum_{p=0}^{n}\ihat_p k_p\frac{\pi}{2}.}
The terms on the right-hand side of (91) are elements of special type, which
proves the second part of the theorem.

\medskip\noindent\textbf{Remarks.}
In the semiexponential representation (72), an element is of special type if
and only if
\eqn{93}{\rho_1=\rho_2=\cdots=\rho_{2^n-1}=0,
\qquad \rho_0\ne0.}
For an element written in the basis
$B=\{\epshat_p,\ihat_0\epshat_p\}_{p=0}^{2^n-1}$,
\eqn{94}{\what=\sum_{p=0}^{2^n-1}r_pe^{\ihat_0\varphi_p}\epshat_p,}
condition (93), taking account of (75) and the form of the matrix
$\sigma H_w$, becomes
\eqn{95}{r_0=r_1=\cdots=r_{2^n-1}\ne0.}

\begin{center}
\renewcommand{\eqprefix}{III.2}
{\Large\bfseries III.2. The Generalized Sobrero Algebra}\addcontentsline{toc}{section}{III.2. The Generalized Sobrero Algebra}\par
\vspace{3mm}
{\large\bfseries III.2.1. Preliminaries: Bases in $S$}\addcontentsline{toc}{subsection}{III.2.1. Preliminaries: Bases}
\end{center}

Let $S$ be the real algebra obtained as the direct product of the complex
field $C$ and the polynomial algebra generated by
$p(\lambda)=\lambda^{n+1}$, $n\in\mathbb N$, called the nilpotent algebra of
order $n+1$ and denoted by $P$.

A basis of $C$ is $\{\onehat,\widehat j\}$, with multiplication table
\eqn{1}{\onehat\widehat j=\widehat j\onehat=\widehat j,
\qquad \onehat^2=-\widehat j^2=\onehat.}
A basis of $P$ is $\{\widehat e_0,\widehat e_1,\ldots,
\widehat e_n\}$, with
\eqn{2}{\widehat e^k\widehat e^p=\widehat e^{k+p},
\qquad k,p\in\{0,1,\ldots,n\},}
and nilpotency condition
\eqn{3}{\widehat e^{n+1}=\widehat0.}
The order of $P$ is $n+1$. The algebra $S=C\otimes P$ is commutative,
associative, and unital, and its order over $\mathbb R$ is $2(n+1)$.

For the bases of $C$ and $P$ given above, a basis of $S$ is
$\{\widehat\varepsilon^k,\widehat j\widehat\varepsilon^k\}
_{k=0}^{n}$, and the multiplication table is
\eqn{4}{\begin{aligned}
\widehat\varepsilon^k\widehat\varepsilon^p
&=\widehat\varepsilon^{k+p},\\
\widehat\varepsilon^k\widehat j\widehat\varepsilon^p
&=\widehat j\widehat\varepsilon^{k+p},\\
\widehat j\widehat\varepsilon^k
\widehat j\widehat\varepsilon^p
&=-\widehat\varepsilon^{k+p}.
\end{aligned}}
In (4), condition (3) is understood. An arbitrary element is
\eqn{5}{\what=\sum_{k=0}^{n}
\left(a_k\widehat\varepsilon^k
+b_k\widehat j\widehat\varepsilon^k\right),
\qquad a_k,b_k\in\mathbb R,}
or
\eqn{6}{\what=\sum_{k=0}^{n}
(a_k+\widehat j b_k)\widehat\varepsilon^k.}
Thus, $S$ may be regarded as a nilpotent $C$-algebra of order $n+1$, and
\eqn{7}{\what=\sum_{k=0}^{n}\xi_k\widehat\varepsilon^k,
\qquad \xi_k\in C.}
The real algebra $S$ is thereby obtained by ``passing to the real form'' [4]
of the nilpotent $C$-algebra of order $n+1$.

For $n=1$, one obtains the Sobrero algebra [18, 140, 148]; hence $S$ will be
called the \emph{generalized Sobrero algebra}.

\begin{center}
{\large\bfseries III.2.2. Distinguished Elements, Zero Divisors, and
Elements of Special Type}\addcontentsline{toc}{subsection}{III.2.2. Distinguished Elements}
\end{center}

Using table (4), equation (7) gives
\eqn{8}{\widehat\varepsilon^k\what=
\xi_0\widehat\varepsilon^k+\xi_1\widehat\varepsilon^{k+1}
+\cdots+\xi_{n-k}\widehat\varepsilon^n,
\qquad k\in\{0,1,\ldots,n\}.}
The matrix canonically associated with element (7) is
\eqn{9}{[\what]=\begin{bmatrix}
\xi_0&\xi_1&\cdots&\xi_n\\
0&\xi_0&\cdots&\xi_{n-1}\\
\vdots&\vdots&\ddots&\vdots\\
0&0&\cdots&\xi_0
\end{bmatrix},\qquad \xi_k\in C.}
The characteristic polynomial is
\eqn{10}{q(\lambda)=\det(\lambda I_{n+1}-[\what])
=(\lambda-\xi_0)^{n+1},\qquad q(\what)=\widehat0.}
Element (7) is a zero divisor in the $C$-algebra $S$ if and only if the
constant term of (10) vanishes, that is,
\eqn{11}{\det[\what]=\xi_0^{n+1}=0.}
In the real algebra, the roots $\xi_0=a_0+\widehat j b_0$ and
$\xi_0^*=a_0-\widehat j b_0$ each have multiplicity $n+1$; therefore,
\eqn{12}{p(\lambda)=
(\lambda^2-2a_0\lambda+a_0^2+b_0^2)^{n+1}.}
The zero divisors are characterized by
\eqn{13}{a_0=b_0=0,}
and have the form
\eqn{14}{\what=\sum_{k=1}^{n}
(a_k+\widehat j b_k)\widehat\varepsilon^k.}
Equation (3) shows that these zero divisors are nilpotent.

More precisely,
\eqn{15}{\what^{n+1}=\widehat0,
\qquad n\in\mathbb N^*\text{ fixat}.}
The exponential form of complex numbers gives the semiexponential
representation
\eqn{16}{\what=\sum_{k=0}^{n}
A_ke^{\widehat j\varphi_k}\widehat\varepsilon^k,
\qquad A_k\in\mathbb R_+,\quad\varphi_k\in[0,2\pi).}
The numbers $A_k,\varphi_k$ are uniquely determined by
\eqn{17}{a_k=A_k\cos\varphi_k,
\qquad b_k=A_k\sin\varphi_k,}
that is,
\eqn{18}{A_k=\sqrt{a_k^2+b_k^2},
\qquad\varphi_k=\operatorname{atan2}(b_k,a_k).}
Here $\operatorname{atan2}$ denotes the unique angle in $[0,2\pi)$ satisfying
(17).

An element is called \emph{of special type} if, in (16),
\eqn{19}{\varphi_0=\varphi_1=\cdots=\varphi_n.}
In representation (6), this condition is equivalent to the existence of a
$c\in\mathbb R$ such that
\eqn{20}{a_k=c\,b_k,\qquad k\in\{0,1,\ldots,n\}.}
By (13) and (18), a zero divisor is characterized in (16) by
\eqn{21}{A_0=0.}
The characteristic polynomial of element (16) is
\eqn{22}{p(\lambda)=
(\lambda^2-2\lambda A_0\cos\varphi_0+A_0^2)^{n+1}.}

\medskip\noindent\textbf{Theorem III.2.2.1.}
Every element of the real algebra $S$ can be written as the sum of two
elements of special type. The product of two elements of special type is an
element of special type.

\noindent\textit{Proof.}
Let $\what\in S$ be arbitrary,
\eqn{23}{\what=\sum_{k=0}^{n}a_k\widehat\varepsilon^k
+\sum_{k=0}^{n}b_k\widehat j\widehat\varepsilon^k.}
Since $\onehat=e^{\widehat j,0}$ and
$\widehat j=e^{\widehat j\pi/2}$,
\eqn{24}{\what=
\left(\sum_{k=0}^{n}a_ke^{\widehat j,0}
\widehat\varepsilon^k\right)
+\left(\sum_{k=0}^{n}b_ke^{\widehat j\pi/2}
\widehat\varepsilon^k\right),}
which proves the first assertion.

Let two elements of special type be
\eqn{25}{\widehat w_1=
\left(\sum_{k=0}^{n}x_k\widehat\varepsilon^k\right)
e^{\widehat j\alpha},\qquad x_k\in\mathbb R_+,}
\eqn{26}{\widehat w_2=
\left(\sum_{k=0}^{n}y_k\widehat\varepsilon^k\right)
e^{\widehat j\beta},\qquad y_k\in\mathbb R_+.}
Using table (4),
\eqn{27}{\widehat w_1\widehat w_2=
\left(\sum_{k=0}^{n}\zeta_k\widehat\varepsilon^k\right)
e^{\widehat j\gamma},}
where
\eqn{28}{\gamma=\alpha+\beta\pmod{2\pi},}
\eqn{29}{\zeta_k=x_0y_k+x_1y_{k-1}+\cdots+x_{k-1}y_1+x_ky_0,
\qquad k\in\{0,1,\ldots,n\}.}
The coefficients $\zeta_k$ are nonnegative real numbers, and (27) is an
element of special type. This completes the proof.

\renewcommand{\eqprefix}{III.3}

\section*{III.3. The Walsh Algebra}
\subsection*{III.3.1. Preliminaries: Bases in $\W$}

Let $U_k$, $k\in\{0,1,\ldots,n-1\}$, be $n$ associative, commutative,
unital algebras of order two over the real field. If
$\{\wh 1_k,\wh u_k\}$ is a basis of the $\R$-algebra $U_k$, its
multiplication table is
\eqn{1}{\wh 1_k\wh u_k=\wh u_k\wh 1_k=\wh u_k,\qquad \wh u_k^{2}=\wh 1_k^{2}=\wh 1_k.}
In (1), $\wh 1_k$ is the identity of $U_k$. The algebras $U_k$ are called
bireal algebras [20].

Consider the algebra $\W=U_0\otimes U_1\otimes\cdots\otimes U_{n-1}$.
It is an associative, commutative, unital $\R$-algebra of order $2^n$ [54].
Adopt the convention
\eqn{2}{\wh u_k^{\,\alpha}=\begin{cases}\wh 1_k,&\alpha=0,\\ \wh u_k,&\alpha=1.
\end{cases}}
The elements of a basis of $\W$ can be written as
\eqn{3}{\wh w_p=\wh u_0^{p_0}\wh u_1^{p_1}\cdots\wh u_{n-1}^{p_{n-1}},\quad p=(p_0p_1\ldots p_{n-1})_2,\quad 0\le p<2^n.}
Clearly, $p_k\in\{0,1\}$. The identity is
$\wh w_0=\wh 1_0\wh 1_1\cdots\wh 1_{n-1}$ and will be denoted by $\one$.
In (3), the identities of the component algebras will be omitted. This
convention considerably simplifies the calculations.

For a second basis element,
\eqn{4}{\wh w_q=\wh u_0^{q_0}\wh u_1^{q_1}\cdots\wh u_{n-1}^{q_{n-1}},\quad q=(q_0,q_1,\ldots,q_{n-1})_2,}
we obtain
\eqn{5}{\wh w_p\wh w_q=\wh u_0^{p_0+q_0}\wh u_1^{p_1+q_1}\cdots\wh u_{n-1}^{p_{n-1}+q_{n-1}},}
and, taking (1) into account,
\eqn{6}{\wh w_p\wh w_q=\wh w_{p\oplus q},\qquad p,q\in\{0,1,\ldots,2^n-1\}.}

Here $p\oplus q$ denotes the carry-free, bitwise modulo-$2$ sum of $p$ and
$q$ [1]. Under multiplication (6), the set
$\{\wh w_p\}_{p=0}^{2^n-1}$ forms the dyadic group of order $2^n$.
Thus, $\W$ is a group algebra [136].

\subsection*{III.3.2. Distinguished Elements, Orthogonal Idempotents, and Zero Divisors}

Consider the elements
\eqn{7}{\wh e_p=\frac1{2^n}\prod_{k=0}^{n-1}\left[\one+(-1)^{p_k}\wh u_k\right],\quad p=(p_0,p_1,\ldots,p_{n-1})_2.}

\textbf{Theorem III.3.2.1.} The elements defined by (7) are $2^n$
orthogonal idempotents in the $\R$-algebra $\W$.

\emph{Proof.} For $p,q\in\{0,1,\ldots,2^n-1\}$, equations (7) and (1) give
\eqn{8}{\wh e_p\wh e_q=\frac1{2^{n+1}}\prod_{k=0}^{n-1}\left\{(-1)^{p_k+q_k}+1+\left[(-1)^{p_k}+(-1)^{q_k}\right]\wh u_k\right\}.}
If $p\ne q$, there exists at least one $k$ such that $p_k\ne q_k$, hence
$\wh e_p\wh e_q=\z$. If $p=q$, then $\wh e_p\wh e_q=\wh e_p$. Therefore,
\eqn{9}{\wh e_p\wh e_q=\begin{cases}\wh e_p,&p=q,\\ \z,&p\ne q.
\end{cases}}
Since the order of the algebra is $2^n$, we have [54]
\eqn{10}{\one=\sum_{p=0}^{2^n-1}\wh e_p.}

\textbf{Theorem III.3.2.2.} The $\R$-algebra $\W$ is polynomial and
semisimple, being the direct sum of $2^n$ algebras isomorphic to the real
field $\R$.

We first prove the relations
\eqn{11}{\wh u_k\wh e_p=(-1)^{p_k}\wh e_p,\quad p\in\{0,\ldots,2^n-1\},\quad k\in\{0,\ldots,n-1\}.}

Indeed, from (7) and (1),
\[
\wh u_k\wh e_p=\wh u_k\frac{\one+(-1)^{p_k}\wh u_k}{2}
\prod_{\substack{m=0\\m\ne k}}^{n-1}\frac{\one+(-1)^{p_m}\wh u_m}{2}
=(-1)^{p_k}\wh e_p.
\]
Summing (11) over $p$ gives
\eqn{12}{\wh u_k\sum_{p=0}^{2^n-1}\wh e_p=\sum_{p=0}^{2^n-1}(-1)^{p_k}\wh e_p.}
From (10) and (12),
\eqn{13}{\wh u_k=\sum_{p=0}^{2^n-1}(-1)^{p_k}\wh e_p,\quad k=0,1,\ldots,n-1.}
Using (9), it follows that
\eqn{14}{\wh u_k^{m}=\sum_{p=0}^{2^n-1}(-1)^{mp_k}\wh e_p,\quad m\in\N.}
For $q=(q_0,q_1,\ldots,q_{n-1})_2$, equations (4) and (14) give
\eqn{15}{\wh w_q=\prod_{k=0}^{n-1}\wh u_k^{q_k}=\sum_{p=0}^{2^n-1}(-1)^{p_0q_0+\cdots+p_{n-1}q_{n-1}}\wh e_p.}
With $\langle p,q\rangle=p_0q_0+p_1q_1+\cdots+p_{n-1}q_{n-1}$,
\eqn{16}{\wh w_q=\sum_{p=0}^{2^n-1}(-1)^{\langle p,q\rangle}\wh e_p.}
Introduce the matrix notation
\eqn{17}{\boldsymbol w=[\wh w_0,\ldots,\wh w_{2^n-1}]^T,\quad \boldsymbol e=[\wh e_0,\ldots,\wh e_{2^n-1}]^T,\quad H=\bigl((-1)^{\langle p,q\rangle}\bigr)_{p,q=0}^{2^n-1}.}
Then
\eqn{18}{\boldsymbol w=H\boldsymbol e.}

The square matrix $H$ is the Hadamard matrix of order $2^n$ [1,151]. It can
be written as the Kronecker power
\eqn{19}{H=\begin{bmatrix}1&1\\1&-1\end{bmatrix}^{\!\otimes n}.}
The Hadamard matrix is symmetric and invertible, and [1,151]
\eqn{20}{H^{-1}=\frac1{2^n}H.}
From (18) and (20),
\eqn{21}{\boldsymbol e=\frac1{2^n}H\boldsymbol w,}
that is,
\eqn{22}{\wh e_p=\frac1{2^n}\sum_{q=0}^{2^n-1}(-1)^{\langle p,q\rangle}\wh w_q.}
Relations (18), (21), together with their componentwise forms (16), (22),
show that the orthogonal idempotents $\wh e_p$ form a basis of $\W$. The
change-of-basis matrix from the ``direct-product'' basis to the orthogonal-
idempotent basis is the Hadamard matrix.

An arbitrary element $\wh f\in\W$,
\eqn{23}{\wh f=\sum_{p=0}^{2^n-1}f_p\wh w_p,\qquad f_p\in\R,}
has the unique expansion
\eqn{24}{\wh f=\sum_{p=0}^{2^n-1}a_p\wh e_p,\qquad a_p\in\R,}
where
\eqn{25}{a_p=\sum_{q=0}^{2^n-1}(-1)^{\langle p,q\rangle}f_q.}
Since $\wh f$ is arbitrary, (24) proves that $\W$ is the direct sum of $2^n$
algebras isomorphic to $\R$. Because the real field is a commutative,
associative, unital, polynomial, and semisimple $\R$-algebra, the theorem is
proved. $\square$

Equation (24) gives
\eqn{26}{\wh e_k\wh f=a_k\wh e_k,\qquad k=0,1,\ldots,2^n-1.}
The matrix corresponding to element (24) under the canonical isomorphism is
\eqn{27}{[\wh f]=\operatorname{diag}(a_0,a_1,\ldots,a_{2^n-1}).}
Consequently,
\eqn{28}{\det[\wh f]=a_0a_1\cdots a_{2^n-1},}
and
\eqn{29}{\det(\lambda I_{2^n}-[\wh f])=\prod_{p=0}^{2^n-1}(\lambda-a_p).}
The zero divisors are characterized by
\eqn{30}{a_0a_1\cdots a_{2^n-1}=0.}
The element $\wh f\in\W$ is called a zero divisor of rank $p$ if its
representation (24) satisfies
\eqn{31}{a_p=0.}
By (25), the equivalent condition in the ``direct-product'' basis is
\eqn{32}{\sum_{q=0}^{2^n-1}(-1)^{\langle p,q\rangle}f_q=0.}
If $\wh f$ is a zero divisor of rank $p$, then (26) gives
$\wh e_p\wh f=\z$. Conversely, this equality and the independence of the
idempotents imply $a_p=0$.

\textbf{Theorem III.3.2.3.} The element $\wh f\in\W$ is a zero divisor of
rank $p$ if and only if
\[
\wh f\wh e_p=\z,\qquad p\in\{0,1,\ldots,2^n-1\}.
\]

\emph{Remark.} An element may simultaneously be a zero divisor of several
ranks. The idempotent $\wh e_p$ is a zero divisor of all ranks
$0,1,\ldots,p-1,p+1,\ldots,2^n-1$.

The characteristic polynomial of $\wh f$ written in the form (24) is
\eqn{33}{p(\lambda)=\prod_{p=0}^{2^n-1}(\lambda-a_p),\qquad p(\wh f)=\z.}
If $\wh f$ is given in the ``direct-product'' basis by (23), then
\eqn{34}{p(\lambda)=\prod_{p=0}^{2^n-1}\left[\lambda-\sum_{q=0}^{2^n-1}(-1)^{\langle p,q\rangle}f_q\right],\qquad p(\wh f)=\z.}

\subsection*{III.3.3. Exponential Form of the Elements: An Extension of the Walsh Algebra}

We investigate the conditions under which an element $\wh f\in\W$ admits the
representation
\eqn{35}{\wh f=\exp\wh g,\qquad \wh g\in\W.}

\textbf{Theorem III.3.3.1.} If the roots of the characteristic polynomial of
$\wh f\in\W$ are strictly positive, then $\wh f$ admits an exponential
representation.

\emph{Proof.} Let
\eqn{36}{\wh f=\sum_{p=0}^{2^n-1}f_p\wh w_p.}
The roots of characteristic polynomial (34) are
\eqn{37}{a_p=\sum_{q=0}^{2^n-1}(-1)^{\langle p,q\rangle}f_q.}
The element can be written as
\eqn{38}{\wh f=\sum_{p=0}^{2^n-1}a_p\wh e_p.}
Under the assumptions of the theorem, $a_p>0$ for every $p$, and
\eqn{39}{\wh f=\sum_{p=0}^{2^n-1}e^{\ln a_p}\wh e_p
=\exp\!\left(\sum_{p=0}^{2^n-1}\ln a_p\,\wh e_p\right)
=\exp\!\left(\sum_{p=0}^{2^n-1}g_p\wh w_p\right),}
where
\eqn{40}{g_p=\frac1{2^n}\sum_{k=0}^{2^n-1}(-1)^{\langle p,k\rangle}\ln a_k.}

\emph{Remarks.} 1$^\circ$ Zero divisors and elements for which at least one
root of the characteristic polynomial is negative do not admit exponential
representations.

2$^\circ$ A computational formula can be given for the exponential of
\eqn{41}{\exp\wh f=\exp\!\left(\sum_{p=0}^{2^n-1}f_p\wh w_p\right).}
Passing from the ``direct-product'' basis to the idempotent basis gives
\eqn{42}{\sum_{p=0}^{2^n-1}f_p\wh w_p=\sum_{p=0}^{2^n-1}a_p\wh e_p,}
where
\eqn{43}{a_p=\sum_{k=0}^{2^n-1}(-1)^{\langle p,k\rangle}f_k.}
Therefore,
\eqn{44}{\exp\wh f=\frac1{2^n}\sum_{p=0}^{2^n-1}\left[\sum_{k=0}^{2^n-1}(-1)^{\langle p,k\rangle}
\exp\!\left(\sum_{q=0}^{2^n-1}(-1)^{\langle k,q\rangle}f_q\right)\right]\wh w_p.}

Since not every element of $\W$ admits an exponential representation, we
extend the algebra. Let $\W^*$ be an associative, commutative, unital
$\C$-algebra of order $2^n$, whose basis satisfies
\eqn{45}{\wh w_p\wh w_q=\wh w_{p\oplus q}.}
An arbitrary element is
\eqn{46}{\wh\xi=\sum_{p=0}^{2^n-1}\xi_p\wh w_p,\qquad \xi_p\in\C.}
The elements
\eqn{47}{\wh e_p=\frac1{2^n}\sum_{k=0}^{2^n-1}(-1)^{\langle p,k\rangle}\wh w_k}
are $2^n$ orthogonal idempotents in $\W^*$, and
\eqn{48}{\wh w_k=\sum_{p=0}^{2^n-1}(-1)^{\langle p,k\rangle}\wh e_p.}

Thus,
\eqn{49}{\wh\xi=\sum_{p=0}^{2^n-1}y_p\wh e_p,\qquad y_p\in\C,}
where
\eqn{50}{y_p=\sum_{k=0}^{2^n-1}(-1)^{\langle p,k\rangle}\xi_k.}

\textbf{Theorem III.3.3.2.} The algebra $\W^*$ is the direct sum of $2^n$
algebras isomorphic to the complex field.

\emph{Remark.} Passing to the real form of the $\C$-algebra $\W^*$ yields an
$\R$-algebra of order $2^{n+1}$. A basis is
\eqn{51}{B=\{\wh w_p,i_0\wh w_p\}_{p=0}^{2^n-1},\qquad i_0^2=-\wh w_0=-\one.}
Also let
\eqn{52}{B^*=\{\wh e_p,i_0\wh e_p\}_{p=0}^{2^n-1}.}
With the notation
\eqn{53}{\boldsymbol w^*=[\wh w_0,\ldots,\wh w_{2^n-1},i_0\wh w_0,\ldots,i_0\wh w_{2^n-1}]^T,}
\eqn{54}{\boldsymbol e^*=[\wh e_0,\ldots,\wh e_{2^n-1},i_0\wh e_0,\ldots,i_0\wh e_{2^n-1}]^T,}
and
\eqn{55}{S=\begin{bmatrix}H&0\\0&H\end{bmatrix},}
we obtain
\eqn{56}{\boldsymbol w^*=S\boldsymbol e^*.}
The matrix $S$ is invertible and
\eqn{57}{S^{-1}=\frac1{2^n}S,}
hence
\eqn{58}{\boldsymbol e^*=\frac1{2^n}S\boldsymbol w^*.}

Therefore, $B^*$ is also a basis of $\W^*$. The algebra $\W^*$ contains a
subalgebra isomorphic to $\W$, generated by the first $2^n$ elements of basis
(51).

In $\W^*$, every element that is not a zero divisor admits an exponential
representation.

Let
\eqn{59}{\wh w=\sum_{p=0}^{2^n-1}(a_p\wh w_p+b_pi_0\wh w_p),\qquad a_p,b_p\in\R.}
We show that there exists a $\wh u\in\W^*$ such that
\eqn{60}{\wh w=\exp\wh u.}
In the basis $B^*$,
\eqn{61}{\wh w=\sum_{p=0}^{2^n-1}(c_p\wh e_p+d_pi_0\wh e_p),}
where
\eqn{62}{c_p=\sum_{k=0}^{2^n-1}(-1)^{\langle p,k\rangle}a_k,}
\eqn{63}{d_p=\sum_{k=0}^{2^n-1}(-1)^{\langle p,k\rangle}b_k.}
With the notation
\eqn{64}{\rho_p=\sqrt{c_p^2+d_p^2},}
\eqn{65}{\cos\theta_p=\frac{c_p}{\rho_p},}
\eqn{66}{\sin\theta_p=\frac{d_p}{\rho_p},}
the properties of the orthogonal idempotents give
\eqn{67}{\wh w=\exp\!\left[\sum_{p=0}^{2^n-1}\bigl(\ln\rho_p\,\wh e_p+\theta_p i_0\wh e_p\bigr)\right].}
This proves (60) for every element with $\rho_p\ne0$ for all $p$, that is,
for every element that is not a zero divisor.

\emph{This construction completes the complex extension of the Walsh algebra.}

\renewcommand{\eqprefix}{III.4}
\section*{III.4. A Characterization of Finite-Dimensional $\R$-Algebras}

Let $A_n$ be a finite-dimensional, commutative, associative, unital
$\R$-algebra of order $n$, and let $B=\{\wh\varepsilon_j\}_{j=1}^n$ be a
basis of $A_n$. We establish conditions under which $A_n$ is polynomial and
determine the change-of-basis matrix from $B$ to a power basis.

\subsection*{III.4.1. General Results}

\textbf{Theorem III.4.1.} If the algebra $A_n$ contains at least one element
whose characteristic polynomial is minimal, then $A_n$ is polynomial.

\emph{Proof.} Let $\wh d\in A_n$ be such an element. Compute
\eqn{1}{\exp(\wh d\,t)=\sum_{j=1}^n f_j(t)\wh\varepsilon_j,\qquad f_j\in C^\infty(\R).}
By Theorem II.2.3, the functions $f_j$ form a fundamental system of solutions
of the linear differential equation with constant coefficients whose
characteristic polynomial is that of $\wh d$. Differentiating $k$ times gives
\eqn{2}{\wh d^{k}\exp(\wh d\,t)=\sum_{j=1}^n f_j^{(k)}(t)\wh\varepsilon_j,\quad  k=0,1,\ldots,n-1.}
At $t=0$,
\eqn{3}{\wh d^{k}=\sum_{j=1}^n f_j^{(k)}(0)\wh\varepsilon_j.}

In matrix form,
\eqn{4}{
\begin{bmatrix}\wh d^0\\\wh d\\\vdots\\\wh d^{n-1}\end{bmatrix}
=\begin{bmatrix}
f_1(0)&f_2(0)&\cdots&f_n(0)\\
f'_1(0)&f'_2(0)&\cdots&f'_n(0)\\
\vdots&\vdots&\ddots&\vdots\\
f_1^{(n-1)}(0)&f_2^{(n-1)}(0)&\cdots&f_n^{(n-1)}(0)
\end{bmatrix}
\begin{bmatrix}\wh\varepsilon_1\\\wh\varepsilon_2\\\vdots\\\wh\varepsilon_n\end{bmatrix}.}
The matrix is nonsingular: its determinant is the Wronskian of the linearly
independent system $f_j$. Hence $B^*=\{\wh d^{k}\}_{k=0}^{n-1}$ is a basis
of $A_n$, so the algebra is polynomial [54]. Relation (4) gives the change of
basis explicitly.

\subsection*{III.4.2. Applications}

The algebras considered in III.1, III.2, and III.3 are polynomial. We now give
their generating polynomials.

\textbf{a) The algebra $\bigotimes \C^{n+1}$.}

Let
\eqn{5}{\wh d=i_0\omega_0+i_1\omega_1+\cdots+i_n\omega_n,\qquad \omega_k\in\R.}
In accordance with III.1, relation (21),
\eqn{6}{i_k=\sum_{p=0}^{2^n-1}(-1)^{p_k}i_0\wh e_p,\quad k=1,\ldots,n.}
Then
\eqn{7}{\wh d=\sum_{p=0}^{2^n-1}\Omega_p i_0\wh e_p,}
where
\eqn{8}{\Omega_p=\omega_0+\sum_{k=1}^n(-1)^{p_k}\omega_k,\qquad p=(p_1p_2\ldots p_n)_2.}

If $\Omega_p\ne0$ for every $p$, the characteristic polynomial of element
(5) is minimal:
\eqn{9}{p(\lambda)=\prod_{p=0}^{2^n-1}(\lambda^2+\Omega_p^2).}
This polynomial follows from (7) by using the properties of the orthogonal
idempotents. By the preceding theorem, the algebra is polynomial, with a
generating polynomial having distinct, purely complex roots.

\textbf{b) The generalized Sobrero algebra.}

Let $\wh d\in S$,
\eqn{10}{\wh d=\alpha j\wh\varepsilon^{0}+\beta\wh\varepsilon,\qquad \alpha,\beta\in\R.}
By III.2.1, relation (12), the characteristic polynomial is
\eqn{11}{p(\lambda)=(\lambda^2+\alpha^2)^{n+1}.}
We have
\[
p(\wh d)=\left[\beta(\beta\wh\varepsilon+2\alpha j\beta)\wh\varepsilon\right]^{n+1}.
\]
If $\alpha\ne0$ and $\beta\ne0$, polynomial (11) is minimal, because
\eqn{12}{\wh\varepsilon^{p}\ne\z\ (p=0,1,\ldots,n),\qquad \wh\varepsilon^{n+1}=\z.}
Thus, $S$ is polynomial, generated by
$p(\lambda)=(\lambda^2+\alpha^2)^{n+1}$, $\alpha\in\R^*$, with roots in
$\C\setminus\R$ of multiplicity $n+1$. For $\alpha=\beta=1$, the element
\eqn{13}{\wh d=j\wh\varepsilon^{0}+\wh\varepsilon}
generates the generalized Sobrero algebra, with polynomial
\eqn{14}{p(\lambda)=(\lambda^2+1)^{n+1}.}

\textbf{c) The Walsh algebra.}

Let
\eqn{15}{\wh d=\sum_{k=0}^{n-1}\alpha_k\wh u_k,\qquad \alpha_k\in\R.}
In accordance with III.3.1, relation (14),
\eqn{16}{\wh d=\sum_{p=0}^{2^n-1}\Omega_p\wh e_p,}
where
\eqn{17}{\Omega_p=\sum_{k=0}^{n-1}(-1)^{p_k}\alpha_k,\qquad p=(p_0,p_1,\ldots,p_{n-1})_2.}
If the real numbers $\Omega_p$ are distinct and nonzero, the polynomial
\eqn{18}{p(\lambda)=\prod_{p=0}^{2^n-1}(\lambda-\Omega_p)}
is minimal for element (15). Thus, the algebra $\W$ is polynomial, generated
by polynomial (18), which has distinct real roots.

\chaptertoc{CHAPTER IV}{Applications of Symbolic Representations}
\renewcommand{\eqprefix}{IV.1.1}

\begin{center}
{\Large\bfseries CHAPTER IV}\vspace{1.2em}

{\large\bfseries Applications of Symbolic Representations}
\end{center}

\section*{IV.1. Symbolic Representation of Signals over Particular $R$-Algebras}
\subsection*{IV.1.1. Symbolic Representations over $\A$}

Let $\A$ be the algebra introduced in III.1. Two bases of this $R$-algebra, denoted by $B$ and $B^*$, were identified:
\eqn{1}{B=\{\wh v_p\}_{p=0}^{2^{n+1}-1},}
where
\eqn{2}{\wh v_p=i_0^{p_0}i_1^{p_1}\cdots i_n^{p_n},\quad
p=(p_0p_1\ldots p_n)_{\rm Gray},\quad i_k^2=-\one.}
Here $p=(p_0p_1\ldots p_n)_{\rm Gray}$ is the representation of the number $p$ in bit-reversed binary Gray code. The elements of the basis $B$ satisfy
\eqn{3}{\wh v_p^{2}=\begin{cases}
\one,&0\le p\le 2^n-1,\\
-\one,&2^n\le p\le 2^{n+1}-1.
\end{cases}}
The second basis is
\eqn{4}{B^*=\{\wh e_p,i_0\wh e_p\}_{p=0}^{2^n-1},\qquad
\wh e_p=\frac1{2^n}\prod_{k=1}^n\left[\one-(-1)^{p_k}i_0i_k\right].}
The elements $\wh e_p$ are the $2^n$ orthogonal idempotents of the algebra.

The changes between the bases $B$ and $B^*$ are given by the matrix relations
\eqn{5}{
\begin{bmatrix}\wh v_0\\\vdots\\\wh v_{2^n-1}\\i_0\wh v_0\\\vdots\\i_0\wh v_{2^n-1}\end{bmatrix}
=\begin{bmatrix}\sigma H_W&0\\0&\sigma'H_W\end{bmatrix}
\begin{bmatrix}\wh e_0\\\vdots\\\wh e_{2^n-1}\\i_0\wh e_0\\\vdots\\i_0\wh e_{2^n-1}\end{bmatrix},}
and, respectively,
\eqn{6}{
\begin{bmatrix}\wh e_0\\\vdots\\\wh e_{2^n-1}\\i_0\wh e_0\\\vdots\\i_0\wh e_{2^n-1}\end{bmatrix}
=\frac1{2^n}\begin{bmatrix}H_W\sigma&0\\0&H_W\sigma'^{T}\end{bmatrix}
\begin{bmatrix}\wh v_0\\\vdots\\\wh v_{2^{n+1}-1}\end{bmatrix}.}
The matrices $\sigma$, $\sigma'$, and $H_W$ were specified in III.1.2, relations (29), (31), and (32). Componentwise,
\eqn{7}{\wh v_p=\sum_{k=0}^{2^n-1}(-1)^{\langle k,p\rangle}\sigma(p)\wh e_k,\quad 0\le p<2^n,}
\eqn{8}{\wh e_p=\frac1{2^n}\sum_{k=0}^{2^n-1}(-1)^{\langle k,p\rangle}\sigma(k)\wh v_k.}
In these formulas,
\[
\sigma(p)=(-1)^{\left\lfloor(p_0+p_1+\cdots+p_n)/2\right\rfloor},
\quad p=(p_0p_1\ldots p_n)_{\rm Gray}.
\]
Let $\wh d\in\A$ be expressed in the basis $B$:
\eqn{9}{\wh d=i_0\omega_0+i_1\omega_1+\cdots+i_n\omega_n,\qquad \omega_k\in\mathbb R.}

In the basis $B^*$, element (9) has the expression
\eqn{10}{\wh d=\sum_{p=0}^{2^n-1}\Omega_p i_0\wh e_p,}
where
\eqn{11}{\Omega_p=\omega_0+\sum_{k=1}^{n}(-1)^{p_k}\omega_k,\qquad p=(p_1p_2\ldots p_n)_2.}
Assuming that the numbers $\Omega_p$ are nonzero and have distinct absolute values, the characteristic polynomial, which is also the minimal polynomial, is
\eqn{12}{p(\lambda)=\prod_{p=0}^{2^n-1}(\lambda^2+\Omega_p^2),\qquad p(\wh d)=\zero.}
According to Theorem II.2.1, operator (9) provides a symbolic representation of every solution of the equation
\eqn{13}{p(D)f=0,\qquad Df=\frac{df}{dt}.}
Let $V$ be the vector space of solutions. To determine a symbolic representation
\eqn{14}{\Psi:V\longrightarrow\A,\qquad \Psi(\dot f)=\wh d\,\Psi(f),}
we compute $\exp(\wh d t)$.

In basis (9), using Euler-type relations,
\eqn{15}{
\begin{aligned}
e^{\wh d t}
&=\prod_{k=0}^{n}e^{i_k\omega_kt}
=\prod_{k=0}^{n}(\cos\omega_kt+i_k\sin\omega_kt)\\
&=\sum_{p=0}^{2^{n+1}-1}
\left[\prod_{k=0}^{n}\cos\!\left(\omega_kt-p_k\frac\pi2\right)\right]\wh v_p,
\quad p=(p_0p_1\ldots p_n)_{\rm Gray}.
\end{aligned}}
Consequently, the representation is uniquely determined by
\eqn{16}{\Psi\!\left[\prod_{k=0}^{n}\cos\!\left(\omega_kt-p_k\frac\pi2\right)\right]=\wh v_p.}

In basis (10), using the orthogonal idempotents,
\eqn{17}{e^{\wh d t}=\sum_{p=0}^{2^n-1}\left(\cos\Omega_pt\,\wh e_p+\sin\Omega_pt\,i_0\wh e_p\right).}
Comparing relations (15) and (17), after using (8), yields the identities
\eqn{18}{
2^n\prod_{k=0}^{n}\cos\!\left(\omega_kt-p_k\frac\pi2\right)
=\sum_{q=0}^{2^n-1}(-1)^{\langle p,q\rangle}
\cos\!\left(\Omega_qt-\frac\pi2\sum_{i=0}^{n}p_i\right).}
Conversely, expressing (15) in the basis $B^*$ gives
\eqn{19}{\cos\Omega_pt=\sum_{k=0}^{2^n-1}(-1)^{\langle k,p\rangle}\sigma(k)\varphi_k(t),}
\eqn{20}{\sin\Omega_pt=\sum_{k=0}^{2^n-1}(-1)^{\langle k,p\rangle}\sigma(2^n+k)\varphi_{2^n+k}(t),}
where
\eqn{21}{\varphi_k(t)=\prod_{j=0}^{n}\cos\!\left(\omega_jt-k_j\frac\pi2\right),\qquad
k=(k_0k_1\ldots k_n)_{\rm Gray}.}

Now let
\eqn{22}{\wh w=\sum_{p=0}^{2^{n+1}-1}w_p\wh v_p,\qquad w_p\in\mathbb R.}
Under the operator $\wh d$, element (22) represents the function $\varphi=\Psi^{-1}(\wh w)$:
\eqn{23}{\varphi=\sum_{p=0}^{2^{n+1}-1}w_p\Psi^{-1}(\wh v_p),}
hence
\eqn{24}{\varphi(t)=\sum_{p=0}^{2^{n+1}-1}w_p\varphi_p(t).}

Explicitly,
\eqn{25}{\varphi_p(t)=\prod_{k=0}^{n}\cos\!\left(\omega_kt-p_k\frac\pi2\right),\qquad
p=(p_0p_1\ldots p_n)_{\rm Gray}.}

\emph{Remark.} In bit-reversed Gray code, $p_0=0$ for even $p$ and $p_0=1$ for odd $p$. Signal (24) is therefore a harmonic oscillation with simultaneous amplitude and phase modulation, with carrier angular frequency $\omega_0$:
\eqn{26}{\varphi(t)=a(t)\cos\omega_0t+b(t)\sin\omega_0t,}
where
\eqn{27}{a(t)=\sum_{p\in P'}w_p\prod_{k=1}^{n}\cos\!\left(\omega_kt-p_k\frac\pi2\right),}
\eqn{28}{b(t)=\sum_{p\in P''}w_p\prod_{k=1}^{n}\cos\!\left(\omega_kt-p_k\frac\pi2\right),}
with $P'=\{0,2,4,\ldots,2^{n+1}-2\}$ and $P''=\{1,3,5,\ldots,2^{n+1}-1\}$. Relation (26) takes the form
\eqn{29}{\varphi(t)=A(t)\cos[\omega_0t-\theta(t)],}
where
\eqn{30}{A(t)=\sqrt{a^2(t)+b^2(t)},\qquad
\cos\theta(t)=\frac{a(t)}{A(t)},\qquad
\sin\theta(t)=\frac{b(t)}{A(t)}.}
The signals $a(t)$ and $b(t)$ are themselves harmonic signals with simultaneous amplitude and phase modulation.

\medskip
\textbf{1$^\circ$ Orthogonal idempotents.} From (8) and the linearity of $\Psi$,
\eqn{31}{\Psi^{-1}(\wh e_p)=\frac1{2^n}\sum_{k=0}^{2^n-1}(-1)^{\langle k,p\rangle}\sigma(k)\Psi^{-1}(\wh v_k).}

From (16), (19), and (31),
\eqn{32}{\Psi^{-1}(\wh e_p)=\frac1{2^n}\cos\Omega_pt.}
Similarly, from the expansion
\eqn{33}{i_0\wh e_p=\frac1{2^n}\sum_{k=0}^{2^n-1}(-1)^{\langle k,p\rangle}\sigma(2^n+k)\wh v_{2^n+k},}
and identity (20), it follows that
\eqn{34}{\Psi^{-1}(i_0\wh e_p)=\frac1{2^n}\sin\Omega_pt.}
Relations (32) and (34) combine into
\eqn{35}{\Psi^{-1}(i_0^{m}\wh e_p)=\frac1{2^n}\cos\!\left(\Omega_pt-m\frac\pi2\right),\qquad m\in\mathbb N.}
The counterpart of relation III.1, (20),
\eqn{36}{i_k\wh e_p=(-1)^{p_k}i_0\wh e_p,}
gives
\eqn{37}{\Psi^{-1}(i_k\wh e_p)=\frac1{2^n}(-1)^{p_k}\sin\Omega_pt,}
or
\eqn{38}{\Psi^{-1}(i_k\wh e_p)=\frac1{2^n}\sin(\Omega_pt+p_k\pi).}
Identities (18) are the counterparts of the algebraic identities
\eqn{39}{\wh v_k=i_0^{k_0}i_1^{k_1}\cdots i_n^{k_n}
=\sum_{p=0}^{2^n-1}(-1)^{\langle k,p\rangle}i_0^{k_0+\cdots+k_n}\wh e_p.}
Relations (32) and (34) also yield the correspondence
\eqn{40}{\Psi^{-1}\!\left(e^{i_0\phi}\wh e_p\right)=\frac1{2^n}\cos(\Omega_pt-\phi).}

Let an arbitrary element be expressed in the basis $B^*$:
\eqn{41}{\wh w=\sum_{p=0}^{2^n-1}A_pe^{i_0\phi_p}\wh e_p,\qquad A_p\ge0,\quad \phi_p\in[0,2\pi).}
From (40),
\eqn{42}{\Psi^{-1}(\wh w)=\frac1{2^n}\sum_{p=0}^{2^n-1}A_p\cos(\Omega_pt-\phi_p),}
which is the general solution of equation (13).

\textbf{2$^\circ$ Zero divisors.} Element (41) is a zero divisor of rank $k$ if
\eqn{43}{\wh w\wh e_k=\zero\quad\Longleftrightarrow\quad A_k=0.}
It represents the signal
\eqn{44}{\Psi^{-1}(\wh w)=\frac1{2^n}\sum_{\substack{p=0\\p\ne k}}^{2^n-1}A_p\cos(\Omega_pt-\phi_p),}
which is the general solution of the equation with characteristic polynomial
\eqn{45}{p_k(\lambda)=\prod_{\substack{p=0\\p\ne k}}^{2^n-1}(\lambda^2+\Omega_p^2),\qquad \deg p_k=2^{n+1}-2.}
If $\wh w$ is a divisor of ranks $k_1,\ldots,k_r$, then
\eqn{46}{A_{k_1}=A_{k_2}=\cdots=A_{k_r}=0,}
and it represents
\eqn{47}{\Psi^{-1}(\wh w)=\frac1{2^n}\sum_{\substack{p=0\\p\notin\{k_1,\ldots,k_r\}}}^{2^n-1}A_p\cos(\Omega_pt-\phi_p),}
with characteristic polynomial
\eqn{48}{p_{k_1\ldots k_r}(\lambda)=\prod_{\substack{p=0\\p\notin\{k_1,\ldots,k_r\}}}^{2^n-1}(\lambda^2+\Omega_p^2),\quad \deg p_{k_1\ldots k_r}=2^{n+1}-2r.}

Thus, the zero divisors of $\A$ symbolically represent degenerate signals, namely solutions of differential equations whose order is lower than the degree of the characteristic polynomial of the differentiation operator. Their characteristic polynomials are among the divisors of the characteristic polynomial of $\wh d$.

\emph{Remarks.} a) The idempotents $\wh e_p$ are zero divisors of every rank other than $p$. The functions they represent satisfy equations with polynomial
\eqn{49}{p(\lambda)=\lambda^2+\Omega_p^2.}
b) These properties follow from the fact that the minimal polynomials of zero divisors have degree strictly lower than that of their characteristic polynomials.

\textbf{3$^\circ$ Special-type elements.} In general, an arbitrary element represents a harmonic signal with simultaneous amplitude and phase modulation. We seek the elements representing signals that are only amplitude-modulated with respect to the carriers $\omega_0,\ldots,\omega_n$:
\eqn{50}{s(t)=A\prod_{k=0}^{n}\cos(\omega_kt-\phi_k),\qquad A>0.}
We use identity (18) for $p=0$:
\eqn{51}{2^n\prod_{k=0}^{n}\cos\alpha_k
=\sum_{p=0}^{2^n-1}\cos\!\left[\alpha_0+\sum_{k=1}^{n}(-1)^{p_k}\alpha_k\right].}
Setting $\alpha_k=\omega_kt-\phi_k$ gives
\eqn{52}{s(t)=\frac{A}{2^n}\sum_{p=0}^{2^n-1}\cos(\Omega_pt-\Phi_p),}
where
\eqn{53}{\Omega_p=\omega_0+\sum_{k=1}^{n}(-1)^{p_k}\omega_k,}
\eqn{54}{\Phi_p=\phi_0+\sum_{k=1}^{n}(-1)^{p_k}\phi_k.}

Using relation (40), equation (52) yields
\eqn{55}{\Psi[s(t)]=A\sum_{p=0}^{2^n-1}e^{i_0\Phi_p}\wh e_p.}
The properties of the orthogonal idempotents give
\eqn{56}{\Psi[s(t)]=A\exp\!\left(\sum_{p=0}^{2^n-1}\Phi_pi_0\wh e_p\right).}
Passing to the basis $B$, analogously to the passage from (9) to (10), gives
\eqn{57}{\Psi[s(t)]=A\exp(i_0\phi_0+i_1\phi_1+\cdots+i_n\phi_n).}
Elements of the form (57) are called \emph{special-type elements}.

We now state their characteristic conditions. Let
\eqn{58}{\wh w=\sum_{p=0}^{2^n-1}A_pe^{i_0\Phi_p}\wh e_p.}
In the basis $B$, it admits the semi-exponential form
\eqn{59}{\wh w=\wh\rho\,e^{\wh\theta},}
where
\eqn{60}{\wh\rho=\sum_{k=0}^{2^n-1}\rho_k\wh v_k,\qquad \rho_k\in\mathbb R,}
\eqn{61}{\wh\theta=\sum_{k=0}^{2^n-1}\theta_k\wh v_{2^n+k},\qquad \theta_k\in[0,2\pi).}
The coefficients are computed from
\eqn{62}{\rho_k=\frac1{2^n}\sum_{p=0}^{2^n-1}(-1)^{\langle k,p\rangle}\sigma(k)A_p,}
\eqn{63}{\theta_k=\frac1{2^n}\sum_{p=0}^{2^n-1}(-1)^{\langle 2^n+k,p\rangle}\sigma(2^n+k)\Phi_p\pmod{2\pi}.}

Element (59) is of special type if
\eqn{64}{\rho_1=\rho_2=\cdots=\rho_{2^n-1}=0,\qquad \rho_0\ne0,}
and
\eqn{65}{\theta_k=0\quad\text{if}\quad \bar k_0+\bar k_1+\cdots+\bar k_n\ne1,
\quad 2^n+k=(\bar k_0\bar k_1\ldots\bar k_n)_{\rm Gray}.}
Thus, $\wh\theta$ has nonzero components only, at most, along the directions $i_0,i_1,\ldots,i_n$. For element (58), condition (64) becomes
\eqn{66}{A_0=A_1=\cdots=A_{2^n-1}\overset{\rm not}=A.}
Under these conditions,
\eqn{67}{\wh w=A\sum_{p=0}^{2^n-1}e^{i_0\Phi_p}\wh e_p.}
It follows from (67) that, under condition (66), $\wh w$ is of special type if and only if
\eqn{68}{\sum_{p=0}^{2^n-1}e^{i_0\Phi_p}\wh e_p
=\exp(i_0\phi_0+i_1\phi_1+\cdots+i_n\phi_n),\qquad \phi_k\in[0,2\pi).}

Using the relation
\eqn{69}{\sum_{k=0}^{n}i_k\phi_k=
\sum_{p=0}^{2^n-1}\left[\phi_0+\sum_{k=1}^{n}(-1)^{p_k}\phi_k\right]i_0\wh e_p,
\qquad p=(p_1p_2\ldots p_n)_2,}
it follows from (68) and (69) that, in representation (58), $\wh w$ is a special-type element if and only if (66) holds and there exist $\phi_k\in\mathbb R$, $k=0,\ldots,n$, such that
\eqn{70}{\Omega_p=\phi_0+\sum_{k=1}^{n}(-1)^{p_k}\phi_k,\qquad p=0,1,\ldots,2^n-1.}
With the matrix notation
\eqn{71}{[\Omega]=[\Omega_0,\Omega_1,\ldots,\Omega_{2^n-1}]^T,}
\eqn{72}{[\phi]=[\phi_0,\phi_1,\ldots,\phi_n]^T,}
\eqn{73}{[\alpha]=(\alpha_{pk})_{\substack{p=0,\ldots,2^n-1\\k=0,\ldots,n}},\qquad
\alpha_{p0}=1,\quad \alpha_{pk}=(-1)^{p_k}\ (k\ge1),}
relation (70) becomes
\eqn{74}{[\Omega]=[\alpha][\phi].}
Since the columns of $[\alpha]$ are orthogonal,
\eqn{75}{[\alpha]^T[\alpha]=2^n I_{n+1},}
and
\eqn{76}{[\phi]=\frac1{2^n}[\alpha]^T[\Omega].}
For (76) to satisfy (74), it is necessary and sufficient that
\eqn{77}{[\Omega]=\frac1{2^n}[\alpha][\alpha]^T[\Omega].}
Of the $2^n$ relations, only $2^n-(n+1)$ are independent. Therefore, element (58) is of special type if
\eqn{78}{A_0=A_1=\cdots=A_{2^n-1}}
and
\eqn{79}{[\Omega]=\frac1{2^n}[\alpha][\alpha]^T[\Omega].}

\renewcommand{\eqprefix}{IV.1.2}
\subsection*{IV.1.2. Symbolic Representations over the Generalized Sobrero Algebra}

Let $S$ be the generalized Sobrero algebra introduced in III.2.1. It is a commutative, associative, unital $R$-algebra of order $2(n+1)$, $n\in\mathbb N$. One of its bases is
\eqn{1}{B=\{\wh\varepsilon^{k},j\wh\varepsilon^{k}\}_{k=0}^{n}.}
The algebra is the direct product of the complex field and the nilpotent algebra of order $n+1$, and the multiplication table of the basis is specified by
\eqn{2}{\wh\varepsilon^{n+1}=\zero,}
\eqn{3}{j^2=-\one.}
The unit $\one$ is $\wh\varepsilon^{0}$. In III.4 it was shown that $S$ is polynomial, generated by the element
\eqn{5}{\wh v=j+\wh\varepsilon,}
whose characteristic and minimal polynomial is
\eqn{4}{p(\lambda)=(\lambda^2+1)^{n+1}.}
Consequently, there is also the power basis
\eqn{6}{B^*=\{(j+\wh\varepsilon)^k\}_{k=0}^{2n+1}.}
The algebra contains the distinguished subalgebras
\eqn{7}{
\begin{aligned}
\mathbb R&=\{\mu\one:\mu\in\mathbb R\},\\
\mathbb C&=\{\mu\one+\lambda j:\mu,\lambda\in\mathbb R\},\\
P&=\left\{\sum_{k=0}^{n}\alpha_k\wh\varepsilon^{k}:\alpha_k\in\mathbb R\right\}.
\end{aligned}}
We shall write simply $\mu$ instead of $\mu\one$, and $j$ instead of $j\wh\varepsilon^0$.

Consider the fixed element
\eqn{8}{\wh d=\delta+\wh\varepsilon+j\omega,\qquad \delta,\omega\in\mathbb R.}

The characteristic polynomial of element (8) is
\eqn{9}{p(\lambda)=\bigl[(\lambda-\delta)^2+\omega^2\bigr]^{n+1},\qquad p(\wh d)=\zero.}
If $\omega\ne0$, it is also the minimal polynomial. According to Theorem II.2.1, operator (8) symbolically represents the solutions of the differential equation
\eqn{10}{\bigl[D^2-2\delta D+(\delta^2+\omega^2)I\bigr]^{n+1}f=0,\qquad Df=\frac{df}{dt}.}
Let $V$ be the solution space. We seek the $R$-linear bijection
\eqn{11}{\Psi:V\longrightarrow S,}
\eqn{12}{\Psi(\dot f)=\wh d\,\Psi(f).}
According to Theorem II.2.3,
\eqn{13}{e^{\wh d t}=e^{\delta t}e^{\wh\varepsilon t}e^{j\omega t}.}
Euler-type relations and nilpotency give
\eqn{14}{e^{\wh\varepsilon t}=\sum_{k=0}^{n}\frac{t^k}{k!}\wh\varepsilon^{k},}
\eqn{15}{e^{j\omega t}=\cos\omega t+j\sin\omega t.}
Upon multiplication,
\eqn{16}{
e^{\wh d t}=\sum_{k=0}^{n}\frac{t^k}{k!}e^{\delta t}\cos\omega t\,\wh\varepsilon^{k}
+\sum_{k=0}^{n}\frac{t^k}{k!}e^{\delta t}\sin\omega t\,j\wh\varepsilon^{k}.}
A symbolic representation is therefore specified by
\eqn{17}{\Psi\!\left(\frac{t^k}{k!}e^{\delta t}\cos\omega t\right)=\wh\varepsilon^{k},\qquad k=0,\ldots,n,}
\eqn{18}{\Psi\!\left(\frac{t^k}{k!}e^{\delta t}\sin\omega t\right)=j\wh\varepsilon^{k},\qquad k=0,\ldots,n.}

Let an arbitrary element be
\eqn{19}{\wh w=\sum_{k=0}^{n}(a_k+b_kj)\wh\varepsilon^{k},\qquad a_k,b_k\in\mathbb R.}
From (17), (18), and the linearity of $\Psi$, one obtains the continuous signal
\eqn{20}{
s(t)=e^{\delta t}\left(\sum_{k=0}^{n}\frac{a_k}{k!}t^k\right)\cos\omega t
+e^{\delta t}\left(\sum_{k=0}^{n}\frac{b_k}{k!}t^k\right)\sin\omega t.}
It is simultaneously amplitude- and phase-modulated. If we set
\eqn{21}{A(t)\cos\theta(t)=e^{\delta t}\sum_{k=0}^{n}\frac{a_k}{k!}t^k,}
\eqn{22}{A(t)\sin\theta(t)=e^{\delta t}\sum_{k=0}^{n}\frac{b_k}{k!}t^k,}
then
\eqn{23}{A(t)=e^{\delta t}\sqrt{\left(\sum_{k=0}^{n}\frac{a_k}{k!}t^k\right)^2+
\left(\sum_{k=0}^{n}\frac{b_k}{k!}t^k\right)^2},}
\eqn{24}{\theta(t)=\operatorname{atan2}\!\left(\sum_{k=0}^{n}\frac{b_k}{k!}t^k,
\sum_{k=0}^{n}\frac{a_k}{k!}t^k\right),}
and
\eqn{25}{s(t)=A(t)\cos[\omega t-\theta(t)].}
Relation (25) explicitly exhibits the character of a harmonic signal with simultaneous amplitude and phase modulation.

Element (19) may also be written in semi-exponential form:
\eqn{26}{\wh w=\sum_{k=0}^{n}A_ke^{j\phi_k}\wh\varepsilon^{k},\qquad A_k\ge0,\quad \phi_k\in[0,2\pi),}
where
\eqn{27}{a_k=A_k\cos\phi_k,\qquad b_k=A_k\sin\phi_k.}

In this case,
\eqn{28}{\Psi^{-1}\!\left(A_ke^{j\phi_k}\wh\varepsilon^{k}\right)
=A_ke^{\delta t}\frac{t^k}{k!}\cos(\omega t-\phi_k).}
The relation follows from (17), (18), and the correspondence
\eqn{29}{\Psi^{-1}\!\left(A_ke^{j\phi_k}\wh\varepsilon^{k}\right)
=\Psi^{-1}(A_k\cos\phi_k\,\wh\varepsilon^{k})
+\Psi^{-1}(A_k\sin\phi_k\,j\wh\varepsilon^{k}).}
From (26),
\eqn{30}{s(t)=\sum_{k=0}^{n}A_ke^{\delta t}\frac{t^k}{k!}\cos(\omega t-\phi_k).}
Signal (30) is the sum of at most $n+1$ signals with the same harmonic carrier, amplitude-modulated according to polynomial--exponential laws.

\textbf{The subalgebra $\mathbb R$.} An element is
\eqn{31}{\wh w=A,\qquad A\in\mathbb R.}
From (17), for $k=0$,
\eqn{32}{s(t)=\Psi^{-1}(A)=Ae^{\delta t}\cos\omega t.}
It represents a harmonic signal with exponential amplitude modulation and zero initial phase, and satisfies the second-order equation with polynomial
\eqn{33}{q(\lambda)=\lambda^2-2\delta\lambda+\delta^2+\omega^2,}
which divides polynomial (9).

\textbf{The subalgebra $\mathbb C$.} An element is written as
\eqn{34}{\wh w=Ae^{j\phi},\qquad A\ge0,\quad \phi\in[0,2\pi).}
From (28), for $k=0$,
\eqn{35}{s(t)=Ae^{\delta t}\cos(\omega t-\phi).}
It satisfies the equation with characteristic polynomial (33).

\textbf{The subalgebra $P$.} An element is written as
\eqn{36}{\wh w=\sum_{k=0}^{n}\alpha_k\wh\varepsilon^{k},\qquad \alpha_k\in\mathbb R.}
From (17),
\eqn{37}{s(t)=e^{\delta t}\left(\sum_{k=0}^{n}\frac{\alpha_k}{k!}t^k\right)\cos\omega t.}
It has polynomial--exponential amplitude modulation and a harmonic carrier, and is a solution of the order-$2(n+1)$ equation with polynomial (9).

\textbf{Zero divisors.} In representation (19), the zero divisors are characterized by
\eqn{38}{a_0=b_0=0,}
and, in representation (26), by
\eqn{39}{A_0=0.}
From (20) and (30), a zero divisor represents a degenerate signal, which is a solution of a differential equation of order at most $2n$:
\eqn{40}{
s(t)=e^{\delta t}\left(\sum_{k=1}^{n}\frac{a_k}{k!}t^k\right)\cos\omega t
+e^{\delta t}\left(\sum_{k=1}^{n}\frac{b_k}{k!}t^k\right)\sin\omega t,}
or, equivalently,
\eqn{41}{s(t)=\sum_{k=1}^{n}A_ke^{\delta t}\frac{t^k}{k!}\cos(\omega t-\phi_k).}

\textbf{Special-type elements.} These were introduced in III.2.1 and are characterized, in representations (19) and (26), respectively, by the conditions
\eqn{42}{a_k=c\,b_k,\qquad c\in\mathbb R,\quad k=0,\ldots,n,}
and, respectively,
\eqn{43}{\phi_0=\phi_1=\cdots=\phi_n=\phi.}
From (30), the special-type element
\eqn{44}{\wh w=\sum_{k=0}^{n}A_ke^{j\phi}\wh\varepsilon^{k}}
represents the signal
\eqn{45}{s(t)=e^{\delta t}\left(\sum_{k=0}^{n}\frac{A_k}{k!}t^k\right)\cos(\omega t-\phi).}
Special-type elements represent signals with a harmonic carrier and pure amplitude modulation of polynomial--exponential type.

An arbitrary element (19) can be decomposed into the sum of two special-type elements:
\eqn{46}{\wh w=\sum_{k=0}^{n}a_ke^{j0}\wh\varepsilon^{k}
+\sum_{k=0}^{n}b_ke^{j\pi/2}\wh\varepsilon^{k}.}
It represents the sum of two signals of the form (45), with carriers in quadrature:
\eqn{47}{s(t)=e^{\delta t}\left(\sum_{k=0}^{n}\frac{a_k}{k!}t^k\right)\cos\omega t
+e^{\delta t}\left(\sum_{k=0}^{n}\frac{b_k}{k!}t^k\right)\sin\omega t.}

\renewcommand{\eqprefix}{IV.2}

\section*{IV.2. Identification of Nonlinear Dynamical Systems Using Walsh Algebras}

Truncated Walsh series, also called Walsh polynomials, form an algebra isomorphic to the algebra introduced in III.3. On the basis of its properties, we propose a procedure for the parametric identification of a class of nonlinear dynamical systems with analytic nonlinearities.

\subsection*{IV.2.1. The Algebra of Truncated Walsh Series}

Rademacher [137] introduced an incomplete system of orthogonal functions in the set of real-valued integrable functions on $[0,1]$:
\eqn{1}{u_n:[0,1]\to\{-1,1\},\qquad u_n(t)=\operatorname{sgn}[\sin(2^n\pi t)],\quad n\in\mathbb N.}
Walsh [164] proposed a complete orthogonal system of integrable functions, defined on $\mathbb R$, periodic with period 1, and taking values in $\{-1,1\}$. This system is expressed in terms of the Rademacher functions [2]. If $p=(p_{n-1}\ldots p_0)_2$, the Walsh function of order $p$, in Hadamard ordering [12], is the periodic extension of the function
\eqn{2}{w_p=u_{n-1}^{p_{n-1}}u_{n-2}^{p_{n-2}}\cdots u_0^{p_0}.}
The Walsh functions, together with ordinary multiplication, form an Abelian group [81], with
\eqn{3}{w_pw_q=w_{p\oplus q},\qquad p,q\in\mathbb N,}
where $p\oplus q$ is the bitwise sum modulo 2, without carry [1]. The set $\{w_p\}_{p=0}^{2^n-1}$ is a subgroup of order $2^n$.

For an integrable 1-periodic function $f$, its generalized Walsh--Fourier series is
\eqn{4}{\widetilde f=\sum_{p=0}^{\infty}c_pw_p,\qquad c_p=\int_0^1f(t)w_p(t)\,dt.}

The series converges uniformly to $f$ at points of continuity, whereas at discontinuities convergence holds in the mean [12]. We consider series truncated to $2^n$ terms, called Walsh polynomials [163]:
\eqn{5}{f=\sum_{p=0}^{2^n-1}c_pw_p,\qquad c_p\in\mathbb R.}
Since $u_k^2=1$, the set $W$ of polynomials (5), endowed with the usual operations, is an algebra isomorphic to the Walsh algebra of III.3. It has order $2^n$ and is associative, commutative, polynomial, and semisimple. Its two bases are
\eqn{6}{B=\{w_p\}_{p=0}^{2^n-1},}
\eqn{7}{B^*=\{e_p\}_{p=0}^{2^n-1}.}
The basis of orthogonal idempotents is
\eqn{8}{e_p=\frac1{2^n}\prod_{k=0}^{n-1}\left[1+(-1)^{p_k}u_k\right],
\qquad p=(p_{n-1}\ldots p_0)_2,}
and satisfies
\eqn{9}{e_pe_q=\begin{cases}0,&p\ne q,\\e_p,&p=q.\end{cases}}
The functions $e_p$ are ``block-pulse'' functions [12,81]. The change of basis is
\eqn{10}{\boldsymbol w=H\boldsymbol e,}
\eqn{11}{\boldsymbol e=\frac1{2^n}H\boldsymbol w,}
where
\eqn{12}{\boldsymbol w=[w_0,w_1,\ldots,w_{2^n-1}]^T,}
\eqn{13}{\boldsymbol e=[e_0,e_1,\ldots,e_{2^n-1}]^T,}
\eqn{14}{H=\begin{bmatrix}1&1\\1&-1\end{bmatrix}^{\!\otimes n}.}

The Kronecker product in (14) yields the Hadamard matrix. Equivalently [1,12],
\eqn{15}{H=(h_{pq})_{p,q=0}^{2^n-1},\qquad h_{pq}=(-1)^{p_0q_0+\cdots+p_{n-1}q_{n-1}}.}
The matrix is symmetric and nonsingular, and
\eqn{15a}{H^{-1}=\frac1{2^n}H.}

\subsection*{IV.2.2. A Symbolic Representation over Walsh Algebras}

For $f\in W$, the unary operation
\eqn{16}{\dot f\overset{\rm not}=\int_0^t f(\tau)\,d\tau}
is well defined in the sense of truncated projection. We seek a representation
\eqn{17}{\Psi:W\to W,\qquad \Psi\!\left[\int_0^t f(\tau)\,d\tau\right]=P\Psi(f),}
where $P$ does not depend on $f$. The following relation is used in the identification of linear and bilinear systems [53,95]:
\eqn{18}{\int_0^t\boldsymbol w(x)\,dx=P_{2^n}\boldsymbol w.}
The operational matrix $P_{2^n}$ is recursively defined by
\eqn{19}{P_1=\frac12,\qquad
P_m=\begin{bmatrix}P_{m/2}&-\dfrac1{2m}I_{m/2}\\[2mm]
\dfrac1{2m}I_{m/2}&0_{m/2}\end{bmatrix},\qquad m=2^n.}
Thus, (17) exists if the matrix of $P$ in the basis $B$ is (19). To simplify powers and inversions, we shall use the block-pulse basis.

The operational matrix in the basis $B^*$ satisfies [152]
\eqn{20}{\int_0^t\boldsymbol e(x)\,dx=Q_{2^n}\boldsymbol e,}
where
\eqn{21}{Q_{2^n}=\frac1{2^n}
\begin{bmatrix}
\frac12&1&1&\cdots&1\\
0&\frac12&1&\cdots&1\\
\vdots&\ddots&\ddots&\ddots&\vdots\\
0&\cdots&0&\frac12&1\\
0&\cdots&\cdots&0&\frac12
\end{bmatrix}.}

From (10), (18), and (20), the matrices are similar:
\eqn{22}{P_{2^n}=H^{-1}Q_{2^n}H,}
hence
\eqn{23}{P_{2^n}=\frac1{2^n}HQ_{2^n}H.}
They have the same characteristic and minimal polynomial,
\[
\chi(\lambda)=\left(\lambda-\frac1{2^{n+1}}\right)^{2^n}.
\]
From (22), for every $k\in\mathbb Z$ for which the operation is defined,
\eqn{24}{P_{2^n}^{k}=H^{-1}Q_{2^n}^{k}H.}
The form of $Q_{2^n}$ makes the algebraic computations particularly simple.

\subsection*{IV.2.3. A Procedure for Identifying Nonlinear Dynamical Systems}

\textbf{1.} Consider the system
\eqn{25}{(S):\quad \ddot q+f(q)=u(t),\qquad q(0)=\dot q(0)=0.}
The known input signal $u(t)$ and output signal $q(t)$ are integrable and periodic with period 1; the unknown function $f$ is analytic at the origin. Equation (25) is equivalent to
\eqn{26}{q(t)+\int_0^t\!dx\int_0^x f(q(y))\,dy
=\int_0^t\!dx\int_0^x u(y)\,dy.}
Let the truncated Walsh series be
\[
q=\sum_{k=0}^{2^n-1}q_kw_k,\qquad u=\sum_{k=0}^{2^n-1}u_kw_k,
\]
and let the vectors be
\eqn{27}{\boldsymbol q=[q_0,q_1,\ldots,q_{2^n-1}]^T,}
\eqn{28}{\boldsymbol u=[u_0,u_1,\ldots,u_{2^n-1}]^T,}
with
\eqn{29}{q_k=\int_0^1q(t)w_k(t)\,dt,}
\eqn{30}{u_k=\int_0^1u(t)w_k(t)\,dt.}

Applying the operational representation twice to the function $u$ gives
\eqn{31}{\int_0^t\!dx\int_0^x u(y)\,dy=\boldsymbol u^TP_{2^n}^{2}\boldsymbol w.}
The analyticity of $f$ and the properties of the idempotents give
\eqn{32}{f(q)=\sum_{k=0}^{2^n-1}f(\boldsymbol q^T\boldsymbol h_k)e_k,}
where $\boldsymbol h_k$ is column $k$ of the Hadamard matrix. By (20),
\eqn{33}{\int_0^t\!dx\int_0^x f(q(y))\,dy
=\boldsymbol f^TQ_{2^n}^{2}\boldsymbol e,}
with
\eqn{34}{\boldsymbol f^T=\left[f(\boldsymbol q^T\boldsymbol h_0),f(\boldsymbol q^T\boldsymbol h_1),\ldots,
f(\boldsymbol q^T\boldsymbol h_{2^n-1})\right].}
From (11) and (22),
\eqn{35}{\int_0^t\!dx\int_0^x f(q(y))\,dy
=\frac1{2^n}\boldsymbol f^TP_{2^n}^{2}\boldsymbol w.}
Integral equation (26) becomes
\eqn{36}{\boldsymbol q^T\boldsymbol w+\frac1{2^n}\boldsymbol f^TP_{2^n}^{2}\boldsymbol w
=\boldsymbol u^TP_{2^n}^{2}\boldsymbol w.}
The independence of the components of $\boldsymbol w$ leads to a linear system that determines the vector $\boldsymbol f$, and hence the values of $f$ at $2^n$ specific points. The function is then estimated by one of the classical procedures [61]. If $f$ is a polynomial of degree at most $2^n-1$, interpolation at the $2^n$ points is exact.

\textbf{2.} The procedure extends to the system
\eqn{37}{q^{(m)}+a_1q^{(m-1)}+\cdots+a_{m-1}\dot q+a_mq+f(q)=u,\qquad m\ge2,}
where $f$ is analytic at the origin. We impose zero initial conditions:
\eqn{38}{q^{(m-1)}(0)=q^{(m-2)}(0)=\cdots=\dot q(0)=q(0)=0.}

Integrating (37) $m$ times gives
\eqn{39}{
q(t)+\sum_{k=1}^{m}a_k\underbrace{\int_0^t\!dx_1\int_0^{x_1}\!dx_2\cdots
\int_0^{x_{k-1}}q(y)\,dy}_{k\ \mathrm{integrations}}
+\underbrace{\int_0^t\!dx_1\cdots\int_0^{x_{m-1}}f(q(y))\,dy}_{m\ \mathrm{integrations}}
=\underbrace{\int_0^t\!dx_1\cdots\int_0^{x_{m-1}}u(y)\,dy}_{m\ \mathrm{integrations}}.}
For the truncated series
\eqn{40}{u=\boldsymbol u^T\boldsymbol w,}
\eqn{41}{q=\boldsymbol q^T\boldsymbol w,}
the symbolic representation gives
\eqn{42}{\underbrace{\int_0^t\!dx_1\cdots\int_0^{x_{k-1}}q(y)\,dy}_{k}
=\boldsymbol q^TP_{2^n}^{k}\boldsymbol w,
\qquad k=0,1,\ldots,m,}
and
\eqn{43}{\underbrace{\int_0^t\!dx_1\cdots\int_0^{x_{m-1}}f(q(y))\,dy}_{m}
=\frac1{2^n}\boldsymbol f^TP_{2^n}^{m}\boldsymbol w.}
Consequently,
\eqn{44}{\boldsymbol q^T\boldsymbol w+a_1\boldsymbol q^TP\boldsymbol w+\cdots
+a_m\boldsymbol q^TP^m\boldsymbol w+\frac1{2^n}\boldsymbol f^TP^m\boldsymbol w
=\boldsymbol u^TP^m\boldsymbol w,}
or, compactly,
\eqn{45}{\boldsymbol q^TF(P)\boldsymbol w+\frac1{2^n}\boldsymbol f^TP^m\boldsymbol w
=\boldsymbol u^TP^m\boldsymbol w,}
where
\eqn{46}{F(P)=I+a_1P+\cdots+a_mP^m.}

\emph{Remark.} The powers of the operational matrix are computed from
\eqn{47}{P^k=H^{-1}Q^kH.}
The matrix $Q$ is written as
\eqn{48}{Q=\frac1{2^{n+1}}I_{2^n}+\Delta,}
where
\eqn{49}{\Delta=\frac1{2^n}
\begin{bmatrix}
0&1&1&\cdots&1\\
0&0&1&\cdots&1\\
\vdots&\ddots&\ddots&\ddots&\vdots\\
0&\cdots&0&0&1\\
0&\cdots&\cdots&0&0
\end{bmatrix}.}
The matrix $\Delta$ is strictly upper triangular and nilpotent of index $2^n$:
\eqn{50}{\Delta^{2^n}=0.}
For $m\in\mathbb N$,
\eqn{51}{Q^m=\sum_{k=0}^{\min(m,2^n-1)}\binom{m}{k}
\left(\frac1{2^{n+1}}\right)^{m-k}\Delta^k.}
The independence of the components of $\boldsymbol w$ transforms (45) into a linear system that determines the $2^n$ values of $f$. If $f$ is a polynomial of degree at most $2^n-1$, identification by interpolation is exact.

\renewcommand{\eqprefix}{IV.3}
\section*{IV.3. Response of Dynamical Systems to Excitations Symbolically Representable over Algebras}
\subsection*{IV.3.1. Statement of the Problem}

Let $A$ be a commutative, associative, unital $R$-algebra of finite order $n$. If $\wh d\in A$ has characteristic and minimal polynomial
\eqn{1}{r(\lambda)=\lambda^n+c_1\lambda^{n-1}+\cdots+c_{n-1}\lambda+c_n,\qquad r(\wh d)=\wh0,}
then, according to Theorem III.4.1, the algebra $A$ is polynomial and generated by $r$.

Let $V^*$ be a real vector space endowed with a linear unary operation, denoted by a dot:
\eqn{2}{\dot{\phantom{x}}:V^*\to V^*,}
\eqn{3}{\dot{(\lambda a+\mu b)}=\lambda\dot a+\mu\dot b.}
If the space $V\subset V^*$ of solutions of the equation
\eqn{4}{f^{(n)}+c_1f^{(n-1)}+\cdots+c_{n-1}\dot f+c_nf=0}
has dimension $n$, then exact symbolic representations exist:
\eqn{5}{\Psi:V\to A,}
\eqn{6}{\Psi(\dot f)=\wh d\,\Psi(f).}
The representation is a vector-space isomorphism satisfying property (6).

\emph{Remarks.} 1$^\circ$ If $V^*$ is the space of indefinitely differentiable real-valued functions on an interval, the operation is ordinary differentiation, and the representations of continuous signals from II.2 are obtained.

2$^\circ$ If $V^*$ is the space of real-valued functions defined on an interval of $\mathbb Z$, and the operation is the right shift, the representations of discrete signals from II.3 are obtained.

\subsection*{IV.3.2. The Transmittance of Dynamical Systems}

\textbf{Definition.} A dynamical system on $V^*$, with input $f$ and output $x$, is the object described by the operational equation
\eqn{7}{
a_0x^{(m)}+a_1x^{(m-1)}+\cdots+a_{m-1}\dot x+a_mx
=b_0f^{(p)}+b_1f^{(p-1)}+\cdots+b_{p-1}\dot f+b_pf,}
where $m>p$, $m,p\in\mathbb N^*$. For differentiable real-valued functions, (7) describes linear time-invariant analog systems [23,151]; for discrete signals and the shift operator, it describes digital filters [1,23].

The general solution is
\eqn{8}{x=x_{\rm om}+x_{\rm par},}
where $x_{\rm om}$ satisfies
\eqn{9}{a_0x^{(m)}+a_1x^{(m-1)}+\cdots+a_{m-1}\dot x+a_mx=0,}
and $x_{\rm par}$ is a particular solution. If $f$ is symbolically represented over $A$ by the operator $\wh d$, we seek $x_{\rm par}$ in the same class. Applying $\Psi$ to equation (7) gives
\eqn{10}{
(a_0\wh d^{m}+a_1\wh d^{m-1}+\cdots+a_m)\Psi(x)
=(b_0\wh d^{p}+b_1\wh d^{p-1}+\cdots+b_p)\Psi(f).}
Let
\eqn{11}{p(\lambda)=\sum_{k=0}^{m}a_{m-k}\lambda^k,}
\eqn{12}{q(\lambda)=\sum_{k=0}^{p}b_{p-k}\lambda^k.}
Then
\eqn{13}{p(\wh d)\Psi(x)=q(\wh d)\Psi(f).}
If $p(\wh d)$ is invertible in $A$,
\eqn{14}{\Psi(x)=\frac{q(\wh d)}{p(\wh d)}\Psi(f).}
The element
\eqn{15}{H(\wh d)=\frac{q(\wh d)}{p(\wh d)}}
is called the \emph{transmittance} of dynamical system (7) in the sense of the symbolic representation. The element
\eqn{16}{x=\Psi^{-1}\!\left[H(\wh d)\Psi(f)\right]}
is a particular solution of equation (7).

Form (16) permits the following cases to be distinguished.

a) If $H(\wh d)=\wh0$, then $x=0$: the system completely filters out the input. This occurs if $q(\lambda)$ is divisible by the characteristic polynomial of $\wh d$.

b) If the minimal polynomial of the element $H(\wh d)\Psi(f)$ has degree strictly less than $n$, the response satisfies an equation of order strictly less than $n$: the system partially filters the input. This occurs, for example, when $H(\wh d)$ is a zero divisor.

c) If $\Psi(f)$ and $H(\wh d)$ are special-type elements, their product is also of special type. The response has the same type as the excitation, and the system does not distort its form.

d) If $H(\wh d)\Psi(f)$ is of special type whereas $\Psi(f)$ is not, the system corrects the distortions of the excitation.

e) If $p(\wh d)$ is zero or a zero divisor, resonance occurs; in these cases $x_{\rm par}$ cannot be represented over $A$ by the operator $\wh d$.

\medskip
\textbf{Theorem IV.3.2.1.} The response of system (7) can be symbolically represented over the polynomial algebra generated by
\eqn{17}{\Pi(\lambda)=p(\lambda)r(\lambda).}

\emph{Proof.} The component $x_{\rm om}$ satisfies homogeneous equation (9), whose characteristic polynomial is $p$. The component $x_{\rm par}$, represented over $A$ by the operator $\wh d$, satisfies an equation with characteristic polynomial $r$. The sum $x=x_{\rm om}+x_{\rm par}$ is annihilated by the product $p(\lambda)r(\lambda)$ and, according to Theorem II.1.2, is representable over the polynomial algebra generated by $\Pi$. $\square$

\emph{Remarks.} 1$^\circ$ If the polynomials $p$ and $r$ are relatively prime, the characteristic polynomial of the element representing the response is minimal. If they have common roots, it is no longer minimal, and the response may be represented over an algebra whose order is smaller than the degree of $\Pi$.

2$^\circ$ For analog systems, the response is unique once the initial conditions
\eqn{18}{x(0)=x_0,\quad \dot x(0)=x_1,\quad\ldots,\quad x^{(m-1)}(0)=x_{m-1}.}
are specified. These conditions select a particular solution $x_{\rm om}$ from the solution class of the homogeneous equation.

3$^\circ$ For discrete systems, the response is unique once the conditions
\eqn{19}{x(0)=x_0,\quad x(1)=x_1,\quad\ldots,\quad x(m-1)=x_{m-1}.}

are specified. 4$^\circ$ In the traditional study of analog and discrete systems of the form (7), the Laplace transfer function and the $Z$ transfer function, respectively, are used [23,151]. These transforms algebraize differential equations and difference equations, respectively. Transmittance (15) has similar properties. Its complete relationship with the Laplace and $Z$ transfer functions will be established after introducing integral transforms over algebras that contain the complex field as a subalgebra.

The Laplace and $Z$ transforms are complex-valued functions of a complex variable; transmittance (15) is a hypercomplex-valued function of a hypercomplex variable.

\chaptertoc{CHAPTER V}{Integral Transforms on \texorpdfstring{$\mathbb R$}{R}-Algebras}
\renewcommand{\eqprefix}{V.1}

\begin{center}{\Large\bfseries CHAPTER V}\vspace{1em}

{\large\bfseries Integral Transforms over $R$-Algebras}\end{center}

\section*{V.1. Statement of the Problem}

The Fourier transform [23,151] associates with a signal in
\[
L_1=\left\{f:\mathbb R\to\mathbb R:\int_{-\infty}^{\infty}|f(t)|\,dt<\infty\right\}
\]
a complex-valued function of a real variable:
\eqn{1}{\mathcal Ff(i\omega)=\int_{-\infty}^{\infty}f(t)e^{-i\omega t}\,dt,
\qquad \omega\in\mathbb R,\quad i^2=-1.}
The mapping $f\mapsto\mathcal Ff$ is a symbolic representation. Under the usual conditions, the inversion formula is
\eqn{2}{\frac{f(t+0)+f(t-0)}2=\frac1{2\pi}\int_{-\infty}^{\infty}\mathcal Ff(i\omega)e^{i\omega t}\,d\omega.}
It expresses the possibility of synthesizing signals in $L_1$ from harmonic signals representable over the complex field by the operator $\wh d=i\omega$.

The Laplace transform [23,151] associates with an original signal $f$ the function
\eqn{3}{\mathcal Lf(s)=\int_0^{\infty}f(t)e^{-st}\,dt,
\qquad s=\delta+i\omega.}

Under specified conditions, the Laplace transform admits the inverse
\eqn{4}{\frac{f(t+0)+f(t-0)}2=\frac1{2\pi i}\int_{\delta-i\infty}^{\delta+i\infty}\mathcal Lf(s)e^{st}\,ds.}
The formula indicates the possibility of synthesizing a broad class of signals from exponentially modulated harmonic signals representable over $\mathbb C$ by the operator $\wh d=\delta+i\omega$.

These considerations motivate the study of the integral transforms
\eqn{5}{\FT f\overset{\rm def}=\int_{-\infty}^{\infty}f(t)e^{-\wh d t}\,dt,}
where $\wh d$ belongs to a commutative, associative, unital algebra $A$ of finite order over $\mathbb R$. Such a transform associates with the signal a hypercomplex-valued function of the hypercomplex variable $\wh d$.

\renewcommand{\eqprefix}{V.2}
\section*{V.2. Integral Transforms over the Algebra $\A$}
\subsection*{V.2.1. Definition. General Properties}

Let $\A$ be the algebra introduced in III.1, with bases
\eqn{6}{B=\{\wh v_p\}_{p=0}^{2^{n+1}-1},\qquad
\wh v_p=i_0^{p_0}i_1^{p_1}\cdots i_n^{p_n},\quad p=(p_0\ldots p_n)_{\rm Gray},}
and
\eqn{7}{B^*=\{\wh e_p,i_0\wh e_p\}_{p=0}^{2^n-1},\qquad
\wh e_p=\frac1{2^n}\prod_{k=1}^{n}\left[\wh1-(-1)^{p_k}i_0i_k\right].}
The algebra contains the subalgebras
\eqn{8}{\mathbb R=\{\lambda\wh1:\lambda\in\mathbb R\},}
\eqn{9}{\mathbb C=\{\lambda\wh1+\mu i_0:\lambda,\mu\in\mathbb R\}.}

Hereafter, we identify the real and complex fields with subalgebras (8) and (9), omitting the unit from the notation.

Let $L_1=L_1(\mathbb R)$, and let $\mathcal H$ be the set of functions of $n+1$ real variables with values in $\A$. For
\eqn{10}{\wh\omega=\sum_{k=0}^{n}\omega_ki_k,\qquad \omega_k\in\mathbb R,}
we define the mapping
\eqn{11}{\FT:L_1\to\mathcal H,\qquad
\FT f=\int_{-\infty}^{\infty}f(t)e^{-\wh\omega t}\,dt.}

\textbf{Theorem V.2.1.} Transform (11) exists for every $f\in L_1$ and is monogenic as a hypercomplex-valued function of the variable $\wh\omega\in\A$.

\emph{Proof.} The properties of the exponential and Euler's formulas give
\eqn{12}{e^{-\wh\omega t}=\prod_{k=0}^{n}e^{-\omega_kti_k}
=\prod_{k=0}^{n}(\cos\omega_kt-i_k\sin\omega_kt).}
For $p=(p_0\ldots p_n)_{\rm Gray}$, let
\eqn{13}{\varphi_p(t)=\prod_{k=0}^{n}\cos\!\left(\omega_kt+p_k\frac\pi2\right),}
so that
\eqn{14}{e^{-\wh\omega t}=\sum_{p=0}^{2^{n+1}-1}\varphi_p(t)\wh v_p.}
Consequently,
\eqn{15}{\FT f=\sum_{p=0}^{2^{n+1}-1}\left[\int_{-\infty}^{\infty}f(t)\varphi_p(t)\,dt\right]\wh v_p.}
Since $|\varphi_p(t)|\le1$, all $2^{n+1}$ integrals converge for $f\in L_1$.

Setting
\eqn{16}{F_p(\omega_0,\ldots,\omega_n)=\int_{-\infty}^{\infty}f(t)\varphi_p(t)\,dt,}
the transform has the expansion
\eqn{17}{\FT f(\wh\omega)=\sum_{p=0}^{2^{n+1}-1}F_p(\omega_0,\ldots,\omega_n)\wh v_p.}
To establish monogenicity, we compute the exterior product [140]
\eqn{18}{\FT f\wedge d\wh\omega=
\int_{-\infty}^{\infty}f(t)\bigl(e^{-\wh\omega t}\wedge d\wh\omega\bigr)\,dt.}
The monogenicity of the exponential implies
\eqn{19}{e^{-\wh\omega t}\wedge d\wh\omega=\wh0,}
whence
\eqn{20}{\FT f\wedge d\wh\omega=\wh0.\qquad\square}

\emph{Remarks.} 1$^\circ$ Monogenicity implies differential relations among the components in (17). According to Theorem 1 in [140, §4.2, p.25],
\eqn{21}{i_0\frac{\partial\FT f}{\partial\omega_k}
=i_k\frac{\partial\FT f}{\partial\omega_0},\qquad k=0,1,\ldots,n.}

2$^\circ$ In the basis $B^*$, variable (10) is written as
\eqn{22}{\wh\omega=\sum_{p=0}^{2^n-1}\Omega_pi_0\wh e_p,}
where
\eqn{23}{\Omega_p=\omega_0+\sum_{k=1}^{n}(-1)^{p_k}\omega_k,\qquad
p=(p_1\ldots p_n)_2.}
The properties of the orthogonal idempotents give
\eqn{24}{e^{-\wh\omega t}=\sum_{p=0}^{2^n-1}e^{-i_0\Omega_pt}\wh e_p.}

By (24), definition (11) becomes
\eqn{25}{\FT f=\sum_{p=0}^{2^n-1}\left[\int_{-\infty}^{\infty}f(t)e^{-i_0\Omega_pt}\,dt\right]\wh e_p.}
Let the complex Fourier transform be
\eqn{26}{\mathcal Ff(i_0\omega)=\int_{-\infty}^{\infty}f(t)e^{-i_0\omega t}\,dt.}
Then the unique expression of transform (11) in the basis $B^*$ is
\eqn{27}{\FT f=\sum_{p=0}^{2^n-1}\mathcal Ff(i_0\Omega_p)\wh e_p.}

\textbf{Theorem V.2.2.} Integral transform (11) is $\mathbb R$-linear and injective, and satisfies
\eqn{28}{\FT(f^{(k)})=\wh\omega^{\,k}\FT f,\qquad k\in\mathbb N,}
for $f,f^{(k)}\in L_1$, as well as
\eqn{29}{\FT(f*g)=\FT f\,\FT g,\qquad f,g\in L_1.}

\emph{Proof.} Using (27) and the properties of the Fourier transform,
\eqn{30}{
\begin{aligned}
\FT(f^{(k)})
&=\sum_{p=0}^{2^n-1}(i_0\Omega_p)^k\mathcal Ff(i_0\Omega_p)\wh e_p\\
&=\left[\sum_{p=0}^{2^n-1}(i_0\Omega_p)^k\wh e_p\right]
\left[\sum_{p=0}^{2^n-1}\mathcal Ff(i_0\Omega_p)\wh e_p\right]
=\wh\omega^{\,k}\FT f.
\end{aligned}}

Similarly,
\eqn{31}{
\begin{aligned}
\FT(f*g)
&=\sum_{p=0}^{2^n-1}\mathcal F(f*g)(i_0\Omega_p)\wh e_p\\
&=\sum_{p=0}^{2^n-1}\mathcal Ff(i_0\Omega_p)\mathcal Fg(i_0\Omega_p)\wh e_p
=\FT f\,\FT g.
\end{aligned}}
Linearity follows from the linearity of the integral, whereas injectivity follows from the injectivity of the complex Fourier transform and the uniqueness of decomposition (27). $\square$

\emph{Remarks.} 1$^\circ$ Theorem V.2.2 shows that (11) is a symbolic representation of the functions in $L_1$, with operator $\wh\omega$.

2$^\circ$ Relation (27) permits the extension of the transform to tempered distributions $\mathcal S'$ [151]. Under this extension,
\eqn{32}{\FT\delta=\wh1,}
\eqn{33}{\FT\delta^{(k)}=\wh\omega^{\,k},\qquad k\in\mathbb N.}
Indeed, the properties of the distributional Fourier transform are
\eqn{34}{\mathcal F\delta=1,}
\eqn{35}{\mathcal F\delta^{(k)}=(i_0\omega)^k.}
Relation (33) shows that the natural powers of the representation operator are the images of the derivatives of the Dirac distribution.

We say that $f\in L_{1,\mathrm{loc}}$ satisfies condition $(*)$ at a point $t_0$ if there exists $\alpha>0$ such that
\eqn{*}{\lim_{\omega\to\infty}\int_0^{\alpha}
[f(t_0+u)-2f(t_0)+f(t_0-u)]\frac{\sin\omega u}{u}\,du=0.}
It is sufficient that $f$ have bounded variation on $[t_0-\alpha,t_0+\alpha]$ and that
$f(t_0)=[f(t_0+)+f(t_0-)]/2$; in particular, the condition is satisfied at points of differentiability.

\textbf{Theorem V.2.3.} If $f$ satisfies condition $(*)$ at the point $t$, transform (11) is invertible at that point, and
\eqn{36}{\widetilde f(t)=\FT^{-1}(\FT f)(t)
=\frac1{2\pi i_0}\int_{\A}\FT f(\wh\omega)e^{\wh\omega t}\,d\wh\omega,}
where $\widetilde f(t)=[f(t+)+f(t-)]/2$.

\emph{Proof.} In the basis $B^*$, using (22), (24), and (27), together with the convention of componentwise integration,
\eqn{37}{
\begin{aligned}
\frac1{2\pi i_0}\int_{\A}\FT f\,e^{\wh\omega t}\,d\wh\omega
&=\frac1{2\pi}\sum_{p=0}^{2^n-1}
\left[\int_{-\infty}^{\infty}\mathcal Ff(i_0\Omega_p)e^{i_0\Omega_pt}\,d\Omega_p\right]\wh e_p.
\end{aligned}}
The Fourier inversion formula gives
\eqn{38}{\widetilde f(t)=\frac1{2\pi}\int_{-\infty}^{\infty}\mathcal Ff(i_0\omega)e^{i_0\omega t}\,d\omega.}
Consequently,
\eqn{39}{\frac1{2\pi i_0}\int_{\A}\FT f\,e^{\wh\omega t}\,d\wh\omega
=\widetilde f(t)\sum_{p=0}^{2^n-1}\wh e_p=\widetilde f(t).\qquad\square}

\emph{Remark.} Formula (36) expresses the synthesis of signals satisfying $(*)$ from quasiperiodic signals with simultaneous amplitude and phase modulation, representable over $\A$ by operator (10).

\subsection*{V.2.2. Energy Aspects}

For expansion (17), we define the conjugate
\eqn{40}{(\FT f)^*=\sum_{p=0}^{2^{n+1}-1}F_p(\omega_0,\ldots,\omega_n)\wh v_p^{\,*},}
where
\eqn{41}{\wh v_p^{\,*}=\begin{cases}\wh v_p,&0\le p<2^n,\\-\wh v_p,&2^n\le p<2^{n+1}.\end{cases}}

\textbf{Theorem V.2.4.} The energy of the signal $f\in L_1\cap L_2$ is completely determined by its transform through Parseval's relation
\eqn{42}{\int_{-\infty}^{\infty}|f(t)|^2\,dt
=\frac1{2\pi i_0}\int_{\A}\FT f\,(\FT f)^*\,d\wh\omega.}

\emph{Proof.} Let $f_1,f_2\in L_1\cap L_2$, $\FT_1=\FT f_1$, and $\FT_2=\FT f_2$. Consider
\eqn{43}{I=\int_{-\infty}^{\infty}f_1(t)f_2(t)\,dt.}
Application of the inversion formula and decomposition with respect to the orthogonal idempotents reduce the identity, component by component, to the complex Parseval relation.

\renewcommand{\eqprefix}{V.2}

Using Theorem V.2.3, we successively obtain
\begin{align}\tag{\eqprefix.44}
\int_{-\infty}^{\infty}f_1(t)f_2(t)\,dt
&=\frac1{2\pi i_0}\int_{\A}(\FT f_1)(\FT f_2)^*\,d\wh\omega.
\end{align}
The following relation is proved analogously:
\eqn{45}{\int_{-\infty}^{\infty}f_1(t)f_2(t)\,dt
=\frac1{2\pi i_0}\int_{\A}(\FT f_1)^*\,\FT f_2\,d\wh\omega.}
For $f_1=f_2=f$, the assertion of the theorem follows from relation (45).

\textbf{Remarks.} 1$^\circ$. Let $L=L_1\cap L_2$, and let $\FT L$ be the space of images of functions in $L$ under integral transform (11). We define
\eqn{46}{H:\FT L\times\FT L\to\mathbb R,\qquad
H(\FT_1,\FT_2)\overset{\rm def}=\frac1{2\pi i_0}\int_{\A}\FT_1(\FT_2)^*\,d\wh\omega.}
It follows from the definition that
\eqn{47}{H(\FT_1,\FT_2)=H(\FT_2,\FT_1),\quad H(\FT_1,\FT_1)\ge0,\quad
H(\FT_1,\FT_1)=0\Longleftrightarrow\FT_1=\wh0.}
Consequently, $H$ is an inner product on $\FT L$. For the induced norm, let
\eqn{48}{\|\FT\|\overset{\rm def}=\sqrt{H(\FT,\FT)}.}
Theorem V.2.4 can be written as
\eqn{49}{\int_{-\infty}^{\infty}|f(t)|^2\,dt=\|\FT f\|^2.}

2$^\circ$. The preceding considerations, relation (45), and Plancherel's theorem [151] yield the following result.

\textbf{Corollary.} The restriction of integral transform (11), $\FT:L\to\FT L$, is a vector-space isomorphism that preserves inner products.

\subsection*{V.2.3. Relationship with the Transmittance of Dynamical Systems}
\textbf{Theorem V.2.5.} The transmittance of a linear dynamical system over the algebra $\A$, with operator (10), is obtained by applying transform (11) to the weighting function of the system.

\emph{Proof.} The input--output relation of a linear invariant dynamical system is of convolution type [151]:
\eqn{50}{x=h*f.}
The function (distribution) $h$ is the response of the system to the Dirac impulse $\delta$ and is called the weighting function. Applying integral transform (11) and Theorem V.2.2 gives
\eqn{51}{\FT x=\FT h\cdot\FT f.}
Taking into account the definition of transmittance and the fact that $\FT$ is a symbolic representation with operator $\wh\omega$, we obtain
\eqn{52}{\wh H=\FT h.}

Let $H(s)$ be the Laplace transfer function. From (52) and (27), it follows that
\eqn{53}{\wh H(\wh\omega)=\sum_{p=0}^{2^n-1}H(i_0\Omega_p)\ep.}
Here $\wh\omega$ is given by (10), whereas $\Omega_p\in\mathbb R$ are defined by (23).

\textbf{Application.} Consider a dynamical system with Laplace transfer function $H(s)$. We determine the conditions under which its steady-state response to the input
\eqn{54}{f(t)=A\cos(\omega_0t-\varphi_0)\cos(\omega_1t-\varphi_1)\cos(\omega_2t-\varphi_2)}
is a signal without parasitic phase modulation, where $A,\omega_k>0$, $\varphi_k\in[0,2\pi)$.

The signal is represented over the algebra $\bigotimes C^3$ by the operator
\eqn{55}{\wh\omega=\omega_0i_0+\omega_1i_1+\omega_2i_2}
by the special-type element
\eqn{56}{\wh f=Ae^{i_0\varphi_0+i_1\varphi_1+i_2\varphi_2}.}
If the transmittance of the system is $\wh H$, the steady-state response is represented by
\eqn{57}{\wh x=\wh H\cdot\wh f.}
The signal $x(t)$ has pure amplitude modulation if and only if $\wh H$ is a special-type element. According to (53),
\eqn{58}{\wh H(\wh\omega)=\sum_{p=0}^{3}H(i_0\Omega_p)\ep,}
where
\eqn{59}{\Omega_0=\omega_0+\omega_1+\omega_2,\quad
\Omega_1=\omega_0+\omega_1-\omega_2,\quad
\Omega_2=\omega_0-\omega_1+\omega_2,\quad
\Omega_3=\omega_0-\omega_1-\omega_2.}

Let
\eqn{60}{H(i_0\Omega_p)=H_pe^{i_0\alpha_p},\qquad H_p\in\mathbb R_+,\quad
\alpha_p\in[0,2\pi),\quad p=0,1,2,3.}
According to IV.1.1, (64)--(65), transmittance (58) is a special-type element if
\eqn{61}{H_0=H_1=H_2=H_3\overset{\rm not}=H,}
\eqn{62}{\alpha_0+\alpha_3=\alpha_1+\alpha_2\overset{\rm not}=2\Psi_0,}
\eqn{63}{\alpha_0-\alpha_2=\alpha_1-\alpha_3\overset{\rm not}=2\Psi_1,}
\eqn{64}{\alpha_0-\alpha_1=\alpha_2-\alpha_3\overset{\rm not}=2\Psi_2.}
Under relations (61)--(64),
\eqn{65}{\wh H=He^{i_0\Psi_0+i_1\Psi_1+i_2\Psi_2}.}
From (56), (57), and (65), it follows that
\eqn{66}{\wh x=HAe^{i_0(\varphi_0+\Psi_0)+i_1(\varphi_1+\Psi_1)+i_2(\varphi_2+\Psi_2)}.}
Therefore, the response of the system is
\begin{equation}\tag{\eqprefix.67}
\begin{split}
x(t)=HA&\cos(\omega_0t-\varphi_0-\Psi_0)
\cos(\omega_1t-\varphi_1-\Psi_1)\\[-2pt]
&\times\cos(\omega_2t-\varphi_2-\Psi_2),\qquad t\in\mathbb R,
\end{split}
\end{equation}
and exhibits no parasitic phase modulation.

\renewcommand{\eqprefix}{V.3}

\section*{V.3. Integral Transforms over the Generalized Sobrero Algebra}
\subsection*{V.3.1. Definition. General Properties}
Let $S$ be the real, commutative, and associative generalized Sobrero algebra of dimension $2(n+1)$, with basis
\eqn{1}{B=\{\eps^{\,k},\jhat\eps^{\,k}\}_{k=0}^{n}.}
The generators satisfy
\eqn{2}{\eps^{\,n+1}=\wh0,\qquad \eps^{\,n}\ne\wh0,}
\eqn{3}{\jhat^{\,2}=-\wh1.}
The unit of the algebra is $\wh1=\eps^{\,0}$. The real and complex fields are identified with the subalgebras
\eqn{4}{\mathbb R=\{\mu\wh1\}_{\mu\in\mathbb R},}
\eqn{5}{\mathbb C=\{\mu\wh1+\lambda\jhat\}_{\lambda,\mu\in\mathbb R}.}
Hereafter, we omit the unit when writing scalars.

Let $\Omega$ denote the space of real-valued, locally integrable, piecewise-continuous functions supported in $[0,+\infty)$ and of at most exponential growth:
\eqn{6}{\Omega=\{f:\mathbb R\to\mathbb R\mid f\in L^1_{\rm loc},\ \operatorname{supp}f\subset[0,+\infty),\ |f(t)|\le Me^{at}\}.}
For
\eqn{7}{\wh\omega=\delta+\eps+\jhat\omega,\qquad \delta,\omega\in\mathbb R,}
we define the integral transform
\eqn{8}{\Lh f(\wh\omega)=\int_0^{\infty}f(t)e^{-\wh\omega t}\,dt.}

\textbf{Theorem V.3.1.} Transform (8) exists for every $f\in\Omega$ of exponential growth and is monogenic in the hypercomplex variable $\wh\omega$.

\emph{Proof.} In the complex subalgebra of $S$, set
\eqn{9}{\wh\omega=s+\eps,\qquad s=\delta+j\omega\in\mathbb C.}
The nilpotency of $\eps$ gives
\eqn{10}{e^{-\wh\omega t}=e^{-st}e^{-\eps t}=e^{-st}\sum_{k=0}^{n}\frac{(-1)^kt^k}{k!}\eps^{\,k}.}
Consequently,
\begin{equation}\tag{\eqprefix.11}\begin{split}
\Lh f&=\int_0^{\infty}\sum_{k=0}^{n}\frac{(-1)^kt^k}{k!}\eps^{\,k}f(t)e^{-st}\,dt\\
&=\sum_{k=0}^{n}\frac{\eps^{\,k}}{k!}\int_0^{\infty}(-1)^kt^kf(t)e^{-st}\,dt.
\end{split}\end{equation}
The classical Laplace transform is
\eqn{12}{\mathcal Lf(s)=\int_0^{\infty}f(t)e^{-st}\,dt.}
For a function of exponential growth $|f(t)|\le Me^{ct}$, its domain of convergence is
\eqn{13}{\Sigma(c)=\{s\in\mathbb C\mid \operatorname{Re}s>c\}.}
In this half-plane,
\eqn{14}{\Sigma(c)\subset C_f,}
\eqn{15}{\mathcal Lf:C_f\to\mathbb C\ \text{is holomorphic},}
and differentiation under the integral sign gives
\eqn{16}{\frac{d^k}{ds^k}\mathcal Lf(s)=\int_0^{\infty}(-1)^kt^kf(t)e^{-st}\,dt.}

With the notation $\mathcal Lf^{(k)}=d^k(\mathcal Lf)/ds^k$, relation (11) becomes
\eqn{17}{\Lh f(\wh\omega)=\sum_{k=0}^{n}\frac{\eps^{\,k}}{k!}\frac{d^k}{ds^k}\mathcal Lf(s).}
This formula proves the existence of the transform and, through the generalized Cauchy--Riemann conditions, its monogenicity.

\textbf{Remark.} Let $a_f$ be the abscissa of convergence of the Laplace transform:
\eqn{18}{a_f=\inf\{a\in\mathbb R\mid \exists s\in C_f, a=\operatorname{Re}s\}.}
It follows from (17) that the same abscissa of convergence applies to $\Lh f$.

For a direct verification of monogenicity, computation of the exterior product yields
\eqn{19}{\Lh f\wedge d\wh\omega
=\int_0^{\infty}f(t)\bigl(e^{-\wh\omega t}\wedge d\wh\omega\bigr)\,dt=\wh0.}

\textbf{Theorem V.3.2.} Integral transform (8) is an $\mathbb R$-linear mapping and has the properties
\eqn{20}{\Lh(f^{(k)})=\wh\omega^{\,k}\Lh f-\sum_{p=0}^{k-1}\wh\omega^{\,k-1-p}f^{(p)}(0),\qquad k\in\mathbb N,}
\eqn{21}{\Lh(f*g)=\Lh f\cdot\Lh g.}

\emph{Proof.} Linearity is immediate. For $f\in\Omega$, relation (17) gives
\eqn{22}{\Lh\dot f=\sum_{k=0}^{n}\frac{\eps^{\,k}}{k!}\frac{d^k}{ds^k}\mathcal L\dot f.}
The classical Laplace transform satisfies
\eqn{23}{\mathcal L\dot f=s\mathcal Lf-f(0).}
Substituting (23) into (22) and using (17), we obtain
\eqn{24}{\Lh\dot f=(s+\eps)\Lh f-f(0)=\wh\omega\Lh f-f(0).}
Relation (20) follows by induction.

For convolution, relation (17) gives
\eqn{25}{\Lh(f*g)=\sum_{k=0}^{n}\frac{\eps^{\,k}}{k!}\frac{d^k}{ds^k}\mathcal L(f*g).}
The classical formula is
\eqn{26}{\mathcal L(f*g)=\mathcal Lf\cdot\mathcal Lg.}
Thus
\eqn{27}{\frac{d^k}{ds^k}\mathcal L(f*g)=\frac{d^k}{ds^k}(\mathcal Lf\cdot\mathcal Lg),}
and Leibniz's formula gives
\eqn{28}{\frac{d^k}{ds^k}(\mathcal Lf\cdot\mathcal Lg)=\sum_{p=0}^{k}\binom{k}{p}(\mathcal Lf)^{(p)}(\mathcal Lg)^{(k-p)}.}

Consequently,
\eqn{29}{\Lh(f*g)=\sum_{k=0}^{n}\frac{\eps^{\,k}}{k!}\frac{d^k}{ds^k}(\mathcal Lf\cdot\mathcal Lg).}
On the other hand, multiplying expansions (17) gives
\begin{align}\tag{\eqprefix.30}
\Lh f\cdot\Lh g
&=\left(\sum_{p=0}^{n}\frac{\eps^{\,p}}{p!}(\mathcal Lf)^{(p)}\right)
\left(\sum_{q=0}^{n}\frac{\eps^{\,q}}{q!}(\mathcal Lg)^{(q)}\right).
\end{align}
Since $\eps^{\,m}=\wh0$ for $m\ge n+1$, collecting terms by powers of $\eps$ and using Leibniz's formula yield
\eqn{31}{\Lh f\cdot\Lh g=\sum_{k=0}^{n}\frac{\eps^{\,k}}{k!}\frac{d^k}{ds^k}(\mathcal Lf\cdot\mathcal Lg).}
Relation (21) follows from (29) and (31).

\textbf{Remark.} Relation (17) permits the extension of transform (8) to distributions. For $s=\alpha+j\omega$, $\alpha>a_f$, the Laplace transform is related to the Fourier transform of $g(t)=e^{-\alpha t}f(t)$ by
\eqn{32}{\mathcal Lf(s)=\mathcal Fg(\omega).}
Under the extension to $\mathcal D'$, one obtains
\eqn{33}{\Lh\delta=\wh1,}
\eqn{34}{\Lh\delta^{(k)}=\wh\omega^{\,k},\qquad k\in\mathbb N.}
Indeed,
\eqn{35}{\mathcal L\delta=1,}
\eqn{36}{\mathcal L\delta^{(k)}=s^k.}

Relations (17), (35), and (36) imply (33) and (34). Thus, the natural powers of the operator $\wh\omega$ are the images of the derivatives of the Dirac impulse.

\textbf{Theorem V.3.3.} Integral transform (8) is invertible for every $f\in\Omega$ of exponential growth, and the inversion formula is
\eqn{37}{f(t)=\frac1{2\pi j}\int_{\alpha-j\infty}^{\alpha+j\infty}\Lh f(\wh\omega)e^{\wh\omega t}\,d\wh\omega,\qquad \alpha>a_f.}

\emph{Proof.} For $\wh\omega=s+\eps$ and $d\wh\omega=ds$, the right-hand side is
\eqn{38}{\frac1{2\pi j}\int_{\alpha-j\infty}^{\alpha+j\infty}\Lh f(s+\eps)e^{(s+\eps)t}\,ds.}
Using (17),
\begin{align}\tag{\eqprefix.39}
\frac1{2\pi j}\int_{\alpha-j\infty}^{\alpha+j\infty}
\left[\sum_{k=0}^{n}\frac{\eps^{\,k}}{k!}(\mathcal Lf)^{(k)}(s)\right]
e^{st}\left[\sum_{q=0}^{n}\frac{t^q}{q!}\eps^{\,q}\right]ds.
\end{align}
The classical inversion formula and its derivatives give
\eqn{40}{(-t)^kf(t)=\frac1{2\pi j}\int_{\alpha-j\infty}^{\alpha+j\infty}(\mathcal Lf)^{(k)}(s)e^{st}\,ds.}

Substituting (40) into (39), the coefficient of $\eps^{\,m}$ is
\[
f(t)\sum_{k=0}^{m}\frac{(-t)^kt^{m-k}}{k!(m-k)!}
=\frac{f(t)t^m}{m!}(1-1)^m,
\]
hence all nilpotent terms vanish, while the zeroth-order term is $f(t)$. Consequently,
\eqn{41}{\frac1{2\pi j}\int_{\alpha-j\infty}^{\alpha+j\infty}\Lh f(s+\eps)e^{(s+\eps)t}\,ds=f(t),}
which proves (37).

\textbf{Remark.} The inversion formula exhibits the possibility of synthesizing signals with pure amplitude modulation and polynomial phase modulation.

\subsection*{V.3.2. Relationship with the Transmittance of Dynamical Systems}
\textbf{Theorem V.3.4.} The transmittance of a linear dynamical system over the algebra $S$, with operator (7), is obtained by applying transform (8) to the weighting function of the system.

\emph{Proof.} The input--output relation is
\eqn{42}{x=h*f.}
Applying (8) and Theorem V.3.2 gives
\eqn{43}{\Lh x=\Lh h\cdot\Lh f.}

Since $\Lh$ is a symbolic representation with operator $\wh\omega$, relation (43) gives
\eqn{44}{\wh H=\Lh h.}
The Laplace transfer function is
\eqn{45}{H=\mathcal Lh.}
From (17), (44), and (45), the relationship between the two transmittances is
\eqn{46}{\wh H(\wh\omega)=\sum_{k=0}^{n}\frac{\eps^{\,k}}{k!}\frac{d^kH}{ds^k}(s),\qquad \wh\omega=s+\eps.}

\textbf{Application.} Consider the input with polynomial phase modulation
\eqn{47}{f(t)=\sum_{k=0}^{n}a_k(t)\cos\!\left(\omega t-\varphi_k(t)\right),}
symbolically represented over the algebra $S$ by a special-type element, denoted by
\eqn{48}{\wh f=a_0e^{\eps\varphi_0}+a_1e^{\eps\varphi_1}+\cdots+a_ne^{\eps\varphi_n}.}
The steady-state response is represented by
\eqn{49}{\wh x=\wh H\cdot\wh f.}

Written in form (46), the transmittance $\wh H$ is a special-type element if and only if the derivatives of the transfer function satisfy the structural conditions
\eqn{50}{\frac{1}{k!}\frac{d^kH}{ds^k}=H\,P_k(\Psi_1,\ldots,\Psi_k),\qquad k=1,\ldots,n,}
where $P_k$ are the polynomials determined by the expansion of the nilpotent exponential. Under these conditions, product (49) remains of special type and the response introduces no parasitic modulation.

\textbf{Remark.} Conditions V.2 (61)--(64) and V.3 (50), imposed on the transfer function of a system, constitute criteria for preserving the modulation type of the input signal.

\chaptertoc{CHAPTER VI}{Symbolic Representation of Vector-Valued Functions}
\renewcommand{\eqprefix}{VI.1}

\thispagestyle{empty}
\begin{center}\vspace*{18mm}{\Large\bfseries CHAPTER VI}\par\vspace{6mm}
{\large\bfseries Symbolic Representation of Vector-Valued Functions}\end{center}
\section*{VI.1. Matrix Symbolic Representation of Vector-Valued Functions}
\subsection*{VI.1.1. Statement of the Problem}
Let $V_3$ be the vector space of free vectors, and let $B=\{\vct e_1,\vct e_2,\vct e_3\}$ be a right-handed orthonormal basis. Every $\vct v\in V_3$ has the unique expression
\eqn{1}{\vct v=v_1\vct e_1+v_2\vct e_2+v_3\vct e_3,\qquad v_k\in\mathbb R, k=1,2,3.}
Let $\mathcal V_3$ be the space of real three-row column matrices. One of its elements is
\eqn{2}{v=[v_1,v_2,v_3]^T,\qquad v_k\in\mathbb R.}
The mapping
\eqn{3}{\Psi:V_3\to\mathcal V_3,\qquad \Psi(\vct v)=v}
is a vector-space isomorphism. For fixed $\vct\omega\in V_3$, consider the mapping $\alpha(\vct v)=\vct\omega\times\vct v$.

Explicitly,
\eqn{4}{\alpha:V_3\to V_3,\qquad \alpha(\vct v)=\vct\omega\times\vct v.}
For a nonzero real square matrix $\mat\omega$, define analogously
\eqn{5}{\beta:\mathcal V_3\to\mathcal V_3,\qquad \beta(v)=\mat\omega v.}
We determine the relationship between the vector $\vct\omega$ and the matrix $\mat\omega$ for which (3) is a symbolic representation. If
\eqn{6}{\vct\omega=\omega_1\vct e_1+\omega_2\vct e_2+\omega_3\vct e_3,}
then
\eqn{7}{\vct\omega\times\vct v=(\omega_2v_3-\omega_3v_2)\vct e_1+(\omega_3v_1-\omega_1v_3)\vct e_2+(\omega_1v_2-\omega_2v_1)\vct e_3.}
Consequently,
\eqn{8}{\Psi(\vct\omega\times\vct v)=
\begin{bmatrix}0&-\omega_3&\omega_2\\\omega_3&0&-\omega_1\\-\omega_2&\omega_1&0\end{bmatrix}
\begin{bmatrix}v_1\\v_2\\v_3\end{bmatrix}.}
With the notation
\eqn{9}{\mat\omega=\begin{bmatrix}0&-\omega_3&\omega_2\\\omega_3&0&-\omega_1\\-\omega_2&\omega_1&0\end{bmatrix},}
we obtain
\eqn{10}{\Psi(\vct\omega\times\vct v)=\mat\omega\,v.}

\textbf{Theorem VI.1.1.} The mapping $\Psi:V_3\to\mathcal V_3$, $\Psi(\vct v)=v$, is an isomorphism with the property
\eqn{11}{\Psi(\vct\omega\times\vct v)=\mat\omega\Psi(\vct v),}
and hence an exact symbolic representation of $V_3$, with respect to the cross-product automorphism, over $\mathcal V_3$, with operator $\mat\omega$.

\textbf{Remark.} The matrix $\mat\omega$ is the representation, in the basis $\{\vct e_p\otimes\vct e_q\}_{p,q=1}^3$, of the second-order antisymmetric tensor associated with the vector $\vct\omega$. We shall write $\vct\omega=\operatorname{vect}\mat\omega$.

\subsection*{VI.1.2. General Results}
The characteristic polynomial of antisymmetric matrix (9) is
\eqn{12}{p(\lambda)=\det(\mat\omega-\lambda I_3)=-\lambda(\lambda^2+\|\vct\omega\|^2).}
Its roots are $0$ and $\pm i\|\vct\omega\|$. Hence:

\textbf{Theorem VI.1.2.1.} If $\vct\omega\ne0$ and $\vct\omega=\operatorname{vect}\mat\omega$, then
\eqn{13}{e^{-\mat\omega t}=\varphi_0I_3+\varphi_1\mat\omega+\varphi_2\mat\omega^2,}
where
\eqn{14}{\varphi_0=1,\qquad \varphi_1=-\frac{\sin(\|\vct\omega\|t)}{\|\vct\omega\|},\qquad
\varphi_2=\frac{1-\cos(\|\vct\omega\|t)}{\|\vct\omega\|^2}.}
By the Cayley--Hamilton theorem,
\eqn{15}{\mat\omega^3=-\|\vct\omega\|^2\mat\omega.}

Denoting by $\Psi^{-1}$ the inverse of correspondence (11), we have
\eqn{16}{\mat\omega v\ \xmapsto{\Psi^{-1}}\ \vct\omega\times\vct v,}
\eqn{17}{\mat\omega^2v\ \xmapsto{\Psi^{-1}}\ \vct\omega\times(\vct\omega\times\vct v),}
\eqn{18}{\mat\omega^3v\ \xmapsto{\Psi^{-1}}\ -\|\vct\omega\|^2(\vct\omega\times\vct v).}
For the spaces of vector-valued functions and column-matrix-valued functions, (11) induces
\eqn{19}{\Psi:V_3^*\to\mathcal V_3^*,}
\eqn{20}{\Psi(\vct\omega\times\vct v)=\mat\omega\Psi(\vct v).}

\textbf{Theorem VI.1.2.2.} For every $\vct r_0\in V_3$ and every function $\vct f:I\to V_3$,
\eqn{21}{e^{-\mat\omega t}r_0\ \xmapsto{\Psi^{-1}}\ \varphi_0\vct r_0+\vct\omega\times(\varphi_1\vct r_0)+\vct\omega\times[\vct\omega\times(\varphi_2\vct r_0)],}
\eqn{22}{e^{-\mat\omega t}*f\ \xmapsto{\Psi^{-1}}\ \varphi_0*\vct f+\vct\omega\times(\varphi_1*\vct f)+\vct\omega\times[\vct\omega\times(\varphi_2*\vct f)].}
In these formulas, $*$ denotes convolution, and the functions $\varphi_k$ are given by (14).

\emph{Proof.} From (13),
\eqn{23}{e^{-\mat\omega t}r_0=\varphi_0r_0+\varphi_1\mat\omega r_0+\varphi_2\mat\omega^2r_0.}
Applying $\Psi^{-1}$ and using (16)--(17) yields (21). Similarly,
\eqn{24}{e^{-\mat\omega t}*f=\varphi_0*f+\varphi_1*\mat\omega f+\varphi_2*\mat\omega^2f,}
which yields (22).

In relation (22), the convolution products are understood componentwise:
\eqn{25}{\varphi_k*\vct f=\sum_{p=1}^{3}(\varphi_k*f_p)\vct e_p
=\sum_{p=1}^{3}\left[\int_0^t\varphi_k(\tau)f_p(t-\tau)\,d\tau\right]\vct e_p,\quad k=0,1,2,}
where $\vct f=\sum_{p=1}^3f_p\vct e_p$.

Consider the Cauchy problems
\eqn{26}{\begin{cases}\dot{\vct r}+\vct\omega\times\vct r=\vct f,\\ \vct r(0)=\vct r_0,
\end{cases}}
\eqn{27}{\begin{cases}\ddot{\vct r}+2\vct\omega\times\dot{\vct r}+\vct\omega\times(\vct\omega\times\vct r)=\vct f,\\
\vct r(0)=\vct r_0,\quad\dot{\vct r}(0)=\vct v_0.
\end{cases}}
These problems model motion in non-inertial frames and in gyroscopic force fields.

\textbf{Theorem VI.1.2.3.} The solution of problem (26) is
\eqn{28}{\vct r=\vct f_0+\vct\omega\times\vct f_1+\vct\omega\times(\vct\omega\times\vct f_2),}
where
\[
\vct f_k=\varphi_k\vct r_0+\varphi_k*\vct f,\qquad k=0,1,2,
\]
and the functions $\varphi_k$ are those in (14).

\emph{Proof.} Applying the correspondence $\Psi$ to equation (26), we obtain the following matrix problem.

The matrix problem is
\eqn{29}{\begin{cases}\dot r+\mat\omega r=f,\\r(0)=r_0.
\end{cases}}
Its solution is
\eqn{30}{r=e^{-\mat\omega t}r_0+e^{-\mat\omega t}*f.}
Relations (21) and (22) give formula (28).

\textbf{Theorem VI.1.2.4.} The solution of problem (27) is
\eqn{31}{\vct r=\vct g_0+\vct\omega\times\vct g_1+\vct\omega\times(\vct\omega\times\vct g_2),}
where
\eqn{32}{\vct g_k=\varphi_k[\vct r_0+(\vct v_0+\vct\omega\times\vct r_0)t]+t\varphi_k*\vct f,
\qquad k=0,1,2.}

\emph{Proof.} The matrix form of the problem is
\eqn{33}{\begin{cases}\ddot r+2\mat\omega\dot r+\mat\omega^2r=f,\\r(0)=r_0,\quad\dot r(0)=v_0.
\end{cases}}
It is equivalent to the cascade
\eqn{34}{\begin{cases}\dfrac d{dt}(\dot r+\mat\omega r)+\mat\omega(\dot r+\mat\omega r)=f,\\
(\dot r+\mat\omega r)_{t=0}=v_0+\mat\omega r_0.
\end{cases}}
The first integration yields
\eqn{35}{\dot r+\mat\omega r=e^{-\mat\omega t}(v_0+\mat\omega r_0)+e^{-\mat\omega t}*f.}

The second integration gives
\eqn{36}{r=e^{-\mat\omega t}r_0+e^{-\mat\omega t}*
\left[e^{-\mat\omega t}(v_0+\mat\omega r_0)+e^{-\mat\omega t}*f\right].}
Using the identity
\eqn{37}{e^{-\mat\omega t}*e^{-\mat\omega t}=t,e^{-\mat\omega t},}
we obtain the simplified form
\eqn{38}{r=e^{-\mat\omega t}[r_0+(v_0+\mat\omega r_0)t]+t,e^{-\mat\omega t}*f.}
Applying $\Psi^{-1}$ and relations (21) and (22) yields (31).

\textbf{Remark.} The symbolic representation permits the direct determination of the solutions of problems (26) and (27), without projecting the equations onto privileged frames. The vector forms are useful both for qualitative analysis and for studying dependence on the initial conditions.

\subsection*{VI.1.3. Applications}
\subsubsection*{VI.1.3.1. Motion in a Non-Inertial Frame under a Uniform Gravitational Field}
The problem for the relative velocity of a particle in the vicinity of the Earth is
\eqn{39}{\begin{cases}\dot{\vct v}+2\vct\omega\times\vct v=\vct g,\\\vct v(0)=\vct v_0,
\end{cases}}
where $\vct\omega$ is the angular velocity of the frame and $\vct g$ is the gravitational acceleration. Differentiation gives the acceleration problem
\eqn{40}{\begin{cases}\dot{\vct a}+2\vct\omega\times\vct a=\vct0,\\
\vct a(0)=\vct g+2\vct v_0\times\vct\omega.
\end{cases}}
According to Theorem VI.1.2.3,
\eqn{41}{\vct a=\varphi_0\vct a(0)+\vct\omega\times[\varphi_1\vct a(0)]
+\vct\omega\times\{\vct\omega\times[\varphi_2\vct a(0)]\}.}

For the operator $2\vct\omega\times$, the functions in (14) are
\eqn{42}{\varphi_0=1,\qquad \varphi_1=-\frac{\sin(2\omega t)}{2\omega},\qquad
\varphi_2=\frac{1-\cos(2\omega t)}{4\omega^2},\qquad \omega=\|\vct\omega\|.}
Setting $\vct A=\vct g+2\vct v_0\times\vct\omega$ and applying (41) to the operator $2\vct\omega\times$, we obtain
\eqn{43}{\vct a=\vct A-\frac{\vct\omega\times\vct A}{\omega}\sin(2\omega t)
+\frac{\vct\omega\times(\vct\omega\times\vct A)}{\omega^2}[1-\cos(2\omega t)].}
By the vector triple-product identity,
\eqn{44}{\vct a=\frac{\vct\omega\cdot\vct g}{\omega^2}[1-\cos(2\omega t)]\vct\omega
+\cos(2\omega t)\vct A-\frac{\sin(2\omega t)}{\omega}\vct\omega\times\vct A.}
An antiderivative is
\eqn{45}{\vct v=\frac{\vct\omega\cdot\vct g}{2\omega^3}[2\omega t-\sin(2\omega t)]\vct\omega
+\frac{\sin(2\omega t)}{2\omega}\vct A
+\frac{\cos(2\omega t)}{2\omega^2}\vct\omega\times\vct A+\vct c.}
The initial condition determines
\eqn{46}{\vct c=\frac{\vct g\times\vct\omega}{2\omega^2}
+\frac{\vct v_0\cdot\vct\omega}{\omega^2}\vct\omega.}
Consequently,
\begin{equation}\tag{\eqprefix.47}\begin{split}
\vct v={}&\frac{\vct g\times\vct\omega}{2\omega^2}+
\frac{\vct v_0\cdot\vct\omega}{\omega^2}\vct\omega+
\frac{\sin(2\omega t)}{2\omega}\vct A\\
&+\frac{\cos(2\omega t)}{2\omega^2}\vct\omega\times\vct A+
\frac{\vct\omega\cdot\vct g}{2\omega^3}[2\omega t-\sin(2\omega t)]\vct\omega.
\end{split}\end{equation}
For the position vector,
\eqn{48}{\vct r(t)=\vct r_0+\int_0^t\vct v(\tau)\,d\tau.}

Integrating relation (47), after elementary computations, yields
\begin{equation}\tag{\eqprefix.49}\begin{split}
\vct r={}&\vct r_0+\vct v_0t+
\frac{\vct g+2\vct v_0\times\vct\omega}{4\omega^2}[1-\cos(2\omega t)]\\
&+\frac{\vct\omega\times(\vct g+2\vct v_0\times\vct\omega)}{4\omega^3}
[\sin(2\omega t)-2\omega t]\\
&+\frac{\vct\omega\cdot\vct g}{4\omega^4}
[2\omega^2t^2+\cos(2\omega t)-1]\vct\omega,
\end{split}\end{equation}
where $\omega=\|\vct\omega\|$. This is the exact solution of the problem
\eqn{50}{\begin{cases}\ddot{\vct r}+2\vct\omega\times\dot{\vct r}=\vct g,\\
\vct r(0)=\vct r_0,\quad\dot{\vct r}(0)=\vct v_0.
\end{cases}}
Traditionally, higher-order terms in $\omega$ are neglected. Using
\eqn{51}{1-\cos(2\omega t)\simeq2\omega^2t^2,}
\eqn{52}{\sin(2\omega t)-2\omega t\simeq-\frac{(2\omega t)^3}{3!}=-\frac43\omega^3t^3,}
one obtains the approximation
\eqn{53}{\vct r\simeq\vct r_0+\vct v_0t+\frac{t^2}{2}\vct g
+(\vct v_0\times\vct\omega)t^2+\frac23\vct\omega\times(\vct v_0\times\vct\omega)t^3
+\frac13(\vct\omega\times\vct g)t^3.}
With $\omega^2\simeq0$,
\eqn{54}{\vct r\simeq\vct r_0+\vct v_0t+\frac{t^2}{2}\vct g
+(\vct v_0\times\vct\omega)t^2-\frac23(\vct\omega\cdot\vct v_0)\vct\omega t^3
+\frac13(\vct\omega\times\vct g)t^3.}
For $\vct v_0=\vct0$,
\eqn{55}{\vct r\simeq\vct r_0+\frac{t^2}{2}\vct g+\frac{t^3}{3}\vct\omega\times\vct g.}

If the centrifugal force is also taken into account, the motion is modeled by
\eqn{56}{\begin{cases}\ddot{\vct r}+2\vct\omega\times\dot{\vct r}
+\vct\omega\times(\vct\omega\times\vct r)=\vct g,\\
\vct r(0)=\vct r_0,\quad\dot{\vct r}(0)=\vct v_0.
\end{cases}}
This is a problem of type (27). Evaluating the convolutions in (31) and (32) gives the exact solution
\begin{equation}\tag{\eqprefix.57}\begin{split}
\vct r={}&\vct r_0+\vct v_0t+
\frac{\vct\omega\cdot\vct g}{2\omega^2}\vct\omega t^2
-\frac{\sin(\omega t)}{\omega}\vct\omega\times(\vct r_0^*+\vct v_0t)\\
&+\frac{1-\cos(\omega t)}{\omega^2}
\vct\omega\times[\vct\omega\times(\vct r_0^*+\vct v_0t)]\\
&+\omega t\cos(\omega t)\frac{\vct\omega\times\vct r_0^*}{\omega}
+\omega t\sin(\omega t)\frac{\vct\omega\times(\vct\omega\times\vct r_0^*)}{\omega^2},
\end{split}\end{equation}
where $\vct r_0^*=\vct r_0+\vct g/\omega^2$. Solution (57) had not been reported in the literature cited in the text.

\subsubsection*{VI.1.3.2. Motion of Charges in a Uniform Electromagnetic Field}
For a particle of charge $q$ and mass $m$, in uniform fields $\vct E$ and $\vct B$, the velocity equation is
\eqn{58}{\begin{cases}m\dot{\vct v}=q(\vct v\times\vct B+\vct E),\\\vct v(0)=\vct v_0.
\end{cases}}
Let
\eqn{59}{\vct\omega\overset{\rm not}=\frac qm\vct B,}
\eqn{60}{\vct\alpha\overset{\rm not}=\frac qm\vct E.}

Problem (58) becomes
\eqn{61}{\begin{cases}\dot{\vct v}+\vct\omega\times\vct v=\vct\alpha,\\\vct v(0)=\vct v_0.
\end{cases}}
Applying general solution (47), with $\omega=qB/m$, gives
\begin{equation}\tag{\eqprefix.62}\begin{split}
\vct v={}&\frac{\vct E\times\vct B}{B^2}+
\frac{\vct v_0\cdot\vct B}{B^2}\vct B+
\frac{\sin(\omega t)}{B}(\vct v_0\times\vct B+\vct E)\\
&+\frac{\cos(\omega t)}{B^2}\vct B\times(\vct v_0\times\vct B+\vct E)
+\frac{\vct E\cdot\vct B}{B^3}[\omega t-\sin(\omega t)]\vct B.
\end{split}\end{equation}
For perpendicular fields, $\vct E\cdot\vct B=0$, this reduces to
\begin{equation}\tag{\eqprefix.63}\begin{split}
\vct v={}&\frac{\vct E\times\vct B}{B^2}+
\frac{\vct v_0\cdot\vct B}{B^2}\vct B+
\frac{\sin(\omega t)}{B}(\vct v_0\times\vct B+\vct E)\\
&+\frac{\cos(\omega t)}{B^2}\vct B\times(\vct v_0\times\vct B+\vct E).
\end{split}\end{equation}

\textbf{Remarks.} If $\vct E\cdot\vct B\ne0$, the velocity hodograph is a helix of variable pitch wound around an elliptic cylinder whose axis is parallel to $\vct B$. If $\vct E\cdot\vct B=0$, the hodograph is a circle in a plane perpendicular to $\vct B$.

Denoting the initial acceleration by $\vct a_0=(q/m)(\vct v_0\times\vct B+\vct E)$, differentiation of (62) gives
\eqn{64}{\vct a=\frac{\vct a_0\cdot\vct B}{B^2}\vct B
+\frac{\sin(\omega t)}B\vct a_0\times\vct B
+\frac{\cos(\omega t)}{B^2}\vct B\times(\vct a_0\times\vct B).}

Relation (64) shows that the acceleration has constant magnitude; its hodograph is a circle lying in a plane perpendicular to the magnetic-field lines. The acceleration vector undergoes uniform precession about the direction $\vct B$, with Larmor angular velocity
\[
\vct\omega_L=-\frac qm\vct B.
\]

Denoting the initial force by $\vct F_0=q(\vct v_0\times\vct B+\vct E)$ and setting $\vct F=m\vct a$, relation (64) gives
\eqn{65}{\vct F=\frac{\vct F_0\cdot\vct B}{B^2}\vct B
+\frac{\sin(\omega t)}B\vct F_0\times\vct B
+\frac{\cos(\omega t)}{B^2}\vct B\times(\vct F_0\times\vct B).}
The magnitude of the resultant force remains constant:
\eqn{66}{\|\vct F(t)\|=\|\vct F_0\|=q\|\vct v_0\times\vct B+\vct E\|,
\qquad t\ge0.}
Its projection onto the direction of the magnetic field is also constant:
\eqn{67}{\operatorname{pr}_{\vct B}\vct F=
\frac{\vct F_0\cdot\vct B}{B^2}\vct B.}
The variable component is
\[
\vct F^{*}=\frac{\sin(\omega t)}B\vct F_0\times\vct B
+\frac{\cos(\omega t)}{B^2}\vct B\times(\vct F_0\times\vct B).
\]

\begin{center}
\includegraphics[width=.42\textwidth]{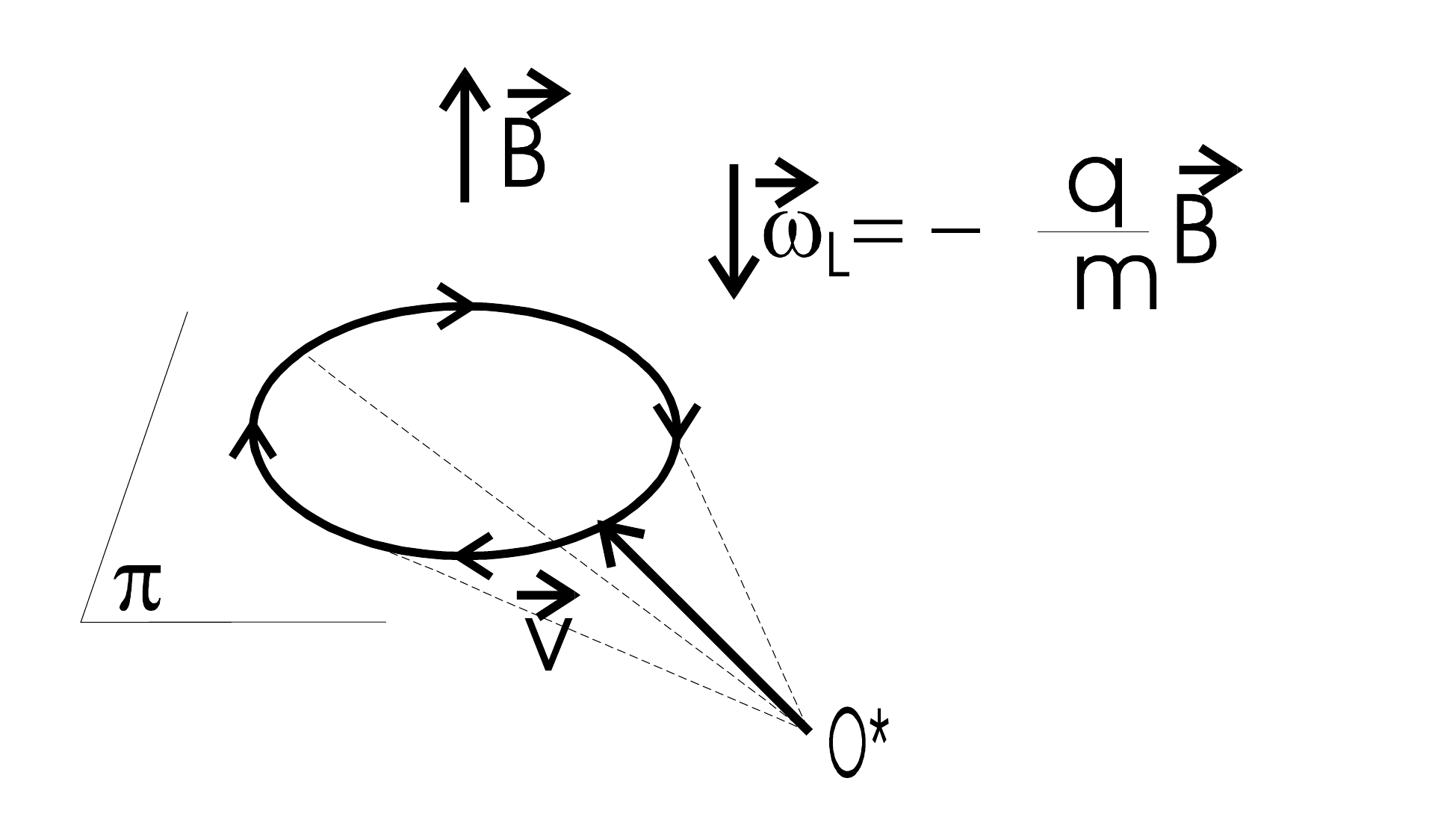}

\small\textnormal{Fig. 1}

\medskip
\begin{minipage}{.47\textwidth}\centering
\includegraphics[width=\linewidth]{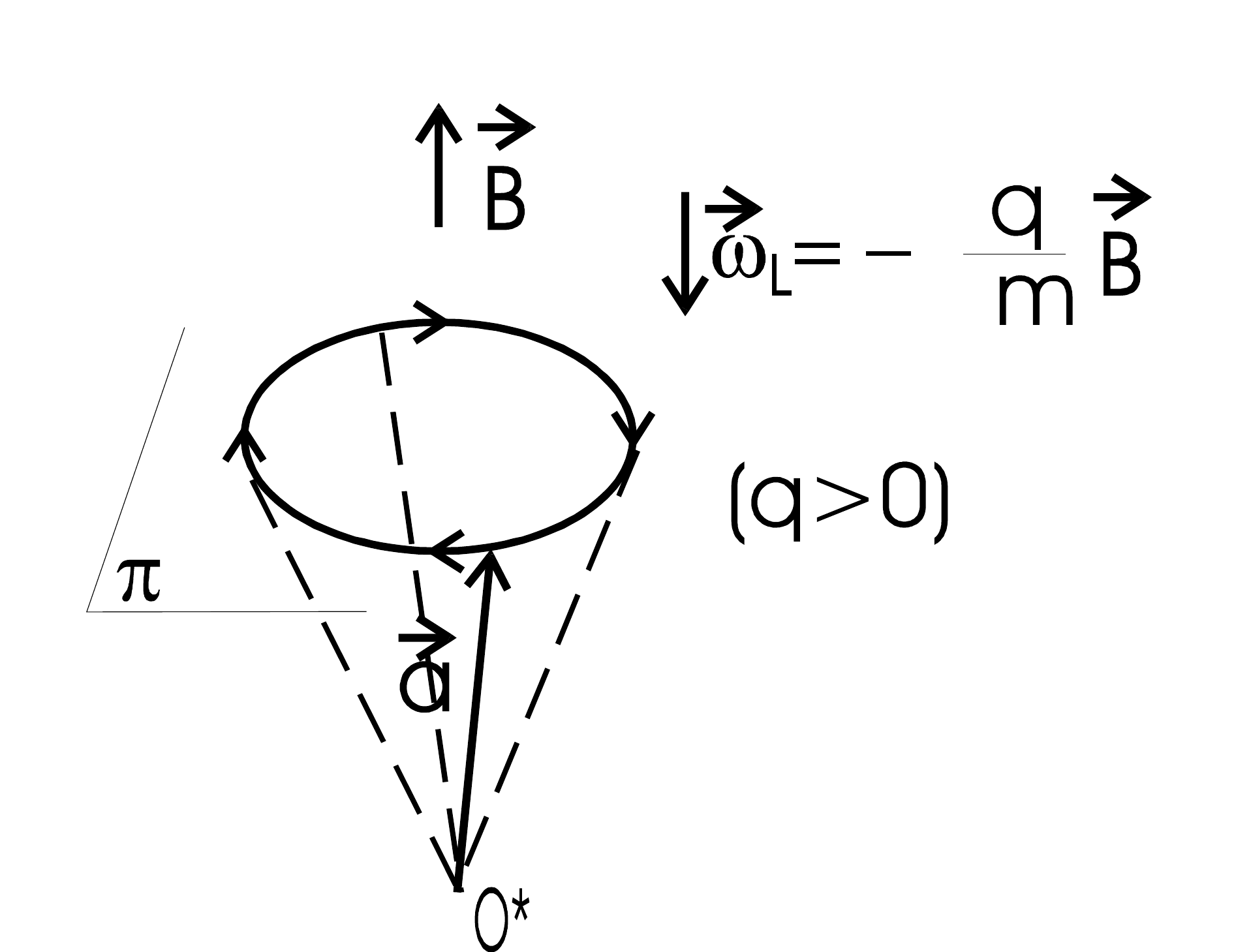}

\small\textnormal{Fig. 2}
\end{minipage}\hfill
\begin{minipage}{.47\textwidth}\centering
\includegraphics[width=\linewidth]{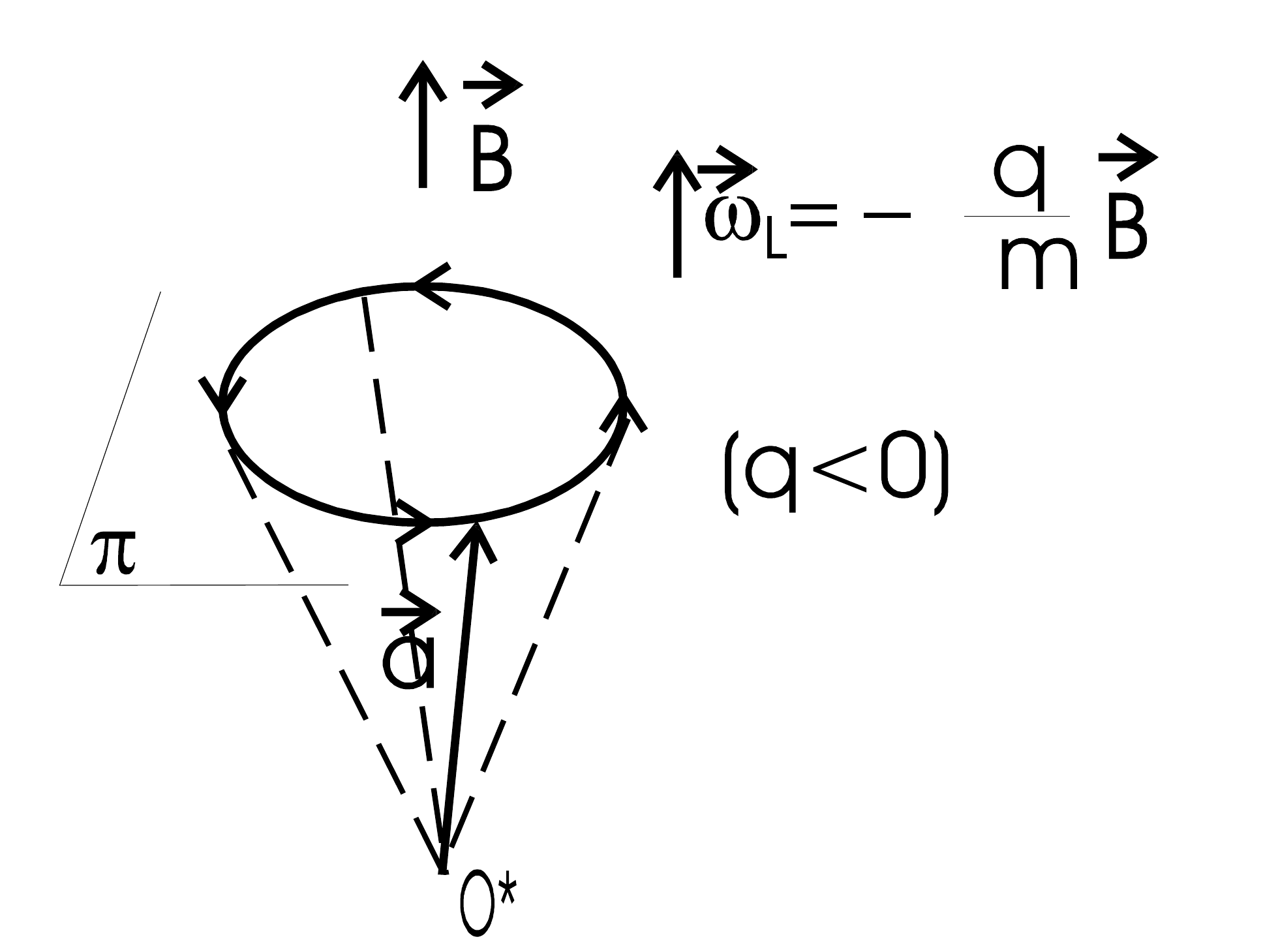}

\small\textnormal{Fig. 3}
\end{minipage}
\end{center}

The variable component of the force, introduced on the preceding page, is a rotating vector with angular velocity
\[
\vct\omega_L=-\frac qm\vct B,
\]
of constant magnitude and with direction perpendicular to the magnetic-field lines. Thus, the resultant force has constant magnitude and undergoes uniform precession about the field lines.

\begin{center}
\begin{minipage}{.47\textwidth}\centering
\includegraphics[width=\linewidth]{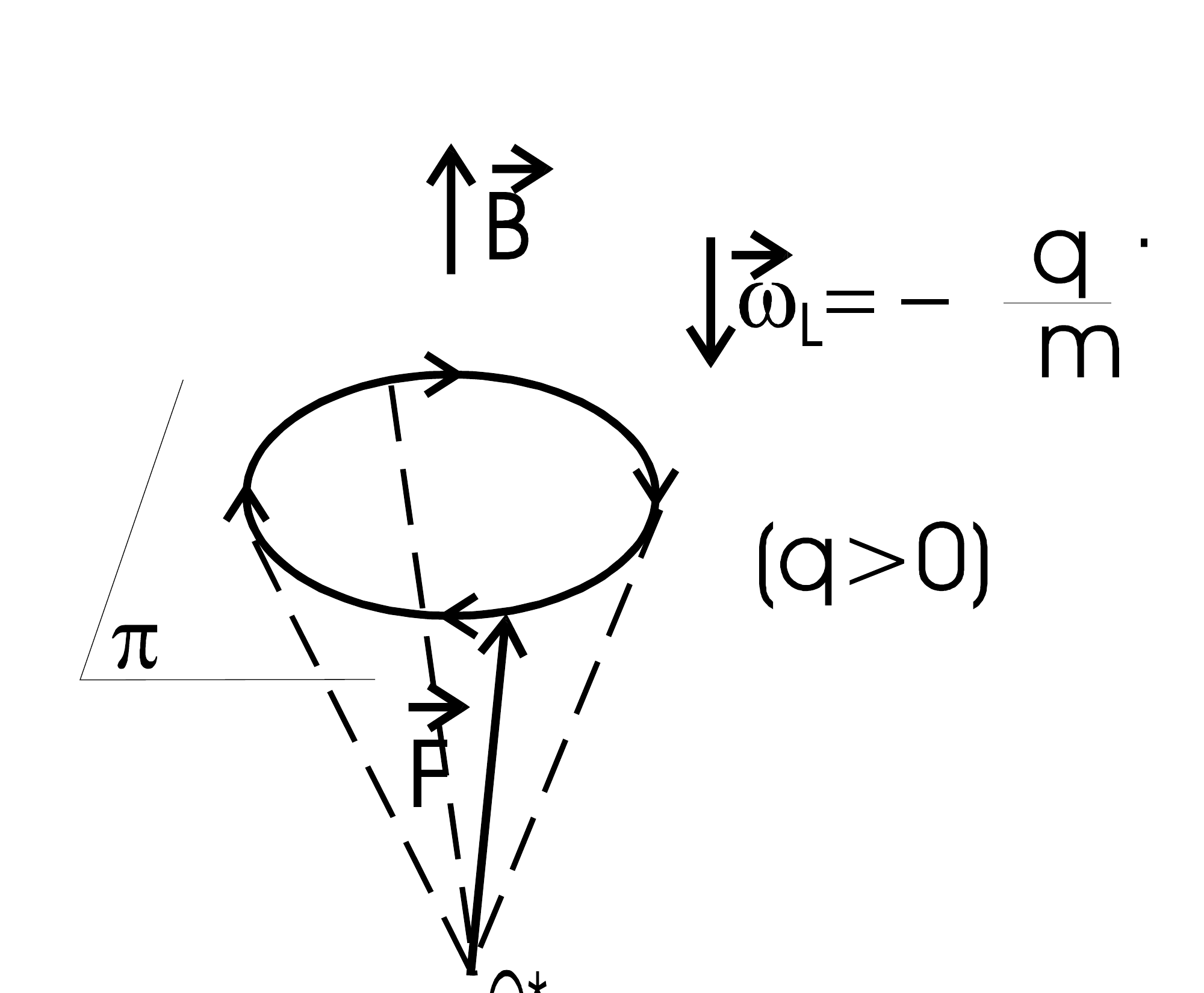}

\small\textnormal{Fig. 4}
\end{minipage}\hfill
\begin{minipage}{.47\textwidth}\centering
\includegraphics[width=\linewidth]{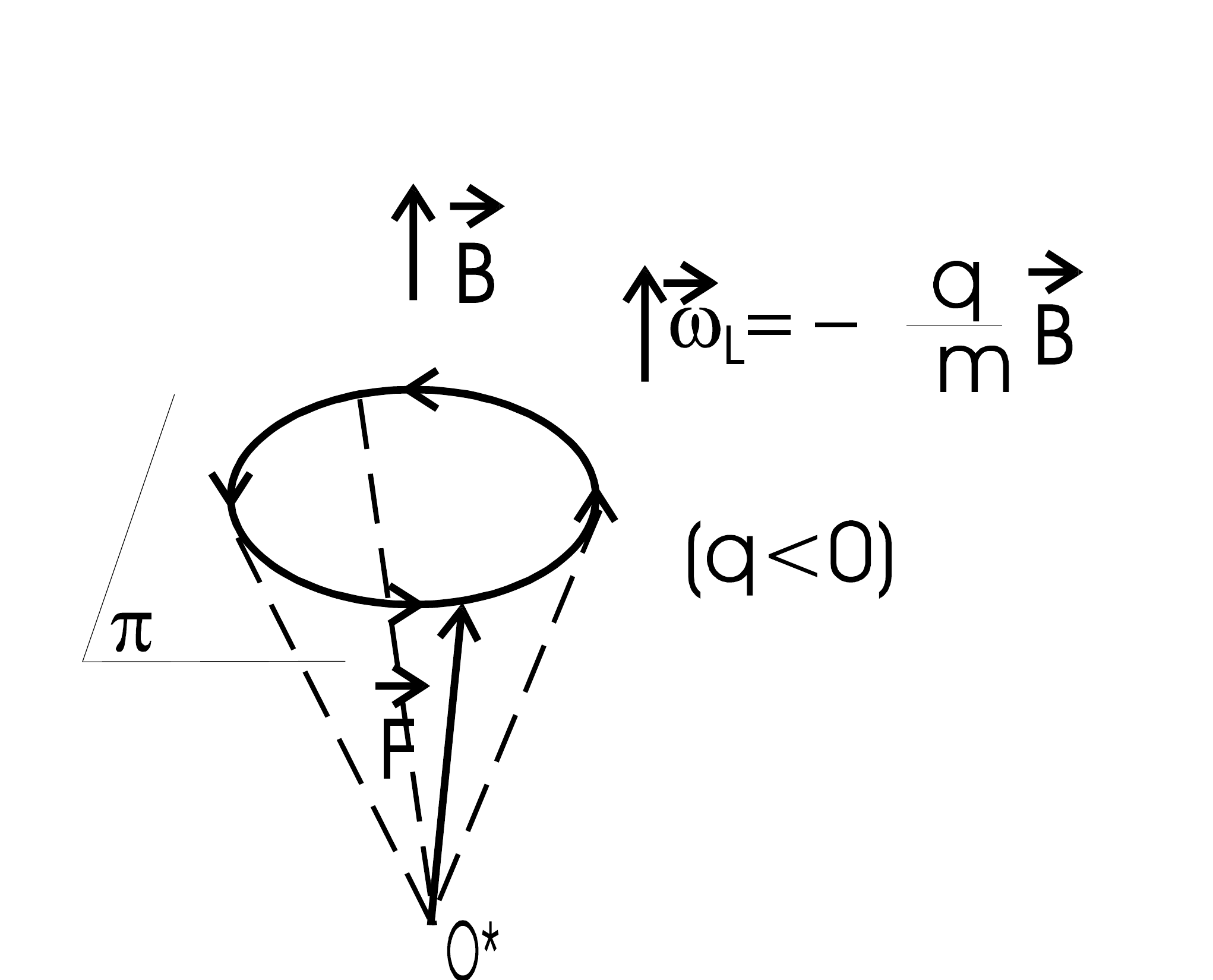}

\small\textnormal{Fig. 5}
\end{minipage}
\end{center}

The electric field $\vct E$ determines the strength of the interaction through $\|\vct F\|=q\|\vct v_0\times\vct B+\vct E\|$, whereas the magnetic field determines the rate $\vct\omega_L=-(q/m)\vct B$. These results generalize Larmor's theorem [99], [103].

\subsubsection*{VI.1.3.3. An Exact Solution for the Foucault Pendulum}
The mathematical model of Foucault's experiment is the Cauchy problem
\eqn{68}{\begin{cases}
\ddot{\vct r}+2\vct\omega\times\dot{\vct r}
+\vct\omega\times(\vct\omega\times\vct r)+\omega_0^2\vct r=\vct0,\\
\vct r(0)=\vct r_0,\quad\dot{\vct r}(0)=\vct v_0,
\end{cases}}
where $\vct r$ is the position relative to the relative equilibrium, $\vct\omega$ is the angular velocity of the Earth, and $\omega_0>0$ is the angular frequency of the oscillator. Applying $\Psi$, we obtain the matrix counterpart
\eqn{69}{\begin{cases}
\ddot r+2\mat\omega\dot r+\mat\omega^2r+\omega_0^2r=0,\\
r(0)=r_0,\quad\dot r(0)=v_0.
\end{cases}}

Consider the complex column matrix
\eqn{70}{\xi=r-\frac{j}{\omega_0}(\dot r+\mat\omega r),\qquad j^2=-1.}
Using (69), one obtains the first-order problem
\eqn{71}{\begin{cases}\dot\xi=\mat\omega^*\xi,\qquad
\mat\omega^*=-\mat\omega+j\omega_0I_3,\\
\xi(0)=r_0-\dfrac{j}{\omega_0}(v_0+\mat\omega r_0).
\end{cases}}
The solution is
\eqn{72}{\xi=e^{\mat\omega^*t}\xi(0).}
The real solution of problem (69) is
\eqn{73}{r=\operatorname{Re}\xi
=\operatorname{Re}(e^{\mat\omega^*t})\operatorname{Re}\xi(0)
-\operatorname{Im}(e^{\mat\omega^*t})\operatorname{Im}\xi(0).}
Since $\mat\omega$ commutes with $I_3$,
\begin{equation}\tag{\eqprefix.74}\begin{split}
e^{\mat\omega^*t}
&=e^{(-\mat\omega+j\omega_0I_3)t}
=e^{-\mat\omega t}e^{j\omega_0t}I_3\\
&=[\cos(\omega_0t)+j\sin(\omega_0t)]e^{-\mat\omega t}.
\end{split}\end{equation}
From (70), (73), and (74),
\eqn{75}{r=\cos(\omega_0t)e^{-\mat\omega t}r_0
+\sin(\omega_0t)e^{-\mat\omega t}\frac{v_0+\mat\omega r_0}{\omega_0}.}
Equivalently,
\eqn{76}{r=e^{-\mat\omega t}\left[r_0\cos(\omega_0t)+
\frac{v_0+\mat\omega r_0}{\omega_0}\sin(\omega_0t)\right].}
Let
\eqn{77}{\psi_0=r_0\cos(\omega_0t)+
\frac{v_0+\mat\omega r_0}{\omega_0}\sin(\omega_0t).}
Then
\eqn{78}{r=e^{-\mat\omega t}\psi_0.}

Applying $\Psi^{-1}$ to (78) and using Rodrigues' formula, we obtain the vector solution
\eqn{79}{\vct r=\vct\psi_0-\frac{\sin(\omega t)}{\omega}\vct\omega\times\vct\psi_0
+\frac{1-\cos(\omega t)}{\omega^2}\vct\omega\times(\vct\omega\times\vct\psi_0),}
where
\eqn{80}{\vct\psi_0=\vct r_0\cos(\omega_0t)+
\frac{\vct v_0+\vct\omega\times\vct r_0}{\omega_0}\sin(\omega_0t),
\qquad \omega=\|\vct\omega\|.}
An equivalent form is
\eqn{81}{\vct r=\frac{\vct\omega\cdot\vct\psi_0}{\omega^2}\vct\omega
-\frac{\sin(\omega t)}\omega\vct\omega\times\vct\psi_0
-\frac{\cos(\omega t)}{\omega^2}\vct\omega\times(\vct\omega\times\vct\psi_0).}

\textbf{Remarks.} The function $\vct\psi_0$ is the solution of the problem
\eqn{82}{\begin{cases}\ddot{\vct r}+\omega_0^2\vct r=\vct0,\\
\vct r(0)=\vct r_0,\quad\dot{\vct r}(0)=\vct v_0+\vct\omega\times\vct r_0.
\end{cases}}
The differential equation in (82) is that in (68) for $\vct\omega=\vct0$.

Relation (81) is the vector counterpart of (78). The matrix-valued function
\eqn{83}{Q(t)=e^{-\mat\omega t}}
is orthogonal, since $\mat\omega^T=-\mat\omega$:
\eqn{84}{Q^TQ=e^{-\mat\omega^Tt}e^{-\mat\omega t}
=e^{\mat\omega t}e^{-\mat\omega t}=I_3.}
Moreover, $\det Q=e^{-\operatorname{tr}(\mat\omega)t}=1$; hence (78) is a proper rotation. Its angular velocity is
\eqn{85}{\vct\omega^*=\operatorname{vect}(\dot QQ^T)
=\operatorname{vect}(-\mat\omega)=-\vct\omega.}
Denoting the operator of this rotation by $F_{-\omega}$, the solution is
\eqn{86}{\vct r=F_{-\omega}(\vct\psi_0).}

The preceding remarks are summarized in the following result.

\textbf{Theorem VI.1.3.3.1.} The solution of the problem
\[
\begin{cases}
\ddot{\vct r}+2\vct\omega\times\dot{\vct r}
+\vct\omega\times(\vct\omega\times\vct r)+\omega_0^2\vct r=\vct0,\\
\vct r(0)=\vct r_0,\quad\dot{\vct r}(0)=\vct v_0,
\end{cases}
\]
is obtained by applying the rotation operator with angular velocity $-\vct\omega$,
\[
F_{-\omega}(\,\cdot\,)=
\frac{\vct\omega\cdot(\,\cdot\,)}{\omega^2}\vct\omega
-\frac{\sin(\omega t)}\omega\vct\omega\times(\,\cdot\,)
-\frac{\cos(\omega t)}{\omega^2}\vct\omega\times[\vct\omega\times(\,\cdot\,)],
\]
to the solution of the nonrotating oscillator
\[
\begin{cases}\ddot{\vct r}+\omega_0^2\vct r=\vct0,\\
\vct r(0)=\vct r_0,\quad\dot{\vct r}(0)=\vct v_0+\vct\omega\times\vct r_0.
\end{cases}
\]

\textbf{Remark.} The hodograph of
\[
\vct\psi_0=\vct r_0\cos(\omega_0t)+
\frac{\vct v_0+\vct\omega\times\vct r_0}{\omega_0}\sin(\omega_0t)
\]
is an ellipse, possibly degenerate, with conjugate diameters in the directions of $\vct r_0$ and $\vct v_0+\vct\omega\times\vct r_0$. The solution of the Foucault pendulum is visualized by rotating the plane of this ellipse with angular velocity $-\vct\omega$.

The tensor relation suggests the direction developed in the following section, VI.2.

\renewcommand{\eqprefix}{VI.2}

\section*{VI.2. Tensorial Symbolic Representation of Vector-Valued Functions}
\subsection*{VI.2.1. General Results}
Let $V_3^*=\{\vct f:I\subset\mathbb R\to V_3\}$, and let $\vct\omega\in V_3^*$ be a continuous function. Denote by $\mat\omega$ the antisymmetric tensor associated with $\vct\omega$. Consider the Cauchy problem
\eqn{1}{\begin{cases}\dot Q=\mat\omega Q,\qquad \mat\omega^T=-\mat\omega,\\Q(0)=I_3.
\end{cases}}
For continuous $\mat\omega(t)$, the problem has a unique solution, and $Q$ is a proper orthogonal tensor. Indeed, for
\eqn{2}{A=QQ^T}
one obtains
\eqn{3}{\begin{cases}\dot A=\mat\omega A-A\mat\omega,\\A(0)=I_3.
\end{cases}}
The unique solution is $A=I_3$, hence $QQ^T=I_3$. Since $\det Q\in\{-1,1\}$ and $\det Q(0)=1$, we have $\det Q=1$.

Let
\eqn{4}{\vct\omega=\operatorname{vect}\mat\omega.}
The solution of (1) is called the rotation tensor corresponding to the angular velocity $\vct\omega$. Define
\eqn{5}{F_{\omega}:V_3^*\to V_3^*,\qquad F_{\omega}(\vct r)=Q\vct r.}
If the angular velocity has a fixed direction,
\eqn{6}{\vct\omega(t)=\omega(t)\vct u,\qquad \|\vct u\|=1,}
then the corresponding tensors commute with one another:
\eqn{7}{\mat\omega(t_1)\mat\omega(t_2)=\mat\omega(t_2)\mat\omega(t_1).}

Under conditions (6), the solution of problem (1) is explicit:
\eqn{8}{Q(t)=\exp\!\left(\int_0^t\mat\omega(\tau)\,d\tau\right).}
Rodrigues' formula gives
\eqn{9}{Q=I_3+\sin\alpha(t)\,\mat u+[1-\cos\alpha(t)]\mat u^{\,2},}
where
\eqn{10}{\alpha(t)=\int_0^t\omega(\tau)\,d\tau.}
Consequently,
\eqn{11}{F_{\omega}(\vct r)=\vct r+\sin\alpha(t)\,\vct u\times\vct r
+[1-\cos\alpha(t)]\vct u\times(\vct u\times\vct r).}
Equivalently,
\[
F_{\omega}(\vct r)=(\vct u\cdot\vct r)\vct u+
\sin\alpha(t)\,\vct u\times\vct r-
\cos\alpha(t)\,\vct u\times(\vct u\times\vct r).
\]
Since $\vct u=\vct\omega/\omega$, it follows that
\eqn{12}{F_{\omega}(\vct r)=
\frac{\vct\omega\cdot\vct r}{\omega^2}\vct\omega+
\frac{\sin\alpha(t)}\omega\vct\omega\times\vct r-
\frac{\cos\alpha(t)}{\omega^2}\vct\omega\times(\vct\omega\times\vct r).}
For $\vct\omega=$ const., $\alpha(t)=\omega t$, and
\eqn{12'}{F_{\omega}(\vct r)=
\frac{\vct\omega\cdot\vct r}{\omega^2}\vct\omega+
\frac{\sin(\omega t)}\omega\vct\omega\times\vct r-
\frac{\cos(\omega t)}{\omega^2}\vct\omega\times(\vct\omega\times\vct r).}

\textbf{Theorem VI.2.1.1.} If $\vct\omega$ has a fixed direction, transformation (5) has the following properties:
\begin{enumerate}
\item $F_\omega$ is $\mathbb R$-linear;
\item $F_\omega(\vct\omega\times\vct r)=\vct\omega\times F_\omega(\vct r)$;
\item $F_\omega(\dot{\vct r}+\vct\omega\times\vct r)=\dfrac d{dt}F_\omega(\vct r)$;
\end{enumerate}

\begin{enumerate}\setcounter{enumi}{3}
\item $F_\omega[\ddot{\vct r}+2\vct\omega\times\dot{\vct r}
+\vct\omega\times(\vct\omega\times\vct r)+\dot{\vct\omega}\times\vct r]
=\dfrac{d^2}{dt^2}F_\omega(\vct r)$;
\item $\|F_\omega(\vct r)\|=\|\vct r\|$;
\item $F_\omega(\vct r)|_{t=0}=\vct r_0$ and
$\dfrac d{dt}F_\omega(\vct r)|_{t=0}=\vct v_0+\vct\omega_0\times\vct r_0$;
\item $F_\omega$ is invertible and $F_\omega^{-1}=F_{-\omega}$.
\end{enumerate}

\emph{Proof.} Linearity follows from the definition. The vector $\vct\omega$ is an eigenvector of the rotation:
\eqn{13}{F_\omega(\vct\omega)=\vct\omega.}
The proper tensor preserves the cross product:
\eqn{14}{F_\omega(\vct\omega\times\vct r)=
F_\omega(\vct\omega)\times F_\omega(\vct r).}
From (13) and (14),
\eqn{15}{F_\omega(\vct\omega\times\vct r)=\vct\omega\times F_\omega(\vct r).}
The matrix counterpart is
\eqn{16}{Q\mat\omega r=\mat\omega Qr,}
whence
\eqn{17}{Q\mat\omega=\mat\omega Q.}

For $F=Qr$, differentiation gives
\eqn{18}{\dot F=\dot Qr+Q\dot r.}

Using (1),
\eqn{19}{\dot F=\mat\omega Qr+Q\dot r.}
With (17),
\eqn{20}{\dot F=Q\mat\omega r+Q\dot r=Q(\dot r+\mat\omega r).}
In tensorial form,
\eqn{21}{\frac d{dt}F_\omega(\vct r)=
F_\omega(\dot{\vct r}+\vct\omega\times\vct r).}
Applying (21) twice gives
\begin{equation}\tag{\eqprefix.22}\begin{split}
F_\omega\bigl[\ddot{\vct r}+2\vct\omega\times\dot{\vct r}
&+\vct\omega\times(\vct\omega\times\vct r)+\dot{\vct\omega}\times\vct r\bigr]\\
&=\frac{d^2}{dt^2}F_\omega(\vct r).
\end{split}\end{equation}
Since it is a proper rotation, $F_\omega$ is an isometry:
\eqn{23}{\|F_\omega(\vct r)\|=\|\vct r\|.}
At $t=0$,
\eqn{24}{F_\omega(\vct r)|_{t=0}=Q(0)r(0)=r_0,}
that is,
\eqn{25}{F_\omega(\vct r)|_{t=0}=\vct r_0.}
Using (21),
\eqn{26}{\frac d{dt}F_\omega(\vct r)\bigg|_{t=0}
=\vct v_0+\vct\omega_0\times\vct r_0.}
Finally, the tensor of $F_\omega^{-1}$ is $Q^T$. Transposing (1) gives
\eqn{27}{\dot Q^T=Q^T\mat\omega^T,}
and antisymmetry yields
\eqn{28}{\dot Q^T=-Q^T\mat\omega.}

Transposing commutation relation (17) and using $\mat\omega^T=-\mat\omega$ gives
\eqn{29}{Q^T\mat\omega=\mat\omega Q^T.}
From (28) and (29),
\eqn{30}{\dot Q^T=-\mat\omega Q^T.}
Thus, the tensor of $F_{-\omega}$ is $Q^T$, which proves
\[
F_\omega^{-1}=F_{-\omega}.
\]

\textbf{Remark.} The transformation $F_\omega:V_3^*\to V_3^*$ is a symbolic representation of the space of vector-valued functions in itself. Property 3 of Theorem VI.2.1.1 shows that the operation
\[
\vct r\longmapsto\dot{\vct r}+\vct\omega\times\vct r
\]
may be regarded as differentiation relative to the angular velocity $\vct\omega$, while the operator of the symbolic representation is ordinary differentiation. The transformation algebraizes a class of vector differential equations for motion in frames undergoing nonuniform rotation about a fixed direction and for motion in gyroscopic fields.

\subsection*{VI.2.2. Applications}
\subsection*{VI.2.2.1. Motion in a Non-Inertial Frame under a Uniform Gravitational Field}
Consider the Cauchy problem
\eqn{31}{\begin{cases}
\ddot{\vct r}+2\vct\omega\times\dot{\vct r}+
\vct\omega\times(\vct\omega\times\vct r)+\dot{\vct\omega}\times\vct r=\vct g,\\
\vct r(0)=\vct r_0,\quad\dot{\vct r}(0)=\vct v_0,
\end{cases}}
where the variable angular velocity has a fixed direction,
\eqn{32}{\vct\omega(t)=\omega(t)\vct u,\qquad \dot{\vct u}=\vct0.}
Define the transformation of the form (12)
\eqn{33}{F_\omega(\vct r)=
\frac{\vct\omega\cdot\vct r}{\omega^2}\vct\omega+
\frac{\sin\alpha(t)}\omega\vct\omega\times\vct r-
\frac{\cos\alpha(t)}{\omega^2}\vct\omega\times(\vct\omega\times\vct r).}

Here
\[
\alpha(t)=\int_0^t\omega(x)\,dx.
\]
We apply transformation (33) to problem (31) and set
\eqn{34}{\vct\psi=F_\omega(\vct r).}
Theorem VI.2.1.1 gives
\eqn{35}{\begin{cases}\ddot{\vct\psi}=F_\omega(\vct g),\\
\vct\psi(0)=\vct r_0,\quad
\dot{\vct\psi}(0)=\vct v_0+\vct\omega_0\times\vct r_0.
\end{cases}}
A first integration gives
\eqn{36}{\dot{\vct\psi}=\vct v_0+\vct\omega_0\times\vct r_0+
\int_0^tF_\omega(\vct g)(x)\,dx.}
Using (33),
\begin{equation}\tag{\eqprefix.37}\begin{split}
\dot{\vct\psi}={}&\vct v_0+\vct\omega_0\times\vct r_0+
\frac{\vct\omega\cdot\vct g}{\omega^2}\vct\omega,t\\
&+\frac{\vct\omega\times\vct g}{\omega}\int_0^t\sin\alpha(x)\,dx
-\frac{\vct\omega\times(\vct\omega\times\vct g)}{\omega^2}
\int_0^t\cos\alpha(x)\,dx.
\end{split}\end{equation}
The second integration yields
\begin{equation}\tag{\eqprefix.38}\begin{split}
\vct\psi_0={}&\vct r_0+[\vct v_0+\vct\omega_0\times\vct r_0]t+
\frac{\vct\omega\cdot\vct g}{2\omega^2}\vct\omega t^2\\
&+\frac{\vct\omega\times\vct g}{\omega}\int_0^t dx\int_0^x\sin\alpha(\tau)\,d\tau\\
&-\frac{\vct\omega\times(\vct\omega\times\vct g)}{\omega^2}
\int_0^t dx\int_0^x\cos\alpha(\tau)\,d\tau.
\end{split}\end{equation}
The solution of problem (31) is
\eqn{39}{\vct r=F_\omega^{-1}(\vct\psi_0)=F_{-\omega}(\vct\psi_0),}
that is,
\eqn{40}{\vct r=\frac{\vct\omega\cdot\vct\psi_0}{\omega^2}\vct\omega-
\frac{\sin\alpha(t)}\omega\vct\omega\times\vct\psi_0-
\frac{\cos\alpha(t)}{\omega^2}\vct\omega\times(\vct\omega\times\vct\psi_0).}

\textbf{Remark.} If $\vct\omega$ is constant, then $\alpha(t)=\omega t$, and (38) and (40) yield
\begin{equation}\tag{\eqprefix.41}\begin{split}
\vct r={}&\vct r_0+\vct v_0t+
\frac{\vct\omega\cdot\vct g}{2\omega^2}\vct\omega t^2
-\frac{\sin(\omega t)}\omega\vct\omega\times(\vct r_0^*+\vct v_0t)\\
&+\frac{1-\cos(\omega t)}{\omega^2}
\vct\omega\times[\vct\omega\times(\vct r_0^*+\vct v_0t)]\\
&+\omega t\cos(\omega t)\frac{\vct\omega\times\vct r_0^*}{\omega}
+\omega t\sin(\omega t)
\frac{\vct\omega\times(\vct\omega\times\vct r_0^*)}{\omega^2},
\end{split}\end{equation}
where $\vct r_0^*=\vct r_0+\vct g/\omega^2$. This coincides with relation VI.1.3.1 (57).

\subsection*{VI.2.2.2. Motion of Electric Charges in a Slowly Varying Electromagnetic Field}
\subsubsection*{a) Variable Magnetic Field of Fixed Direction}
The velocity equation of a particle of charge $q$ and mass $m$ is
\eqn{42}{\begin{cases}m\dot{\vct v}=q\vct v\times\vct B,\\\vct v(0)=\vct v_0.
\end{cases}}
Setting
\eqn{43}{\vct\omega=\frac qm\vct B,}
the problem becomes
\eqn{44}{\begin{cases}\dot{\vct v}+\vct\omega\times\vct v=\vct0,\\
\vct v(0)=\vct v_0.
\end{cases}}
If $\vct B(t)$ is variable but has a fixed direction, then $\vct\omega$ has the same property. Applying $F_\omega$ and setting $\vct\psi=F_\omega(\vct v)$, we obtain
\eqn{45}{\begin{cases}\dot{\vct\psi}=\vct0,\\\vct\psi(0)=\vct v_0.
\end{cases}}
Consequently,
\eqn{46}{\vct\psi_0=\vct v_0,}
and
\eqn{47}{\vct v=F_{-\omega}(\vct v_0).}

Explicitly,
\eqn{48}{\vct v=\frac{\vct\omega\cdot\vct v_0}{\omega^2}\vct\omega-
\frac{\sin\alpha(t)}\omega\vct\omega\times\vct v_0-
\frac{\cos\alpha(t)}{\omega^2}\vct\omega\times(\vct\omega\times\vct v_0),}
where
\[
\alpha(t)=\int_0^t\omega(x)\,dx=\frac qm\int_0^tB(x)\,dx.
\]
In physical variables, this becomes
\eqn{49}{\vct v=\frac{\vct v_0\cdot\vct B}{B^2}\vct B+
\frac{\sin\alpha(t)}B\vct v_0\times\vct B+
\frac{\cos\alpha(t)}{B^2}\vct B\times(\vct v_0\times\vct B).}
The exact solution had not been reported in the literature cited by the author.

The law of motion is the solution of the problem
\eqn{50}{\begin{cases}m\ddot{\vct r}=q\dot{\vct r}\times\vct B,
\qquad \vct B\times\dot{\vct B}=\vct0,\\
\vct r(0)=\vct r_0,\quad\dot{\vct r}(0)=\vct v_0.
\end{cases}}
The condition $\vct B\times\dot{\vct B}=\vct0$ is necessary and sufficient for $\vct B(t)$ to have a fixed direction. Integration of (49) gives
\eqn{51}{\vct r=\vct r_0+\int_0^t\vct v(x)\,dx,}
that is,
\begin{equation}\tag{\eqprefix.52}\begin{split}
\vct r={}&\vct r_0+
\frac{\vct v_0\cdot\vct B}{B^2}\vct B,t+
\frac{\vct v_0\times\vct B}{B}\int_0^t\sin\alpha(x)\,dx\\
&+\frac{\vct B\times(\vct v_0\times\vct B)}{B^2}
\int_0^t\cos\alpha(x)\,dx.
\end{split}\end{equation}

\textbf{Remark.} The velocity hodograph is a circle in a plane perpendicular to the magnetic-field lines. The velocity vector has constant magnitude and undergoes nonuniform precession about $\vct B$, with instantaneous angular velocity $\vct\omega_L=-(q/m)\vct B$.

\begin{center}
\begin{minipage}{.47\textwidth}\centering
\includegraphics[width=\linewidth]{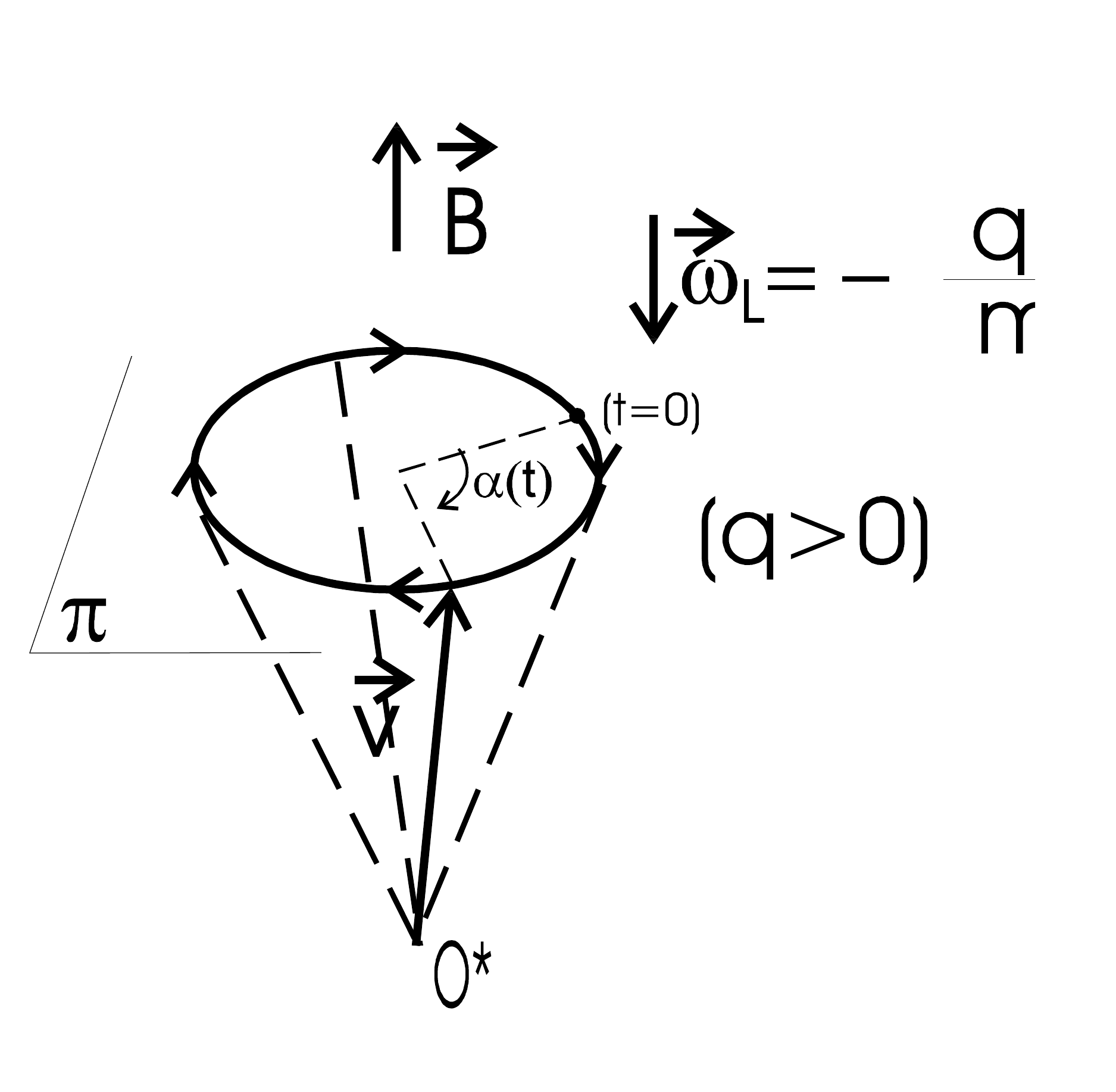}

\small\textnormal{Fig. 6}
\end{minipage}\hfill
\begin{minipage}{.47\textwidth}\centering
\includegraphics[width=\linewidth]{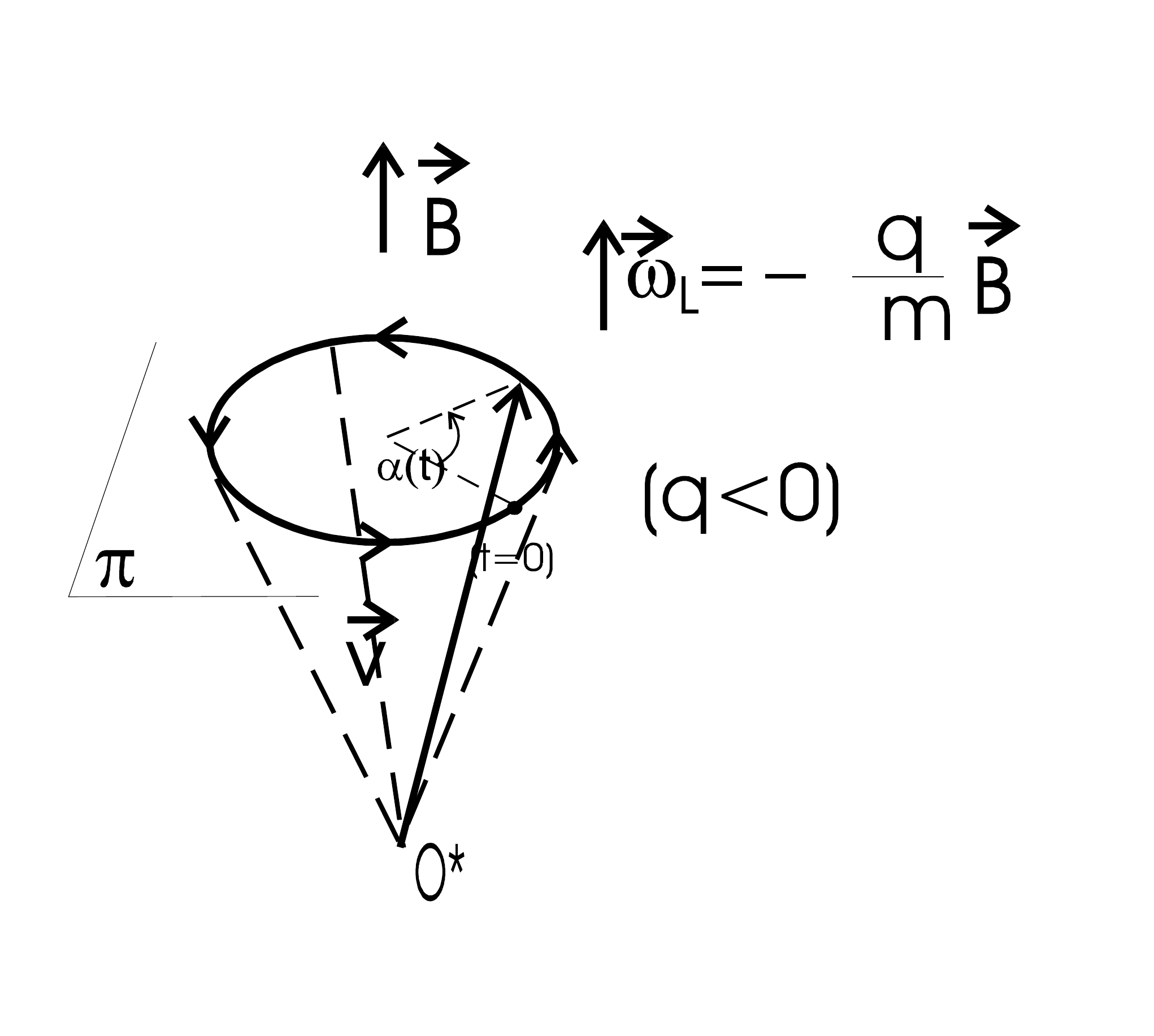}

\small\textnormal{Fig. 7}
\end{minipage}
\end{center}

The acceleration of the particle follows by differentiating relation (49):
\[
\vct a=\dot{\vct v}=\dot\alpha(t)\cos\alpha(t)
\frac{\vct v_0\times\vct B}{B}
-\dot\alpha(t)\sin\alpha(t)
\frac{\vct B\times(\vct v_0\times\vct B)}{B^2}.
\]
In carrying out the differentiation, we used the fact that $\vct B/B$ is constant when the direction of the magnetic induction remains fixed. Since $\vct F=m\vct a$ and $\dot\alpha=(q/m)B$, it follows that
\[
\vct F=\cos\alpha(t)(q\vct v_0\times\vct B)
-\sin\alpha(t)\frac{\vct B\times(q\vct v_0\times\vct B)}{B}.
\]
The force remains perpendicular to the magnetic-field lines, has variable magnitude
$\|\vct F\|=q\|\vct v_0\times\vct B\|$, and is a rotating vector with nonuniform angular velocity
\[
\vct\omega_L=-\frac qm\vct B.
\]
These observations extend Larmor's theorem to a variable magnetic field of fixed direction.

\subsubsection*{b) Electromagnetic Field with Slowly Varying Magnetic Induction}
Consider a slowly varying magnetic field $\vct B(t)$ of fixed direction and a uniform electric field $\vct E$. The velocity problem is
\eqn{53}{\begin{cases}
m\dot{\vct v}=q(\vct v\times\vct B+\vct E),
\qquad \dot{\vct B}\times\vct B=\vct0,\quad \vct E=\mathrm{const.},\\
\vct v(0)=\vct v_0.
\end{cases}}

With the notation
\eqn{54}{\vct\omega=\frac qm\vct B,}
\eqn{55}{\vct\alpha=\frac qm\vct E,}
problem (53) becomes
\eqn{56}{\begin{cases}
\dot{\vct v}+\vct\omega\times\vct v=\vct\alpha,
\qquad \dot{\vct\omega}\times\vct\omega=\vct0,\quad
\vct\alpha=\mathrm{const.},\\
\vct v(0)=\vct v_0.
\end{cases}}
Applying the transformation $F_\omega$ and setting $\vct\psi=F_\omega(\vct v)$ gives
\eqn{57}{\begin{cases}\dot{\vct\psi}=F_\omega(\vct\alpha),\\
\vct\psi(0)=\vct v_0.
\end{cases}}
The solution is
\eqn{58}{\vct\psi_0=\vct v_0+\int_0^tF_\omega(\vct\alpha)(x)\,dx.}
Using the explicit expression of the transformation and the fact that $\vct\alpha$ is constant, we obtain
\begin{equation}\tag{\eqprefix.59}\begin{split}
\vct\psi_0={}&\vct v_0+
\frac{\vct\omega\cdot\vct\alpha}{\omega^2}\vct\omega\,t
+\frac{\vct\omega\times\vct\alpha}{\omega}
\int_0^t\sin\alpha(x)\,dx\\
&-\frac{\vct\omega\times(\vct\omega\times\vct\alpha)}{\omega^2}
\int_0^t\cos\alpha(x)\,dx,
\end{split}\end{equation}
where the scalar phase is
\[
\alpha(t)=\int_0^t\omega(\tau)\,d\tau.
\]
The solution of problem (56) is
\eqn{60}{\vct v=F_{-\omega}(\vct\psi_0),}
that is,
\eqn{61}{\vct v=
\frac{\vct\omega\cdot\vct\psi_0}{\omega^2}\vct\omega-
\frac{\sin\alpha(t)}\omega\vct\omega\times\vct\psi_0-
\frac{\cos\alpha(t)}{\omega^2}\vct\omega\times(\vct\omega\times\vct\psi_0).}

In formula (61), $\vct\psi_0$ is given by (59). For a more readable notation, set
\[
I_s(t)=\int_0^t\sin\alpha(x)\,dx,\qquad
I_c(t)=\int_0^t\cos\alpha(x)\,dx.
\]
After expanding $F_{-\omega}(\vct\psi_0)$, the solution takes the form
\begin{equation}\tag{\eqprefix.62}\begin{split}
\vct v={}&\frac{\vct v_0\cdot\vct\omega}{\omega^2}\vct\omega
+\frac{\vct\alpha\times\vct\omega}{\omega^2}
[\cos\alpha(t)+\omega\sin\alpha(t)I_c-\omega\cos\alpha(t)I_s]\\
&+\frac{\vct\alpha+\vct v_0\times\vct\omega}{\omega}\sin\alpha(t)
+\frac{\vct\omega\times(\vct\alpha+\vct v_0\times\vct\omega)}{\omega^2}\cos\alpha(t)\\
&+\frac{\vct\omega\cdot\vct\alpha}{\omega^2}
[t-\sin\alpha(t)I_s-\cos\alpha(t)I_c]\vct\omega\\
&+[\sin\alpha(t)I_s+\cos\alpha(t)I_c-\tfrac{\sin\alpha(t)}\omega]\vct\alpha.
\end{split}\end{equation}
Here $\vct\omega=(q/m)\vct B$, $\vct\alpha=(q/m)\vct E$, and the scalar phase is $\alpha(t)=\int_0^t\omega(\tau)\,d\tau$.

Returning to the physical quantities gives
\begin{equation}\tag{\eqprefix.63}\begin{split}
\vct v={}&\frac{\vct v_0\cdot\vct B}{B^2}\vct B
+\frac{\vct E\times\vct B}{B^2}
[\cos\alpha(t)+\omega\sin\alpha(t)I_c-\omega\cos\alpha(t)I_s]\\
&+\frac{\vct E+\vct v_0\times\vct B}{B}\sin\alpha(t)
+\frac{\vct B\times(\vct E+\vct v_0\times\vct B)}{B^2}\cos\alpha(t)\\
&+\frac{\vct E\cdot\vct B}{B^3}
[\omega t-\omega\sin\alpha(t)I_s-\omega\cos\alpha(t)I_c]\vct B\\
&+[\omega\sin\alpha(t)I_s+\omega\cos\alpha(t)I_c-\sin\alpha(t)]\frac{\vct E}{B}.
\end{split}\end{equation}
This is the solution for $\vct B=B(t)\vct u$, $\dot{\vct u}=\vct0$.

For constant $\vct B$, $\alpha(t)=\omega t$, and the integrals can be evaluated in elementary form. Relation (63) reduces to
\begin{equation}\tag{\eqprefix.64}\begin{split}
\vct v={}&\frac{\vct v_0\cdot\vct B}{B^2}\vct B+
\frac{\vct E\times\vct B}{B^2}+
\frac{\vct E+\vct v_0\times\vct B}{B}\sin(\omega t)\\
&+\frac{\vct B\times(\vct E+\vct v_0\times\vct B)}{B^2}\cos(\omega t)
+\frac{\vct E\cdot\vct B}{B^3}[\omega t-\sin(\omega t)]\vct B.
\end{split}\end{equation}
This coincides with solution VI.1.3.2 (62).

\subsection*{VI.2.2.3. Motion in a Non-Inertial Frame under Positional Force Fields}
The procedure of VI.2.1 makes it possible to solve motion problems in a frame undergoing nonuniform rotation about a fixed direction, under central force fields of the form
\[
\vct F(\vct r)=f(r)\frac{\vct r}{r},\qquad r=\|\vct r\|.
\]

\textbf{Theorem VI.2.2.3.1.} The solution of the problem
\eqn{65}{\begin{cases}
\ddot{\vct r}+2\vct\omega\times\dot{\vct r}+
\vct\omega\times(\vct\omega\times\vct r)+
\dot{\vct\omega}\times\vct r+f(r)\dfrac{\vct r}{r}=\vct0,\\
\dot{\vct\omega}\times\vct\omega=\vct0,\qquad
\vct r(0)=\vct r_0,\quad\dot{\vct r}(0)=\vct v_0,
\end{cases}}
is obtained by applying the rotation operator
\begin{equation}\tag{\eqprefix.66}\begin{split}
F_{-\omega}(\,\cdot\,)={}&
\frac{\vct\omega\cdot(\,\cdot\,)}{\omega^2}\vct\omega
-\frac{\sin\alpha(t)}\omega\vct\omega\times(\,\cdot\,)\\
&-\frac{\cos\alpha(t)}{\omega^2}
\vct\omega\times[\vct\omega\times(\,\cdot\,)],
\end{split}\end{equation}
where $\vct\omega=\omega(t)\vct u$, $\dot{\vct u}=\vct0$, and
$\alpha(t)=\int_0^t\omega(x)\,dx$, to the solution of the inertial problem
\eqn{67}{\begin{cases}
\ddot{\vct\psi}+f(\psi)\dfrac{\vct\psi}{\psi}=\vct0,\\
\vct\psi(0)=\vct r_0,\quad
\dot{\vct\psi}(0)=\vct v_0+\vct\omega_0\times\vct r_0.
\end{cases}}
The functions $\vct\omega$ and $f$ are assumed sufficiently regular to ensure existence and uniqueness of the solutions.

\emph{Proof.} Applying $F_\omega$ to problem (65) and setting $\vct\psi=F_\omega(\vct r)$, Theorem VI.2.1.1 yields
\eqn{68}{\begin{cases}
\ddot{\vct\psi}+f(\psi)\dfrac{\vct\psi}{\psi}=\vct0,\\
\vct\psi(0)=\vct r_0,\quad
\dot{\vct\psi}(0)=\vct v_0+\vct\omega_0\times\vct r_0.
\end{cases}}
Here we have used the implication
\eqn{69}{\vct\psi=F_\omega(\vct r)\Longrightarrow
F_\omega\!\left[f(r)\frac{\vct r}{r}\right]
=f(\psi)\frac{\vct\psi}{\psi},}
which follows from the linearity, homogeneity, and isometry of $F_\omega$.

Let $\vct\psi_0$ be the solution of problem (68). Since $F_\omega$ is invertible and $F_\omega^{-1}=F_{-\omega}$, the solution of (65) is
\[
\vct r=F_{-\omega}(\vct\psi_0).
\]

\textbf{Remark.} The theorem reduces the problem in the non-inertial frame to the corresponding problem in the inertial frame, with the initial velocity condition modified according to the composition law. The trajectory described by the solution of (65) is obtained by rotating, with angular velocity $-\vct\omega$, the trajectory described by the solution of (67).

\textbf{Application: the Foucault pendulum.} For an elastic field, $f(r)=\omega_*^2r$, $\omega_*>0$, problem (65) becomes
\eqn{70}{\begin{cases}
\ddot{\vct r}+2\vct\omega\times\dot{\vct r}+
\vct\omega\times(\vct\omega\times\vct r)+
\dot{\vct\omega}\times\vct r+\omega_*^2\vct r=\vct0,\\
\dot{\vct\omega}\times\vct\omega=\vct0,\qquad
\vct r(0)=\vct r_0,\quad\dot{\vct r}(0)=\vct v_0.
\end{cases}}
The reduced problem (68) is
\eqn{71}{\begin{cases}
\ddot{\vct\psi}+\omega_*^2\vct\psi=\vct0,\\
\vct\psi(0)=\vct r_0,\quad
\dot{\vct\psi}(0)=\vct v_0+\vct\omega_0\times\vct r_0,
\qquad \vct\omega_0=\vct\omega(0).
\end{cases}}

The solution of problem (71) is
\eqn{72}{\vct\psi_0=\vct r_0\cos(\omega_*t)+
\frac{\vct v_0+\vct\omega_0\times\vct r_0}{\omega_*}\sin(\omega_*t).}
Its hodograph is an ellipse, possibly degenerate, with conjugate diameters along the directions of $\vct r_0$ and $\vct v_0+\vct\omega_0\times\vct r_0$. By Theorem VI.2.2.3.1,
\eqn{73}{\vct r=F_{-\omega}\!\left[\vct r_0\cos(\omega_*t)+
\frac{\vct v_0+\vct\omega_0\times\vct r_0}{\omega_*}\sin(\omega_*t)\right].}
The trajectory of the solution is obtained by rotating the plane of the ellipse with the nonuniform angular velocity $-\vct\omega$. Denoting the expression in brackets by $\vct\psi_0$ and setting
\[
\alpha(t)=\int_0^t\omega(x)\,dx,\qquad \vct\omega=\omega(t)\vct u,\quad\dot{\vct u}=\vct0,
\]
the explicit form is
\eqn{74}{\vct r=\frac{\vct\omega\cdot\vct\psi_0}{\omega^2}\vct\omega-
\frac{\sin\alpha(t)}\omega\vct\omega\times\vct\psi_0-
\frac{\cos\alpha(t)}{\omega^2}\vct\omega\times(\vct\omega\times\vct\psi_0).}
This exact solution had not been reported in the literature cited by the author.

We also consider the traditional reduced form of the Foucault-pendulum model:
\eqn{75}{\begin{cases}
\ddot{\vct r}+2\vct\omega\times\dot{\vct r}+\omega_0^2\vct r=\vct0,
\qquad \vct\omega=\mathrm{const.},\quad\omega_0>0,\\
\vct r(0)=\vct r_0,\quad\dot{\vct r}(0)=\vct v_0.
\end{cases}}
This form neglects the centrifugal contribution $\vct\omega\times(\vct\omega\times\vct r)$. If $\vct r_0\times\vct\omega=\vct0$ and $\vct v_0\times\vct\omega=\vct0$, the solution remains parallel to the axis and satisfies
\eqn{76}{\ddot{\vct r}+\omega_0^2\vct r=\vct0.}

\vspace*{8mm}
The solution of equation (76), with the initial data in (75), is
\eqn{77}{\vct r=\vct r_0\cos(\omega_0t)+\frac{\vct v_0}{\omega_0}\sin(\omega_0t).}

Now assume $\vct r_0\cdot\vct\omega=0$ and $\vct v_0\cdot\vct\omega=0$, and seek solutions satisfying $\vct r\cdot\vct\omega=0$. Let
\eqn{78}{\vct\psi=F_\omega(\vct r)=
\frac{\vct\omega\cdot\vct r}{\omega^2}\vct\omega+
\frac{\sin(\omega t)}\omega\vct\omega\times\vct r-
\frac{\cos(\omega t)}{\omega^2}\vct\omega\times(\vct\omega\times\vct r).}
From $\vct\omega\cdot\vct r=0$ it follows that $\vct\omega\cdot\vct\psi=0$. By Theorem VI.2.1.1,
\begin{equation}\tag{\eqprefix.79}\begin{split}
\ddot{\vct\psi}
&=F_\omega[\ddot{\vct r}+2\vct\omega\times\dot{\vct r}
+\vct\omega\times(\vct\omega\times\vct r)]\\
&=F_\omega(\ddot{\vct r}+2\vct\omega\times\dot{\vct r})
+\vct\omega\times(\vct\omega\times\vct\psi).
\end{split}\end{equation}
Therefore,
\eqn{80}{F_\omega(\ddot{\vct r}+2\vct\omega\times\dot{\vct r})
=\ddot{\vct\psi}-\vct\omega\times(\vct\omega\times\vct\psi).}
Applying $F_\omega$ to problem (75), we obtain
\eqn{81}{\begin{cases}
\ddot{\vct\psi}-\vct\omega\times(\vct\omega\times\vct\psi)+
\omega_0^2\vct\psi=\vct0,\\
\vct\psi(0)=\vct r_0,\quad
\dot{\vct\psi}(0)=\vct v_0+\vct\omega\times\vct r_0.
\end{cases}}
Under the condition $\vct\psi\cdot\vct\omega=0$, this becomes
\eqn{82}{\begin{cases}
\ddot{\vct\psi}+(\omega^2+\omega_0^2)\vct\psi=\vct0,\\
\vct\psi(0)=\vct r_0,\quad
\dot{\vct\psi}(0)=\vct v_0+\vct\omega\times\vct r_0.
\end{cases}}
The solution is
\eqn{83}{\vct\psi_* =\vct r_0\cos(\omega_*t)+
\frac{\vct v_0+\vct\omega\times\vct r_0}{\omega_*}\sin(\omega_*t),
\qquad \omega_*^2=\omega^2+\omega_0^2.}
It also satisfies $\vct\psi_*\cdot\vct\omega=0$. The solution of problem (75) is
\eqn{84}{\vct r=F_{-\omega}(\vct\psi_*).}

Since $\vct\psi_*\cdot\vct\omega=0$, relation (84) reduces to
\eqn{85}{\vct r=-\frac{\sin(\omega t)}\omega\vct\omega\times\vct\psi_*
+\cos(\omega t)\vct\psi_*.}

Decompose the initial data into components parallel and perpendicular to $\vct\omega$:
\eqn{86}{\begin{cases}
\vct r_{01}=\dfrac{\vct\omega\cdot\vct r_0}{\omega^2}\vct\omega,
\qquad
\vct v_{01}=\dfrac{\vct\omega\cdot\vct v_0}{\omega^2}\vct\omega,
\end{cases}}
\eqn{87}{\begin{cases}
\vct r_{02}=\dfrac{\vct\omega\times(\vct r_0\times\vct\omega)}{\omega^2},
\qquad
\vct v_{02}=\dfrac{\vct\omega\times(\vct v_0\times\vct\omega)}{\omega^2}.
\end{cases}}
The corresponding solutions are
\eqn{88}{\vct r_1=\vct r_{01}\cos(\omega_0t)+
\frac{\vct v_{01}}{\omega_0}\sin(\omega_0t)
=\frac{\vct\omega\cdot\vct\psi_0}{\omega^2}\vct\omega,}
where
\eqn{89}{\vct\psi_0=\vct r_0\cos(\omega_0t)+
\frac{\vct v_0+\vct\omega\times\vct r_0}{\omega_0}\sin(\omega_0t),}
and
\eqn{90}{\vct r_2=-\frac{\sin(\omega t)}\omega\vct\omega\times\vct\psi_2
+\cos(\omega t)\vct\psi_2,}
with
\eqn{91}{\vct\psi_2=\vct r_{02}\cos(\omega_*t)+
\frac{\vct v_{02}+\vct\omega\times\vct r_{02}}{\omega_*}\sin(\omega_*t)
=\frac{\vct\omega\times(\vct\psi_*\times\vct\omega)}{\omega^2}.}
It follows that
\eqn{92}{\vct r_2=-\frac{\sin(\omega t)}\omega\vct\omega\times\vct\psi_*
+\frac{\cos(\omega t)}{\omega^2}\vct\omega\times(\vct\psi_*\times\vct\omega).}
Since $\vct r_0=\vct r_{01}+\vct r_{02}$ and $\vct v_0=\vct v_{01}+\vct v_{02}$, linearity implies $\vct r=\vct r_1+\vct r_2$.

\par\medskip
\textbf{Theorem VI.2.2.3.2.} The solution of the Cauchy problem
\eqn{93}{\begin{cases}
\ddot{\vct r}+2\vct\omega\times\dot{\vct r}+\omega_0^2\vct r=\vct0,
\qquad \vct\omega=\mathrm{const.},\quad\omega_0>0,\\
\vct r(0)=\vct r_0,\quad\dot{\vct r}(0)=\vct v_0,
\end{cases}}
is
\eqn{94}{\vct r=\frac{\vct\omega\cdot\vct\psi_0}{\omega^2}\vct\omega-
\frac{\sin(\omega t)}\omega\vct\omega\times\vct\psi_*-
\frac{\cos(\omega t)}{\omega^2}\vct\omega\times(\vct\omega\times\vct\psi_*),}
where
\eqn{95}{\vct\psi_0=\vct r_0\cos(\omega_0t)+
\frac{\vct v_0+\vct\omega\times\vct r_0}{\omega_0}\sin(\omega_0t),}
\eqn{96}{\vct\psi_*=\vct r_0\cos(\omega_*t)+
\frac{\vct v_0+\vct\omega\times\vct r_0}{\omega_*}\sin(\omega_*t),
\qquad \omega_*^2=\omega^2+\omega_0^2.}
The solution is exact and applies, in particular, to models of vibrating elements in gyroscopic instruments.

\subsubsection*{b) First Integrals in Central Positional Force Fields}
Theorem VI.2.2.3.1 yields the following result.

\par\medskip
\textbf{Theorem VI.2.2.3.3.} For the problem
\eqn{97}{\begin{cases}
\ddot{\vct r}+2\vct\omega\times\dot{\vct r}+
\vct\omega\times(\vct\omega\times\vct r)+
\dot{\vct\omega}\times\vct r+f(r)\dfrac{\vct r}{r}=\vct0,\\
\dot{\vct\omega}\times\vct\omega=\vct0,\qquad
\vct r(0)=\vct r_0,\quad\dot{\vct r}(0)=\vct v_0,
\end{cases}}
there exist the first integrals
\eqn{98}{\frac{\|\dot{\vct r}\|^2}{2}+(\vct\omega,\vct r,\dot{\vct r})+
\frac12\|\vct\omega\times\vct r\|^2+\int f(r)\,dr=\mathrm{const.},}
\eqn{99}{\vct r\times(\dot{\vct r}+\vct\omega\times\vct r)=
\frac{\vct\omega\cdot\vct\Omega_0}{\omega^2}\vct\omega-
\frac{\sin\alpha(t)}\omega\vct\omega\times\vct\Omega_0-
\frac{\cos\alpha(t)}{\omega^2}\vct\omega\times(\vct\omega\times\vct\Omega_0),}
where $\vct\Omega_0=\vct r_0\times(\vct v_0+\vct\omega_0\times\vct r_0)$ and
$\alpha(t)=\int_0^t\omega(x)\,dx$.

The transformed problem is
\eqn{100}{
\begin{cases}
\ddot{\vct\psi}+f(\psi)\dfrac{\vct\psi}{\psi}=\vct0,\\
\vct\psi(0)=\vct r_0,\\
\dot{\vct\psi}(0)=\vct v_0+\vct\omega_0\times\vct r_0.
\end{cases}}

Problem (100) admits the first integrals of energy and angular momentum [87], [89]:
\eqn{101}{\frac{\|\dot{\vct\psi}\|^2}{2}+\int_0^{\psi}f(x)\,dx=\mathrm{const.},}
\eqn{102}{\vct\psi\times\dot{\vct\psi}
=\vct r_0\times(\vct v_0+\vct\omega_0\times\vct r_0).}

By Theorem VI.1.2.1, we have
\eqn{103}{\dot{\vct\psi}=F_{\omega}(\dot{\vct r}+\vct\omega\times\vct r).}

Since the transformation $F_{\omega}$ is an isometry, (103) gives
\eqn{104}{\|\dot{\vct\psi}\|^2
=\|\dot{\vct r}+\vct\omega\times\vct r\|^2
=\|\dot{\vct r}\|^2+2(\vct\omega\times\vct r)\!\cdot\!\dot{\vct r}
+\|\vct\omega\times\vct r\|^2,}
\eqn{105}{\psi=\|\vct\psi\|=\|F_{\omega}(\vct r)\|=\|\vct r\|=r.}

Using relations (104) and (105), the first integral (101) can be written as
\eqn{106}{\frac{\|\dot{\vct r}\|^2}{2}
+\dot{\vct r}\!\cdot\!(\vct\omega\times\vct r)
+\frac12\|\vct\omega\times\vct r\|^2
+\int f(r)\,dr=\mathrm{const.}}

Using the properties of the scalar triple product, (106) takes the equivalent form
\eqn{107}{\frac{\|\dot{\vct r}\|^2}{2}
+(\vct\omega,\vct r,\dot{\vct r})
+\frac12\|\vct\omega\times\vct r\|^2
+\int f(r)\,dr=\mathrm{const.}}

Setting $\vct\Omega_0=\vct r_0\times(\vct v_0+\vct\omega_0\times\vct r_0)$,
relation (102) becomes
\begin{equation}\tag{\eqprefix.108}
\begin{aligned}
\vct\psi\times\dot{\vct\psi}=\vct\Omega_0
&\iff F_{\omega}(\vct r)\times\dot F_{\omega}(\vct r)=\vct\Omega_0\\
&\iff F_{\omega}(\vct r)\times
F_{\omega}(\dot{\vct r}+\vct\omega\times\vct r)=\vct\Omega_0\\
&\iff F_{\omega}[\vct r\times(\dot{\vct r}+\vct\omega\times\vct r)]=\vct\Omega_0\\
&\iff \vct r\times(\dot{\vct r}+\vct\omega\times\vct r)
=F_{-\omega}(\vct\Omega_0).
\end{aligned}
\end{equation}
Consequently,
\eqn{109}{
\vct r\times(\dot{\vct r}+\vct\omega\times\vct r)
=\frac{\vct\omega\cdot\vct\Omega_0}{\omega^2}\vct\omega
-\frac{\sin\alpha(t)}{\omega}\vct\omega\times\vct\Omega_0
-\frac{\cos\alpha(t)}{\omega^2}\vct\omega\times(\vct\omega\times\vct\Omega_0),
\quad t\geq0.}

The chain of equivalences (108) uses properties of the transformation $F_{\omega}$
established by Theorem VI.1.2.1.

\medskip
\noindent\underline{Remark}

The first integral (107) proves the existence, in the force field consisting of
the central positional force $\vct F=f(r)\vct r/r$ and the inertial forces, of a
function playing the role of a generalized potential energy:

\vspace*{6mm}
\eqn{110}{V(t,\vct r,\dot{\vct r})
=(\vct\omega,\vct r,\dot{\vct r})
+\frac12\|\vct\omega\times\vct r\|^2+\int f(r)\,dr.}

Using the first integral (99), the term that depends explicitly on the velocity
$\dot{\vct r}$ in relation (110) can be eliminated. Thus, taking the scalar
product of relation (109) with $\vct\omega$, we obtain
\eqn{111}{(\vct\omega,\vct r,\dot{\vct r})
+(\vct\omega,\vct r,\vct\omega\times\vct r)
=\vct\omega\cdot\vct\Omega_0.}

Relation (111) takes the equivalent form
\eqn{112}{(\vct\omega,\vct r,\dot{\vct r})
+\|\vct\omega\times\vct r\|^2
=\vct\omega\cdot\vct\Omega_0.}

Using relation (112), the potential energy (110) becomes
\eqn{113}{V(t,\vct r)=\vct\omega\cdot\vct\Omega_0
-\frac12\|\vct\omega\times\vct r\|^2+\int f(r)\,dr,}
or, equivalently, taking into account that
$\vct\Omega_0=\vct r_0\times(\vct v_0+\vct\omega_0\times\vct r_0)$:
\eqn{114}{V(t,\vct r)
=(\vct v_0+\vct\omega_0\times\vct r_0)\!\cdot\!(\vct\omega\times\vct r_0)
-\frac12\|\vct\omega\times\vct r\|^2+\int f(r)\,dr.}

In relation (114), $V$ depends explicitly on time through the function
$\vct\omega=\vct\omega(t)$. In the particular case $\vct\omega=\mathrm{const.}$,
the force field is positional, and (114) yields
\eqn{115}{V(\vct r)=-\frac12\|\vct\omega\times\vct r\|^2+\int f(r)\,dr.}

Result (115), valid in the particular case $\vct\omega=\mathrm{const.}$,
is the only one reported in the literature [102], [153].

\renewcommand{\eqprefix}{VI.3}

\renewcommand{\eqprefix}{VI.3}
\section*{VI.3. An Extension of the Tensorial Symbolic Representation of Vector-Valued Functions}
\subsection*{VI.3.1. General Results}

The results of VI.2 can be established under considerably more general conditions.
Let $\vct\omega\in V_3^*$ be a vector-valued function of a real variable,
differentiable on an interval $I\subseteq\mathbb R$. Let $\tens\omega$ be the
second-order skew-symmetric tensor corresponding to the vector-valued function
$\vct\omega$ [97]:
\eqn{116}{\vct\omega=\operatorname{vect}\tens\omega.}

As shown in VI.2.1, the second-order tensor Cauchy problem
\eqn{117}{\begin{cases}
\dot Q=-\tens\omega Q,\\ Q(0)=I_3,
\end{cases}\qquad \tens\omega^{T}=-\tens\omega,
\quad \dot Q=\dfrac{dQ}{dt}}
considered in the set of second-order tensor-valued functions of a real variable,
has a unique solution, described for every $t\geq0$ by a proper orthogonal tensor.

The transformation
\eqn{118}{F_{-\omega}:V_3^*\to V_3^*,\qquad
F_{-\omega}(\vct r)\stackrel{\rm def}{=}Q\vct r,}
where $Q$ is the solution of problem (117), was termed the rotation tensor with
angular velocity $-\vct\omega$.

Now consider the transformation
\eqn{119}{F_{\omega}^{*}:V_3^*\to V_3^*,\qquad
F_{\omega}^{*}(\vct r)\stackrel{\rm def}{=}Q^{T}\vct r,}
where $Q$ is the solution of problem (117). The essential properties of
transformation (119) are stated in the following theorem.

\medskip\noindent\textbf{Theorem VI.3.1.1.}

If $\vct\omega\in V_3^*$ is a vector-valued function of a real variable,
differentiable throughout its domain, transformation (119) has the following
properties:

\noindent$1^\circ$ $F_{\omega}^{*}$ is $\mathbb R$-linear;

\vspace*{5mm}
\noindent$2^\circ$
\[
F_{\omega}^{*}(\dot{\vct r}+\vct\omega\times\vct r)
=\dot F_{\omega}^{*}(\vct r),\qquad \forall\vct r\in V_3^*,
\]
for every differentiable function $\vct r$;

\noindent$3^\circ$
\[
F_{\omega}^{*}[\ddot{\vct r}+2\vct\omega\times\dot{\vct r}
+\vct\omega\times(\vct\omega\times\vct r)+\dot{\vct\omega}\times\vct r]
=\ddot F_{\omega}^{*}(\vct r),
\]
$\forall\vct r\in V_3^*$, with $\vct r$ twice differentiable;

\noindent$4^\circ$ $\|F_{\omega}^{*}(\vct r)\|=\|\vct r\|$,
$\forall\vct r\in V_3^*$;

\noindent$5^\circ$
\[
\left. F_{\omega}^{*}(\vct r)\right|_{t=0}=\vct r_0,
\qquad
\left.\dot F_{\omega}^{*}(\vct r)\right|_{t=0}
=\vct v_0+\vct\omega_0\times\vct r_0,
\]
where $\vct r_0=\vct r(0)$, $\vct v_0=\dot{\vct r}(0)$, and
$\vct\omega_0=\vct\omega(0)$;

\noindent$6^\circ$ $F_{\omega}^{*}$ is invertible and
$(F_{\omega}^{*})^{-1}=F_{-\omega}$;

\noindent$7^\circ$
$F_{\omega}^{*}(\vct r_1\times\vct r_2)
=F_{\omega}^{*}(\vct r_1)\times F_{\omega}^{*}(\vct r_2)$.

\medskip\noindent\underline{Proof}

$1^\circ$ The $\mathbb R$-linearity follows directly from definition (119).

$2^\circ$ We prove the matrix counterpart of the relation in item $2^\circ$. Let
\eqn{120}{\vct F=Q^T\vct r,}
where $Q$ is the orthogonal matrix corresponding to the proper orthogonal tensor
that solves tensor equation (117), and $\vct r$ is the column matrix associated
with the vector $\vct r$. Differentiating (120) with respect to $t\in\mathbb R$
gives
\eqn{121}{\dot{\vct F}=\dot Q^T\vct r+Q^T\dot{\vct r}.}
Transposing (117) shows that $Q^T$ is the solution of the Cauchy problem
\eqn{122}{\begin{cases}\dot Q^T=Q^T\tens\omega,\\Q^T(0)=I_3.
\end{cases}}
Using (122), relation (121) becomes
\eqn{123}{\dot{\vct F}=Q^T\tens\omega\vct r+Q^T\dot{\vct r}
=Q^T(\dot{\vct r}+\tens\omega\vct r).}
In tensorial language, this is equivalent to

\vspace*{5mm}
\eqn{124}{\dot F_{\omega}^{*}(\vct r)
=F_{\omega}^{*}(\dot{\vct r}+\vct\omega\times\vct r),
\qquad\forall\vct r\in V_3^*.}

$3^\circ$ Applying relation (124) twice gives successively
\begin{equation}\tag{\eqprefix.125}
\begin{aligned}
F_{\omega}^{*}[\ddot{\vct r}+2\vct\omega\times\dot{\vct r}
+\vct\omega\times(\vct\omega\times\vct r)+\dot{\vct\omega}\times\vct r]
&=F_{\omega}^{*}\!\left[\frac d{dt}(\dot{\vct r}+\vct\omega\times\vct r)
+\vct\omega\times(\dot{\vct r}+\vct\omega\times\vct r)\right]\\
&=\dot F_{\omega}^{*}(\dot{\vct r}+\vct\omega\times\vct r)
=\ddot F_{\omega}^{*}(\vct r).
\end{aligned}
\end{equation}

$4^\circ$ Transformation (119) is a proper rotation and is therefore an isometry:
\eqn{126}{\|F_{\omega}^{*}(\vct r)\|=\|\vct r\|,
\qquad\forall\vct r\in V_3^*.}

$5^\circ$ Setting $t=0$ in (120) gives
\eqn{127}{\left.\vct F\right|_{t=0}=Q^T(0)\vct r(0).}
Using (122)$_2$, we have $Q^T(0)=I_3$, and (127) yields
\eqn{128}{\left.\vct F\right|_{t=0}=I_3\vct r(0)=\vct r(0).}
With the notation
\eqn{129}{\vct r_0\stackrel{\rm not}{=}\vct r(0),
\qquad\forall\vct r\in V_3^*,}
the tensorial counterpart of relation (128) is
\eqn{130}{\left.F_{\omega}^{*}(\vct r)\right|_{t=0}=\vct r_0.}
Using (130) and relation (124), we obtain
\eqn{131}{\left.\dot F_{\omega}^{*}(\vct r)\right|_{t=0}
=\left.F_{\omega}^{*}(\dot{\vct r}+\vct\omega\times\vct r)\right|_{t=0}
=\vct v_0+\vct\omega_0\times\vct r_0,}
where $\vct r_0=\vct r(0)$, $\vct v_0=\dot{\vct r}(0)$, and
$\vct\omega_0=\vct\omega(0)$.

$6^\circ$ This follows from the orthogonality of the tensor $Q$,
$QQ^T=Q^TQ=I_3$, and from the definitions of the transformations $F_{-\omega}$
and $F_{\omega}^{*}$ in relations (118) and (119), respectively.

$7^\circ$ This follows from the fact that $F_{\omega}^{*}$ is a proper orthogonal
tensor, which preserves orientation and metric [5], and hence also preserves
vector products.

\subsection*{VI.3.2. Non-Inertial Motion in Central Positional Fields. First Integrals}

Transformation (119), with the general properties stated in VI.3.1.1, makes it
possible to address a particularly difficult problem: the study of relative
motion in central positional force fields with respect to non-inertial frames
undergoing arbitrary rotational motion. The principal results are summarized
in the following theorems.

\medskip\noindent\textbf{Theorem VI.3.2.1.}

The solution of the Cauchy problem
\eqn{132}{\begin{cases}
\ddot{\vct r}+2\vct\omega\times\dot{\vct r}
+\vct\omega\times(\vct\omega\times\vct r)
+\dot{\vct\omega}\times\vct r+f(r)\dfrac{\vct r}{r}=\vct0,\\
\vct r(0)=\vct r_0,\qquad \dot{\vct r}(0)=\vct v_0,
\end{cases}\quad \vct\omega\in V_3^*}
is obtained by applying the rotation tensor with angular velocity
$-\vct\omega$, $F_{-\omega}$, to the solution of the problem
\eqn{133}{\begin{cases}
\ddot{\vct r}+f(r)\dfrac{\vct r}{r}=\vct0,\\
\vct r(0)=\vct r_0,\\
\dot{\vct r}(0)=\vct v_0+\vct\omega_0\times\vct r_0,
\quad \vct\omega_0\stackrel{\rm not}{=}\vct\omega(0).
\end{cases}}
The function $f:\mathbb R_+^*\to\mathbb R$ satisfies the conditions ensuring
existence and uniqueness of the solution of Cauchy problem (133).

\medskip\noindent\underline{Proof}

Let
\eqn{134}{\vct\psi=F_{\omega}^{*}(\vct r).}
By item $3^\circ$ of Theorem VI.3.1.1,
\eqn{135}{\ddot{\vct\psi}=F_{\omega}^{*}[\ddot{\vct r}
+2\vct\omega\times\dot{\vct r}+\vct\omega\times(\vct\omega\times\vct r)
+\dot{\vct\omega}\times\vct r].}
By the homogeneity of transformation $F_{\omega}^{*}$ and property $4^\circ$,
\eqn{136}{F_{\omega}^{*}\!\left[f(r)\frac{\vct r}{r}\right]
=f(\psi)\frac{\vct\psi}{\psi}.}

Taking into account property $5^\circ$ of Theorem VI.3.1.1 and relations
(135), (136), applying transformation $F_{\omega}^{*}$ to Cauchy problem (132)
shows that the vector-valued function $\vct\psi$ in (134) is the solution of
\eqn{137}{\begin{cases}
\ddot{\vct\psi}+f(\psi)\dfrac{\vct\psi}{\psi}=\vct0,\\
\vct\psi(0)=\vct r_0,\\
\dot{\vct\psi}(0)=\vct v_0+\vct\omega_0\times\vct r_0.
\end{cases}}
Let $\vct\psi_0$ be the unique solution of Cauchy problem (137). From (134),
by property $6^\circ$ of Theorem VI.3.1.1,
\eqn{138}{\vct r=F_{-\omega}(\vct\psi_0),}
which proves the assertion.

\medskip\noindent\textbf{Theorem VI.3.2.2.}

The Cauchy problem
\eqn{139}{\begin{cases}
\ddot{\vct r}+2\vct\omega\times\dot{\vct r}
+\vct\omega\times(\vct\omega\times\vct r)
+\dot{\vct\omega}\times\vct r+f(r)\dfrac{\vct r}{r}=\vct0,
\quad\vct\omega\in V_3^*,\\
\vct r(0)=\vct r_0,\qquad\dot{\vct r}(0)=\vct v_0
\end{cases}}
admits the following first integrals:
\eqn{140}{\frac{\|\dot{\vct r}\|^2}{2}
+(\vct\omega,\vct r,\dot{\vct r})
+\frac12\|\vct\omega\times\vct r\|^2
+\int f(r)\,dr=\mathrm{const.},}
\eqn{141}{\vct r\times(\dot{\vct r}+\vct\omega\times\vct r)
=F_{-\omega}(\vct\Omega_0),\quad
\vct\Omega_0\stackrel{\rm not}{=}\vct r_0\times
(\vct v_0+\vct\omega_0\times\vct r_0),\quad
\vct\omega_0\stackrel{\rm not}{=}\vct\omega(0).}
Under the hypotheses of the theorem, the function
$f:\mathbb R_+^*\to\mathbb R$ is assumed to ensure existence and uniqueness for
problem (139) when $\vct\omega=\vct0$ and to possess antiderivatives. The symbol
$F_{-\omega}$ denotes the rotation tensor with angular velocity $-\vct\omega$.

\noindent\underline{Proof.}
By Theorem VI.3.2.1, if $\vct r$ is the solution of problem (139), then
$\vct\psi=F_{\omega}^{*}(\vct r)$ is the solution of
\eqn{142}{\begin{cases}
\ddot{\vct\psi}+f(\psi)\dfrac{\vct\psi}{\psi}=\vct0,\\
\vct\psi(0)=\vct r_0,\\
\dot{\vct\psi}(0)=\vct v_0+\vct\omega_0\times\vct r_0,
\quad\vct\omega_0\stackrel{\rm not}{=}\vct\omega(0).
\end{cases}}
Problem (142) admits the first integrals of energy and angular momentum [58], [87]:
\eqn{143}{\frac{\|\dot{\vct\psi}\|^2}{2}+\int f(\psi)\,d\psi=\mathrm{const.},}
\eqn{144}{\vct\psi\times\dot{\vct\psi}
=\vct r_0\times(\vct v_0+\vct\omega_0\times\vct r_0)
\stackrel{\rm not}{=}\vct\Omega_0.}

\mbox{}\par\vspace{8mm}
Using properties $2^\circ$ and $4^\circ$ of Theorem VI.3.1.1, we have
\eqn{145}{\|\dot{\vct\psi}\|^2
=\|F_{\omega}^{*}(\dot{\vct r}+\vct\omega\times\vct r)\|^2
=\|\dot{\vct r}+\vct\omega\times\vct r\|^2
=\|\dot{\vct r}\|^2+2(\vct\omega,\vct r,\dot{\vct r})
+\|\vct\omega\times\vct r\|^2.}
Using (145), relation (143) becomes
\eqn{146}{\frac{\|\dot{\vct r}\|^2}{2}
+(\vct\omega,\vct r,\dot{\vct r})
+\frac12\|\vct\omega\times\vct r\|^2
+\int f(r)\,dr=\mathrm{const.}}
Using properties $2^\circ$ and $7^\circ$, we obtain successively
\begin{equation}\tag{\eqprefix.147}\begin{aligned}
\vct\psi\times\dot{\vct\psi}
&=F_{\omega}^{*}(\vct r)\times\dot F_{\omega}^{*}(\vct r)\\
&=F_{\omega}^{*}(\vct r)\times
F_{\omega}^{*}(\dot{\vct r}+\vct\omega\times\vct r)\\
&=F_{\omega}^{*}[\vct r\times(\dot{\vct r}+\vct\omega\times\vct r)].
\end{aligned}\end{equation}
The invertibility of transformation $F_{\omega}^{*}$ and (147) give
\eqn{148}{\vct r\times(\dot{\vct r}+\vct\omega\times\vct r)
=F_{-\omega}(\vct\psi\times\dot{\vct\psi}).}
From (148) and (144), we obtain
\eqn{149}{\vct r\times(\dot{\vct r}+\vct\omega\times\vct r)
=F_{-\omega}(\vct\Omega_0).}
This completes the proof.

\subsection*{VI.3.3. Applications}
For particular forms of the force field $\vct F=f(r)\vct r/r$, Theorems
VI.3.2.1 and VI.3.2.2 yield consequences with important practical implications.
For an elastic force field, the following result holds.

\noindent\textbf{Consequence VI.3.3.1.}

The solution of the Cauchy problem
\eqn{150}{\begin{cases}
\ddot{\vct r}+2\vct\omega\times\dot{\vct r}
+\vct\omega\times(\vct\omega\times\vct r)
+\dot{\vct\omega}\times\vct r+\omega_*^2\vct r=\vct0,
\quad\vct\omega\in V_3^*,\ \omega_*\in\mathbb R_+^*,\\
\vct r(0)=\vct r_0,\qquad\dot{\vct r}(0)=\vct v_0
\end{cases}}
is given by the vector-valued function
\eqn{151}{\vct r=F_{-\omega}\!\left(\vct r_0\cos\omega_*t
+\frac{\vct v_0+\vct\omega_0\times\vct r_0}{\omega_*}\sin\omega_*t\right).}

\mbox{}\par\vspace{8mm}
Here $F_{-\omega}$ is the rotation tensor with angular velocity
$-\vct\omega$, where $\vct\omega\in V_3^*$ is differentiable, and
$\vct\omega_0\stackrel{\rm not}{=}\vct\omega(0)$.

\noindent\underline{Proof.}
Problem (133) in Theorem VI.3.2.1 becomes
\eqn{152}{\begin{cases}
\ddot{\vct r}+\omega_*^2\vct r=\vct0,\\
\vct r(0)=\vct r_0,\\
\dot{\vct r}(0)=\vct v_0+\vct\omega_0\times\vct r_0,
\end{cases}}
with solution
$\vct r^{*}=\vct r_0\cos\omega_*t+
(\vct v_0+\vct\omega_0\times\vct r_0)\sin(\omega_*t)/\omega_*$,
$t\geq0$. Result (151) follows by applying Theorem VI.3.2.1 to the solution
$\vct r^{*}$.

\noindent\underline{Remark.}
Consequence VI.3.3.1 generalizes the properties established in VI.1.3.3 and
VI.2.2.3, relation (70), for the Foucault pendulum when the angular velocity of
the non-inertial frame varies in both magnitude and direction.

For the Newtonian force field, the following result holds.

\noindent\textbf{Consequence VI.3.3.2.}

The solution of the Cauchy problem
\eqn{153}{\begin{cases}
\ddot{\vct r}+2\vct\omega\times\dot{\vct r}
+\vct\omega\times(\vct\omega\times\vct r)
+\dot{\vct\omega}\times\vct r+\dfrac{k}{r^3}\vct r=\vct0,
\quad\vct\omega\in V_3^*,\ k\in\mathbb R_+^*,\\
\vct r(0)=\vct r_0,\qquad\dot{\vct r}(0)=\vct v_0
\end{cases}}
is bounded if and only if $E<0$, where
\eqn{154}{E=\frac{\|\vct v_0\|^2}{2}
+(\vct\omega_0,\vct r_0,\vct v_0)
+\frac12\|\vct\omega_0\times\vct r_0\|^2-\frac{k}{r_0},
\quad\vct\omega_0\stackrel{\rm not}{=}\vct\omega(0).}

\noindent\underline{Proof.}
In this case, problem (133) becomes
\eqn{155}{\begin{cases}
\ddot{\vct r}+\dfrac{k}{r^3}\vct r=\vct0,\\
\vct r(0)=\vct r_0,\\
\dot{\vct r}(0)=\vct v_0+\vct\omega_0\times\vct r_0.
\end{cases}}
Problem (155) is classical in theoretical mechanics [87], [89]. The hodograph
of the solution is a conic, possibly degenerate, lying in the plane determined
by the vectors $\vct r(0)$ and $\dot{\vct r}(0)$, with one focus at the origin.

\mbox{}\par\vspace{8mm}
With respect to this focus, the hodograph is traversed with constant areal
velocity. The type of conic is determined by the sign of the expression [87], [89]
\eqn{156}{E=\frac{\|\dot{\vct r}(0)\|^2}{2}-\frac{k}{r_0}.}
If $E<0$, the conic is elliptic; if $E=0$, it is parabolic; and if $E>0$, it is
hyperbolic. Using the initial conditions in (155), the expression becomes
\eqn{157}{E=\frac{\|\vct v_0\|^2}{2}
+(\vct\omega_0,\vct r_0,\vct v_0)
+\frac12\|\vct\omega_0\times\vct r_0\|^2-\frac{k}{r_0}.}
Let $\vct r^{*}$ be the solution of problem (155). By Theorem VI.3.2.1, the
solution of problem (153) is $\vct r=F_{-\omega}(\vct r^{*})$. Since
$F_{-\omega}$ is an isometry and the solution of problem (155) is bounded only
for $E<0$, the assertion of Consequence VI.3.3.2 follows. For $E\geq0$, the
solution is unbounded.

The expression determining whether the solution is bounded or unbounded is the
value of the first integral (140), evaluated for the Newtonian force field.

\noindent\underline{Remark.}
The trajectory of the solution of problem (153) is obtained by rotating, with
angular velocity $-\vct\omega$, the plane of the conic solving problem (155).
This permits a qualitative study based on the geometric properties of
transformation $F_{-\omega}$ and of the classical solution, an approach not
developed in the present work.

\noindent\textbf{Consequence VI.3.3.3.}

The velocity hodograph in the motion of a point electric charge $P(m,q)$ in an
arbitrary magnetic field of induction $\vct B$ is a spherical curve obtained by
rotating the velocity vector $\vct v$, $\|\vct v\|=v_0$, with angular velocity
$\vct\omega_L=-q\vct B/m$ relative to the inertial frame.

\noindent\underline{Proof.}
The velocity problem is described by the Cauchy problem
\eqn{158}{\begin{cases}
m\dfrac{d\vct v}{dt}=q\vct v\times\vct B,\\
\vct v(0)=\vct v_0.
\end{cases}}

With the notation $\vct\omega\stackrel{\rm not}{=}q\vct B/m$ and
$\dot{\vct v}\stackrel{\rm not}{=}d\vct v/dt$, problem (158) becomes
\eqn{159}{\begin{cases}
\dot{\vct v}+\vct\omega\times\vct v=\vct0,\\
\vct v(0)=\vct v_0.
\end{cases}}

Let $\vct\psi=F_{\omega}^{*}(\vct v)$. By properties $2^\circ$ and
$5^\circ$ of Theorem VI.3.1.1, the vector-valued function $\vct\psi$ satisfies
the Cauchy problem
\eqn{160}{\begin{cases}
\dot{\vct\psi}=\vct0,\\
\vct\psi(0)=\vct v_0.
\end{cases}}

The unique solution of problem (160) is
\eqn{161}{\vct\psi^{*}=\vct v_0.}

By property $6^\circ$ of Theorem VI.3.1.1,
$\vct v=F_{-\omega}(\vct v_0)$. The solution of problem (158), which is also
the equation of the velocity hodograph, is
\eqn{162}{\vct v=F_{-\frac qm\vct B}(\vct v_0),}
which proves Consequence VI.3.3.3.

The results of this consequence extend the well-known Larmor theorem [103],
[99] to the motion of electric charges in an arbitrary magnetic field.

% D6, D66--D67: original pages 181--204, editable recovery
\subsectionstar{VI.3.4. The Two-Body Problem in Non-Inertial Frames}

Consider a mechanical system consisting of two interacting material particles
$P_1$ and $P_2$, of masses $m_1$ and $m_2$, respectively. We study the motion
of the system relative to a non-inertial frame undergoing arbitrary rotational
motion, with angular velocity $\vct\omega\in V_3^*$, a differentiable
vector-valued function of a scalar variable.

We assume that the interaction forces satisfy the principle of mutual actions.

\begin{figure}[htbp]
\centering
\begin{minipage}{.46\textwidth}\centering\figeight\par\small\textnormal{Fig. 8}\end{minipage}\hfill
\begin{minipage}{.46\textwidth}\centering\fignine\par\small\textnormal{Fig. 9}\end{minipage}
\end{figure}

Denote by $\vct r_k^0$ and $\vct v_k^0$, $k=1,2$, the position vectors and,
respectively, the velocities of the particles at the initial instant $t=0$
(Fig. 9). Their motion is described by the Cauchy problems

\eqn{163}{\begin{cases}
m_1\ddot{\vct r}_1+2m_1\vct\omega\times\dot{\vct r}_1
+m_1\vct\omega\times(\vct\omega\times\vct r_1)
+m_1\dot{\vct\omega}\times\vct r_1=\vct F_{12},\\
\vct r_1(0)=\vct r_1^0,\\
\dot{\vct r}_1(0)=\vct v_1^0.
\end{cases}}

\eqn{164}{\begin{cases}
m_2\ddot{\vct r}_2+2m_2\vct\omega\times\dot{\vct r}_2
+m_2\vct\omega\times(\vct\omega\times\vct r_2)
+m_2\dot{\vct\omega}\times\vct r_2=\vct F_{21},\\
\vct r_2(0)=\vct r_2^0,\\
\dot{\vct r}_2(0)=\vct v_2^0.
\end{cases}}

In (163) and (164), $\vct r_k$ is the position vector, at a given instant, of
particle $P_k$, $k=1,2$, relative to the frame in which the motion is studied.
The symbol $\vct F_{ij}$ denotes the force acting on particle $P_i$ due to its
interaction with particle $P_j$, $i\ne j$, $i,j\in\{1,2\}$. Since the
interaction forces satisfy the principle of mutual actions, the following
relations hold:

\eqn{165}{\vct F_{12}+\vct F_{21}=\vct0,}
\eqn{166}{\vct F_{12}\times(\vct r_2-\vct r_1)=\vct0.}

\subsubsectionstar{VI.3.4.1. The Two-Body Problem in Inertial Frames}

If $\vct\omega=\vct0$, the frame relative to which the motion is studied is
inertial. The classical form of the ``two-body problem'' [19], [59], [160] is
obtained. In this case, problems (163) and (164) become

\eqn{167}{\begin{cases}
m_1\ddot{\vct r}_1=\vct F_{12},\\
\vct r_1(0)=\vct r_1^0,\\
\dot{\vct r}_1(0)=\vct v_1^0,
\end{cases}}
\eqn{168}{\begin{cases}
m_2\ddot{\vct r}_2=\vct F_{21},\\
\vct r_2(0)=\vct r_2^0,\\
\dot{\vct r}_2(0)=\vct v_2^0,
\end{cases}}
subject to the conditions given by relations (165) and (166).

In the classical treatment, the study of the motion of the two interacting
point bodies decouples into the study of two separate particle-type motions.
The procedure is summarized below.

a) Let \emph{$C$} be the \textbf{center of mass} of the system, with position
vector (Fig. 8)

\eqn{169}{{\vct r}_{c} = \frac{m_{1}{\vct r}_{1} + m_{2}{\vct r}_{2}}{m_{1} + m_{2}}\text{.}}
From (167) and (168), using condition (165), the motion of the center of mass
of the system is described by the Cauchy problem

\eqn{170}{
\begin{cases}
(m_1+m_2)\ddot{\vct r}_c=\vct0,\\
\vct r_c(0)=\vct r_c^{0},\\
\dot{\vct r}_c(0)=\vct v_c^{0}.
\end{cases}}

In (170), \({\vct r}_{c}^{0}\) and \({\vct v}_{c}^{0}\) denote the position
vector and velocity of the center of mass at the initial instant $t=0$ (Fig. 9):

\({\vct r}_{c}^{0} = \frac{m_{1}{\vct r}_{1}^{0} + m_{2}{\vct r}_{2}^{0}}{m_{1} + m_{2}}\),

\eqn{171}{{\vct v}_{c}^{0} = \frac{m_{1}{\vct v}_{1}^{0} + m_{2}{\vct v}_{2}^{0}}{m_{1} + m_{2}}\text{.}}
The solution of problem (170) is

\eqn{172}{{\vct r}_{c} = {\vct r}_{c}^{0} + {\vct v}_{c}^{0}t}; $t\in\mathbb R$,

hence the \textbf{center of mass of the system moves uniformly along a straight line}.

Problem (170) models the \textbf{inertial motion of a particle whose mass equals}
\((m_{1} + m_{2})\)\textbf{, the total mass of the system}.

b) We now study the motion of particle P$_2$ relative to P$_1$. Set

\eqn{173}{\vct r = {\vct r}_{2} - {\vct r}_{1},}
the position vector of particle P$_2$ (m$_2$) relative to P$_1$ (m$_1$). From
(167) and (168), using condition (165), one obtains the differential equation
satisfied by the vector-valued function (173):

\eqn{174}{\frac{m_{1} \cdot m_{2}}{m_{1} + m_{2}}\ddot{\vct r} = {\vct F}_{\text{21}}}.

Using condition (166), it follows that the force \({\vct F}_{\text{21}}\) is of
``central type'' with respect to P$_1$:

\eqn{175}{{\vct F}_{\text{21}} = F\frac{\vct r}{r}},

where $F$ denotes a scalar function. If we set

\eqn{176}{m = \frac{m_{1} \cdot m_{2}}{m_{1} + m_{2}}}
the ``reduced mass'' of the system [14], [59], then the vector-valued function
(173) is the solution of the Cauchy problem

\eqn{177}{\begin{array}{r}
\left\{ m\ddot{\vct r} = F\frac{\vct r}{r} \middle| \right.\ \left\{ \vct r(0) = {\vct r}^{0} \middle| \right.\ 
\end{array}}

where we have set

\eqn{178}{\begin{array}{r}
{\vct r}^{0} = {\vct r}_{2}^{0} - {\vct r}_{1}^{0} \\
{\vct v}^{0} = {\vct v}_{2}^{0} - {\vct v}_{1}^{0}\text{.}
\end{array}}

Problem (177) models the \textbf{motion of a particle whose mass equals the
``reduced mass'' of the system,} \(m = \frac{m_{1} \cdot m_{2}}{m_{1} + m_{2}}\)
\textbf{, in a central force field}.

\paragraph{Remarks.}

1$^\circ$ Knowledge of the solutions of the particle-type problems (170) and
(177) completely determines the \textbf{laws of motion} of the point bodies
P$_1$ (m$_1$) and P$_2$ (m$_2$). Indeed, (169) and (173) give

\eqn{179}{\begin{array}{r}
{\vct r}_{1} = {\vct r}_{c} - \frac{m_{2}}{m_{1} + m_{2}}\vct r, \\
{\vct r}_{2} = {\vct r}_{c} + \frac{m_{1}}{m_{1} + m_{2}}\vct r\text{.}
\end{array}}

Relations (179) give the laws of motion of particles P$_1$ (m$_1$) and P$_2$
(m$_2$) relative to an inertial frame with origin at $O$.

The solutions of problems (170) and (177) also determine the \textbf{global
characteristics} of the mechanical system under consideration.

The \textbf{linear momentum} of the mechanical system,

\eqn{180}{\vct P = \sum_{i = 1}^{2}{m_{i}{\dot{\vct r}}_{i}}}
is expressed in terms of the solution of problem (170), as follows from (179)
and (180):

\eqn{181}{\vct P = (m_{1} + m_{2}){\dot{\vct r}}_{c},}
and is conserved during the motion, by (172) and (181):

\eqn{182}{\vct P = (m_{1} + m_{2}){\vct v}_{c}^{0} = \text{const}\text{.}}
The \textbf{angular momentum} of the mechanical system about the pole $O$,

\eqn{183}{{\vct K}_{0} = \sum_{i = 1}^{2}{{\vct r}_{i}\times m_{i}{\dot{\vct r}}_{i}}}
is expressed in terms of the solutions of problems (170) and (177), as follows
from (179) and (183):

\eqn{184}{{\vct K}_{0} = \vct r\times m\dot{\vct r} + {\vct r}_{c}\times (m_{1} + m_{2}){\dot{\vct r}}_{c}},

and is conserved, in accordance with (170) and (177):

\eqn{185}{{\vct K}_{0} = {\vct r}^{0}\times m\vct v + {\vct r}_{c}^{0}\times (m_{1} + m_{2}){\vct v}_{c}^{0} = \text{const}\text{.}}
The \textbf{kinetic energy} of the mechanical system,

\eqn{186}{E_{\text{cin}} = \frac{1}{2}\sum_{i = 1}^{2}{m_{i}{\dot{\vct r}}_{i}^{2}}}
is expressed in terms of the solutions of problems (170) and (177), as follows
from (179) and (186):

\eqn{187}{E_{\text{cin}} = \frac{1}{2}m{\dot{\vct r}}^{2} + \frac{1}{2}(m_{1} + m_{2}){\dot{\vct r}}_{c}^{2}}.

If the central force field in (177) is positional, that is, $F=F(r)$, problem
(177) admits the \textbf{first integral of energy}:

\eqn{188}{\frac{m{\dot{\vct r}}^{2}}{2} + V(r) = E = \text{const}\text{.},}
where \emph{\textbf{V}} denotes the potential energy of the central positional
force field,

\eqn{189}{V\text{=-}\int_{}^{}{F(r)\text{dr}},}
and $E$ denotes the total energy:

\eqn{190}{E = \frac{m{\vct v}^{\text{02}}}{2} + V(r^{0})\text{.}}
Under these hypotheses, the \textbf{energy of the mechanical system is conserved}:

\eqn{191}{E_{\text{cin}} + V\text{= const,}}
as follows from (187), (188), and (182).

2$^\circ$ The motion of the point bodies P$_1$ (m$_1$) and P$_2$ (m$_2$)
relative to the center of mass of the system, whose motion is known from (172),
can likewise be described completely from the solutions of problem (177).

Let \({\vct r}_{1}^{'}\) and \({\vct r}_{2}^{'}\) be the position vectors of
particles P$_1$ (m$_1$) and P$_2$ (m$_2$), respectively, relative to the center
of mass \emph{$C$} of the system. From the relations

\eqn{192}{\begin{array}{r}
\left\{ m_{1}{\vct r}_{1}^{'} + m_{2}{\vct r}_{2}^{'} = 0 \middle| \right.\ 
\end{array}}

it follows that

\eqn{193}{\begin{array}{r}
\left\{ {\vct r}_{1}^{'}\text{=-}\frac{m_{2}}{m_{1} + m_{2}}\vct r \middle| \right.\ 
\end{array}}

Relations (193), with \(\vct r\) the solution of Cauchy problem (177), represent
the \textbf{laws of motion} of particles P$_1$ (m$_1$) and P$_2$ (m$_2$)
relative to an inertial frame whose origin is at the center of mass.

Relative to this inertial frame, called the \textbf{center-of-mass frame}, the
\textbf{global characteristics} of the mechanical system are as follows.

The \textbf{linear momentum}, \({\vct P}^{c} = \sum_{i = 1}^{2}{m_{i}{\dot{\vct r}}_{i}^{'}}\):

\eqn{194}{{\vct P}^{c} = 0}.

The \textbf{angular momentum}, \({\vct K}_{c} = \sum_{i = 1}^{2}{{\vct r}_{i}^{'}\times m_{i}{\dot{\vct r}}_{i}^{'}}\):

\eqn{195}{{\vct K}_{c} = \vct r\times m\dot{\vct r}}.

The \textbf{kinetic energy}, \(E_{\text{cin}}^{c} = \frac{1}{2}\sum_{i = 1}^{2}{m_{i}{\dot{\vct r}}_{i}^{}}\):

\eqn{196}{E_{\text{cin}}^{c} = \frac{1}{2}m{\dot{\vct r}}^{2}},

as follows, after elementary calculations, from relations (193).

The \textbf{angular momentum} of the mechanical system relative to the
center-of-mass frame is conserved, as follows from (177), whose right-hand side
is a central force, and (195):

\eqn{197}{{\vct K}_{c} = {\vct r}^{0}\times {\vct v}^{0} = \text{const}\text{.}}
This fact has a special significance, which will become apparent below.

The Cauchy problems describing the motion of particles P$_1$ (m$_1$) and P$_2$
(m$_2$) relative to the center-of-mass frame follow from (193) and (177):

\eqn{198}{\begin{cases}
m_1\ddot{\vct r}'_1=-F\dfrac{\vct r'_1}{r'_1},\\
\vct r'_1(0)=\vct r'_{01},\\ \dot{\vct r}'_1(0)=\vct v'_{01}.
\end{cases}}

\eqn{199}{\begin{cases}
m_2\ddot{\vct r}'_2=F\dfrac{\vct r'_2}{r'_2},\\
\vct r'_2(0)=\vct r'_{02},\\ \dot{\vct r}'_2(0)=\vct v'_{02}.
\end{cases}}

In (198) and (199), \(\vct r_{0k}^{'}\) and \(\vct v_{0k}^{'}\) denote,
respectively, the position vector and velocity of particle P$_k$ (m$_k$),
\(k=1,2\), relative to the center-of-mass frame at the initial instant $t=0$.

The trajectories of the motions described by (198) and (199) are plane curves,
because the forces acting on the particles are central. Since (193) gives

\eqn{200}{{\vct r}_{1}^{'}\text{=-}\frac{m_{2}}{m_{1}}{\vct r}_{2}^{'},}
the trajectories of the two particles in the center-of-mass frame are
homothetic curves with ratio \(-\frac{m_2}{m_1}\). Hence the planes of the two
trajectories coincide with the plane through the center of mass $C$ of the
system, whose normal is given by the angular-momentum vector (197). If
\(\vct r^{0}\times {\vct v}^{0}=0\), the trajectories of the two particles in
the center-of-mass frame are straight lines; the homothety property (200)
continues to hold in this case.

If the interaction forces between particles P$_1$ (m$_1$) and P$_2$ (m$_2$)
depend only on their relative distance, the energy of the mechanical system
relative to the center-of-mass frame is conserved, as follows from (188) and
(196):

\eqn{201}{E_{\text{cin}}^{c} + V = \text{const},}
where notation (189) has been used for the potential energy
\emph{\textbf{V}} of the interaction forces.

\subsubsectionstar{VI.3.4.2. The Two-Body Problem in Non-Inertial Frames with Arbitrary Rotation}

Using the results of VI.2, we shall show that, even when the two-body problem is studied in a non-inertial frame rotating with angular velocity \(\vct\omega\in V_3^*\) (where \(\vct\omega\) is an arbitrary differentiable vector-valued function of a scalar variable), the motion can be described by considering two ``decoupled'' particle-type problems.

\fourthsectionstar{VI.3.4.2.1. Motion of the Center of Mass}

Let \emph{$C$} be the center of mass of the system, with position vector \({\vct r}_{c}\) given by relation (169). From (163) and (164), using condition (165), it follows that the motion of the center of mass \emph{$C$} is described by the Cauchy problem

\eqn{202}{\begin{cases}
\ddot{\vct r}_c+2\vct\omega\times\dot{\vct r}_c
+\vct\omega\times(\vct\omega\times\vct r_c)+\dot{\vct\omega}\times\vct r_c=\vct0,\\
\vct r_c(0)=\vct r_c^0,\\ \dot{\vct r}_c(0)=\vct v_c^0.
\end{cases}}

In (202), \({\vct r}_{c}^{0}\) and \({\vct v}_{c}^{0}\) denote, respectively, the position vector and velocity of the center of mass at the initial instant, $t=0$, as specified by relations (171).

The solution of problem (202) is given by Theorem VI.3.2.1.

First consider the Cauchy problem

\eqn{203}{\begin{cases}
\ddot{\vct r}=\vct0,\\ \vct r(0)=\vct r_c^0,\\
\dot{\vct r}(0)=\vct v_c^0+\vct\omega_0\times\vct r_c^0.
\end{cases}}

whose solution is

\eqn{204}{\vct r = {\vct r}_{c}^{0} + ({\vct v}_{c}^{0} + {\vct\omega}_{0}\times {\vct r}_{c}^{0})t;{t\in\mathbb R}_{+}}.

According to Theorem VI.3.2.1, the solution of problem (202) is obtained by applying the rotation tensor with angular velocity \(-\vct\omega\) to the vector-valued function (204):

\eqn{205}{{\vct r}_{c} = F_{- \vct\omega}\lbrack{\vct r}_{c}^{0} + ({\vct v}_{c}^{0} + {\vct\omega}_{0}\times {\vct r}_{c}^{0})t\rbrack;{t\in\mathbb R}_{+}}.

\paragraph{Remarks.}

1$^\circ$ Since the transformation \(F_{- \vct\omega}\colon V_3^*\to V_3^*\) is also an isometry, (205) yields

\eqn{206}{\left| {\vct r}_{c} \right| = \left| {\vct r}_{c}^{0} + (\vct v_c^0+\vct\omega_0\times\vct r_c^0)t \right|;{t\in\mathbb R}_{+}\text{.}}
It follows from (206) that, at any fixed time $t\in\mathbb R_+$, the center of mass of the system lies on the sphere centered at $O$ and having radius \(R = \left| {\vct r}_{c}^{0} + ({\vct v}_{c}^{0} + {\vct\omega}_{0}\times {\vct r}_{c}^{0})t \right|\).

2$^\circ$ If the initial conditions in (163) and (164) are such that

\eqn{207}{{\vct v}_{c}^{0} + {\vct\omega}_{0}\times {\vct r}_{c}^{0} = 0},

the hodograph of the vector-valued function describing the motion of the center of mass is a spherical curve obtained by rotating the vector \(\vct r_{c}^{0}\) with angular velocity \(-\vct\omega\):

\eqn{208}{{\vct r}_{c} = F_{- \vct\omega}({\vct r}_{c}^{0})}.

3$^\circ$ If the angular-velocity vector \(\vct\omega\) has a fixed direction,

\eqn{209}{\vct\omega = w\vct u},

where $w\colon\mathbb R_+\to\mathbb R$ is a differentiable real-valued function of a real variable and \(\vct u\) is a constant unit vector, the rotation tensor \(F_{- \vct\omega}\) is written explicitly as in (12):

\eqn{210}{\begin{array}{r}
F_{- \vct\omega}:V_{3}^{} \rightarrow V_{3}^{} \\
F_{- \vct\omega}(\vct r) = \frac{\vct\omega \cdot \vct r}{{\vct\omega}^{2}}\vct\omega - \sin\alpha(t)\frac{\vct\omega\times \vct r}{w} - \cos\alpha(t)\frac{\vct\omega\times (\vct\omega\times \vct r)}{{\vct\omega}^{2}};

\end{array}}

where \(\alpha(t) = \int_{0}^{t}{w(t)dt}\), $t\ge 0$.

The solution (205) of problem (202) can now be written explicitly.

In this case, \textbf{the hodograph of the vector-valued function describing the motion of the center of mass is a curve lying on the ruled surface} obtained by rotating, with angular velocity \(-\vct\omega\), the vector-valued function (204), whose hodograph is a half-line. If \(\vct\omega = \text{const}\), the ruled surface is conical or cylindrical, and the hodograph of the position vector of the center of mass is a conical or cylindrical helix. \textbf{The trajectory of the center of mass is then one of the geodesics of these surfaces}.

In particular, if condition (207) is satisfied, the trajectory of the center of mass is circular. The magnitude of its position vector is fixed at $t=0$ and remains constant throughout the motion. Under the conditions given by (207) and (209), the vector \({\vct r}_{c}\) undergoes precession about the fixed direction of \(\vct\omega\), with angular velocity \(-\vct\omega\), as follows from relation (208).

4$^\circ$ By Theorem VI.3.2.2, problem (202) admits two first integrals: a \textbf{vector first integral} given by

\eqn{211}{{\vct r}_{c}\times ({\dot{\vct r}}_{c} + \vct\omega\times {\vct r}_{c}) = F_{- \vct\omega}({\vct\Omega}_{0});{\vct\Omega}_{0}\overset{\text{not}}{=}{\vct r}_{c}^{0}\times ({\vct v}_{c}^{0} + {\vct\omega}_{0}\times {\vct r}_{c}^{0}),}
and a \textbf{scalar first integral} having an energy interpretation:

\eqn{212}{\frac{{\dot{\vct r}}_{c}^{2}}{2} + (\vct\omega,{\vct r}_{c},{\dot{\vct r}}_{c}) + \frac{1}{2}(\vct\omega\times {\vct r}_{c})^{2} = \text{const}}.

As noted in VI.2.2.3, relation (113), if the vector-valued function \(\vct\omega\) has a fixed direction, the first integral (212) takes the form

\eqn{213}{\frac{{\dot{\vct r}}_{c}^{2}}{2} + \vct\omega{\vct\Omega}_{0} - \frac{1}{2}(\vct\omega\times {\vct r}_{c})^{2} = \text{const},}
and, under the assumption \(\vct\omega=\text{const}\), it becomes

\eqn{214}{\frac{{\dot{\vct r}}_{c}^{2}}{2} - \frac{1}{2}(\vct\omega\times {\vct r}_{c})^{2} = \text{const}\text{.}}
\fourthsectionstar{VI.3.4.2.2. Relative Motion}

The motion of particle P$_2$ (m) relative to P$_1$ (m) is described by the vector-valued function

\eqn{215}{\vct r = {\vct r}_{2} - {\vct r}_{1}}
In relation (215), \({\vct r}_{1}\) and \({\vct r}_{2}\) are the solutions of the Cauchy problems (163) and (167). Using condition (165) from (163) and (164), it follows that the vector-valued function defined by (215) satisfies the differential equation

\eqn{216}{\ddot{\vct r} + 2\vct\omega\times \dot{\vct r} + \vct\omega\times (\vct\omega\times \vct r) + \dot{\vct\omega}\times \vct r = \frac{{\vct F}_{\text{21}}}{m}}
In (216), \(m = \frac{m_{1} \cdot m_{2}}{m_{1} + m_{2}}\) denotes the ``reduced mass'' of the system, as in (176). By (166), \({\vct F}_{\text{21}}\) is a force of ``central type'' with respect to P$_1$:

\eqn{217}{{\vct F}_{\text{21}} = F\frac{\vct r}{r}},

where $F$ denotes a scalar function.

Consequently, using the notation (178), the vector-valued function \(\vct r\) in (215) is the solution of the Cauchy problem

\eqn{218}{\begin{array}{r}
\left\{ \ddot{\vct r} + 2\vct\omega\times \dot{\vct r} + \vct\omega\times (\vct\omega\times \vct r) + \dot{\vct\omega}\times \vct r = \frac{F}{m}\frac{\vct r}{r}, \middle| \right.\ \left\{ \vct r(0) = {\vct r}^{0}, \middle| \right.\ 
\end{array}}

Problem (218) describes \textbf{the motion of a particle whose mass equals the reduced mass of the mechanical system, relative to the non-inertial frame in which the two-body problem is studied.}

The solution of problem (218) is described by Theorem VI.3.2.1. Consider the problem

\eqn{219}{\begin{array}{r}
\left\{ \ddot{\vct r} = \frac{F}{m}\frac{\vct r}{r} \middle| \right.\ \left\{ \vct r(0) = {\vct r}^{0} \middle| \right.\ 
\end{array}}

If the solution of this problem is denoted by \({\vct r}^{}\), then, by Theorem VI.3.2.1, the solution of problem (218) is obtained by applying the rotation tensor with angular velocity \(-\vct\omega\) to the solution of problem (219):

\eqn{220}{\vct r = F_{- \vct\omega}({\vct r}^{})}
\paragraph{Remarks.}

1$^\circ$ By Theorem VI.3.2.2, problem (218) admits the \textbf{vector first integral}

\eqn{221}{\vct r\times (\dot{\vct r} + \vct\omega\times \vct r) = F_{- \vct\omega}({\vct\Omega}_{0})},

where \(\vct\Omega_{0}{= \vct r}^{0}\times ({\vct v}^{0} + {\vct\omega}_{0}\times {\vct r}^{0})\).

Relation (221) shows that the hodograph of the vector-valued function \(\vct\Omega = \vct r\times (\dot{\vct r} + \vct\omega\times \vct r)\) is a spherical curve obtained by rotating its initial value with angular velocity \(-\vct\omega\). In particular, (221) implies that the magnitude of this vector quantity is conserved during the motion:

\eqn{222}{\left| \vct r\times (\dot{\vct r} + \vct\omega\times \vct r) \right| = \left| \vct\Omega_{0} \right|}
2$^\circ$ If the angular-velocity vector \(\vct\omega\) has a fixed direction, then, taking (210) into account, the first integral (221) takes the form

\eqn{223}{\begin{array}{r}
\vct r\times (\dot{\vct r} + \vct\omega\times \vct r) = \frac{\vct\omega \cdot \vct\Omega}{{\vct\omega}^{2}}\vct\omega - \sin\alpha(t)\frac{\vct\omega\times \vct\Omega}{w} - \cos\alpha(t)\frac{\vct\omega\times (\vct\omega\times \vct\Omega)}{{\vct\omega}^{2}}, \\
\text{                                                                                          }{t\in\mathbb R}_{+},
\end{array}}

where $w$ is given by relation (209) and \(\alpha(t) = \int_{0}^{t}{w(t)dt}\).

In this case, the hodograph of the vector-valued function \(\vct\Omega\) is a circle (or circular arc) lying in a plane perpendicular to the fixed direction of \(\vct\omega\). During the motion, the vector \(\vct\Omega\)

sweeps out the surface of a right circular cone, undergoing precession with angular velocity \(-\vct\omega\) about the fixed direction of \(\vct\omega\) (Fig.~10).

\begin{figure}[htbp]
\centering\figten\par\small\textnormal{Fig. 10}
\end{figure}

3$^\circ$ If the interaction forces between the particles of the system depend only on their relative distance, then $F=F(r)$ and, by Theorem VI.3.2.2, problem (218) admits a \textbf{scalar first integral} having an energy interpretation:

\eqn{224}{\frac{m{\dot{\vct r}}^{2}}{2} + m(\vct\omega,\vct r,\dot{\vct r}) + \frac{m}{2}(\vct\omega\times \vct r)^{2} - \int_{}^{}{F(r)\text{dr}} = \text{const},}
The first integral (224) reveals the existence of a function having the role of a generalized potential energy:

\eqn{225}{V(t,\vct r,\dot{\vct r}) = m(\vct\omega,\vct r,\dot{\vct r}) + \frac{m}{2}(\vct\omega\times \vct r)^{2} - \int_{}^{}{F(r)\text{dr}},}
in terms of which the first integral (224) may be written as

\eqn{226}{\frac{m{\dot{\vct r}}^{2}}{2} + V(t,\vct r,\dot{\vct r}) = \text{const}\text{.}}
As noted in the remarks of VI.2.2.3, if the vector-valued function \(\vct\omega\) has a fixed direction, the generalized potential energy does not depend on \(\dot{\vct r}\). One obtains

\eqn{227}{V(t,\vct r) = m\,\vct\omega\cdot\vct\Omega_0 - \frac{m}{2}(\vct\omega\times \vct r)^{2} - \int_{}^{}{F(r)\text{dr}},}
with the property

\eqn{228}{\frac{m{\dot{\vct r}}^{2}}{2} + V(t,\vct r) = \text{const}\text{.}}
If \(\vct\omega = \text{const}\text{.}\), one obtains the classical result [102], [153]:

\eqn{229}{V(\vct r)\text{=-}\frac{m}{2}(\vct\omega\times \vct r)^{2} - \int_{}^{}{F(r)\text{dr}},}
with the property

\eqn{230}{\frac{m{\dot{\vct r}}^{2}}{2} + V(\vct r) = \text{const}\text{.}}
\fourthsectionstar{VI.3.4.2.3. Laws of Motion in the Non-Inertial Frame}

Knowledge of the solutions of the particle-type Cauchy problems (202) and (218) completely resolves the determination of the \textbf{laws of motion} of the point bodies P$_1$ (m$_1$) and P$_2$ (m$_2$) relative to the non-inertial frame under consideration. Indeed, according to relations (179),

\eqn{231}{\begin{array}{r}
{\vct r}_{1} = {\vct r}_{c} - \frac{m_{2}}{m_{1} + m_{2}}\vct r, \\
{\vct r}_{2} = {\vct r}_{c} + \frac{m_{1}}{m_{1} + m_{2}}\vct r\text{.}
\end{array}}

If \(\vct r_{c}^{}\) and \(\vct r^{}\) denote the solutions of the Cauchy problems (203) and (219), respectively, relations (205) and (220) give

\eqn{232}{\begin{array}{r}
{\vct r}_{c} = F_{- \vct\omega}(r_{c}^{}) \\
\vct r = F_{- \vct\omega}(r^{})
\end{array}}

By the $\mathbb R$-linearity of the rotation tensor \(F_{-\vct\omega}\) with angular velocity \(-\vct\omega\), relations (231) and (232) yield

\eqn{233}{\begin{array}{r}
{\vct r}_{1} = F_{- \vct\omega}\left( {\vct r}_{c}^{} - \frac{m_{2}}{m_{1} + m_{2}}{\vct r}^{} \right), \\
{\vct r}_{2} = F_{- \vct\omega}\left( {\vct r}_{c}^{} + \frac{m_{1}}{m_{1} + m_{2}}{\vct r}^{} \right)\text{.}
\end{array}}

Relations (233) and the results of VI.3.4.1 suggest the following technique for solving the two-body problem in an arbitrarily rotating non-inertial frame $\{R\}$:

I. The non-inertial frame $\{R\}$ is ``frozen'' at time $t=0$, thereby producing an inertial frame denoted by $\{R^0\}$. The set of initial conditions in $\{R^0\}$ is modified accordingly: the position vectors of the particles P$_k$ (m$_k$), \(k=1,2\), retain their values from the frame \(\{R\}\), namely \({\vct r}_{k}(0)={\vct r}_{k}^{0}\), whereas the initial transport velocities \({\vct\omega}_{0}\times {\vct r}_{k}^{0}\) are added to their initial velocities in $\{R\}$, so that \({\dot{\vct r}}_{k}(0)={\vct v}_{k}^{0}+{\vct\omega}_{0}\times {\vct r}_{k}^{0}\), $k=1,2$.

\begin{figure}[htbp]
\centering
\begin{minipage}{.46\textwidth}\centering
\figeleven\par\smallskip
\small\textnormal{Fig. 11. Initial-condition configuration\\in the non-inertial frame $\{R\}$}
\end{minipage}\hfill
\begin{minipage}{.46\textwidth}\centering
\figtwelve\par\smallskip
\small\textnormal{Fig. 12. Initial-condition configuration\\in the inertial frame $\{R^0\}$}
\end{minipage}
\end{figure}

II. The two-body problem is solved in the inertial frame $\{R^0\}$ as described in VI.3.4.1.

Let \(\vct r_{1}^{}\) and \(\vct r_{2}^{}\) be the laws of motion, in the frame $\{R^0\}$, of the two interacting point bodies P$_k$ (m$_k$), \(k=1,2\).

III. The laws of motion of the particles P$_k$ (m$_k$), \(k=1,2\), in the non-inertial frame $\{R\}$ are obtained by applying the rotation tensor \(F_{-\vct\omega}\), with angular velocity \(-\vct\omega\), to their laws of motion in the inertial frame $\{R^0\}$:

\eqn{235}{\begin{array}{r}
{\vct r}_{1} = F_{- \vct\omega}({\vct r}_{1}^{}), \\
{\vct r}_{2} = F_{- \vct\omega}({\vct r}_{2}^{})\text{.}
\end{array}}

\textbf{The configuration of the motion in the non-inertial frame} $\{R\}$ \textbf{is obtained by applying the rotation tensor with angular velocity} \(-\vct\omega\) \textbf{to the configuration of the motion in the inertial frame} $\{R^0\}$.

\paragraph{Remark.}

The behavior of vector-valued functions of the form

\eqn{236}{\vct r=F_{-\vct\omega}(\vct r^0),}
is sometimes surprising. For example, suppose that \(\vct r^{}\) is the law of a uniform circular motion with velocity \(\vct v^{}\), while \(\vct\omega\) is a constant vector lying in the plane of the circular motion. The hodograph of the vector-valued function (236) then densely fills [4] the surface of the sphere of radius \(\left|\vct r^{}\right|\) if \(\frac{w}{w'}\notin\mathbb R\), where \(w^{}=\frac{\left|{\vct v}^{}\right|}{\left|\vct r^{}\right|}\).

This remark shows that, under certain initial conditions, Cauchy problems of the form (218) have \textbf{solutions with semi-chaotic behavior} [161]. Surprisingly, this situation may occur even when the differential equation (218) is linear.

\fourthsectionstar{VI.3.4.2.4. Global Mechanical Characteristics in the Non-Inertial Frame}

The expressions for \textbf{linear momentum, angular momentum}, and \textbf{kinetic energy}, as defined by relations (180), (183), and (186), are unsuitable for studying motion in non-inertial frames. In the two-body problem relative to a non-inertial frame, these quantities are not conserved during the motion and do not provide first integrals, as they do for motion in inertial frames.

We shall introduce a set of vector and scalar quantities that can serve as \textbf{global characteristics} of the mechanical system under consideration.

The \textbf{modified linear momentum} of the mechanical system is

\eqn{237}{\vct H\overset{\text{def}}{=}\sum_{i = 1}^{2}{m_{i}({\dot{\vct r}}_{i} + \vct\omega\times {\vct r}_{i})}\text{.}}
Using relation (169), equation (237) yields

\eqn{238}{\vct H = (m_{1} + m_{2})({\dot{\vct r}}_{c} + \vct\omega\times {\vct r}_{c})\text{.}}
Direct calculation from (238) gives

\eqn{239}{\dot{\vct H} + \vct\omega\times \vct H = (m_{1} + m_{2})\lbrack{\dot{\vct r}}_{c} + 2\vct\omega\times {\dot{\vct r}}_{c} + \vct\omega\times (\vct\omega\times {\vct r}_{c}) + \dot{\vct\omega}\times {\vct r}_{c}\rbrack\text{.}}
Equations (202) and (239) imply

\eqn{240}{\dot{\vct H} + \vct\omega\times \vct H = 0\text{.}}
Relation (240) is of the form (159), with solution

\eqn{241}{\vct H = F_{- \vct\omega}({\vct H}^{0})\text{.}}
In (241), \(\vct H^{0}\) denotes the value of the modified linear momentum at the initial instant, $t=0$. By (241) and the properties of the tensor \(F_{-\vct\omega}\), \textbf{the magnitude of the modified linear momentum of the system is conserved during the motion:}

\eqn{242}{\left| \vct H \right| = \left| {\vct H}^{0} \right|\text{.}}
The \textbf{modified angular momentum} of the mechanical system about the pole $O$ is

\eqn{243}{{\vct L}_{0}\overset{\text{def}}{=}\sum_{i = 1}^{2}{m_{i}{\vct r}_{i}\times ({\dot{\vct r}}_{i} + \vct\omega\times {\vct r}_{i})}\text{.}}
Using relations (231), direct calculation from (243) gives

\eqn{244}{{\vct L}_{0} = m\vct r\times (\dot{\vct r} + \vct\omega\times \vct r) + (m_{1} + m_{2}){\vct r}_{c}\times ({\dot{\vct r}}_{c} + \vct\omega\times {\vct r}_{c})\text{.}}
Using the first integrals (211) and (221), equation (244) yields

\eqn{245}{{\vct L}_{0} = F_{- \vct\omega}({\vct L}_{0}^{0})\text{.}}
In (245), \(\vct L_{0}^{0}\) denotes the value of the modified angular momentum \(\vct L_{0}\) at the initial instant, $t=0$.

It follows from relation (245) that, during the motion, \textbf{the magnitude of the modified angular-momentum vector of the system about the pole $O$ is conserved:}

\eqn{246}{\left| {\vct L}_{0} \right| = \left| {\vct L}_{0}^{0} \right|,}
and it satisfies, at every instant, the relation

\eqn{247}{{\dot{\vct L}}_{0} + \vct\omega\times {\vct L}_{0} = 0\text{.}}
The \textbf{modified kinetic energy} of the mechanical system is

\eqn{248}{T\overset{\text{def}}{=}\sum_{i = 1}^{2}{\frac{m_{i}}{2}({\dot{\vct r}}_{i} + \vct\omega\times {\vct r}_{i})^{2}}\text{.}}
Using relations (231), direct calculation from (248) gives

\eqn{249}{T = \frac{m}{2}(\dot{\vct r} + \vct\omega\times \vct r)^{2} + \frac{m_{1} + m_{2}}{2}({\dot{\vct r}}_{c} + \vct\omega\times {\vct r}_{c})^{2}\text{.}}
The first integrals (212) and (224) can be written in the form

\eqn{250}{({\dot{\vct r}}_{c} + \vct\omega\times {\vct r}_{c})^{2} = \text{const},}
\eqn{251}{\frac{m}{2}(\dot{\vct r} + \vct\omega\times \vct r)^{2} - \int_{}^{}{F(r)\text{dr}} = \text{const}\text{.}}
Consequently, if the interaction forces depend only on the distance between the particles of the system, equations (249)--(251) yield the following \textbf{first integral} having an energy interpretation:

\eqn{252}{T + V = \text{const},}
where \(V\text{=-}\int_{}^{}{F(r)\text{dr}}\) denotes the interaction potential energy of the system.

\paragraph{Remarks.}

1$^\circ$ The quantities defined by (237), (243), and (248) are, in fact, the \textbf{linear momentum, angular momentum}, and \textbf{kinetic energy}, respectively, of the system in the inertial frame $\{R^0\}$ obtained by ``freezing'' the non-inertial frame $\{R\}$.

The \textbf{linear momentum, angular momentum} about the pole $O$, and \textbf{kinetic energy} of the mechanical system consisting of the interacting point bodies P$_k$ (m$_k$), \(k=1,2\), relative to the non-inertial frame $\{R\}$ are defined by

\eqn{253}{\vct P = \sum_{i = 1}^{2}{m_{i}{\dot{\vct r}}_{i}},}
\eqn{254}{{\vct K}_{0} = \sum_{i = 1}^{2}{m_{i}{\vct r}_{i}\times {\dot{\vct r}}_{i}},}
\eqn{255}{E_{\mathrm{cin}}=\sum_{i=1}^{2}\frac12m_i\|\dot{\vct r}_i\|^2.}
Direct calculation establishes the following relations between the quantities defined by (237), (243), and (248) and those given by (253), (254), and (255):

\eqn{256}{\vct H = \vct P + (m_{1} + m_{2})\vct\omega\times {\vct r}_{c},}
\eqn{257}{{\vct L}_{0} = {\vct K}_{0} + I_{0}\vct\omega,}
\eqn{258}{T = E_{\text{cin}} + V_{c}} .

In (256), \({\vct r}_{c}\) denotes the position vector of the center of mass of the system; in (257), $I_0$ is the inertia operator (tensor) in the frame $\{R\}$:

\eqn{259}{\begin{aligned}
I_0&:V_3^*\longrightarrow V_3^*,\\
I_0\vct\omega&\overset{\mathrm{def}}{=}
\sum_{i=1}^{2}m_i\,\vct r_i\times
(\vct\omega\times\vct r_i),
\qquad \forall\,\vct\omega\in V_3^*.
\end{aligned}}

and, in (258),

\eqn{260}{V_{c} = \sum_{i = 1}^{2}\left\lbrack m_{i}(\vct\omega,{\vct r}_{i},{\dot{\vct r}}_{i}) + m_{i}\frac{(\vct\omega\times {\vct r}_{i})^{2}}{2} \right\rbrack\text{.}}
Taking the preceding relations into account, the first integrals (242), (246), and (252), respectively, become

\eqn{261}{\left| \vct P + (m_{1} + m_{2})\vct\omega\times {\vct r}_{c} \right| = \text{const},}
\eqn{262}{\left| {\vct K}_{0} + I_{0}\vct\omega \right| = \text{const},}
\eqn{263}{E_{\text{cin}} + V = \text{const}\text{.}}
In (263), $V$ denotes the function having the role of a generalized potential energy:

\eqn{264}{V(t,{\vct r}_{1},{\vct r}_{2},{\dot{\vct r}}_{1},{\dot{\vct r}}_{2}) = \sum_{i = 1}^{2}\left\lbrack m_{i}(\vct\omega,{\vct r}_{i},{\dot{\vct r}}_{i}) + m_{i}\frac{(\vct\omega\times {\vct r}_{i})^{2}}{2} \right\rbrack - \int_{}^{}{F(r)\text{dr}}\text{.}}
According to the remarks following Theorem VI.3.2.2, if the angular-velocity vector \(\vct\omega\) has a fixed direction, the function (264) does not depend on the velocities \(\dot{\vct r}_i\), $i=1,2$:

\eqn{265}{\begin{array}{r}
V(t,{\vct r}_{1},{\vct r}_{2}) = \sum_{i = 1}^{2}\left\lbrack m_i\,\vct\omega\cdot\vct\Omega_{0i} - m_{i}\frac{(\vct\omega\times {\vct r}_{i})^{2}}{2} \right\rbrack - \int_{}^{}{F(r)\text{dr}}, \\
{\vct\Omega}_{0i}\overset{\text{not}}{=}{\vct r}_{0i}\times ({\vct v}_{0i} + {\vct\omega}_{0}\times {\vct r}_{0i});i = 1,2,
\end{array}}

and, if \(\vct\omega=\text{const}\), one has

\eqn{266}{V(\vct r_1,\vct r_2)=-\frac12\sum_{i=1}^{2}m_i\|\vct\omega\times\vct r_i\|^2-\int F(r)\,dr,}
where one recognizes the expression of the ``centrifugal potential'' [102]: \(V_{\text{cen}}\text{=-}\frac{1}{2}\sum_{i = 1}^{2}{m_{i}(\vct\omega\times {\vct r}_{i})^{2}}\).

2$^\circ$ If the angular-velocity vector \(\vct\omega\) has a fixed direction, the vector \textbf{first integrals} (241) and (245) take the form

\eqn{267}{\vct H=\frac{\vct\omega\cdot\vct H^0}{w^2}\vct\omega-\frac{\sin\alpha(t)}w\vct\omega\times\vct H^0-\frac{\cos\alpha(t)}{w^2}\vct\omega\times(\vct\omega\times\vct H^0),}
\eqn{268}{\vct L_0=\frac{\vct\omega\cdot\vct L_0^0}{w^2}\vct\omega-\frac{\sin\alpha(t)}w\vct\omega\times\vct L_0^0-\frac{\cos\alpha(t)}{w^2}\vct\omega\times(\vct\omega\times\vct L_0^0),}
where $\alpha(t)=\int_0^t w(\tau)\,d\tau$, $t\ge0$, and $w$ is defined by (209).

In this case, the vectors \(\vct H\) and \(\vct L_{0}\) undergo synchronous precessional motions with angular velocity \(-\vct\omega\) about the fixed direction of \(\vct\omega\):

\begin{figure}[htbp]
\centering
\begin{minipage}{.46\textwidth}\centering\figthirteen\par\small\textnormal{Fig. 13}\end{minipage}\hfill
\begin{minipage}{.46\textwidth}\centering\figfourteen\par\small\textnormal{Fig. 14}\end{minipage}
\end{figure}

\fourthsectionstar{VI.3.4.2.5. Motion in the Center-of-Mass Frame}

The motion of the point bodies P$_k$ (m$_k$), \(k=1,2\), relative to the center of mass of the system, whose motion is known from (205), can be described completely from the solutions of problem (218).

Let \(\vct r_{1}^{'}\) and \(\vct r_{2}^{'}\) be the position vectors of particles P$_1$ (m$_1$) and P$_2$ (m$_2$), respectively, relative to the center of mass \emph{$C$} of the system (Fig.~15).

\begin{figure}[htbp]
\centering\figfifteen\par\small\textnormal{Fig. 15}
\end{figure}

From the relation

\eqn{269}{m_1\vct r'_1+m_2\vct r'_2=\vct0.}

it follows that

\eqn{270}{\vct r'_1=-\frac{m_2}{m_1+m_2}\vct r,\qquad \vct r'_2=\frac{m_1}{m_1+m_2}\vct r.}

Relations (270), with \(\vct r\) the solution of Cauchy problem (218), represent the \textbf{laws of motion} of particles P$_1$ (m$_1$) and P$_2$ (m$_2$) relative to a frame whose origin is at the center of mass, denoted by $\{R_c\}$.

The Cauchy problems describing the motion of particles P$_k$ (m$_k$), \(k=1,2\), relative to this frame follow from (270) and (218):

\eqn{271}{\begin{cases}
m_1\ddot{\vct r}'_1+2m_1\vct\omega\times\dot{\vct r}'_1
+m_1\vct\omega\times(\vct\omega\times\vct r'_1)+m_1\dot{\vct\omega}\times\vct r'_1=-F\dfrac{\vct r'_1}{r'_1},\\
\vct r'_1(0)=\vct r'_{01},\\ \dot{\vct r}'_1(0)=\vct v'_{01}.
\end{cases}}

\eqn{272}{\begin{cases}
m_2\ddot{\vct r}'_2+2m_2\vct\omega\times\dot{\vct r}'_2
+m_2\vct\omega\times(\vct\omega\times\vct r'_2)+m_2\dot{\vct\omega}\times\vct r'_2=F\dfrac{\vct r'_2}{r'_2},\\
\vct r'_2(0)=\vct r'_{02},\\ \dot{\vct r}'_2(0)=\vct v'_{02}.
\end{cases}}

In (271) and (272), \(\vct r_{0k}^{'}\) and \(\vct v_{0k}^{'}\) denote the position vector and velocity of particle P$_k$ (m$_k$), \(k=1,2\), relative to the frame $\{R_c\}$.

Examination of relations (271) and (272) shows that \textbf{the center-of-mass frame} $\{R_c\}$ \textbf{is non-inertial} and is characterized by the same angular velocity \(\vct\omega\) as the frame $\{R\}$. We shall say that the frames $\{R\}$ and $\{R_c\}$ are synchronous.

In the center-of-mass frame $\{R_c\}$, the \textbf{global characteristics} corresponding to those defined by relations (253), (254), and (255) are as follows.

The modified \textbf{linear momentum} is

\eqn{273}{{\vct H}^{c} = \sum_{i = 1}^{2}{m_{i}({\dot{\vct r}}_{i}^{'} + \vct\omega\times {\vct r}_{i}^{'})}\text{.}}
The modified \textbf{angular momentum} about the pole \emph{$C$} is

\eqn{274}{{\vct L}_{c} = \sum_{i = 1}^{2}{m_{i}{\vct r}_{i}^{'}\times ({\dot{\vct r}}_{i}^{'} + \vct\omega\times {\vct r}_{i}^{'})}\text{.}}
The modified \textbf{kinetic energy} is

\eqn{275}{T^{c} = \sum_{i = 1}^{2}{\frac{m_{i}}{2}({\dot{\vct r}}_{i}^{'} + \vct\omega\times {\vct r}_{i}^{'})^{2}}\text{.}}
Using relations (270), equations (273)--(275) give

\eqn{276}{{\vct H}^{c} = 0,}
\eqn{277}{{\vct L}_{c} = m\vct r\times (\dot{\vct r} + \vct\omega\times \vct r),}
\eqn{278}{T^{c} = \frac{m}{2}(\dot{\vct r} + \vct\omega\times \vct r)^{2},}
where \(m = \frac{m_{1}m_{2}}{m_{1} + m_{2}}\) denotes the reduced mass of the system and \(\vct r\) is the solution of Cauchy problem (218).

Relations (221) and (277) yield the vector \textbf{first integral}

\eqn{279}{{\vct L}_{c} = F_{- \vct\omega}({\vct L}_{c}^{0}),}
where \(\vct L_{c}^{0}\) is the value of the modified \textbf{angular momentum} about the center of mass at the initial instant, $t=0$.

Equation (279) shows that the hodograph of the vector-valued function \(\vct L_{c}\) is a spherical curve; the magnitude of this quantity is fixed at $t=0$ and remains constant throughout the motion. It also follows that \(\vct L_{c}\) satisfies, at every instant, the vector equation

\eqn{280}{{\dot{\vct L}}_{c} + \vct\omega\times {\vct L}_{c} = 0\text{.}}
If the angular-velocity vector \(\vct\omega\) has a fixed direction, \(\vct L_{c}\) undergoes precession with angular velocity \(-\vct\omega\) about that direction (Fig.~14).

Equations (277) and (270) yield

\eqn{281}{\vct r'_k\cdot\vct L_c=0,\qquad k=1,2.}
and therefore the two particles always lie in a plane, generally moving, that passes through the center of mass of the system and is normal to the variable direction of \(\vct L_{c}\). Since relation (270) gives \({\vct r}_{1}^{'}\text{=-}\frac{m_{2}}{m_{1}}{\vct r}_{2}^{'}\), \textbf{the trajectories of the two particles in the center-of-mass frame are, in general, space curves homothetic with respect to the point} \emph{$C$}, \textbf{with homothety ratio} \(-\frac{m_{2}}{m_{1}}\).

In accordance with form (220) of the solution of Cauchy problem (218), \textbf{the configuration of the motion in the non-inertial center-of-mass frame} $\{R_c\}$ \textbf{is obtained by applying the rotation tensor with angular velocity} \(-\vct\omega\) \textbf{to the configuration of the motion in the inertial frame} \(\{R_{c}^{0}\}\), obtained by ``freezing'' the frame $\{R_c\}$ at time $t=0$.

\paragraph{Remark.}

Under certain initial conditions, planar trajectories relative to the center-of-mass frame may occur when the direction of the angular-velocity vector \(\vct\omega\) is invariant. Under these conditions, the trajectories are planar if and only if \(\vct L_{c}\times\vct\omega=0\) throughout the motion.

Let \(\vct u\) be the unit vector along the fixed direction of \(\vct\omega\). It follows from (279) that the direction of \(\vct L_{c}\) is invariant if and only if

\eqn{282}{\vct u\times {\vct L}_{c}^{0} = 0\text{.}}
Since \(\vct L_{c}^{0}=m{\vct r}^{0}\times({\vct v}^{0}+{\vct\omega}_{0}\times{\vct r}^{0})\), relation (282) yields the following necessary and sufficient condition for the trajectories to be planar:

\eqn{283}{(\vct u \cdot {\vct v}^{0}){\vct r}^{0} = (\vct u \cdot {\vct r}^{0})({\vct v}^{0} + {\vct\omega}_{0}\times {\vct r}^{0})}.

By (281), in this case the plane of the trajectories is normal to the direction of \(\vct u\), hence \(\vct u\cdot{\vct r}^{0}=0\); equation (283) then gives \(\vct u\cdot{\vct v}^{0}=0\). Condition (282) is equivalent to

\eqn{284}{\vct r^0\cdot\vct u=0,\qquad \vct v^0\cdot\vct u=0.}

If the interaction forces between particles P$_1$ (m$_1$) and P$_2$ (m$_2$) depend only on their relative distance, then $F=F(r)$ and one can write a \textbf{scalar first integral} having an energy interpretation:

\eqn{285}{T^c+V=\mathrm{const}.}

Relation (285) follows from (224) and (278), where, as usual, \(V\text{=-}\int_{}^{}{F(r)\text{dr}}\) denotes the interaction potential energy.

Equivalently, (285) can be written as

\eqn{286}{\frac{m{\dot{\vct r}}^{2}}{2} + m(\vct\omega,\vct r,\dot{\vct r}) + \frac{m}{2}(\vct\omega\times \vct r)^{2} - \int_{}^{}{F(r)\text{dr}} = \text{const}\text{.}}
Recognizing in (286) that \(E_{\text{cin}}^{c}=\frac{m{\dot{\vct r}}^{2}}{2}=\sum_{i=1}^{2}{\frac{m_i}{2}{\dot{\vct r}}_i^{}}\) is the \textbf{kinetic energy} of the system in the non-inertial center-of-mass frame $\{R_c\}$, and setting

\eqn{287}{V(t,\vct r,\dot{\vct r}) = m(\vct\omega,\vct r,\dot{\vct r}) + \frac{m}{2}(\vct\omega\times \vct r)^{2} - \int_{}^{}{F(r)\text{dr}},}
as a function having the role of a generalized potential energy, equation (286) becomes

\eqn{288}{E_{\text{cin}}^{c} + V(t,\vct r,\dot{\vct r}) = \text{const}\text{.}}
The preceding relation is read as follows: \textbf{for the two-body problem relative to a non-inertial frame whose origin is at the center of mass and which rotates arbitrarily with angular velocity} \(\vct\omega\)\textbf{, the first integral} (288) \textbf{holds, with the generalized potential energy given by relation} (287)\textbf{.}

If the angular velocity \(\vct\omega\) has a fixed direction (\(\vct\omega\times\dot{\vct\omega}=0\)), the first integral (288) takes the form

\eqn{289}{E_{\text{cin}}^{c} + V(t,\vct r) = \text{const}\text{.},}
where

\eqn{290}{\begin{array}{r}
V(t,\vct r) = m\,\vct\omega\cdot\vct\Omega_0 - \frac{m}{2}(\vct\omega\times \vct r)^{2} - \int_{}^{}{F(r)\text{dr}}, \\
{\vct\Omega}_{0}\overset{\text{not}}{=}{\vct r}^{0}\times ({\vct v}^{0} + {\vct\omega}_{0}\times {\vct r}^{0})\text{.}
\end{array}}

In the particular case \(\vct\omega=\text{const}\), one recovers a result due to Jacobi [102]:

\eqn{291}{E_{\text{cin}}^{c} + V(\vct r) = \text{const}\text{.},}
where

\(V(\vct r)\text{=-}\frac{m}{2}(\vct\omega\times \vct r)^{2} - \int_{}^{}{F(r)\text{dr}}\text{.}\)

\paragraph{Remark.}

Using the same technique, similar results can be obtained in the study of the ``\textbf{$n$-body problem}'' relative to an arbitrarily rotating non-inertial frame.

\chaptertoc{CHAPTER VII}{Conclusions}
% Chapter VII recovered from D7.DOC (Word 6.0, 1995)
% The chapter opening reproduces the graphical format of Chapters I--VI.
\thispagestyle{empty}
\begin{center}
\vspace*{18mm}
{\Large\bfseries CHAPTER VII}\par
\vspace{6mm}
{\large\bfseries Conclusions}
\end{center}
\vspace{12mm}

\par
\noindent\textbf{1.} The notion of a symbolic representation of a vector
space $V$ endowed with an associative, $\mathbb R$-linear unary operation
on another vector space, called the \textbf{representation space}, is
introduced for the first time in the literature. The representation space
is endowed with an automorphism called the \textbf{operator of the symbolic
representation}. A symbolic representation is an $\mathbb R$-linear operator
between the two vector spaces under which the unary operations correspond.

The idea of symbolic representation consists in transferring certain
operator equations from one space to another that offers specific advantages
in applications. This transfer operation has been called the
``algebraization'' of the operator equation on the vector space $V$.
The representation of harmonic oscillation by a complex number, the
representation of certain classes of modulated oscillations on
$\mathbb R$-algebras, and integral transforms of Fourier, Laplace, and
$\mathbb Z$ type are symbolic representations in the sense introduced in
this work.

\par
\noindent\textbf{2.} If the representation space is a finite-dimensional,
associative, commutative, unital $\mathbb R$-algebra and $\widehat d$ is a
fixed element of the algebra, the symbolic representation is specialized as
follows: an algebra element $\widehat w$ symbolically represents the element
$f$ of the vector space $V$ if the algebra element
$\widehat d\,\widehat w$ symbolically represents the element of the vector
space obtained by applying the unary operation to $f\in V$, for every
$f\in V$.

\par
\noindent\textbf{3.} It is proved that, on a finite-dimensional,
associative, commutative, unital $\mathbb R$-algebra, every solution of the
linear operator equation on $V$ whose characteristic polynomial coincides
with the minimal polynomial of the algebra element $\widehat d$ is
symbolically represented with operator $\widehat d$.

\par
\noindent\textbf{4.} It is proved that the vector space of solutions of a
linear operator equation defined on $V$ is symbolically represented on the
polynomial algebra whose generating polynomial coincides with the
characteristic polynomial of the operator equation. The element generating
the cyclic basis of the polynomial algebra plays the role of the operator of
the symbolic representation. It is shown that polynomial algebras suffice to
illustrate the proposed symbolic representation.

\par
\noindent\textbf{5.} By specializing the preceding results, it is shown that
finite-dimensional $\mathbb R$-algebras can symbolically represent continuous
signals that solve linear differential equations with constant coefficients,
as well as discrete signals that solve linear difference equations. Theorems
are given that specify, from the standpoint of symbolic representation, the
relationship between continuous signals and the discrete signals obtained by
their uniform sampling.

\par
\noindent\textbf{6.} The structural properties of the direct-product algebra
of $(n+1)$ algebras isomorphic to the complex field, denoted by
$\bigotimes C^{n+1}$, are studied for the first time. The direct-sum character
of this structure as an algebra of complex components is exhibited. The
passage from the direct-product basis to the basis displaying the direct-sum
character is performed by Walsh--Hadamard transforms, for which fast
computational algorithms exist. Distinguished algebra elements are
characterized: zero divisors, orthogonal idempotents, and special-type
elements. Exponential and semi-exponential forms are presented for an
arbitrary element of the algebra.

\par
\noindent\textbf{7.} The generalized Sobrero algebra, denoted by $S$, is
introduced as the direct product of the second-order complex algebra and the
nilpotent algebra of order $(n+1)$. Its distinguished elements, including zero
divisors and special-type elements, are characterized. Distinguished
subalgebras of $S$ are identified. Exponential and semi-exponential forms are
presented for an arbitrary element of the algebra.

\par
\noindent\textbf{8.} The Walsh algebra, denoted by $W$, is introduced as the
direct product of \emph{n} bireal algebras. The structural properties of this
algebra are studied for the first time, and its character as a direct sum of
$2^n$ real algebras is exhibited. The passage from the direct-product basis to
the basis displaying the direct-sum character is performed by Hadamard-type
transforms, for which fast computational algorithms exist. Distinguished
algebra elements are characterized: zero divisors, orthogonal idempotents,
and special-type elements. Exponential and semi-exponential forms are given
for an arbitrary element of the algebra. It is proved that the direct product
of this algebra with the second-order complex algebra is an algebra isomorphic
to $\bigotimes C^{n+1}$.

\par
\noindent\textbf{9.} A theorem is presented that permits the recovery of
polynomial algebras within the class of finite-dimensional, commutative,
associative, unital algebras. Thus, if the algebra contains an element whose
characteristic polynomial is also its minimal polynomial, then the algebra is
polynomial. The matrix of the basis transformation from an arbitrary basis to
the power basis that explicitly displays the polynomial-algebra structure is
given.

It is proved that the algebra $\bigotimes C^{n+1}$ is polynomial, with a
generating polynomial having $2^{n+1}$ distinct simple purely complex roots.
The algebra $S$ is polynomial, with a generating polynomial having two purely
complex roots of multiplicity $(n+1)$. The algebra $W$ is polynomial, with a
generating polynomial having $2^n$ distinct real roots.

For each algebra studied, procedures are given for determining the generating
polynomial and the element that generates the cyclic basis of the algebra.

\par
\noindent\textbf{10.} It is proved that the algebra
$\bigotimes C^{n+1}$ symbolically represents continuous quasiperiodic signals
with a harmonic carrier and simultaneous amplitude and phase modulation of
the ``product'' type of harmonic signals. These signals are solutions of
constant-coefficient differential equations of order $2^{n+1}$ whose
characteristic polynomial has distinct purely complex roots. Zero divisors
symbolically represent ``degenerate'' signals, which are solutions of linear
constant-coefficient differential equations of order strictly less than
$2^{n+1}$. The orthogonal idempotents of the algebra symbolically represent
harmonic signals. Special-type elements symbolically represent quasiperiodic
signals with a harmonic carrier and pure amplitude modulation.

\par
\noindent\textbf{11.} On the algebra $S$, continuous signals with
amplitude-modulated harmonic carriers in quadrature are symbolically
represented. These signals have the character of polynomial-exponential
oscillations simultaneously modulated in amplitude and phase, and are
solutions of linear constant-coefficient differential equations of order
$2(n+1)$ whose characteristic polynomial has two complex-conjugate roots of
multiplicity $(n+1)$. Zero divisors symbolically represent degenerate
continuous signals, which are solutions of linear constant-coefficient
differential equations of order strictly less than $2(n+1)$.

Special-type elements symbolically represent continuous signals with a
harmonic carrier and pure amplitude modulation.

\par
\noindent\textbf{12.} A symbolic representation of integrable periodic
continuous signals on the Walsh algebra is introduced. This representation
``algebraizes'' a particular class of integral equations with degenerate
kernels. Using specific properties of Walsh algebras proved in this work, a
novel procedure is proposed for the parametric identification of a class of
nonlinear dynamical systems with analytic nonlinearities. For dynamical
systems with polynomial nonlinearities, the identification procedure is
exact.

\par
\noindent\textbf{13.} The notion of a dynamical system on a vector space
endowed with an associative, $\mathbb R$-linear unary operation is introduced
through an operator equation. This notion encompasses, in a unified manner,
both analog linear dynamical systems and digital filters. If the vector space
supporting the dynamical system is symbolically represented on a
finite-dimensional, associative, commutative, unital $\mathbb R$-algebra, the
properties of the symbolic representation permit the introduction of an
algebra element characterizing the system, called the \textbf{transmittance
of the dynamical system}. In practice, the transmittance is a rational
expression in the operator of the symbolic representation. It is proved that
the transmittance provides a concise characterization of the behavior of the
dynamical system in every situation.

\par
\noindent\textbf{14.} A Fourier-type integral transform is introduced on the
algebra $\bigotimes C^{n+1}$. Thus, every continuous signal in $L_1$ is
associated with a hypercomplex-valued function of a hypercomplex variable.
It is proved that this association is a symbolic representation on
$\bigotimes C^{n+1}$ when the usual differentiation of real-valued functions
of a real variable is taken as the unary operation on $L_1$.

The transform is proved to be monogenic as a function of the hypercomplex
variable and invertible. Its invertibility shows that any signal in $L_1$ can
be synthesized from signals symbolically representable on
$\bigotimes C^{n+1}$, having the character of harmonic-carrier oscillations
simultaneously modulated in amplitude and phase. Relations between this
transform and the usual Fourier transform of functions in $L_1$ are
established. The possibility of extending the proposed transform to tempered
distributions is indicated. Under this extension, the image of the
distributional derivative of order $k\in\mathbb N$ of the Dirac impulse
distribution is shown to be the $k$th power of the operator of the symbolic
representation.

\par
\noindent\textbf{15.} A Laplace-type integral transform is introduced on the
generalized Sobrero algebra $S$. This transform associates every signal in
the class of original signals with a hypercomplex-valued function of a
hypercomplex variable. The resulting association is a symbolic
representation on $S$ when the usual differentiation of real-valued
functions of a real variable is taken as the unary operation. The transform
is proved to be monogenic as a function of the hypercomplex variable and
invertible. Its invertibility shows that any signal can be synthesized from
signals symbolically representable on $S$, having the character of
sinusoidal-carrier oscillations simultaneously modulated in amplitude and
phase. The relation between the proposed transform and the usual Laplace
transform of original functions, defined on a subalgebra of $S$, is
established. This relation permits the extension of the transform to the
broader space of distributions. Under this extension, the image of the
distributional derivative of order $k\in\mathbb N$ of the Dirac impulse
distribution is the $k$th power of the operator of the symbolic
representation.

\par
\noindent\textbf{16.} For both integral transforms, defined on
$\bigotimes C^{n+1}$ and on the generalized Sobrero algebra $S$, the
following property is proved: the transmittance of a linear, invariant analog
dynamical system on either algebra is obtained by applying the corresponding
integral transform to the weighting function of the system. This result
establishes the relation between the transmittance of a dynamical system on an
$\mathbb R$-algebra and its Laplace transfer function. Conditions imposed on
the transmittance of a system are thereby reflected in the corresponding
transfer function.

If the transmittance of a dynamical system is a special-type element of the
algebra, the system introduces no parasitic phase distortion. The resulting
conditions on the Laplace transfer function are less restrictive than those
presented in the literature. For the specified signal types, the conditions
given in this work are necessary and sufficient, whereas the literature gives
sufficient conditions for a broader class of continuous signals.

\par
\noindent\textbf{17.} A matrix symbolic representation of vector-valued
functions of a real variable is introduced. A vector-valued function is
associated with the column matrix containing the three real-valued components
of that function in a right-handed orthonormal basis. The cross product with a
fixed vector is taken as the unary operation on the space of vector-valued
functions. This cross product is associated with the product of a suitably
chosen skew-symmetric matrix and the column matrix corresponding to the
vector-valued function. A symbolic representation in the sense defined in
this work is thereby obtained. The resulting transformation ``algebraizes''
certain vector differential equations involving cross products.

\par
\noindent\textbf{18.} Using the proposed symbolic representation, exact
vector solutions are obtained for two classes of Cauchy problems arising in
the dynamics of mechanical systems moving relative to uniformly rotating
non-inertial frames. As a first application, exact vector solutions are
presented for three classical problems of theoretical mechanics for which
only approximate solutions are traditionally given. These are the motion in
uniform fields relative to a uniformly rotating non-inertial frame, the
motion of electric charges in crossed electric and magnetic fields, and the
celebrated Foucault-pendulum problem. The exact vector solutions of these
problems are reported for the first time in this work.

\par
\noindent\textbf{19.} The problems of motion in a uniform gravitational
field relative to an Earth-fixed non-inertial frame and of the Foucault
pendulum have served as models for experiments demonstrating the Earth's
motion about its poles without recourse to astronomical observations. The
existence of exact solutions for these problems may suggest new experiments.

\par
\noindent\textbf{20.} The exact vector solution for the motion of electric
charges in a uniform electromagnetic field suggests several novel
considerations concerning the mechanism of particle--field interaction. The
resultant force acting on the moving particle is shown to have constant
magnitude, fixed at the initial instant. During the motion, the resultant
force undergoes uniform precession about the magnetic-field lines, with a
precession angular velocity depending only on the magnitude of the magnetic
induction vector and on the specific charge of the particle. These results
extend the well-known Larmor theorem.

\par
\noindent\textbf{21.} A tensor symbolic representation of vector-valued
functions of a real variable is introduced. If $\boldsymbol\omega$ is a
differentiable vector-valued function, it is proved that the rotation-tensor
operator with angular velocity $\boldsymbol\omega$ is uniquely determined by
the condition that it coincide with the identity tensor at the initial
instant.

A symbolic representation of the vector space $V_3^*$ of vector-valued
functions of a real variable into itself is defined by means of the rotation
tensor with angular velocity $\boldsymbol\omega$, when
$\boldsymbol\omega$ has a fixed direction. The relative derivative
corresponding to the rotational angular velocity $\boldsymbol\omega$ is
chosen as the unary operation, while the usual derivative of vector-valued
functions with respect to the real variable is chosen as the automorphism,
that is, as the operator of the symbolic representation. The transformation
is proved to be invertible, its inverse being described by the rotation tensor
with angular velocity $-\boldsymbol\omega$. The correspondence law of this
operator is given explicitly in terms of $\boldsymbol\omega$.

The proposed transformation ``algebraizes'' differential equations modeling
motion relative to non-inertial frames undergoing nonuniform rotation about a
fixed direction, as well as motion relative to inertial frames in certain
gyroscopic force fields.

\par
\noindent\textbf{22.} As a first application, exact vector solutions are
given for a class of problems not addressed in the literature, even from the
standpoint of approximate solutions. The following are treated completely:
motion in uniform force fields relative to non-inertial frames undergoing
nonuniform rotation about a fixed direction; the motion of electric charges
in a nonuniform magnetic field of fixed direction; and the motion of electric
charges in an electromagnetic field with a slowly varying magnetic field of
fixed direction and a uniform electric field. An exact vector solution is
also given for the motion of the Foucault pendulum in non-inertial frames
undergoing nonuniform rotation about a fixed direction.

For the first time, an exact solution is presented for a reduced
Foucault-pendulum problem found in most treatises on theoretical mechanics,
for which approximate solutions based on small-parameter methods are
commonly given. The problem also models the motion of the vibrating element
in a class of gyroscopic instruments.

\par
\noindent\textbf{23.} The problem of motion in positional central-force
fields relative to non-inertial frames undergoing nonuniform rotation about a
fixed direction is solved. It is proved that the solution of the Cauchy
problem describing such motion is obtained by applying the rotation tensor
with angular velocity $-\boldsymbol\omega$ to the solution of the
differential equation obtained from the initial problem by setting
$\boldsymbol\omega=0$, with suitably modified initial conditions.

Two novel first integrals are exhibited for motion in positional
central-force fields relative to non-inertial frames undergoing nonuniform
rotation about a fixed direction.

\par
\noindent\textbf{24.} The symbolic representation of the vector space
$V_3^*$ into itself is extended to arbitrary angular velocities
$\boldsymbol\omega$, expressed by arbitrary differentiable vector-valued
functions of a real variable. The unary operation on $V_3^*$ is the relative
derivative with respect to the angular velocity $\boldsymbol\omega$, while
the operator of the symbolic representation is the usual derivative of
vector-valued functions of a real variable.

The symbolic representation so defined is proved to be invertible, its
inverse being expressed by the rotation tensor with angular velocity
$-\boldsymbol\omega$.

The proposed transformation ``algebraizes'' differential equations modeling
phenomena involving motion relative to non-inertial frames undergoing
arbitrary rotation with respect to an inertial frame.

\par
\noindent\textbf{25.} Two novel first integrals are exhibited for motion in
positional central fields relative to non-inertial frames rotating with
angular velocity $\boldsymbol\omega$, where $\boldsymbol\omega$ is an
arbitrary differentiable vector-valued function.

\par
\noindent\textbf{26.} Although, in the general case, the rotation operator
with angular velocity $-\boldsymbol\omega$ cannot be written explicitly, the
geometric properties of the transformation permit a qualitative study of the
motion. Using this methodology, theorems are given concerning the motion of
the Foucault pendulum, criteria concerning the boundedness and form of
solutions in Newtonian force fields, and a generalization of Larmor's theorem
for the motion of electric charges in an arbitrarily varying magnetic field.

\par
\noindent\textbf{27.} Several directions for extending the results obtained
in this work are envisaged:

\par\noindent a) the study of the symbolic representation of solutions, in
the sense of distribution theory, for operator equations on specified vector
spaces;

\par\noindent b) the study of symbolic representation on nonassociative
Lie-type algebras;

\par\noindent c) the study of integral transforms on polynomial algebras;

\par\noindent d) the study of the symbolic representation of random signals;

\par\noindent e) the study of sampling for signals that solve
constant-coefficient differential equations but do not satisfy the
hypotheses of the Whittaker--Kotelnikov--Shannon theorem on the sampling of
continuous signals;

\par\noindent f) the study of quaternionic and spinorial rotation operators
in the context of the tensor symbolic representation of vector-valued
functions of a real variable, with applications in robotics;

\par\noindent g) the qualitative study of the behavior of solutions of
Cauchy problems modeling the motion of mechanical systems in Newtonian
fields relative to arbitrarily rotating non-inertial frames, with
applications in space dynamics;

\par\noindent h) the study of integral transforms suggested by the tensor
symbolic representation of vector-valued functions of a real variable, with
applications in dynamical-systems theory;

\par\noindent i) the study of the effect of symbolic signal representation
on $\mathbb R$-algebras upon the response of nonlinear dynamical systems
whose input--output behavior is described by Volterra series or by Fliess-type
series in noncommuting variables.

\clearpage
\setchapterheader{Bibliography}

\section*{Bibliography}
\begin{enumerate}[label={[\arabic*]},leftmargin=*,itemsep=.28em,parsep=0pt]
\item Ahmed, N., Rao K.R. - Orthogonal Transforms for Digital Signal Processing, Springer-Verlag, Berlin - Heidelberg-New York, 1975.
\item Albert, A. Adrian - Structure of Algebras, American Mathematical Society, Colloquium Publications, Vol. XXIV, Providence, 1961.
\item Antoulas, A.C. - " A New Approach to Synthesis Problems , Linear System Theory", " IEEE Transactions on Automatic Control", Vol. AC-30, N. 5, May 1985.
\item Arnold, V.I. - Equations différentielles ordinaires, Ed. Mir, Moscou, 1974.
\item Arnold, V.I. - Metodele matematice ale mecanicii clasice, Editura Științifică și Enciclopedică, București, 1980.
\item Atherton, D.P. - " A Survey on Nonlinear Oscillations" , " Int. J. Control", Vol. 31, n. 6, 1980, pp. 1041-1105.
\item Balas, M.J. - " Feedback Control of Flexible Systems" , " IEEE Transactions on Automatic Control", Vol. AC-23, N. 4, Aug. 1978.
\item Barbu, Viorel - Ecuații diferențiale, Editura Junimea, Iași, 1985.
\item Baskakov, S.I. - Signals and Circuits, Mir Publishers, Moscow, 1986.
\item Bauman, W.T., Hugh, W.J. - " Feedback Control of Nonlinear Systems by Extended Linearization" , IEEE Transactions on Automatic Control", Vol. AC-31, N. 1, Jan. 1986.
\item Bauman, W.T., Herdman, T.L., Stalford, H., Suchomel, C.F. - " Recent Work Using Volterra Series as a Methodology to Analyze Nonlinear Aircraft Dynamic Properties" , " Math. Comput. Modelling", Vol. 11, Pergamon Press, G.B. 1988, pp. 883-888.
\item Beauchamp, K.G. - Walsh Functions and their Applications, " Techniques of Physics", Academic Press, London-New York, San Francisco, 1975.
\item Bolton, A.G. - " Inverse Nyquist Design Using Region Algebra" , " Int. J. Control", 1981, Vol. 33, N.3, pp. 575-584
\item Borș, I. - Lecții de mecanică, Universitatea "Al.I. Cuza", Iași, Vol. I-II, 1987.
\item Boyd, S.P., Chua, L.O. - " Fading Memory and the Problems of Approximating Nonlinear Operators with Volterra Series" , " IEEE Transactions on Circuits and Systems", Vol. CAS-32, N. 11, Nov. 1985.
\item Boyd, S.P., Chua, L.O. - " Uniqueness of Circuits and Systems Containing one Nonlinearity", IEEE Transactions on Automatic Control", Vol. AC-30, N. 7, July 1985.
\item Boțan, C. - " Contribuții la studiul parametrilor unor sisteme automate de urmărire (de curent alternativ)", summary of the doctoral thesis, Iași, 1973.
\item Braier, A. - " New Algebraic Method in Dynamic System", Lectures held at the Cosserat Session of C.I.M., Udine, 1972.
\item Braier, A. - Mecanica pentru ingineri, Vol. I " Cinematica", Editura Institutului Politehnic, Iași, 1975.
\item Braier, A. - " Utilizarea algebrelor pe corpul real în mecanica vibrațiilor", summary of the doctoral thesis, Iași, 1968.
\item Braier, A., Condurache, D. - " On the Instantaneous Angular Velocity Vector", Bul., Inst. Polit. Iași, XXX (XXXIV), 1-4, s. I, 1984, 25-28.
\item Can, S., Unal, A. - " Transfer Functions for Nonlinear System via Fourier-Borel Transforms" , ISCAS’88, pp. 2445-2449.
\item Cartianu, Gh. et al. - Semnale, circuite și sisteme, Editura Didactică și Pedagogică, București, 1980.
\item Celle, E., Gauthier, J.P., Milani, E. - " Existence of Realizations of Nonlinear Analytic Input-Output Maps" , IEEE Transactions on Automatic Control", Vol. AC-31, N. 4, April 1986.
\item Chang, R.Y., Wang, M.L. - " Parameter Identification via Shifted Legendre Polynomials" Int. J.Syst. Sci., Vol. 13, 1982, pp. 1125-1135.
\item Cheng, C., Hsiao, C.H. - " Time Domain Synthesis Walsh Function", Proc. IEEE, Vol. 122, 1975, pp. 565-570.
\item Chen, Wen-Liang, Shih, Yen-Ping - " Shift Walsh Matrix and Delay Differential Equations" , " IEEE Transactions on Automatic Control", Vol. AC-23, N. 6, Dec. 1978.
\item Cheng, B., Hsu, N.S. - " Single-input, Single-output Systems Identification via Block-Pulse Functions", Int. J. Syst. Sci., Vol. 12, pp. 643-647, 1981.
\item Cheng, D.K., Hung, J.K. - " Applications on Quaternion in Robot Control" , Lecture Notes in Control and Information Sciences, Springer-Verlag, 1984, pp. 269-283.
\item Cheng, R., Tarn, T., Isidori, A. - " Global External Linearization of Nonlinear Systems Via Feedback" , IEEE Transactions on Automatic Control", vol. Ac-30, N. 8, Aug. 1985.
\item Cheng, D.K., Shanker, A.U. - " Walsh - Function Representation and Noise Analysis of Linear Sequential Circuits" , IEEE Transactions on Circuits and Systems", vol. Cas-32, N. 3, Mar. 1985.
\item Chow, J.H., Kokotovic, P.V. - " Two-Time Scale Feedback Design of a Class of Nonlinear Systems" , IEEE Transactions on Automatic Control, Vol. AC-32, N.6, June 1978.
\item Chung-Kwong, Y. - " Upper Bounds on Walsh Transforms" , " IEEE Transactions on Computers", Dec. 1972.
\item Claude, D. - " Linéarisation par diffeomorphisme et immersion des systèmes" , Lecture Notes in Control and Information Sciences, Springer-Verlag, 1984, pp. 340-351.
\item Clement, R.R. - " Laguerre Functions in Signal Analysis and Parameter Identification", J. Franklin Inst., Vol. 313, 1982, pp. 85-95.
\item Condurache, D. - " Symbolic Representation of Signals on Hyperspaces" (I), Bul. Inst. Polit. Iași, XXVII(XXXI), 1-2, s. III, 1981, pp. 33-42.
\item Condurache, D. - " Symbolic Representation of Signals on Hyperspaces" (II), Bul. Inst. Polit. Iași, XXVII(XXXI), 3-4, s. III, 1981, pp. 49-56.
\item Condurache, D. - " On the Symbolic Representation of Discrete Signals", Bul. Inst. Polit. Iași, XXVII(XXXI), 1-2, s. III, 1981, pp. 55-59.
\item Condurache, D. - " Polynomial Algebras and Symbolic Representation of Discrete Signals", Bul. Inst. Polit. Iași, XXVII(XXXI) 3-4, s. I, 1981, pp. 41-46.
\item Condurache, D. - " Integral Transforms on Algebras" (I), Bul. Inst. Polit. Iași, XXVIII(XXXII), 1-4, s.III, 1982, pp. 27-36.
\item Condurache, D. - " Integral Transforms on Algebras" (II), Bul. Inst. Polit. Iași, XXX(XXXIV), 1-4, s.III, 1984, pp. 19-23.
\item Condurache, D. - " Integral Transforms of Fourier-Type on Commutative Algebras", Bul. Inst. Polit. Iași, XXX(XXXII), 1-4, s.I, 1984, pp. 13-17.
\item Condurache, D. - " About a $2^{n+1}$ Order Algebra", Bul. Inst. Polit. Iași, XXXI(XXXII), 1-4, s.I, 1985, pp. 35-39.
\item Condurache, D. - " Identificarea exactă a sistemelor fără memorie utilizînd funcții Walsh", a IV-a Conferință națională de fiabilitate în construcția de mașini" , Iași, 1986.
\item Condurache, D. - " Integral Transforms on Sobrero Algebra", Bul. Inst. Polit. Iași, XXXIII(XXXVII), 1-4, s.I, pp. 5-8, 1987.
\item Condurache, D. - " Asupra soluțiilor matriceale ale unor ecuații vectoriale", " Creația tehnică și fiabilitatea în construcții de mașini", Iași, 22-25 Nov. 1985, pp. 90-94.
\item Condurache, D. - " About some Algebraic Properties of The Walsh Functions", Bul. Inst. Polit. Iași, XXXIII(XXXVII), 1-4, s.III, 1987, pp. 19-23.
\item Condurache, D. - " Proprietățile algebrice ale seriilor Walsh trunchiate în identificarea sistemelor neliniare cu memorie", Conferința națională de Matematică aplicată și Mecanică, 20-23 Oct. 1988, pp. 387-390.
\item Condurache, D. - " Un procedeu simbolic în studiul dinamicii relative a particulei materiale", în Concepție, tehnologie și management în construcția de mașini, Inst. Polit. Iași, Iași, 1992, pp. 121-123.
\item Condurache, D. - " O metodă directă de integrare în mișcarea particulelor în cîmp electromagnetic", în Concepție, tehnologie și management în construcția de mașini, Inst. Polit. Iași, pp. 124-126, Iași, 1992.
\item Condurache, D., Braier, A. - " A Method for the Direct Integration in the Study of Some Classical Problems of Theoretical Mechanics", Bul. Inst. Polit. Iași, XLI(XLV), 1-2, s.I, 1995.
\item Corduneanu, A. - Ecuații diferențiale cu aplicații în electrotehnică, Editura Facla, Timișoara, 1981.
\item Corrington, S. Murlan - " Solution of Differential and Integral Equations with Walsh Function" , " IEEE Transactions and Circuit Theory", Sept. 1973.
\item Creangă I., Enescu, I. - Algebre, Editura Tehnică, 1973.
\item Davis, Artice. M. - " Almost Periodic Extension of Band-Limited Functions and its Application to Nonuniform Sampling" , " IEEE Transactions on Circuits and Systems", Vol. CAS-33, N. 10, Oct. 1986.
\item Dieudonné, J. - Fondements de l’analyse moderne, Gauthier-Villars, Paris, 1963.
\item Diktine, V., Proudnikov, A. - Calcul operationnel, Editions Mir, Moscou, 1979.
\item Dragoș, L. - Principiile mecanicii mediilor continue, Editura Tehnica, București, 1983.
\item Dragoș, L. - Principiile mecanicii analitice, Editura Tehnica, București, 1976.
\item Ekeland, I. - " Oscillation de systèmes hamiltoniens nonlinéaires", " Bulletin de la Société Mathématique de France", 109, 1981, pp. 297-330.
\item Eykhoff, P. - Identificarea sistemelor, Editura Tehnică, București, 1977.
\item Fish Jr., A.J. - " An Algebraically Derived Nonlinear Control Theory" in "Journal of Dynamic Systems, Measurement and Control", Vol. 105, June 1983, pp. 83-91.
\item Fliess, M. - " Fonctionnelles causales nonlinéaires et indéterminées non commutatives" , "Bull. Soc. Math. France", 109, 1981, pp. 3-40.
\item Fliess, M. - " On a possible Connection between Volterra Series and Nonlinear Optimal Control", Proc. 7th Conf. Informat. Sci. Systems, Baltimore, 1983, pp. 402-407.
\item Fliess, M., Lamnabhi, M., Lamnabhi-Lagarrigue, F. - " An Algebraic Approach to Nonlinear Functional Expansions", " IEEE Transactions on Circuits and Systems", Vol. CAS-30, N. 8, Aug. 1983.
\item Fliess, M. - " Lie Brackets and Optimal Nonlinear Feedback Regulation", Lecture Notes in Control and Information Sciences, Springer-Verlag, 1984, pp. 78-83.
\item Fliess, M. - " Esquisses pour une théorie des systèmes non linéaires en temps discret," proc. Conf Linear nonlinear Math. Control Theory, Torino, 1986.
\item Fliess, M. - " Nonlinear Control Theory and Differential Algebra", Proc. I.I.A.S.A. Conf. Modelling Adaptive Control, Sopron, Ungaria, July 1986.
\item Fornasini, E., Marchesini, G. - " Some Connections between Algebraic Properties of pairs of Matrices and 2D Systems Realization", Lecture Notes in Control and Information Sciences, Springer-Verlag, 1984, pp. 117-189.
\item Friedland, B., Hutton, M. - " Theory and Error Analysis of Vibrating-Member Gyroscope", " IEEE Transactions on Automatic Control", Vol. AC-23, N. 4, Aug. 1978.
\item Fukushima, M., Yamamoto, Y. - " A Second-Order Algorithm for Continuous-Time Nonlinear Optimal Control Problems", " IEEE Transactions on Automatic Control", Vol. AC-31, N. 7, July 1986.
\item Gantmakher, F.R. - The Theory of Matrices, I, II, N.Y. Chelsea, 1959.
\item Ghelfond, A.O. - Calculul cu diferențe finite, Editura Tehnica, București, 1956.
\item Gheorghiu, N. - Introducere în analiza funcțională, Editura Academiei, București, 1974.
\item Gonorovsky, I. - Radio Circuits and Signals, Mir Publishers, Moscova, 1981.
\item Goras, L. - Semnale, circuite și sisteme, Editura Gh.Asachi, Iași, 1994.
\item Godunov, S.K., Reabenki, V.S. - Scheme de calcul cu diferențe finite, Editura Tehnica, București, 1977.
\item Goussard, Y., Krenz, W., Stark, L. - " An Improvement of Lee and Schetzen Cross-Correlation Method", " IEEE Transactions on Automatic Control", Vol. AC-30, N. 9, Sept. 1985.
\item Grizzle, J.W., Marcus, S.I. - " The Structure of Nonlinear Control Systems Possessing Symmetries", " IEEE Transactions on Automatic Control", Vol. AC-30, N. 3, March 1985.
\item Grizzle, J.W. - " Controlled Invariance for Discrete - Time Nonlinear Systems with an Application to the Disturbance Decoupling Problems", " IEEE Transactions on Automatic Control", Vol. AC-30, N. 9, Sept. 1985.
\item Harmuth, F. Henning - Sequency Theory, " Foundations and Applications", Academic Press, N.Y., San Francis­co, London, 1977.
\item Hopfield, J.J., Tank, D.W. - " Neural Computation of Decisions in Optimization problems", Biol. Cybern., 1985, 52, pp. 141-152.
\item Herdla, J. - " New Integral Transform between Analogue and Discrete Systems", " Electronics Letters", May 1988, Vol. 24, n. 10.
\item Hsu, I., Kaszkurewicz, E. - " Linearization of Nonlinear Control systems by Recursive transformations: An Algebraic Approach", " Preprints of the 9th World Congress of the International Federation of Automatic Control", Budapest, July 1984, Vol. V, pp. 89-94.
\item Hwang, C., Guo, Y. - " Transfer Function Matrix Identification in MIMO Systems via Shifted Legendre Polynomials", Int. J. Control, Vol. 39, 1984, pp. 807-814.
\item Hwang, C., Shih, Y.P. - " Laguerre Operatorial Matrices for Fractional Calculus and Applications", Int. J. Control, vol. 34, 1981, pp. 577-584.
\item Iacob, Caius - Mecanica teoretică, Editura Didactică și Pedagogică, București, 1980.
\item Irimiciuc, N., Deleanu, Sp., Vieru, D., Budei, R. - Curs și culegere de probleme de mecanică, Vol. I, II, I.P., Iași, 1983.
\item Irimiciuc, N. et al. - Curs de mecanică pentru uzul studenților de la Facultatea de Mecanică, I.P. Iași, 1993
\item Incertis, F. - " A Faster Method of Computing Matrix Pythagorean Sums", " IEEE Transactions on Automatic Control", Vol. AC-30, N. 3, March 1985.
\item Isidori, A. - " The Matching of a Prescribed Linear Input-Output Behaviour in Nonlinear Systems", IEEE Transactions on Automatic Control, Vol. AC-30, N. 3, Mar. 1985.
\item Izvercian, P.N., Popescu, C. - " The Structure of Finite-Dimensional Spaces of Functions with are Representable on Algebras", Simpozionul de matematici și aplicații, Timișoara, 1987, pp. 183-192.
\item Izvercian, P.N., Hadnagy, A. - " Symbolic Representation of Function Spaces by Elements of an Algebra", An Univ. Timișoara, seria " At. Matematice", Vol. XXI, 1-2, 1983.
\item Kao, Huang Y., Huang, C.J. - " Routes to Chaos in Duffing Oscillator with one Potential Well".
\item Karanam, V.R., Frick, P.A. - " Bilinear System Identification by Walsh Functions", " IEEE Transactions on Automatic Control", Vol. AC-23, N. 4, Aug. 1978.
\item Kasumasa, H., Sawai, N. - " A General Criterion for Jump Resonance of Nonlinear Control Systems", " IEEE Transaction on Automatic Control", Vol. AC-23, N. 5, Oct. 1973.
\item Kocin, N.E. - Calculul vectorial și introducere în calculul tensorial, Editura Tehnica, București, 1954.
\item Kostrikin, A., - Introduction a l’Algèbre, Éditions Mir, Moscou, 1981.
\item Kotkine, G., Serbo, V. - Recueil de problèmes de mécanique classique, Éditions Mir, Moscou, 1981.
\item Krasnov, M.L., Kiselev, A.I., Makarenko, G.I. - Vector Analysis, Mir Publishers, Moscow, 1981.
\item Kundert, K.S., Sangiovanni-Vincentelli, A. - " Simulations of Nonlinear Circuits in the Frequency Domain", " IEEE Transactions on Computer-Aided Design", Vol. CAD-5, N. 4, Oct. 1986.
\item Landau, L., Lifchitz, E. - Mécanique, Éditions Mir, Moscou, 1981.
\item Landau, L., Lifchitz, E. - Théorie des champs, Éditions Mir, Moscou, 1970.
\item Landau, I.D. et al. - " Outils et modèles mathématiques pour l’automatique, l’analyse de systèmes et le traitement du signal", Éditions du Centre National de la Recherche Scientifique, Paris, 1983.
\item Lamnabhi-Lagarrigue, P. - " Sur les conditions nécessaires d’optimalité du deuxième et troisième ordre dans les problèmes de commande optimale singulière", Lecture Notes in Control and Information Sciences, Springer-Verlag, 1984, pp. 535-541.
\item Lau, S.L., Cheung, Y.K., Yu, S.Y. - " Incremental Harmonic Balance Method with Multiple Time Scales for Aperiodic Vibration of Nonlinear Systems", " Journal of Applied Mechanics", Vol. 50, Dec. 1983, pp. 871-876.
\item Lee, H. - " The Orthogonal Transformations", " IEEE Transactions on Circuits and Systems", Vol. CAS-32, N. 11, Nov. 1985.
\item Legros R., Martin, A.V.J. - Transform Calculus for Electrical Engineers, Prentice Hall, New Jersey, 1961.
\item Lejeune, R., Hugh, W.J. - " Linearization of Nonlinear Systems about Constant Operating Points", " IEEE Transactions on Automatic Control", Vol. AC-30, N. 8, Aug. 1985.
\item Lesiak, C., Krener, A.J. - " The Existence of Uniqueness of Volterra Series for Nonlinear Systems", " IEEE Transactions on Automatic Control", Vol. AC-23, N. 6, Dec. 1978.
\item Loparo, K.A., Blankenship, G.L. - " Estimating the Domain of Attraction of Nonlinear Feedback Systems" , " IEEE Transactions on Automatic Control", Vol. AC-23, N. 4, Aug. 1978.
\item Marciuk, G.I., Saidurov, V.V. - Creșterea preciziei soluțiilor în scheme cu diferențe, Editura Academiei, București, 1981.
\item Mangeron, D., Irimiciuc, N. - Mecanica rigidului cu aplicații în inginerie, Vol. I-III, Editura Tehnica, București, 1978, 1980, 1981.
\item Mangeron, D., Oguztöreli, M.N., Jasulionis, A.I. - " Application de la méthode matricielle harmonique à l’intégration des équations différentielles linéaires d’ordre n à coefficients périodiques", " Bulletin de la classe des sciences" 5\textasciicircum{}(e) serie, T.LVI, Académie Royale de Belgique, Bruxelles, 1970.
\item Marmarelis, V.Z. - " Identification of Nonlinear Systems by Use of Nonstationary White-Noise Inputs", " Appl. Math., Modelling", 1980, Vol. 4, April 1980, pp 117-124.
\item Maslov, V.P. - Operational Methods, Mir Publishers, Moscow, 1976.
\item Mathis, W., Voigt, B. - " Applications on Lie Series Averaging in Nonlinear Oscillation", " IEEE Transactions on Automatic Control", N. 5, Aug. 1987.
\item Mihalcea, A. - Prelucrarea optimă a semnalelor în sisteme informaționale, Editura Militară, București, 1987.
\item Moisil, Gr. - " Sur une classe de systèmes d’équations aux dérivées partielles de la physique mathématique", in Opera Matematica, Vol. II, Editura Academiei, 1980.
\item Moisil, Gr. - " Sur les petits mouvements des corps élastiques", in Opera Matematica, Vol. II, Editura Academiei, 1980.
\item Moisil, Gr. - " Metoda funcțiilor de variabilă hipercomplexă în hidrodinamica plană a lichidelor vîscoase incompresibile", in Opera Matematica, Vol. II, Editura Academiei, 1980.
\item Moisil, Gr. - " Metoda funcțiilor analitice în hidrodinamica lichidelor vîscoase", in Opera Matematica, Vol. II, Editura Academiei, 1980.
\item Monaco, S., Normand-Cyrot, D. - " Sur la commande non-interactive des systèmes non-linéaires en temps discret", Lecture Notes in Control and Information Sciences, Springer-Verlag, 1984, pp. 364-377.
\item Monaco, S., Normand-Cyrot, D. - " Approximation entrée-sortie d’un système non-linéaire continu par un système discret", Lecture Notes in Control and Information Sciences, Springer-Verlag, 1984, pp. 354-367.
\item Moreno, A., Lagunas, M. - " Envelope/Phase Representation in Signal Modeling", in Signal Processing III. Theories and Applications, Elsevier Science Publishers B.V., Holland, 1986.
\item Naimark, M., Stern, A. - Théorie des représentations des groupes, Éditions Mir, Moscou, 1975.
\item Nijmeijer, H., Schumacher, J.M. - " Zeros at Infinity for Affine Nonlinear Control Systems", " IEEE Transactions on Automatic Control", Vol. AC-30, N. 6, June 1985.
\item Olbort, A., Klein, B. - " Genericity of Open Set Reachability for some Classes of Discrete-Time Nonlinear Systems", " IEEE Transactions on Automatic Control", Vol. AC-31, N. 3, March, 1986.
\item Paraskevopoulos, P.N. - " Chebyshev Series Approach to System Identification, Analysis and Optimal Control", in " Journal of the Franklin Institute", Vol. 361, 1983, pp. 135-157.
\item Paraskevopoulos, P.N., Kekkeris, G. Th. - " Hermite Series Approach to System Identification, Analysis and Optimal Control", Measurement and Control Conference, Atena, Vol. I, 1983, pp. 146-149.
\item Paraskevopoulos, P.N. - " System Analysis and Synthesis via Orthogonal Polynomial Series and Fourier Series", " Mathematics and Computers in Simulation", N. 27, North-Holland, Amsterdam, 1985, pp. 453-469.
\item Pop, E. et al. - Metode în prelucrarea numerică a semnalelor, Editura Facla, Timișoara, 1986.
\item Postnikov, M. - Leçons de géométrie; Algèbre linéaire et géométrie différentielle, Éditions Mir, Moscou, 1981.
\item Postnikov, M. - Leçons de géométrie; Groupes et algèbres de Lie, Éditions Mir, Moscou, 1985.
\item Prasada Rao, G., Tzafestas, S.G. - " A Decade of Piecewise Constant Orthogonal Functions in Systems and Control", " Mathematics and Computers in Simulation", N. 27, 1985, pp. 389-407.
\item Purdea, I. - Tratat de algebră modernă, Editura Academiei, Vol. I, 1977, Vol. II. 1982.
\item Rademacher, H. - " Einige Sätze über Reihen von allgemeinen Orthogonalfunktionen", Math. Ann., Vol. 87, 1922, pp. 112-138.
\item Radu, A. - Probleme de mecanica continuumului, Universitatea " Al. I. Cuza" Iași, 1982.
\item Rao, G.P., Sivakumar, L. - " Transfer Functions Matrix Identification in MIMO Systems via Walsh Functions", Proceedings of IEEE, Vol. 69, Sept. 1981.
\item Roșculeț, M.N. - Funcții monogene pe algebre comutative, Editura Academiei, București, 1975.
\item Rugh, W.J. - " An Extended Linearization Approach to Nonlinear System Inversion", " IEEE Transactions on Automatic Control", Vol. AC-31, N. 8, Aug. 1986.
\item Rugh, W.J. - Nonlinear System Theory: The Volterra-Wiener Approach, Johns Hopkins Univ. Press, Baltimore, MD, 1981.
\item Sage, A., Melsa, J.L. - System Identification, Academic Press, New York, 1971.
\item Schetzen, M. - The Volterra and Wiener Theories of Nonlinear Systems, John Wiley and Sons Inc., New York, 1980.
\item Saito, T. - " A Chaos Generator Based on a Quasi-Harmonic Oscillator", " IEEE Transactions on Circuits and Systems", Vol. CAS-32, N. 4, April 1985.
\item Shankar-Sastry, S. - " The Effects of Small-Noise on Nonlinear Dynamical Systems", " Preprints of the 9th World Congress of the IFAC", Budapest, Vol. 5, 1984, pp. 129-132.
\item Shima, M., Isurugi, Y., Singh, S.W. - " Adaptive Model Following Control Nonlinear Robotic Systems", in " IEEE Transactions on Automatic Control", Vol. AC-30, N. 11, Nov. 1985.
\item Sobrero, L. - Algebra delle Funzioni ipercomplesse e sue applicazioni alla Teoria matematica dell’elasticità, Accad. d’Italia, Memorie della classe di Sci. fiz. mat. e naturali, Roma, N. 6, 1934.
\item Skar, S.J. - " Odd Harmonic Periodic Solutions of Systems with Odd Nonlinearities", " IEEE Transactions on Circuits and Systems", Vol. CAS-32, N. 12, Dec. 1985.
\item Su, R., Hunt, L.R. - " A Canonical Expansion for Nonlinear Systems", " IEEE Transactions on Automatic Control", Vol. AC-31, N. 7, July 1986.
\item Stanomir, D., Stanasila, O. - Metode matematice în teoria semnalelor, Editura Tehnica, București, 1980.
\item Teodorescu, D. - Sisteme automate deterministe, Editura Tehnica, București, 1984.
\item Teodorescu, P.P. - Sisteme mecanice; Modele clasice, Vol. I, II, Editura Tehnica, 1984, 1988.
\item Terashima, K., Akashi, H. - " Explicit Solution of Bilinear Systems with and without Time Delays by Lie Algebraic Approach and its Application to Controllability", " IEEE Transactions on Automatic Control", Vol. AC-31, N. 7, pp. 95-100, July 1986.
\item Vladimirov, V.S. - Ecuațiile fizicii matematice, Editura Științifică și Enciclopedică, București, 1980.
\item Vladimirov, V.S. - Generalized Functions in Mathematical Physics, Mir Publishers, Moscow, 1979.
\item Vladimirov, V.S. et al. - Culegere de probleme de ecuațiile fizicii matematice, Editura Științifică și Enciclopedică, București, 1981.
\item Voicu, M. - Sisteme automate multivariabile. Proiectarea prin metoda frecvențială, Editura " Gh. Asachi", Iași, 1993.
\item Voicu, M. - Teoria sistemelor, Vol. I, II, Institutul Politehnic " Gh. Asachi", Iași, 1980.
\item Voinea, R., Voiculescu, D., Ceaușu, V. - Mecanica, Editura Didactică și Pedagogică, București, 1983.
\item Voinea, R., Stroe, I. - Sisteme dinamice, Institutul Politehnic, București, 1993.
\item Yuan, J.S.C. - " Dynamic Decoupling of a remote Manipulator System", " IEEE Transactions on Automatic Control", Vol. AC-23, N. 4, Aug. 1978.
\item Yuen, C.K. - Upper Bounds on Walsh Transforms, IEEE Transactions on Computers, Vol. C-21, Dec. 1972.
\item Walsh, J.L. - " A closed-set of Normal Orthogonal Functions", Amer. J. Math., Vol. 45, 1923.
\end{enumerate}

\end{document}